\documentclass[aos,reqno]{imsart}

\RequirePackage{amsthm,amsmath,amsfonts,amssymb}
\RequirePackage[numbers]{natbib}
\RequirePackage[colorlinks,citecolor=blue,urlcolor=blue]{hyperref}
\RequirePackage{graphicx}
\usepackage{algorithm2e}
\usepackage{subfiles}
\usepackage{titletoc}
\usepackage{pdfpages}
\usepackage{bbm}
\usepackage{multirow}

\usepackage{mathtools}
\usepackage{algorithm2e}
\usepackage{booktabs}
\RestyleAlgo{ruled}
\SetKwComment{Comment}{/* }{ */}
\usepackage{hyperref}
\usepackage{url}
\usepackage{breqn}
\usepackage{graphicx}
\usepackage{tipa}
\usepackage{bm}
\usepackage{paralist}
\usepackage{calrsfs}
\usepackage{subfiles}
\usepackage{xcolor}
\usepackage{mwe}
\usepackage{dutchcal}
\usepackage[labelfont=bf]{caption}
\usepackage[inline]{showlabels}

\startlocaldefs
\theoremstyle{plain}

\newtheorem{theorem}{Theorem}[section]
\newtheorem{lemma}[theorem]{Lemma}
\newtheorem{remark}{Remark}
\newtheorem{corollary}{Corollary}[theorem]
\newtheorem{proposition}{Proposition}
\theoremstyle{remark}
\newtheorem{definition}[theorem]{Definition}

\newtheorem*{fact}{Fact}

\newcommand{\bb}[1]{\mathbb{#1}}

\newcommand{\scr}[1]{\mathscr{#1}}

\newcommand{\wh}[1]{\widehat{#1}}

\newcommand{\bfa}[1]{\boldsymbol{#1}}
\newcommand{\fk}[1]{\mathfrak{#1}}
\DeclareMathAlphabet{\pazocal}{OMS}{zplm}{m}{n}
\newcommand{\ca}[1]{\pazocal{#1}}
\newcommand{\mca}[1]{\mathcal{#1}}

\newcommand{\nnb}{\nonumber}

\newcommand{\sign}{\text{sign}}

\newcommand{\wha}[1]{\wh{\bfa #1}}
\newcommand{\rank}{\mathtt{rank}}
\newcommand{\lef}{\left}
\newcommand{\rig}{\right}

\newcommand{\tda}[1]{\tilde{\bfa{#1}}}

\newcommand{\tdx}{\tilde x}

\newcommand{\pta}{\partial}

\newcommand{\mbbm}[1]{\mathbbm{#1}}

\newcommand{\tde}{\tilde}

\newcommand{\la}{\langle}
\newcommand{\ra}{\rangle}
\newcommand{\diag}{\text{diag}}

\newcommand{\fron}[1]{\Vert#1\Vert_F}

\newcommand{\ve}[1]{\vect \lef (#1\rig )} 

\newcommand{\ub}[1]{\underbrace{#1}}

\newcommand{\tr}{\text{Tr}}

\newcommand{\hks}{\ca H_{k,S_0}}
\newcommand{\bl}{{\bigg (}}
\newcommand{\br}{{\bigg )}}

\newcommand{\sm}[1]{\small{#1}\normalsize}
\newcommand{\tny}[1]{\footnotesize{#1}\normalsize}
\newcommand{\ttny}[1]
{\scriptsize{#1}\normalsize}
\DeclareMathOperator*{\argmax}{arg\,max}
\DeclareMathOperator{\sech}{sech}
\DeclareMathOperator*{\argmin}{arg\,min}

\endlocaldefs
\pdfoutput=1
\begin{document}

\begin{frontmatter}
\title{Hidden Clique Inference in Random Ising Model I:
the planted random field Curie-Weiss model}
\runtitle{Hidden Clique Inference in the pRFCW Model}

\begin{aug}
\author[A]{\fnms{Yihan}~\snm{He}\ead[label=e1]{yihan.he@princeton.edu}},
\author[B]{\fnms{Han}~\snm{Liu}\ead[label=e2]{hanliu@northwestern.edu}}
\and
\author[A]{\fnms{Jianqing}~\snm{Fan}\ead[label=e3]{jqfan@princeton.edu}}
\address[A]{Princeton University\printead[presep={,\ }]{e1,e3}}

\address[B]{Northwestern University\printead[presep={,\ }]{e2}}
\end{aug}

\begin{abstract}

We study the problem of testing and recovering the hidden $k$-clique Ferromagnetic correlation in the planted Random Field Curie-Weiss model (a.k.a. the pRFCW model). The pRFCW model is a random effect Ising model that exhibits richer phase diagrams both statistically and physically than the standard Curie-Weiss model. Using an alternative characterization of parameter regimes as 'temperatures' and the mean values as 'outer magnetic fields,' we establish the minimax optimal detection rates and recovery rates. The results consist of $7$ distinctive phases for testing and $3$ phases for exact recovery. Our results also imply that the randomness of the outer magnetic field contributes to countable possible convergence rates, which are not observed in the fixed field model. As a byproduct of the proof techniques, we provide two new mathematical results: (1) A family of tail bounds for the average magnetization of the Random Field Curie-Weiss model (a.k.a. the RFCW model) across all temperatures and arbitrary outer fields. (2) A sharp estimate of the information divergence between RFCW models. These play pivotal roles in establishing the major theoretical results in this paper. Additionally, we show that the mathematical structure involved in the pRFCW hidden clique inference problem resembles a 'sparse PCA-like' problem for discrete data. The richer statistical phases than the long-studied Gaussian counterpart shed new light on the theoretical insight of sparse PCA for discrete data.
\end{abstract}



\end{frontmatter}

\section{Introduction}\label{sect1}
We study the problem of testing and recovering positive correlation in the random field Curie-Weiss (or RFCW) model that forms a clique. Letting $\bfa\sigma \in \{-1, +1\}^n$ be $n$ binary-valued random variables,
the RFCW model specifies the distribution of $\bfa\sigma$ as
\begin{align}\label{skdirect}
\bb P(\bfa\sigma)\propto {\exp\bl  \frac{\theta_1}{2n} \sum_{i,j\in[ n]}\sigma_i\sigma_j+\sum_{i=1}^nh_i\sigma_i\br},
\end{align}
where $\{h_i\}_{i\in[n]}$ are realizations of independent random variables and $\theta_1>0$ is a deterministic parameter with the physical meaning of {\emph{inverse temperature}}.  In the following, we refer to the set of $n$ random variables $\bfa\sigma$ as the set of $n$ {spins}. The RFCW model belongs to the class of spin systems in physics characterizing the joint distribution of spins according to their energy function (or {Hamiltonian}) $\mca H^{RFCW}_{\theta_1,\bfa h}:\{-1,+1\}^n\to\bb R$
\begin{align}\label{rfcw-hamiltonian}
    \mca H^{RFCW}_{\theta_1}(\bfa\sigma,\bfa h):=- \frac{\theta_1}{2n} \sum_{i,j\in [n]}\sigma_i\sigma_j-\sum_{i\leq n}h_i\sigma_i,\quad\theta_1 > 0,
\end{align}
which depends on the inverse temperature $\theta_1$ and the random magnetic field $\bfa{h}\in\bb R^n$.  In addition, $\bfa h$ governs the magnitude of the inclination of a certain spin towards positive or negative. The random field in the model makes it semi-parametric and we hope to derive valid statistical procedure that works for an uncountable class of distributions or even without the knowledge of what the actual random field looks like.

The distribution \eqref{skdirect} can be written as a random Gibbs measure depending further on measure $\mu$ supported on $\bb R$ whose $n$-fold product is denoted by $\mu^{\otimes n}$,
\begin{align*}
    \bb P(\bfa\sigma|\bfa h)=\frac{\exp(-\mca H_{\theta_1}^{RFCW}(\bfa\sigma,\bfa h))}{\sum_{\bfa\sigma}\exp(-\mca H_{\theta_1}^{RFCW}(\bfa\sigma,\bfa h))}, \text{ with }\quad \bfa h\sim\mu^{\otimes n}.
\end{align*}
Statistically, the random field Curie-Weiss model is an uncountable mixture of Ising models. Compared with the original Curie-Weiss model, where $h_i$s are deterministic, the random field Curie-Weiss model exhibits more interesting statistical behavior given by the randomness in the field.

For an unknown subset $S\subset[n]$ with $|S|=k$, we define the planted RFCW (or pRFCW) model by its Hamiltonian as
\begin{align}\label{cwhamilt}
    \mca H^{pRFCW}_{\theta_1}(\bfa\sigma,\bfa h):=- \frac{\theta_1}{2k} \sum_{i,j\in S}\sigma_i\sigma_j-\sum_{i\leq n}h_i\sigma_i,\quad\theta_1 > 0.
\end{align}
We study the parameter region of $\theta_1$ and denote it as the inverse temperature of the pRFCW model. The high-temperature regime corresponds to smaller $\theta_1$, and the low-temperature regime corresponds to larger $\theta_1$. We study the coupling coefficient under the scaling of $\frac{1}{k}$, as the phase transition phenomenon happens here \citep{friedli2017statistical}, which makes it the hardest. We later show that the same results obtained in this work under other scalings can be trivially obtained using the identical techniques in this work. Specifically, our results in the high temperature regime can cover the $\theta_1$ vanishing case (with an additional scaling factor of $\theta_1$). And our results for the low temperature regime can cover the $\theta_1$ diverging case.

We denote $\ca G_0(n)$ as the null hypothesis where the Hamiltonian  is given by
\begin{align*}
    \mca H_{0}=-\sum_{i\leq n}h_i\sigma_i,\text{ with } \quad\bfa h\sim\mu^{\otimes n}.
\end{align*}
and $\ca G_1(\theta_1, k, n)$ as the alternative hypothesis of the pRFCW model with Hamiltonian defined by \eqref{cwhamilt} and the set $S$ unknown. Let $\{\bfa\sigma^{(1)}, \ldots, \bfa\sigma^{(m)}\}$ be $m$ independent observations. We aim to test between the hypotheses of
\begin{align}\label{mainproblem}
     \text{Null}: \bfa\sigma^{(i)}\sim \ca G_0(n) \qquad \text{vs}\qquad \text{Alternative}:\bfa\sigma^{(i)}\sim \ca G_1(\theta_1,k,n).
\end{align}
The statistical diagram for testing involves the existence of asymptotic powerful tests, characterized as follows. And the problem of hypothesis testing is fundamental in the history statistics, we refer to \citep{lehmann1986testing} for a historical exposition.
\begin{definition}[Asymptotic Power of Tests]\label{def1}\label{def2}
Let $\bb P_{0}$  be the probability measure under the null. Let $S_0$ be the set of vertices of the hidden clique with size $|S_0|=k$ and define $\bb P_{S_0}$ as the probability measure under the alternative. Let $\bb P_{0,m}$ and $\bb P_{S,m}$ be their product of $m$ measures respectively. We define a sequence of tests $\psi:\{\bfa\sigma^{(i)}\}_{i\in[m]}\to\{0,1\}$ depending only on $m$ i.i.d. samples to be asymptotically powerful if we have 
  \begin{align*}
       \lim_{k\to+\infty} \bigg[\bb P_{0,m(k)}\lef(\psi=1\rig)+\sup_{S:|S|=k}\bb P_{S,m(k)}\lef(\psi=0\rig)\bigg]=0
    .\end{align*}
And we say all tests are asymptotically powerless if test statistics sequence depending on $m$ i.i.d. samples satisfy
\begin{align*}
    \lim_{k\to+\infty}\inf_{\psi}\bigg[\bb P_{0,m(k)}(\psi=1)+\sup_{S:|S|=k}\bb P_{S,m(k)}(\psi=0)\bigg]=1.
\end{align*}
\end{definition}

We emphasize that, unlike classical asymptotic theory where we let $m \rightarrow \infty$, the asymptotic setup in Definition \ref{def2} is driven by $k\rightarrow \infty$. This is because we are considering correlated random variables and in certain regimes even with a single sample $m=1$ is enough for the test (More details are presented in Section \ref{sect3}).

In addition to the testing result, we also provide theoretical results for almost exact recovery and the exact recovery defined below.
\begin{definition}[Recovery Guarantees]\label{recoguarantee}
Let $\Delta$ be the symmetric difference between two sets and $S\subset[n]$ be the index set of a clique. Let $\wh S$ be the estimated index set, we define:
\begin{itemize}
    \item  Exact Recovery:  if  $\bb P(|\wh S\Delta S|=0)=1-o(1)$;
    \item  Almost Exact Recovery:  if $\bb P(|\wh S\Delta S|=o(k))=1-o(1)$.
\end{itemize}
\end{definition}
\subsection{Motivations} 
This work is motivated by the requirement to understand the inference performance on discrete data, which abounds in real-world data analysis. Some typical examples include the voting data from presidential elections and treatment vectors in the causal inference literature.  In particular, we are interested in the statistical phenomenon given by the subset selection and testing procedure under the high dimensional setup. Under this setup we are interested in finding the subset of spins with larger correlations or testing its existence. This problem also shares fundamental similarities to the problem of probablistic PCA, where we are given a graphical model whose projection on a subspace has larger variance than the rest.  At a more fundamental level, we hope to answer the questions: \emph{ Does the algorithm for real-valued data and continuous distribution assumption continue to have solid theoretical guarantees given the discrete data? What are some new phenomena for the discrete distributions in high dimensional inference problems?} 

To achieve the above goal, we use the pRFCW model as a stylized semiparametric model to characterize the unknown dependence among the binary data distributions. Apart from the classical parametric Ising model, the random field in the pRFCW model makes this model to have expressiveness for real-world distributions. Moreover, it can also model a joint binary distribution with an unknown but latent mean. Our results also showcased that the performance guarantee does not rely on the exact distribution of the random field but some typical moment statistics. We are interested in the asymptotics of $n\to\infty$ and $k\to\infty$ with $n$ at different rates, being of typical regimes in classical high dimensional statistics.

From the theoretical side, we show that a unique phase transition phenomenon appears in the discrete graphical model as the inverse temperature  $\theta_1$  approaches some threshold where most of the covariates take the same sign. This is often referred to as the `simultaneous magnetization' in statistical physics literature \citep{ellis2006entropy}. In particular, during the phase transition, the covariates move from being `almost independent' to `ultra-correlated'. Such a phenomenon never appears in Gaussian distributions, and its effect on the high dimensional statistical inference problem, including the PCA and subset selection problems, is largely unknown. Hence, this work set the first step toward analyzing its effect on the minimax rates compared with the well-studied multivariate Gaussian variants \citep{amini2009high,berthet2013complexity}. In particular, at the center of the many new phenomena obtained in this work is the multiple phase transitions at the interplay between the 'simultaneous magnetization' phenomenon and the asymptotic log-ratio between the sparsity $k$ versus the dimension $n$.

\subsection{Contributions}

Our major contribution is establishing sharp upper and lower bound results of the statistical rates for the test and recovery problem in \eqref{mainproblem}. In detail, we show that a diverse phenomenon of multi-level phase transition appears in the test of a hidden clique in the pRFCW model.

\begin{itemize}
\item \textbf {(Top level phase transition characterized by $\theta_1$)} There exists a ciricital parameter $\theta_c$ partitioning the space of $\theta_1$ into three regimes: $\theta_1 < \theta_c$ (high-temperature regime), $\theta_1 = \theta_c$ (critical temperature reigme), $\theta_1 > \theta_c$ (low-temperature regime), and the minimax sample complexity differ significantly across these regions.
\item \textbf{(Intermediate level phase transition characterized by $k$)} We observe that a `mountain climbing' phase transition characterized by $k$ appears at all temperature regimes: when $k\gtrsim n^{\beta}$ for some $\beta\in(0,1)$, the optimal sample complexity is achieved by global tests taking all the spins as input; when $k=o(n^{\beta})$, the optimal complexity is achieved by a class of local scan tests. Moreover, $\beta$ differs across the temperature regimes and can take countable possible values at the critical temperature.
\item \textbf {(Bottom level phase transition characterized by $h$)} Depending on the tail heaviness of  $h$, the phase diagram at the critical temperature can vary significantly, where the optimal rate can take countable values. Moreover, the critical temperature represents an intermediate state between the high and low temperatures where we identify a co-existence of statistical phases in the high and low temperatures.
 
\end{itemize}

To prove the minimax rates for testing, we extend the ideas of (multivariate) Laplace approximation of exponential integral \citep{bolthausen1986laplace} and the transfer principle in \citep{ellis1980limit} to random measures and give concise proofs on the limiting theorem for the average magnetization across all temperature regimes. We also prove a novel all-temperature tail bound with a sharp rate. Based on this method, we give a sharp estimate of the information divergence between the random Gibbs measures, providing an optimality guarantee for the proposed tests. We also propose another method based on the construction of \emph{fake measure} as a sharpening tool for the information divergence of the close-to-critical temperature, completing the phase diagrams.

We utilize the local optimal tests to derive algorithms for almost exact recovery (or weak model selection consistency) to obtain the optimal complexity for exact recovery. Then, we utilize a screening procedure to obtain the exact recovery (or almost sure model selection consistency). The lower bounds are constructed using a leave-one-out prior and Fano's inequality. To derive a sharp estimate of information divergence, we utilize the steepest descent method and a projection procedure.
\smallbreak
\emph{Organization.} The rest of this paper is organized as follows: First, we summarize all the necessary notations; Section \ref{sect2} contains a discussion and review of related works; Section \ref{sect3} states the major results under centered $h$, including the upper and lower bounds at different parameter regimes and the algorithms achieving the optimal rates for both testing and exact recovery; Section \ref{sect4} complements section \ref{sect3} with the results under non-centered $h$; Section \ref{sect5} contains the statement and proof of the CLT for the RFCW model; Section \ref{sect6} discusses the computational feasible algorithms and the computational statistical gaps. We delay extended proofs and technical details to the supplementary material. 
\smallbreak
\emph{Notations.} 
The following notations are used throughout this work. We use $:=$ as the notation for \emph{defining}. We denote $[n]:=\{1,\ldots,n\}$ and $[i:j]:=\{i,i+1,\ldots, j\}$ for $i<j$. For a vector denoted by $\bfa v=(v_1,\ldots v_n)\in\bb R^n$ we denote its $\ell_p$ norm by $\Vert \bfa v\Vert_p=\lef(\sum_{i=1}^p v_i^p\rig)^{1/p}$ for all $p\in[1,\infty)$. Denote $\Vert \bfa v\Vert_\infty=\sup_{i\in[n]}|v_i|$. For a matrix $A\in\bb R^{n\times m}$ with $m,n\in\bb N$ we denote $\Vert A\Vert_\infty = \sup_{i,j}|A_{ij}|$, $\Vert A\Vert_F=\lef(\sum_{i,j}A_{ij}^2\rig)^{1/2}$ and $\Vert A\Vert_p = \sup_{\bfa v:\Vert \bfa v\Vert_p=1}\Vert A\bfa v\Vert_p$ for all $p\in[1,\infty)$. For a vector $\bfa v\in\bb R^n$ and set $A\subset[n]$, we denote $\bfa v_{-A}$ to be the vector constrained to $A^c$. We denote $\mbbm 1_{B}$ for some event $B$ as the indicator function of $B$. Moreover, for some set $A\subset [n]$, we denote $\bfa v=\mbbm 1_A\in\bb R^{n}$ if $v_i=\mbbm 1_{i\in A}$. For a set $A\subset \Omega$ we denote $A^c=\Omega\setminus A$ where $\backslash$ is the notation for set minus. For another set $B\setminus \Omega$, we denote  $A \Delta B=(A\cup B)\setminus(A\cap B)$ to be the symmetric difference between $A$ and $B$. For a function $f(x)$ and $\tau\in\bb N$, we define $f^{\tau}(x)$ to be its $\tau$-th derivative at $x$. Throughout this work, we use $\mca i$ as the notation for the imaginary unit. Let $\bb P$ be a probability measure and $\bb P^{\otimes n}$ be the $n$-th order product measure of $\bb P$. Given two sequences $a_n$ and $b_n$, we denote $a_n\lesssim b_n$ or $a_n = O(b_n)$ if $\limsup_{n\to\infty}\lef|\frac{a_n}{b_n}\rig|<\infty$ and $a_n=o(b_n)$ if $\limsup_{n\to\infty}\lef|\frac{a_n}{b_n}\rig|=0$. Similarly, we denote $a_n\gtrsim b_n$ or $a_n =\Omega(b_n)$ if $b_n = O(a_n)$ and $a_n=\omega(b_n)$ if $b_n=o(a_n)$. We denote $a_n\asymp b_n$ or $a_n=\Theta(b_n)$ if $b_n\lesssim a_n$ and $a_n\gtrsim b_n$ both hold. For two sequence of measurable functions $f_n,g_n$ with $n\in\bb N$, we denote $f_n= O_p(g_n)$ if for all $\epsilon>0$ there exists $C>0$ such that $\limsup_n\bb P\lef (|f_n|>C|g_n|\rig )\leq \epsilon$  and $f_n = o_p(g_n)$ if for all $\delta>0$ $\limsup_{n}\bb P\lef(|f_n|>\delta |g_n|\rig) = 0$. We denote all $z$ in this work as standard Gaussians. Regarding convergence, we denote $\overset{d}{\to}$ to be convergence in distribution. We denote $X\perp Y$ if two random variables are independent. If $\psi$ is a monotonic nondecreasing, convex function with $\psi(0)=0$, the Orlicz norm of an integrable random variable $X$ with respect to a function $\psi:\bb R\to\bb R^+$ is given by $\Vert X\Vert_{\psi}=\sup\lef\{u>0:\bb E\lef[\psi\lef(\frac{|X|}{u}\rig)\rig]\leq 1\rig\}$. In particular, for $\theta\in\bb R^+$ we use the notation of $\psi_{\theta}(x):=\exp(x^\theta)-1$. And we introduce the $O_{\psi_\theta}$ notation as follows: $A-B=O_{\psi_{\theta}}(C)\Leftrightarrow \Vert A-B\Vert_{\psi_\theta}\lesssim \Vert C\Vert_{\psi_{\theta}}$. Similarly, we define $A-B=o_{\psi_{\theta}}(C)\Leftrightarrow \Vert A-B\Vert_{\psi_\theta}=o\lef( \Vert C\Vert_{\psi_{\theta}}\rig)$. We denote $\mca i$ to be the imaginary unit. For a function $f:\bb R\to\bb R$ we denote $f^{(i)}$ to be the $i$-th derivative of $f$. Finally, all the constants denoted by $C$ in this work are ad hoc, and we do not attempt to optimize them.

\section{Related Work}\label{sect2}

A rich literature studies statistical problems on Ising models.  \citep{neykov2019property,cao2022high} study the problem of detecting hidden combinatorial structures in the Ferromagnetic Ising model. In \citep{neykov2019property}, the authors studied the detection of exact combinatorial structure (i.e., the existence versus non-existence), including cliques. However, their work did not achieve optimal sample complexity. \citep{cao2022high} gave the statistical upper and lower bounds for testing general properties (e.g., bi-cliques, $k$-stars, etc.) in the Ferromagnetic Ising model. (We recall that ferromagnetism refers to having only a positive correlation parameter between spins.) They point out that a structure that requires special focus in testing general graph properties is the $k$-clique where the author provides a lower bound for the general graph by the lower bound for testing the smallest clique containing it. This implies that clique is a bottleneck test for general property tests.  However, their result is not optimal in sample complexity. Moreover, both \citep{neykov2019property, cao2022high} considered only the high-temperature region of the Ising Ferromagnetic model, while our results cover all temperature regimes, exhibiting a new phase transition pattern. Moreover, our result gives optimal sample complexity for the planted SK model, which takes the clique test in \citep{cao2022high} as a special case of $\bfa h=0$.  The method of analysis in \citep{cao2022high,neykov2019property} does not hold in our problem since a few essential theoretical results in Ferromagnetic Ising models do not hold when a random field presents. \citep{berthet2019exact} studied the problem of exact recovery of two latent mean field groups in the Ising model. \citep{bresler2020learning,boix2022chow} studied algorithms for Ising tree reconstruction.  In particular, they pointed out that prediction may not necessarily require the exact reconstruction. 

Another line of work on mean field Ising Ferromagnetic models studied the tests against the existence of sparse outer-magnetic field \citep{mukherjee2018global,deb2020detecting, bhattacharya2021sharp}. Though their problem settings are different from ours, a key intuition in both their and our works is that the phase transition of spin systems results in the phase transition of statistical inference, making an optimal test in one phase diagram regime sub-optimal or even invalid in another regime.

For the statistical analysis of random Gibbs measure, \citep{chatterjee2007estimation} studied the pseudo likelihood estimator for the inverse temperature in a few spin glass models.  Our problem differs from theirs since they focus on parameter estimation, whereas we focus on testing and recovering the planted structures.  Their analysis does not apply to our setting since his estimation method requires the correlation matrix among all spins to be known. 

Beyond the literature with Gibbs measure, a rich line of work also focuses on the detection of correlations in Gaussian models. For example, \citep{addario2010combinatorial} studied the testing procedure of specific subsets of components in a Gaussian vector. \citep{arias2012detection} consider the correlation graph of a Gaussian random vector and study the problem of detecting certain classes of fully connected cliques. However, the Gibbs measure based models are more subtle to analyze and the mathematical tools are largely different.

From the technical perspective, \citep{amaro1991fluctuations} attempts to study the weak convergence of average magnetization of the RFCW model. However, the result in \citep{amaro1991fluctuations} has a counterexample when $\theta_1=0$. Therefore, our result also serves as a correction for their results. On the other hand, our method gives the tail bound whereas their method does not. Despite that for fixed effect Curie-Weiss model, \citep{chatterjee2007stein,chatterjee2010applications} give tail bound based on the Chatterjee-Stein's method of exchangeable pair, their method cannot be directly applied to random Gibbs measures.  Hence, to the best of the author's knowledge, the tail bound does not exist in the previous literature for the RFCW model. On the other hand, previous work \citep{mukherjee2018global} studied the non-existence of an asymptotically powerful test based on the Neyman-Pearson Lemma, which is also not directly generalizable to random Gibbs measure. Our work follows another path inspired by \citep{neykov2019property,cao2022high} where the \emph{all or nothing} type of results are provided. This result implies the non-existence of any test that has an asymptotic minimax Type I + Type II error smaller than $1$, which also implies the nonexistence of powerful tests in their work. To achieve this result, we utilized the same method proving the CLT type of result for the RFCW model. Moreover, this method suffers from vacuity when it comes to the close-to-critical temperature. To address this issue we provide a \emph{fake measure method} that sharpens the divergence.

A rich literature of works studying the sparse PCA problems \citep{amini2009high, krauthgamer2013semidefinite, berthet2013complexity, berthet2019exact, vu2013minimax, wang2016statistical,cai2013sparse} under the Gaussian/constrained sub-Gaussian variants. We note that mathematically the classical zero field Curie-Weiss model shares some mathematical similarities with the single spike $k$-sparse PCA model despite that the Ising Gibbs measures are more delicate to analyze than regular Gaussians. Comparing our results with theirs, we show that the rate of ${k\log n}$ that is prevalent in most Sparse PCA literature is only the optimal rate at the high-temperature regime. At the critical and low temperature, this optimal rate can be much smaller and the pRFCW model demonstrates a richer statistical phase diagram than the classical Sparse PCA literature. This also implies an interesting phenomenon when the random variables become binary or categorical. 

Another stream of work in the Sparse PCA problem considers the observation to be a single sample sampled from a spiked random matrix model \citep{onatski2013asymptotic,montanari2015limitation,deshpande2014information,perry2016optimality}, with the Spiked Wigner model being the simplest one. This line of work is under a different formulation where only a single sample from the random matrix model is observed with entry-wise independence. Despite sharing similarities in the name and the problem formulation, the mathematical model is largely different from the one considered in this paper since we consider statistical estimation problem based on multiple samples from a graphical model.

\section{Major Results}\label{sect3}
For presentation clarity, this section only presents theoretical results  when the distribution  of $\tanh(h)$ is centered at $0$:
$\bb E[\tanh(h)]=0$.
We present the more challenging case of non-centered $\tanh(h)$ in section \ref{sect4}.
The major results of this section are summarized in table \ref{table1} for the testing problem and table \ref{table3} for the recovery problem. 

\begin{table}[htp!]\small
\captionsetup{labelfont=bf}
\caption{\textbf{The minimax sample complexity of testing when $\tanh(h)$ is centered.} Shown in the table is the minimum sample size m, as a function of planted click size $k$ and the number of spins $n$. At all temperature regimes, we observe an `ascending-descending' phenomenon of the minimax rate as the clique size $k$ gets larger: Fixing $n$, the optimal sample complexity goes through a monotonically increasing phase followed by a monotonically decreasing phase as $k$ gets larger. And the complexity peaks at a middle point.  This unique phenomenon results from the interaction of two contradicting factors: (1) A larger $k$ makes the clique more observable and reduces the complexity; (2) In the meantime, a larger $k$ also makes the spins in the clique less correlated (due to the existence of the scaling factor $1/k$) and increases the complexity. We also observe that the critical temperature regime undergoes $3$ phases rather than $2$ at the high/low-temperature regimes. This is because the critical temperature is an intermediate phase between the high/low temperatures and exhibits a phase diagram mixing both.}
  \renewcommand*{\arraystretch}{2.4}
\makebox[0pt]{
\begin{tabular}{cl|l|ll}
\toprule
\multicolumn{2}{c|}{\textbf{Centered} $\tanh(h)$ Testing}                                                                                                  & Small Clique Regime                                                              & \multicolumn{2}{c}{Large Clique Regime}                                                                                                            \\ \midrule
\multicolumn{2}{l|}{\textbf{High Temperature}}                                                                                & $k=o\lef(n^{\frac{2}{3}}\rig)$                                                   & \multicolumn{2}{l}{$n^{\frac{2}{3}}\lesssim k\leq n$}                                                                                              \\ \hline
\multicolumn{1}{l|}{\multirow{2}{*}{$\theta_1\in\lef(0,\frac{1}{2\bb E[\sech^2(h)]}\rig)$}}                            & UBs & $O\lef(k\log n\rig)$                                                    & \multicolumn{2}{l}{$O\lef(\frac{n^2}{k^2}\rig)$}                                                                                                   \\ \cline{2-5} 
\multicolumn{1}{l|}{}                                                                                          & LBs & $\Omega\lef(k\log n\rig)$                                               & \multicolumn{2}{l}{$\Omega\lef(\frac{n^2}{k^2}\rig)$}                                                                                              \\ \hline
\multicolumn{1}{l|}{\multirow{2}{*}{$\theta_1\in\lef[\frac{1}{2\bb E[\sech^2(h)]},\frac{1}{\bb E[\sech^2(h)]}\rig)$}} & UBs & $O\lef(k\log n\rig)$                                                    & \multicolumn{2}{l}{$O\lef(\frac{n^2}{k^2}\rig)$}                                                                                                   \\ \cline{2-5} 
\multicolumn{1}{l|}{}                                                                                          & LBs & $\Omega\lef(\frac{k}{\log k}\log n\rig)$                                & \multicolumn{2}{l}{$\Omega\lef(\frac{n^2}{k^2}\rig)$}                                                                                              \\ \hline
\multicolumn{2}{l|}{\textbf{Critical Temperature}}                                                                            & $k=o\lef(n^{\frac{4\tau-2}{8\tau-5}}\rig)$                                       & \multicolumn{1}{c|}{$n^{\frac{4\tau-2}{8\tau-5}}\lesssim k\lesssim n^{\frac{2\tau-1}{4\tau-3}}$} & $k=\omega\lef(n^{\frac{2\tau-1}{4\tau-3}}\rig)$ \\ \hline
\multicolumn{2}{l|}{Upper Bounds}                                                                                             & $O(k^{\frac{1}{2\tau-1}}\log n)$                                        & \multicolumn{1}{l|}{$O\lef(n^2k^{-\frac{2(4\tau-3)}{2\tau-1}}\rig)$}                                    & $O(1)$                                          \\ \hline
\multicolumn{2}{l|}{Lower Bounds}                                                                                             & $\Omega\lef(\lef(\frac{k}{\log k}\rig)^{\frac{1}{2\tau-1}} \log n\rig)$ & \multicolumn{1}{l|}{$\Omega\lef(n^2k^{-\frac{2(4\tau-3)}{2\tau-1}}\rig)$}                               & $\Omega(1)$                                     \\ \hline
\multicolumn{2}{l|}{\textbf{Low Temperature}}                                                                                 & $k=o\lef(n^{\frac{1}{2}}\rig)$                                                   & \multicolumn{1}{l|}{$k\asymp n^{\frac{1}{2}}$}                                                         & $k=\omega\lef(n^{\frac{1}{2}}\rig)$             \\ \hline
\multicolumn{2}{l|}{Upper Bounds}                                                                                             & $O\lef(\log n\rig)$                                                     & \multicolumn{1}{l|}{$O(1)$}                                                                      & \multicolumn{1}{l}{1}                           \\ \hline
\multicolumn{2}{l|}{Lower Bounds}                                                                                             & $\Omega\lef(\log n\rig)$                                                & \multicolumn{1}{l|}{$\Omega(1)$}                                                                 & \multicolumn{1}{l}{1}                           \\ \bottomrule
\end{tabular}}\label{table1}
\end{table}

\smallbreak
\emph{Testing.} Our first main result characterizes a family of multi-level phase transition diagrams of the sample complexity in testing a hidden clique in the pRFCW model: \textbf{At the top level}, let $h\sim\mu$, we characterize a critical parameter $\theta_c:=\frac{1}{\bb E[\sech^2(h)]}$ to partition the space of $\theta_1$ into three regimes: $\theta_1 < \theta_c$ (high-temperature regime), $\theta_1 = \theta_c$ (critical temperature regime), $\theta_1 > \theta_c$ (low-temperature regime), and show that the optimal sample complexity differs significantly across these regions.  \textbf{At the middle level}, we keep the low and high-temperature regimes untouched but further partition the critical temperature regime into infinitely countable sub-regimes based on the flatness of the RFCW distribution with Hamiltonian characterized in \eqref{rfcw-hamiltonian}. More specifically, the flatness parameter $\tau$ depends on the order of the first non-zero derivative at $0$ of a univariate characteristic function, which we discuss in \eqref{taylorcond}.
\begin{align*}
    H(x):=\frac{1}{2}x^2-\bb E[\log\cosh(\sqrt{\theta_1}x+h)].
\end{align*}
We delay more detailed definitions in section \ref{sect33}. \textbf{At the bottom level}, we further partition the obtained regimes into more subregimes according to the scaling of the hidden clique size $k$ concerning $n$. In particular, for a fixed $n$, we partition the high temperature regime into two subregimes by $k\asymp n^{\frac{2}{3}}$; We partition the low temperature regime into two subregimes by $k\asymp n^{\frac{1}{2}}$; We partition the critical temperature regime into three subregimes by $k\asymp n^{\frac{4\tau-2}{8\tau-5}}$ and $k\asymp n^{\frac{2\tau-1}{4\tau-3}}$. The intuition underlying the extra partition at the critical temperature comes from its mixture of phases original to the high and low temperatures. 

To obtain the upper bounds we provide two types of tests: local vs. global. The local tests construct $\binom{n}{k}$ test statistics for each $k$-subset and enumerate all the $k$ subsets to optimize a criterion. The global tests construct test statistics using all the spins at once. We show that they achieve optimal sample complexity for the large and small clique regimes stated in table \ref{table1}. Without loss of generality, we assume the temperature regimes are known. Otherwise, we can construct a simple adaptive procedure carrying out all the tests together and reject the null when any of them are rejected.

\smallbreak
\emph{Recovery.} Table \ref{table3} summarizes our main results in exact recovery. The almost exact recovery results are also proved in sections \ref{sect31}, \ref{sect32}, \ref{sect33}. In particular, our results imply that the scan tests used to attain the complexity upper bounds of testing small cliques and provide almost exact recovery algorithms. We further propose a screening procedure in section \ref{exactreco} to boost the almost exact recovery algorithm to the exact recovery. This implies that for the small clique regimes, there exists no test-recovery gap.

\begin{table}[htp!]
\caption{\textbf{The minimax sample complexity of exact recovery.} Compared with the `ascending-descending' phenomenon in table \ref{table1}, the minimax sample complexity for exact recovery is monotonically increasing across all temperature regimes. In particular, the sample complexity matches the `ascending' phase of the testing but differs significantly from the `descending' phase. This implies that the difficulty of testing is the same as exact recovery for small cliques, whereas for large cliques, recovery is significantly larger than testing.}
  \renewcommand*{\arraystretch}{3.5}
\centering
\begin{tabular}{ll|l|l}
\toprule
\multicolumn{2}{l|}{Exact Recovery}                                                 & Upper Bounds                     & Lower Bounds                          \\ \midrule
\multicolumn{1}{l|}{\multirow{3}{*}{$\tanh(h)$ is centered}} & High Temperature     & $O(k\log  n)$                     & $\Omega(k\log  n)$                     \\ \cline{2-4} 
\multicolumn{1}{l|}{}                                        & Low Temperature      & $O(\log  n)$                      & $\Omega(\log  n)$                      \\ \cline{2-4} 
\multicolumn{1}{l|}{}                                        & Critical Temperature & $O(k^{\frac{1}{2\tau-1}}\log  n)$ & $\Omega(k^{\frac{1}{2\tau-1}}\log  n)$ \\ \hline
\multicolumn{2}{l|}{$\tanh(h)$ is non-centered}                                     & $O(\log n)$                      & $\Omega(\log n)$                      \\ \bottomrule
\end{tabular}\label{table3}
\end{table}

\smallbreak
\emph{Organization.} The rest of this section is organized as follows: section \ref{sect30} presents a technical overview of our proof strategy; section \ref{sect31} presents the optimal test rates for the high-temperature regime; section \ref{sect32} presents the minimax rates for the low-temperature regime; section \ref{sect33} presents the minimax rates at the critical temperature; section \ref{exactreco} presents the minimax rates for exact recovery.

\subsection{Overview of  Techniques}\label{sect30}
Here, we give a brief technical overview of the major methods in this work. Our framework for the upper bounds is based on the concentration inequalities for the average magnetizations in the RFCW model under different temperature regimes. Our framework for the lower bounds is based on Le Cam's method and Fano's inequality, which needs a sharp estimate of information divergences between the null and alternative hypotheses. Despite these frameworks being standard, the tail bounds for the RFCW model and the information divergences for the mixture Gibbs measures (including the RFCW measure) remain challenging problems. Existing literature gives tail bound results for the classical Curie-Weiss model (for example, \citep{chatterjee2010applications} used Stein's exchangeable pair method). Moreover, their method cannot work for the mixture Gibbs measure. For the information divergence, previous works \citep{neykov2019property,cao2022high} propose a few methods, including the polynomial expansion of the non-mixture Ising models. However, the existence of magnetic fields invalidates these methods. Also, their methods only work for extremely high temperatures, and the derivation of the results for all temperature regimes poses another challenge.

In summary, two major technical problems are solved in this work. (1) Existing literature does not give any tail bound for the average magnetization of the RFCW model, and the prevalent method for the non-random field Ising model does not work here \citep{chatterjee2010spin}. (2) The estimate of information divergence of discrete random Gibbs measures like the RFCW model is difficult, and no existing literature covers the results in this work. We discuss our strategies to overcome these barriers.
\smallbreak
\emph{Hubbard-Stratonovich Transformation.} Our centered strategy to solve the above two problems utilizes asymptotic integral expansion along with the Hubbard–Stratonovich (H-S) transformation. For example, we consider a random Gibbs average of a function $f:\{-1,1\}^n\to\bb R$ concerning the RFCW measure of the following form.
\begin{align}\label{gibbsaverageform}
    \bb E[f]&=\bb E\bigg[\frac{\sum_{\bfa\sigma}f(\bfa\sigma)\exp(\frac{\theta_1}{2k}\sum_{i,j\in[n]}\sigma_i\sigma_j+\sum_{i=1}^n\sigma_ih_i)}{\sum_{\bfa\sigma}\exp(\frac{\theta_1}{2k}\sum_{i,j\in[n]}\sigma_i\sigma_j+\sum_{i=1}^n\sigma_ih_i)}\bigg]
\end{align}
This Gibbs average form \eqref{gibbsaverageform} appears naturally when we hope to obtain the moment generating function and the information divergences.

Then, the H-S transformation utilizes the Gaussian moment generating and characteristic function given by
\small
\begin{align*}
    \int_{\bb R}\frac{1}{\sqrt{2\pi}}\exp\bl-\frac{x^2}{2}+yx\br dx=\exp\bl\frac{y^2}{2}\br, \quad\int_{\bb R}\frac{1}{\sqrt{2\pi}}\exp\bl-\frac{x^2}{2}+\mca iyx\br=\exp\bl-\frac{y^2}{2}\br.
\end{align*}
\normalsize
Then we notice that $\sum_{i,j=1}^n\sigma_i\sigma_j=(\sum_{i=1}^n\sigma_i)^2$, which implies that
\begin{align*}
    \bb E[f]&=\bb E\bigg[\frac{\int_{\bb R}\sum_{\bfa\sigma}f(\bfa\sigma)\exp(-\frac{x^2}{2}+\sum_{i=1}^n(\sqrt{\frac{\theta_1}{k}}x+h_i)\sigma_i)dx}{\int_{\bb R}\sum_{\bfa\sigma}\exp(-\frac{x^2}{2}+\sum_{i=1}^n(\sqrt{\frac{\theta_1}{k}}x+h_i)\sigma_i)dx}\bigg].
\end{align*}
And, if $f$ is a 1-degree polynomial of $\bfa\sigma$, the sum of individual spins in the denominator and numerator can be computed explicitly. After this transformation, the quadratic form in the exponential disappears, and the sum becomes simple. This implies that the mean Gibbs average finally boils down to understanding the following integral for some functions $\mca H_{0,n}(x,\bfa h)$ and $\mca H_{1,n}(x,\bfa h)$ as
\begin{align}\label{transformedgibbsaverage}
    \bb E[f]=\bb E\bigg[\frac{\int_{\bb R}\exp(-n\mca H_{0,n}(x,\bfa h))dx}{\int_{\bb R}\exp(-n\mca H_{1,n}(x,\bfa h))dx}\bigg].
\end{align}
We show that the convexity of $\mca H_{1,n}$ concerning $x$ is closely connected to the different temperature regimes that we consider.
To estimate the value of \eqref{transformedgibbsaverage}, we need to utilize asymptotic integral expansion in the univariate and multivariate cases.
\smallbreak
\emph{Asymptotic Integral Expansion.}
 The Laplace method gives the following asymptotic equivalence condition considering a convex function $g$ with unique global minimum $x^*$ and $g^{(2)}(x^*)>0$:
\begin{align*}
    \int_{\bb R}\exp(-ng(x))dx = \sqrt{\frac{2\pi}{ng^{(2)}(x^*)}}\exp(-ng(x^*))\bl 1+O\bl\frac{1}{n}\br\br.
\end{align*}
And we also encounter the integral with $g(x)$ being a function of complex variables. Then, we need a more general form of the Laplace method that accounts for the complex integral.

In the appendix, we also provide a version where this univariate integral expansion is further generalized to the multivariate and stochastic setting. The randomness in the field $h$ results in the randomness of $x^*$ in the above integral. For example, we see that the global minimum point $x^*$ of $\mca H_{1,n}(x,\bfa h)$ in \eqref{transformedgibbsaverage} is, in fact, a fixed point of a random function:
\begin{align}\label{fixedpointanalysis}
    x^*=\frac{\sqrt{\theta_1}}{n}\sum_{i=1}^n\tanh(h_i+\sqrt{\theta_1}x^*),
\end{align}
with $h_i$ i.i.d. Then we use the theory of $Z$-estimators. However, we need stronger control of the tail since the asymptotic expansion of $Z$-estimators only gives weak convergence results, and we need convergence of moment generating function (or m.g.f.). Hence we introduce a notation $o_{\psi_2}$ as the notation for terms with vanishing sub-Gaussian norm. ( a.k.a. Orlicz norm with $\psi_2(x):=\exp(x^2)-1$) This corresponds to stronger control of the $o_p$ terms in the classical linearization of $Z$-estimator \citep{van1996empirical} where instead of having $o_p(1)$ terms we have $o_{\psi_2}(1)$ terms. To account for this, we prove a lemma for the \emph{Strong Linearization of Z-estimators}. Such a term helps us obtain the asymptotics of the moment generating function.

This is pivotal to derive the upper bound on the m.g.f. of the average magnetization in section \ref{sect5}, which gives a tail bound for the RFCW measure. However, the above analysis can only cover high and critical temperature regimes; low temperature regimes need the transfer principle as an extra technical tool. 
\smallbreak
\emph{The Transfer Principle.} It is known \citep{ellis2006entropy} that at the low-temperature regime, the Curie-Weiss model demonstrates a spontaneous magnetization phenomenon where the correlation between spins is extremely strong. Instead of concentrating on a single point, the unique phenomenon of concentration on two separate points appears. This translates to the two global minima of both the functions in the numerator and denominator of \eqref{transformedgibbsaverage} that invalidate the expansion analysis. Therefore, we need a transfer principle as the new analysis method. This idea is first used by  \citep{ellis1980limit} to prove the weak limit of the meta-stable state in the standard Curie-Weiss model. Their results imply that a conditional sum in the spins can be transferred to a partial integral region of the $x$ after the H-S transform. Formally, their results imply that for an event $\ca C:=\{\sum_{i=1}^n\sigma_i>0\}$, there exists $C>0$ such that
\begin{align*}
    \frac{\sum_{\bfa\sigma\in\ca C}f(\bfa\sigma)\exp(\frac{\theta_1}{2n}\sum_{i,j\in[n]}\sigma_i\sigma_j)}{\sum_{\bfa\sigma\in\ca C}\exp(\frac{\theta_1}{2n}\sum_{i,j\in[n]}\sigma_i\sigma_j)}=\frac{\int_{B_0}\exp(-n\mca H_{0,n}(x))dx}{\int_{B_1}\exp(-n\mca H_{1,n}(x))dx}+O(\exp(-nC)),
\end{align*}
where $B_0,B_1\subset \bb R$ such that there exist single global minimum points $x_0^*$ and $x_1^*$ of $\mca H_{0,n}$ and $\mca H_{1,n}$ satisfying $x_0^*\in B_0$ and $x_1^*\in B_1$. The additional global minimum is dropped.
In this work, we extend their results to the RFCW measure (We notice that this is a mixture Gibbs measure and is not covered by \citep{ellis1980limit}.) and obtain that there exists $C>0$ with
\begin{align*}
    \bb E\bigg[\frac{\sum_{\bfa\sigma\in\ca C}f(\bfa\sigma)\exp(\frac{\theta_1}{2n}\sum_{i,j\in[n]}\sigma_i\sigma_j+\sum_{i=1}^nh_i\sigma_i)}{\sum_{\bfa\sigma\in\Sigma}\exp(\frac{\theta_1}{2n}\sum_{i,j\in[n]}\sigma_i\sigma_j+\sum_{i=1}^nh_i\sigma_i)}\bigg]&=\bb E\bigg[\frac{\int_{B_0}\exp(-n\mca H_{0,n}(x))dx}{\int_{B_1}\exp(-n\mca H_{1,n}(x))dx}\bigg]\\
    &+O(\exp(-nC)).
\end{align*}
\smallbreak
\emph{The Local Tail Bounds.}
The upper bounds of testing with multiple samples are reduced to studying the tail bounds of test statistics. The transfer principle method leads to the convergence of the moment generating function (or mgf), for all $t\in\bb R$, pointwise,
\begin{align*}
   \bl\frac{1}{\sqrt n}\sum_{i=1}^n(\sigma_i-\mu)\bigg|\sum_{i=1}^n\sigma_i>0\br\overset{mgf}{\to} N(0,\ca V),\qquad \ca V,\mu>0.
\end{align*}
which also leads to the convergence of moments by the Lipchitzness of mgf for compact intervals of $t$ \citep{billingsley2017probability}. However, to give a tail estimate, the convergence in moment generating function is insufficient to achieve uniform control over the tail in $\bb R^+$. For example, using the standard Chernoff bound, \emph{only when $t=o(f(n))$} for some increasing function $f$ of $n$, we have for all $t\in\bb R$,
\begin{align*}
    \bb P\bl\frac{1}{\sqrt n}\sum_{i=1}^n(\sigma_i-\mu)\geq t\br\leq\inf_{\lambda}\exp\bl\frac{1}{2}\ca V\lambda^2-\lambda t\br(1+o_n(1))\leq (1+o_n(1))\exp\lef(-\frac{1}{2}\ca Vt^2\rig).
\end{align*}
We term this by the \emph{locality} of tail bounds since for the sum of i.i.d. random variables, the standard tail bound does not require this extra condition on the relationship between $t$ and $n$.

To drop the condition of $t=o(f(n))$, we apply the boundedness of $X:=\frac{1}{\sqrt n}\sum_{i=1}^n(\sigma_i-\mu)$ and use the Cauchy-Schwartz inequality to obtain a uniform upper bound for all $t=o(n^{1/2})$. Then, we use the moment computation to derive upper bounds for $\bb E[|X|^p]^{\frac{1}{p}}$ all $p\in\bb N$. This finally leads to the sub-Gaussian norm being bounded according to \citep{vershynin2018high} Proposition 2.5.2. and gives the uniform upper bound on $t\in\bb R^+$, akin to the case of i.i.d. random variables. This method does not lead to sharp control over the constant factor in the exponential term as the price paid for uniformity.

A few challenges arise from deriving lower bounds, which require additional techniques.
\smallbreak
\emph{The Fake Measure Method.}
 A standard method in proving the lower bound is Le Cam's method, which states that it is sufficient to obtain an upper bound of the TV distance between the null and alternative measures, given by (recall the notations given in definition \ref{def1}.)
 \begin{align*}
     \inf_{\psi:\{\bfa\sigma^{(i)}\}_{i\in[m]}\to\{0,1\}}\bb P_{0,m}(\psi=1)&+\sup_{S:|S|=k}\bb P_{S,m}(\psi=0)\geq 1-\frac{1}{2}\bigg\Vert\bb P_0^{\otimes m}-\int_{\ca S}\bb P_S^{\otimes m}\pi(dS)\bigg\Vert_{TV}\\
     &\geq 1-\frac{1}{2}\bb E\bigg[\bigg\Vert\bb P_0^{\otimes m}(\cdot|\bfa h)-\int_{\ca S}\bb P^{\otimes m}_S(\cdot|\bfa h)\pi(dS)\bigg\Vert_{TV}\bigg].
 \end{align*}
 with $\pi$ being an arbitrary prior on the distributions of alternative hypotheses indexed by $\ca S:=\{S:|S|=k,S\subset[n]\}$. However, it is not easy to control the TV distance for the mixture Gibbs measure, we upper bound it with the chi-square divergence between the null and a mixture of alternative hypotheses.
 Due to the phase transition phenomenon of the Gibbs measure, when $\theta_1$ is close to the inverse critical temperature $\frac{1}{\bb E[\sech^2(h)]}$ unboundedness issue arises for the chi-square divergence (discussed with equation \eqref{nonfake}). To get an intuition, we let the distribution $\bb P_S$ be the pRFCW measure with the clique planted in set $S$ and $\bb P_0$ be the non-planted measure. Then we consider the following quantity that naturally arises from the quadratic form in the chi-square divergence \emph{regardless} of the prior we choose,
\begin{align}\label{nonfake}
    \bb E\bigg[\sum_{\bfa\sigma}\frac{(\bb P_S(\bfa\sigma|\bfa h))^2}{\bb P_0(\bfa\sigma|\bfa h)}\bigg]&=\bb E\Bigg[\frac{\sum_{\bfa\sigma}\exp(\sum_{i,j\in[k]}\frac{\theta_1}{k}\sigma_i\sigma_j+\sum_{i=1}^kh_i\sigma_i)\sum_{\bfa\sigma}\exp(\sum_{i=1}^kh_i\sigma_i)}{\lef(\sum_{\bfa\sigma}\exp(\sum_{i,j\in[k]}\frac{\theta_1}{2k}\sigma_i\sigma_j+\sum_{i=1}^kh_i\sigma_i)\rig)^2}\Bigg]\nnb\\
    &=\bb E\bigg[\frac{\int_{\bb R}\exp(-k\mca H_{0,k}(x))dx\prod_{i=1}^k \cosh(h_i)}{(\int_{\bb R}\exp(-k\mca H_{1,k}(x))dx)^2}\bigg],
\end{align}
with \small$$\mca H_{0,k}(x):=-\frac{x^2}{2}+\frac{1}{k}\sum_{i=1}^k\log\cosh(\sqrt{2\theta_1}x+h_i),\quad \mca H_{1,k}(x):=-\frac{x^2}{2}+\frac{1}{k}\sum_{i=1}^k\log\cosh(\sqrt{\theta_1}x+h_i).$$ 
\normalsize
And we note that the numerator corresponds to double the inverse temperature $\theta_1$ compared with the term in the denominator. This results in the numerator arriving at the low-temperature regime while the denominator remains at the high-temperature regime.  Specifically, notice that the Laplace method yields the value of the exponential integral to roughly be the value evaluated at the sequence of stationary points of $\mca H_{1,k}$ and $\mca H_{0,k}$. When they are not converging to $0$ together, the numerator and denominator have a difference of $\exp(kC)$ that diverges. This is pessimistic and leads to trivial lower bounds of testing. Our strategy is instead considering a \emph{fake measure} in the form of $\mbbm 1_{A}\bb P_S$ where $A$ is a constrained high probability set (depends also on $S$) that screens out the configurations of $\bfa\sigma$ making \eqref{nonfake} blow up. (It is termed fake since it is not a probability measure that integrates to $1$). It is checked that interpolating this fake measure into the TV distance does not lose too much,

\begin{align*}
    &\sup_{S}\Vert \bb P^{\otimes}_S-\bb P^{\otimes m}_S\mbbm 1_{A(S)}\Vert_{TV}=o(1)\Rightarrow\quad \bigg\Vert\int_{\ca S}\bb P^{\otimes m}_S\pi(dS)-\int_{\ca S}\bb P^{\otimes m}_S\mbbm 1_{A(S)}\pi(dS)\bigg\Vert_{TV}=o(1).
\end{align*}
But, the resulting chi-square between the null and this fake measure is tighter.
\footnotesize
\begin{align*}
    \bb E\bigg[\sum_{\bfa\sigma}\frac{(\bb P_S(\bfa\sigma|\bfa h))^2\mbbm 1_{A(S)}}{\bb P_0(\bfa\sigma|\bfa h)}\bigg]=\bb E\Bigg[\frac{\sum_{\bfa\sigma\in A(S)}\exp(\sum_{i,j\in[k]}\frac{\theta_1}{k}\sigma_i\sigma_j+\sum_{i=1}^kh_i\sigma_i)\sum_{\bfa\sigma}\exp(\sum_{i=1}^kh_i\sigma_i)}{\lef(\sum_{\bfa\sigma}\exp(\sum_{i,j\in[k]}\frac{\theta_1}{2k}\sigma_i\sigma_j+\sum_{i=1}^kh_i\sigma_i)\rig)^2}\Bigg].
\end{align*}
\normalsize
We also note that this method is used with the transfer principle given in the last paragraph to arrive at the sharp estimate of the information divergence, which translates to a sharp estimate of lower bounds.

\subsection{High Temperature}\label{sect31}

In this section, we first give the fundamental limits of testing in the high-temperature regime. Then we match the obtained lower bounds with an adaptive composite test. We see that the rate of optimal sample complexity at the high-temperature regime is akin to that of the previous literature on Sparse PCA under continuous distribution. However, it is shown that the derivations of upper and lower bounds for the pRFCW model are largely different from the Gaussian case due to the discreteness.
\subsubsection{Lower Bounds}\label{321}
Our results give exact minimax lower bounds for half of the high-temperature regime when $k=o(\sqrt n)$ with $\theta_1<\frac{1}{2}\theta_c$. And when $\theta_c>\theta_1\geq\frac{1}{2}\theta_c$ we have optimality up to a $\log k$ factor. For the region of $k=\omega(\sqrt n)$ we achieve the exact order of sharpness. The missing logarithmic factor is due to the different methods of proof designated to overcome the difficulty incurred by upper bounding the information divergence between two random Gibbs measures. We show that there exists an `elbow' effect in the minimax rate. This phenomenon is universal and appears in all temperature regimes. This is the direct result of two contradicting factors: (1) The larger clique makes the test easier since it is much easier to find it. (2) The larger clique makes the correlation between vertices weaker, hence making the test harder to perform.
\begin{theorem}[Lower Bounds for High Temperature]\label{minimaxlb}
  Assume that $\tanh(h)$ is centered.  Assume $\theta_1<\frac{1}{\bb E[\sech^2(h)]}$.
   Then the region of sample complexity $m$ making all tests asymptotic powerless is given by:
    \begin{enumerate}
        \item If $k=o( n^{\frac{2}{3}} )$, $\theta_1<\frac{1}{2\bb E[\sech^2(h)]}$, and $m=o\lef(k\log n\rig)$;
        \item If $k=o( n^{\frac{2}{3}})$, $\theta_1\geq\frac{1}{2\bb E[\sech^2(h)]}$, and $m=o\lef(k\frac{\log n}{\log k}\rig)$;
        \item If $k=\Omega( n^{\frac{2}{3}})$ and $m=o\lef(\frac{n^2}{k^2}\rig)$.
    \end{enumerate}
\end{theorem}
We note that the lower bound results are much more difficult to obtain for the pRFCW model as compared with the Gaussian models due to the discreteness of variable support. Moreover, the pRFCW model is also a mixture distribution so the Neyman-Pearson uniform powerful tests that have been utilized in \citep{mukherjee2018global} can no longer be used to obtain the lower bound. To give intuitions on how to prove this lower bound, we provide a proof sketch in the appendix \ref{proofsketch3.1}. The complete proof is delayed to the supplementary materials.

\subsubsection{Upper Bounds}
To attain the fundamental limits given by the lower bound in the past section, we propose a composite test. The test consists of two parts : (1) The local parts that match the fundamental limits of small cliques. (2) The global parts that match the fundamental limits of large cliques. 

The idea of the local part in \ref{alg:two} is scanning over all possible subsets and competing for the most likely sets. Then we reject if the maximal one lies within a typical region. However, as suggested by the lower bounds, this test is not always optimal.We recall that under the high-temperature regime, the difference between null and alternative is in the order of $O(\sqrt k)$.. And, as $k$ becomes larger, the contribution of the spins within the clique to the mean magnetization significantly improves. This implies that when $k$ is large, there is the possibility of having better algorithms. This leads to our design of the global parts of the test, which is just computing the total mean correlation and rejecting if it exceeds a pre-specified threshold. We note that in the previous literature, there exist other methods that can perform testing on the small clique regime of sparse PCA like computing the $k$-sparse principle eigenvalue proposed in \citep{berthet2013optimal} and we believe they can be used interchangeably with our proposed method with more technical difficulties in the analysis. Moreover, we believe that our method is simpler to analyze when it comes to the recovery guarantees than the spectral methods in the literature since the spectral properties of the empirical covariance matrices for Ising models remain open.

\begin{algorithm}[htbp]
\caption{High Temperature Test}\label{alg:two}
\KwData{$\{\bfa\sigma^{(i)}\}_{i\in[m]}$ with $\bfa\sigma\in\{-1,1\}^n$}
\eIf{$k=o(n^{\frac{2}{3}})$}{Compute empirical correlation matrix$\wh {\bb E}[\bfa\sigma\bfa\sigma^\top]= \frac{1}{m}\sum_{i=1}^m\bfa\sigma^{(i)}\bfa\sigma^{(i)\top}$\;
Going over all subset $S\subset[n]$ with $|S|=k$.  Compute $\phi_S = \frac{1}{k}\mbbm 1_{S}^\top\wh {\bb E}[\bfa\sigma\bfa\sigma^\top]\mbbm 1_{S}$\; 
Reject Null if $\phi_1 = \sup_{S:|S|=k}\phi_S\geq \tau_\delta$ where $\tau_\delta\in\lef(0, \frac{1-\theta_1\bb E[\sech^2(h)]^2}{(1-\theta_1\bb E[\sech^2(h)])^2}\rig)$ \;}
{Compute empirical correlation $\phi_2= \frac{1}{mk}\sum_{i=1}^m\bfa\sigma^{(i)\top}\bfa\sigma^{(i)}-1$\;
Reject Null if $\phi_2>\tau_\delta$, with $\tau_\delta\in\lef(0,\frac{2\theta_1\bb E[\sech^2(h)]-\theta_1(1+\theta_1)(\bb E[\sech^2(h)])^2}{(1-\theta_1\bb E[\sech^2(h)])^2}\rig)$;}
\end{algorithm}
\begin{theorem}\label{thm1}
Assume that $\tanh(h)$ is centered.
    Assume that $0<\theta_1<\theta_c$ where $\theta_c=\frac{1}{\bb E[\sech^2(h)]}$ is the critical temperature. Then the algorithm \ref{alg:two} is asymptotically powerful when
    \begin{enumerate}
        \item $k=o(n^{2/3})$ and $m\gtrsim k\log n $;
        \item  $ n^{2/3}\lesssim k\leq n$ and $m\gtrsim \frac{n^2}{k^2}$.
    \end{enumerate}
    
\end{theorem}

Another interesting question is how good is the set $S$ achieved by scanning. Intuition tells us that the set returned by the local part of \ref{alg:two} should have a large overlap with the hidden clique. However, to prove this fact we need to understand the limiting distribution of only a part of the spins in the clique. This is more complicated than considering the clique as a whole due to the correlation structure within the set that contains the clique. The idea to prove this is first confirmed by a tilting method in the large deviation theory that appears in the supplementary material \citep{he2023supp} and then rigorously proved using the Laplace method. Our final result is given by the following corollary, which confirms that the set returned by the local part of \ref{alg:two} guarantees almost exact recovery.
\begin{corollary}\label{hightemperaturerec}
  Assume that $\tanh(h)$ is centered. For arbitrary $\delta>0$, there exists $C>0$ such that with sample size $m\geq Ck\log n$ and the same regular conditions in theorem \ref{thm1}, the following holds for the set targeted by the local part of algorithm \ref{alg:two}, $S_{\max}:=\argmax_{S:|S|=k}\phi_S$ :
$         \bb P\lef(|S\Delta S_{max}|\geq\delta k\rig)=o(1).
$\end{corollary}

\subsection{Low Temperature}\label{sect32}

This section presents the result for the low-temperature regime. We show that the clique is much easier testable for the low-temperature regime than that for the high-temperature regime.  Intuitively, the low-temperature regime is characterized by the correlation between each pair of spins in the clique being so strong that they are much more likely to be positive or negative together.   For the special case of $\tanh h$ having a symmetric distribution with respect to $0$, this induces a unique phenomenon where the concentration of the spins is on two symmetric points around $0$ instead of one at the high temperature. For the general asymmetric $\tanh(h)$ but centered case the results are simpler where only one point out of the two local minima is global. Hence it is enough to discuss the symmetric case since the asymmetric but centered case deserves no special treatment. This phenomenon is also referred to as `spontaneous magnetization' in physics where we observe much larger clusters consisting of all positive/negative spins in the Ferromagnetic material. Under the low-temperature regime, the hidden clique problem shares fundamental similarities with the planted clique problem, where the Erdos-Renyi graph generates a $O(\log n)$ clique almost surely. Therefore, finding an $\omega(\log n)$ clique under the null is almost impossible. Similar to the high-temperature regime, we first present the lower bounds of this problem, then followed by the algorithm achieving the limit.

\subsubsection{Lower Bounds}
Our results for the lower bound cover the region of $o(\sqrt n)$ only, due to the intuition that the larger ones can be either tested with constant sample size or even a single sample.

\begin{theorem}\label{lbltsmallclique}
   Assume that $\tanh(h)$ is centered.  If $\theta_1>\frac{1}{\bb E[\sech^2(h^\prime)]}$ , and $k=o(\sqrt n)$ then all tests are asymptotically powerless if $m\leq C\log\frac{n}{k}$ for some $C>0$.
\end{theorem}
\subsubsection{Upper Bounds}
To attain the fundamental limits given above we present test \ref{alg:three}. For the global part of the test, we further divide it into $k\asymp n$ and $k=\omega(k)$ two cases, where the rejection regions are different.  However, instead of computing the variances, we compute the absolute values, which gives better convergence guarantees and attain the lower bounds at both the local and the global parts. Interestingly, instead of requiring $\frac{n}{k}$ samples during the high temperature, at low temperature we only need $1$ sample to test the large clique of order $\omega(\sqrt n)$.
The underlying intuition is that fluctuation in the null model is in the order of $\frac{1}{\sqrt n}$, which is surpassed by the simultaneous magnetization in the clique that is in the order of $\frac{k}{n}$. This $1$ sample testability never appears in the Gaussian Sparse PCA problem.
\begin{algorithm}[htpb]
\caption{Low Temperature Test}\label{alg:three}
\KwData{$\{\bfa\sigma^{(i)}\}_{i\in[m]}$ with $\bfa\sigma\in\{-1,1\}^n$}
\uIf{$k=o(\sqrt n)$}{Going over all subset $S\subset[n]$ with $|S|=k$.  Compute $\phi_S = \frac{1}{m}\sum_{j=1}^m\lef|\frac{1}{k}\sum_{i\in S}\sigma^{(j)}_i\rig|$\; 
Reject Null if $\phi_3 = \sup_{S:|S|=k}\phi_S\geq \tau_\delta$ with $\tau_\delta\in\lef(0, x\rig)$ with $x$ defined by the positive solution to $x=\bb E[\tanh(\theta_1x+h)] $ \;}
\Else{Compute statistics $\phi_4:=\frac{1}{m}\sum_{j=1}^m\lef|\frac{1}{k}\sum_{i=1}^n\sigma^{(j)}_i\rig|$\;
\uIf{$k\asymp n$}{
Reject Null if $\phi_4 > \tau_\delta$  for $\tau_\delta\in\lef(\sqrt{\frac{2n}{\pi k^2}},\frac{\sqrt n}{k}\sqrt{\frac{2}{\pi}}\exp\lef(-\frac{x^2k^2}{2n}\rig)+x\lef[1-2\Phi\lef(-\frac{xk}{\sqrt n}\rig)\rig]\rig)$ with $\Phi$ being the cumulative distribution function of the standard Gaussian random variable\;}
\Else{Reject Null if $\phi_4>\tau_\delta$ for $\tau_\delta\in (0,x)$ with $m$ defined by the positive solution to with $x$ defined by the positive solution to $x=\bb E[\tanh(\theta_1x+h)] $\;}
}
\end{algorithm}[htbp]
\begin{theorem}\label{thmlowtempfindclique}
   Assume that $\tanh(h)$ is centered and $\theta_1>\frac{1}{\bb E[\sech^2(h)]}$. Then the test given by algorithm \ref{alg:three} is asymptotically powerful when
    \begin{enumerate}
        \item $k=o(\sqrt n)$ and $m\gtrsim \log\frac{n}{k}$.
        \item $k\asymp n^{1/2}$ and $m=\omega(1)$.
        \item $k=\omega(n)$ and $m=1$.
    \end{enumerate}
\end{theorem}
The following result is the analogous, almost exact recovery guarantee for the $o(k)$ clique in the low-temperature regime, where the algorithm outputs the clique selected in algorithm \ref{alg:three}. 
\begin{corollary}\label{recoveryguaranteelowtemp}
    Assume that $\tanh(h)$ is centered. For arbitrary $\epsilon>0$, sample size $m\gtrsim k^{1/(2\tau-1)}\log\frac{n}{k}$, consider $S_{\max}\in\argmax_{S:|S|=k}\phi_{S}$ returned by the local part of algorithm \ref{alg:three}, we have for all $\delta>0$,
    \begin{align*}
         \bb P\lef(|S\Delta S_{\max}|\geq\delta k\rig)=o(1).
    \end{align*}
\end{corollary}

Compared with the results in the high-temperature regime, it is not hard to see that at the low temperature, the sample complexity is much smaller. We believe this phenomenon is also universal in a large class of discrete statistical models, where after the signal (here the parameter $\theta_1$ represents the magnitude of correlation) surpasses a certain constant threshold, the problem becomes significantly easier.

\subsection{Critical Temperature}\label{sect33}
At the critical temperature, the statistical behavior is more delicate, where the rate relies on the tail property of $h$. The critical temperature is characterized by having a flat landscape around zero with the curvature or the second order derivative $0$. In this regime, the behavior at the critical temperature is decided by the following function 
\begin{align}\label{taylorcond}
    H(x)=\frac{1}{2}x^2-\bb E[\log\cosh(\sqrt{\theta_1}x+h)].
\end{align}

\begin{definition}[Flatness of Local Optimum]\label{flatnessdefine}
    We call the local minimum and maximum $x^*$ of \eqref{taylorcond} is $\tau$-flat for $\tau\in\bb N\setminus\{1\}$ if $ H(x)=\frac{H^{(2\tau)}(x^*)}{(2\tau)!}(x-x^*)^{2\tau-1}+O((x-x^*)^{2\tau})$ with $H^{(2\tau)}(x_1^*)>0$ and $H^{(2\tau)}(x_1^*)<0$, respectively. 
\end{definition}

Then, depending on the positivity of $H^{(2\tau)}$ we might have the following two different cases:
\begin{enumerate}
    \item $H^{(2\tau)}(0)>0$: We encounter the intermediate phase of critical temperature.
    \item $H^{(2\tau)}(0)<0$: We skip the intermediate phase of critical temperature and directly go to the low-temperature phase.
\end{enumerate}
In this section, we focus on the first case since the second one can be treated analogously as the low-temperature regime.
Then we let $\tau$ be the order of flatness at the critical temperature. It is checked that this quantity depends only on the moments of $h$. This phenomenon is also quite interesting if compared with the zero field Curie-Weiss model, where \citep{ellis1978limit} establishes that its convergence rate of the average magnetization at the critical temperature is ${n^{-1/4}}$. The random effect in the field then yields much richer statistical behavior than the zero field model. We show that, under the critical temperature, the statistical rate will demonstrate a mixture effect of both the high-temperature phases and the low-temperature phases, causing 2 transitions to happen as $k/n$ gets larger. Moreover, we show that both the statistical minimax rates and the thresholds of transitions can take countable different values according to $\tau$. These unique phenomena were never observed in any other statistical models in the past, which implies the random Ising models can have much richer statistical landscapes than the long-established Gaussian models.

\subsubsection{Lower Bounds}
From the fundamental limits' perspective, we show that the statistical optimal rate experiences two phase transitions happening at $k\asymp n^{\frac{4\tau-2}{8\tau-5}}$ and $k\asymp n^{\frac{2\tau-1}{4\tau-3}}$ respectively. Moreover, we observe that this `double-elbow' effect splices the large clique regimes of the high and low temperatures together.
\begin{theorem}\label{minimaxlbct}
Assume that $\tanh(h)$ is centered, $\theta_1=\frac{1}{\bb E[\sech^2(h)]}$, and the flatness parameter of the global minimum of \eqref{taylorcond} is $\tau$.
    Then the region of sample complexity $m$ making all tests asymptotic powerless is given by:
    \begin{enumerate}
        \item $k=o( n^{\frac{4\tau-2}{8\tau-5}} )$ and $m=o\lef(\lef(\frac{k}{\log k}\rig)^{\frac{1}{2\tau-1}} \log n\rig)$;
        \item $n^{\frac{4\tau-2}{8\tau-5}}\lesssim k\lesssim n^{\frac{2\tau-1}{4\tau-3}}$ and $m=o\lef(n^2k^{-\frac{2(4\tau-3)}{2\tau-1}}\rig)$;
        \item $k=\omega\lef(n^{\frac{2\tau-1}{4\tau-3}}\rig)$ and $m\leq C$ for some $C\geq 1$.
    \end{enumerate}
\end{theorem}

\subsubsection{Upper Bounds}
We use test \ref{alg:five} to attain the lower bounds. And we notice that the global part of the test attains the minimax optimality on two of the large clique phases, which is separated by $k\asymp n^{\frac{2\tau-1}{4\tau-3}}$. Before the threshold, we need an infinite number of samples, whereas after the threshold, it is possible to perform successful testing with constant samples. This appears to be an intermediate state between the high and low-temperature regimes since the pre-transition phase is analogous to the large clique phase of the high temperature and the post-transition phase to the large clique phase of low temperature. To understand this is rather straightforward if we pick $\tau$ to be $1$ and $\infty$, we will recover the high and low temperature regimes, respectively.

\begin{algorithm}[htbp]
\caption{Critical Temperature Test }\label{alg:five}
\KwData{$\{\bfa\sigma^{(i)}\}_{i\in[m]}$ with $\bfa\sigma\in\{-1,1\}^n$}
\uIf{$k=o\lef(n^{\frac{4\tau-2}{8\tau-5}}\rig)$}{Compute Scaled Empirical Correlation $\wh{\bb E}[\bfa\sigma\bfa\sigma^\top]=\frac{1}{m}\sum_{j=1}^m\bfa\sigma^{(j)}\bfa\sigma^{(j)\top}$\;
Go over all subset $S\subset[n]$ with $|S|=k$ and  compute $\phi_{S} = k^{-(4\tau-3)/(2\tau-1)}\lef(\mbbm 1_{S}^\top\wh {\bb E}[\bfa\sigma\bfa\sigma^\top]\mbbm 1_{S}\rig)$\; 
Reject Null if $\phi_{5} = \sup_{S:|S|=k}\phi_S\geq \tau_\delta$  for $\tau_\delta\in\lef(0,\pi^{-\frac{1}{2}}(2\ca V(\tau))^{\frac{1}{2\tau-1}}\Gamma(\frac{2\tau+1}{4\tau-2})\rig)$ where $\ca V(\tau):=\frac{((2\tau)!)^2\bb V(\tanh(h))(\bb E[\sech^2(h)])^{4\tau-2}}{2^{2\tau-1}\lef(\bb E\lef[(1+\tanh(h))\sum_{k=0}^{2\tau-1}\frac{k!}{2^k}S(2\tau-1,k)(\tanh(h)-1)^{k}\rig]\rig)^2}$ and $S(n,k)$ is the second type of Stirling numbers\;}
\Else{Compute Scaled Correlation $\phi_6=m^{-1}k^{-\frac{4\tau-3}{2\tau-1}}\sum_{j=1}^m\lef(\lef(\sum_{i=1}^n\sigma^{(j)}_i\rig)^2-n\rig)$\;
Reject Null if $\phi_6\geq\tau_\delta$ for $\tau_\delta\in\lef(0,\pi^{-\frac{1}{2}}(2\ca V(\tau))^{\frac{1}{2\tau-1}}\Gamma\lef(\frac{2\tau+1}{4\tau-2}\rig)\rig)$ \;}

\end{algorithm}

\begin{theorem}\label{guaranteecriticaltest}
    Assume that $\tanh(h)$ is centered and $\theta_1=\frac{1}{\bb E[\sech^2(h)]}$. Assume that $0$ is a global minimum of \eqref{taylorcond} with flatness $\tau$. Then the test \ref{alg:five} is asymptotically powerful when
    \begin{enumerate}
        \item $k=o\lef(n^{\frac{4\tau-2}{8\tau-5}}\rig)$ and $m=\omega\lef(k^{\frac{1}{2\tau-1}}\log n\rig)$;
        \item $n^{\frac{4\tau-2}{8\tau-5}}\lesssim k\lesssim n^{\frac{2\tau-1}{4\tau-3}}$ and $m=\omega\lef(n^2k^{-\frac{2(4\tau-3)}{2\tau-1}}\rig)$;
        \item $k=\omega\lef( n^{\frac{2\tau-1}{4\tau-3}}\rig)$ and $m=\omega(1)$.
    \end{enumerate}
\end{theorem}
And similar to the high/low temperature we can derive the almost exact recovery guarantee for algorithm \ref{alg:five} given by the following corollary.
\begin{corollary}\label{recoverycrtical}
    Assume the hidden clique is indexed by $S$ at the critical temperature. For arbitrary $\epsilon>0$, sample size $m\gtrsim k^{1/(2\tau-1)}\log n$, consider $S_{\max}\in\argmax_{S:|S|=k}\phi_{S}$ returned by the local part of algorithm \ref{alg:five}, we  have for all $\delta>0$,
$         \bb P\lef(|S\Delta S_{\max}|\geq\delta k\rig)=o(1).
$\end{corollary}

\subsection{The Exact Recovery}\label{exactreco}
In the above three sections, we discuss the almost exact recovery guarantees (as defined in \ref{recoguarantee}) given by the local parts of testing algorithms under the high/critical and low-temperature regimes. Here, we prove that the almost exact recovery can be boosted to exact recovery and also provide matching lower bounds to corroborate their optimality guarantees. Similar to the testing, we show that the exact recovery rate also has significant differences across temperature regimes. At the high temperature, this rate is $k^2\log n$, the same as \citep{amini2009high}. However, the critical and low-temperature regimes deviate significantly from the standard rate of the Gaussian Sparse PCA problem. On the other hand, we show that the derivation of the recovery lower bounds is also more difficult than the testing lower bounds.
\subsubsection{Lower Bounds}
This subsection discusses the statistical barriers to the exact recovery. At the center of the proof is the Fano's inequality. In our case, considering a random variable $S$ uniformly taking values in the set $$\ca S:=\begin{cases}
    \{S:S_{[k-1]}=[k-1], S_k\in[k:n]\}&\text{ when }k\leq (1-\delta)n\text{ for small }\delta>0,\\
    \{S^c:S_{[1,n-k-1]}=[k+1:n], S_{n-k}\in[k]\}&\text{otherwise}.
\end{cases}$$ Essentially, we show that the above prior construction gives the precise order of $\log(n)$ rather than the order of $\log(n-k)$ appears in \citep{amini2009high} for Gaussian SPCA, which is obtained using only the first part of priors. A careful observation of the proof in \citep{amini2009high} suggests that their results can also be improved for vanishing $n-k$ using our strategy of the second part of priors.
We define the measure induced by $S$ as $\mu_S$. Let $\wh S$ as an estimator based on $\{\bfa\sigma^{(i)}\}_{i\in[m]}$. Define the distribution $\bb P_{\bar S}(\bfa\sigma):=\frac{1}{|\ca S|}\sum_{S\in\ca S}\bb P(\bfa\sigma|S)$. Let $I(P,Q)$ be the mutual information between measures $P$ and $Q$, then we  have by the information processing inequality,
\begin{align*}
    \bb P(\wh S\neq S)\geq 1-\frac{I(S;\wh S)+\log 2}{\log |\ca S|}\geq 1-\frac{I(\mu_S^{\otimes m};\bb P_{\bar S}^{\otimes m})+\log 2}{\log |\ca S|}.
\end{align*}
The difficulty in this proof is the estimation of mutual information. In the standard setup, this is often estimated using the maximum entropy distributions like the Gaussians. However, this tool does not work for the pRFCW measures considered in this work since it is both a mixture and is a discrete distribution. Bounding this divergence is an algebraic problem involving the computation of complex asymptotic integrals and the corresponding higher order fixed point analysis illustrated in \eqref{fixedpointanalysis}. Moreover, here we encounter a path integral on the complex plane. This results in the failure of the former method of $Z$-estimator analysis used in the derivation of testing lower bounds. Instead, we analyze the complex `stationary points'. We have to use two technical steps to obtain the proper lower bounds: (1) A distortion step on the integral path in $\bb C$ as required by the method of the steepest descent. (2) A projection step to simplify the algebraic manipulations of the complex fixed point analysis. Further details are delayed to the supplementary materials.
\begin{theorem}\label{lowerboundforrecovery}
   Assume that $\tanh(h)$ is centered. Let the set consisting of all possible priors on the position of clique set $S$ be $\ca P$. Let $\wh S$ be the estimate of $S$ given the $m$ sets of i.i.d. sampled data $\{\bfa\sigma^{(i)}\}_{i\in[m]}$, then according to the different temperature regimes, we  have:
    \begin{enumerate}
        \item At the high temperature regime, $\inf_{\wh S}\sup_{\bb P\in \ca P}\bb P(S\neq\wh S)\geq 1-O\lef(\frac{m}{k\log (n)}\rig)$;
        \item At the low temperature regime, $\inf_{\wh S}\sup_{\bb P\in\ca P}\bb P(S\neq \wh S)\geq 1-O\lef(\frac{m}{\log (n)}\rig)$;
        \item At the critical temperature regime, $\inf_{\wh S}\sup_{\bb P\in\ca P}\bb P(S\neq \wh S)\geq 1-O\lef(\frac{m}{k^{\frac{1}{2\tau-1}}\log (n)}\rig)$.
    \end{enumerate}
\end{theorem}
\subsubsection{Upper Bounds}
  To improve this result to the exact recovery guarantee, we can apply the following set screening procedure in algorithm \ref{alg:screen} to do the model selection. The underlying intuition is to use the almost exact support to denoise, since if we reduce the total entries in the sum statistics from $n$ to $k$, the fluctuation given by the spins in set $S^c$ is significantly reduced. This helps us to screen out the noisy part within $S$ and achieve better guarantees.

\begin{algorithm}[htbp]
\caption{Set Screening}\label{alg:screen}
\KwData{$\{\bfa\sigma^{(i)}\}_{i\in[m]}$ such that $\bfa\sigma^{(i)}\in\{-1,1\}^n$, an almost exact solution $S^\prime$ returned by algorithm \ref{alg:two}, \ref{alg:three}, or algorithm \ref{alg:five} in the high/low/critical temperature regimes respectively.}
Compute the statistics $\phi_i=\begin{cases}
     m^{-1}\sum_{\ell=1}^m\sum_{j\in S^\prime,j\neq i}\sigma^{(\ell)}_i\sigma^{(\ell)}_j&\text{ at the high temperature regime}\\
     k^{-\frac{2\tau-2}{2\tau-1}}m^{-1}\sum_{\ell=1}^m\sum_{j\in S^\prime,j\neq i}\sigma^{(\ell)}_i\sigma^{(\ell)}_j &\text{ at the critical temperature regime with flatness } \tau\\
     k^{-1}m^{-1}\sum_{\ell=1}^m\lef|\sum_{j\in S^\prime,j\neq i}\sigma^{(\ell)}_i\sigma^{(\ell)}_j\rig|&\text{ at the low temperature }
\end{cases}$\;
Then we rank $\phi_i$ and pick $S^{\prime\prime}$ by the set achieving top $k$ values of $\phi_i$.
\end{algorithm}
Then, we  prove theorem \ref{asscreening}, which shows that the screening procedure gives exact recovery guarantees.
\begin{theorem}\label{asscreening}
    Assume that $\tanh(h)$ is centered. The set $S^{\prime\prime}$ returned by algorithm \ref{alg:screen} satisfy $\bb P(S^{\prime\prime}= S)=1-o(1)$ when
    \begin{enumerate}
        \item $m=\omega(k\log (n))$ at the high temperature regime;
        \item $m=\omega( k^{\frac{1}{2\tau-1}}\log (n))$ at critical temperature regime with flatness $\tau$;
        \item $m=\omega( \log (n))$ at the low temperature.
    \end{enumerate}
\end{theorem}

\section{Non-centered Random Field}\label{sect4}

This section discusses the result when $\tan(h)$ is non-centered in the sense that $\bb E[\tanh(h)]\neq 0$. Our results are summarized in table \ref{table2}. 
\begin{table}[htp!]
\caption{\textbf{The Minimax Sample Complexity of Powerful Testing when $\tanh (h)$ is Non-centered} The agnostic case refers to constructing test statistics using no information of $h$, and the oracle refers to constructing test statistics using moment information of $h$. We note from the table that the agnostic case does not match with the oracle case when $\sqrt n\lesssim k\lesssim n^{\frac{2\tau-1}{4\tau-3}}$. It remains an open problem if there exists an agonistic test matching the lower bound or if there exists a certificate for the non-existence.}
\renewcommand*{\arraystretch}{2.0}
\centering
\begin{tabular}{ll|l|l|l}
\toprule
\multicolumn{2}{l|}{\textbf{Non-centered} $\tanh (h)$ Testing}                        & $k=o(\sqrt n)$             & $\sqrt n\lesssim k\lesssim n^{\frac{2\tau-1}{4\tau-3}}$ & $k=\omega(n^{\frac{2\tau-1}{4\tau-3}})$ \\ \midrule
\multicolumn{1}{l|}{\multirow{2}{*}{Upper Bounds}} & Oracle & $O(\log n)$      & $O(1)$                                        & $O(1)$                                  \\ \cline{2-5} 
\multicolumn{1}{l|}{}                              & Agnostic     & $O(\log n)$      & $O(\log n)$                         & $O(1)$                                  \\ \hline
\multicolumn{2}{l|}{Lower Bounds}                                 & $\Omega(\log n)$ & $\Omega(1)$                                   & $\Omega(1)$                             \\ \bottomrule
\end{tabular}
\label{table2}
\end{table}

A difference of non-centered $\tanh (h)$ is that we do not separate the phases according to temperature regimes. Here, the point of measure concentration for the mean magnetization is determined by 
\begin{align}\label{globalminimumcf}
   m^*:=\argmin_{m\in[-1,1]} \frac{\theta_1m^2}{2}-\bb E[\log\cosh(\theta_1m+h)].
\end{align}
It is checked that the above equation only has a single global optimal point and up to three stationary points that depend on the value of $\theta_1$. Specifically, this implies that the low-temperature phase in section \ref{sect3}  disappears, and the two symmetric global maxima reduce to one local maximum and one global maximum. Therefore, instead of using the separation of temperatures in the symmetric case, the test algorithm is based on the same parameter of \emph{flatness} of the global minimum in definition \ref{flatnessdefine}.

A key difference between the non-centered and centered field problem is the information on $h$. We consider the following two cases representing the amount of information we obtained from the distribution of $h$.
\begin{enumerate}
    \item \textbf{Agnostic Case:} The situation when we have zero knowledge on $h$.
    \item \textbf{Oracle Case:} The situation when we have oracle information about $h$, in particular, the mean of $\tanh(h)$.
\end{enumerate}

\subsection{Lower Bounds}\label{sect42}
To give statistical barriers for the tests, we give the following minimax lower bounds. Comparing the rate with the centered case, it observed that the barrier looks more like the low-temperature regime than the high/critical temperature regimes. This is due to the mean shift of the magnetization contributed by the non-centered $h$. From the physics point of view, this corresponds to an outer magnetic field forcing all the spins toward a certain pole, making the ferromagnets polarized.

\begin{theorem}\label{lbnoncentered}
    Assume that the flatness of the global minimum of \ref{taylorcond}, denoted by $x^*$ is $\tau$. Then the region of sample complexity $m$ making all tests asymptotic powerless is given by:
    \begin{enumerate}
        \item If $k=o(\sqrt n)$ and $m=o\lef(\log n\rig)$;
        \item If $k=\Omega(\sqrt n)$ and $m\leq C$ for some $C\geq 1$.
    \end{enumerate}
\end{theorem}

\subsection{Upper Bounds}\label{sect41}
We present the test \ref{alg:seven} to account for the agnostic case and \ref{alg:nine} for the oracle case. The test \ref{alg:seven} gives matching upper bounds for the $k=o(n^{\frac{1}{2}})$ and $k=\omega(n^{\frac{2\tau-1}{4\tau-3}})$ with the lower bounds given. The region of $\sqrt n\lesssim k\lesssim k^{\frac{2\tau-1}{4\tau-3}}$ misses a $\log n$ factor from the lower bound. However, this is further filled by test \ref{alg:nine} under the oracle condition. One intuition for this discrepancy between the oracle and agnostic rate comes from the sample complexity for estimating $\bb E[\tanh(h)]$. 
\smallbreak
\emph{Agnostic Test.} To develop the agnostic tests, we apply a self-comparison strategy. The test is also composed of two parts, corresponding to the local and global parts that appear in the non-centered case. For the local part, we utilize the \emph{Racing Against the Crowd} strategy. The intuition of this strategy is to use the gap between the local mean for a set of spins containing part of the clique with the population mean of all the spins. This procedure is particularly robust since it accounts for any order of flatness of $x^*$. For the global part, we propose the test based on the idea of symmetrization with independent copies. This strategy has the advantage of highlighting the higher-order moment differences when we have no access to the means. This symmetrization idea is also common in the theory of empirical processes.

\begin{algorithm}[htbp]
\caption{Agnostic Non-centered Test}\label{alg:seven}
\KwData{$\{\bfa\sigma^{(i)}\}_{i\in[m]}$ with $\bfa\sigma\in\{-1,1\}^n$}
\uIf{$k=o\lef(n^{\frac{2\tau-1}{4\tau-3}}\rig)$}{Compute the empirical mean for all subsets $S\subset[n]$ with $|S|=k$ and $\phi_S= \frac{1}{mk}\sum_{i=1}^m\mbbm 1_{S}\bfa\sigma^{(i)}$\;
Compute the average of all spins $\xi=\frac{1}{mn}\sum_{j=1}^{m}\sum_{i=1}^n\sigma_i^{(j)}$\;
Reject $H_0$ if $\phi_7^{\max} = \sup_{S:|S|=k}\phi_S> \xi+\delta$ or $\phi_7^{\min}=\inf_{S:|S|=k}\phi_S<\xi-\delta$ for some small constant $\delta$\;}
\Else{Compute statistics $\phi_8:=m^{-1}k^{-\frac{4\tau-3}{2\tau-1}}\sum_{j=1}^m\lef(\sum_{i=1}^n(\sigma_i^{(2j-1)}-\sigma_i^{(2j)})\rig)^2$\;
Reject $H_0$ if $\phi_8\geq\tau_\delta$ for some small constant $\tau_\delta>0$\;}
\end{algorithm}
\begin{theorem}\label{thm4.2}
Consider $\tanh(h)$ be non-centered and the flatness parameter is $\tau$ in \eqref{globalminimumcf}. Then the test given by algorithm \ref{alg:seven} is asymptotically powerful when
\begin{enumerate}
    \item $k=O\lef(n^{\frac{2\tau-1}{4\tau-3}}\rig)$ and $m\gtrsim\log n$;
    \item $k=\omega\lef(n^{\frac{2\tau-1}{4\tau-3}}\rig)$ and $m=\omega(1)$.
\end{enumerate}\end{theorem}
We remark that an interesting phenomenon arises when comparing the guarantee of algorithm \ref{alg:seven} with the non-centered case. Instead of having a centered $\tanh(h)$, where we need different algorithms to achieve the optimality for different temperature regions, here our procedure is adaptive to the flatness of $x^*$.
\smallbreak
\emph{Oracle Test.}
Our oracle test is given by \ref{alg:nine}, which uses explicitly the information on $\bb E[\tanh(h)]$. This matches the information-theoretic lower bound but    induce
\begin{algorithm}[htbp]
\caption{Oracle Non-centered Global Test }\label{alg:nine}
\KwData{$\{\bfa\sigma^{(i)}\}_{i\in[m]}$ such that $\bfa\sigma^{(i)}\in\{-1,1\}^n$}
Compute Statistics $$\phi_9:=m^{-1}k^{-2}\sum_{j=1}^m\bl\bl\sum_{i=1}^n\sigma_i^{(j)}-n\bb E[\tanh(h)]\br^2-n(1-\bb E[\tanh(h)]^2))\br$$
Reject $H_0$ if $\phi_9\geq\tau_\delta$ for some small constant $\tau_\delta>0$\;
\end{algorithm}
\begin{theorem}\label{oracletestguarantee}
   When $k=\Omega(\sqrt n)$, algorithm \ref{alg:nine} is asymptotically powerful given $m=\omega(1)$.
\end{theorem}
\begin{remark}
Being agnostic, algorithm \ref{alg:seven} works in a different setting from the lower bounds. The authors conjecture that no agnostic algorithm can match the oracle case between $\sqrt n\lesssim k\lesssim n^{\frac{2\tau-1}{4\tau-3}}$. However, to rigorously verify this claim, it is necessary to use a new lower bound framework for the class of agnostic tests that accounts for $h$. This is out of the scope of the present paper and might be an interesting future working direction.
\end{remark}
\smallbreak
\emph{Exact Recovery.} The exact recovery results of the non-centered $\tanh(h)$ can be analogously derived using the technique given by section \ref{exactreco}, which is stated as follows.
\begin{theorem}
     Assume that $\tanh(h)$ is non-centered. Let the set consisting of all possible priors on the position of clique $S$ be $\ca P$. Let $\wh S$ be the estimate of $S$ given the $m$ sets of i.i.d. samples $\{\bfa\sigma^{(i)}\}_{i\in[m]}$, then we have
     \begin{enumerate}
         \item $m=\omega(\log n)$ and $\inf_{\wh S}\bb P(\wh S=S)=1-o(1)$;
         \item  $m=o(\log n)$ and $\inf_{\wh S}\sup_{\bb P\in\ca P}\bb P(S\neq\wh S)=1-o(1)$.
     \end{enumerate}
\end{theorem}
\begin{remark}
We show that the non-centered random field has the same statistical minimax rate as the low-temperature regimes, no matter how flat the global maximum of \eqref{taylorcond}. However, the unique $1$ sample test
\end{remark}

\section{The Limiting Theorems}\label{sect5}

In this section, we give proof of the limiting theorems for the RFCW model. All of the results presented in the previous sections are based on this theorem, which characterizes the mgf convergence of average magnetization. Compared with Gaussians, the concentration arguments are much more complicated, where multiple modes of sub-Weilbull distributions coexist with various convergence rates. This is also the reason why discrete SPCA demonstrates much richer statistical rates. In the previous literature, \citep{chatterjee2010applications} develop a novel changeable pair method to prove similar results for the classical zero field Ising model. However, their method is not directly generalizable to mixture measures, which is the unique characteristic of random Gibbs measures like the RFCW.

To prove it, we provide a generalization of the transfer principle in \citep{ellis1980limit} and the Laplace integral approximation (see \citep{bolthausen1986laplace} ) to the random measures. Despite previous work \citep{amaro1991fluctuations} attempts to derive the weak limit of the same quantity, their limiting variance is not consistent with the zero-field result that appears in \citep{ellis1980limit}. Here we not only prove stronger results of m.g.f.s convergence and simultaneously give tail bound control, but also correct their variance results.

Our results demonstrated a unique behavior at the critical temperature compared with the results for the Curie-Weiss model with zero field. Specifically, the convergence rates of the average magnetization can take countable values according to the flatness parameters $\tau$. This phenomenon directly results in the rich statistical minimax rate of this work. Another interesting phenomenon compared with the standard zero field Curie-Weiss model is that the convergence rate is always slower than the $n^{-1/4}$ rate, given $h$ having constant variance. This phenomenon implies that the fluctuations given by $h$ dominate the `latent' fluctuations of the zero field CW models, which is also observed in \citep{amaro1991fluctuations}'s weak convergence results.

\begin{theorem}[Limiting Theorem for the Random Field Curie-Weiss Model with Symmetric $h$]\label{cltrfcr}
Assume that $h_i\sim \mu(h)$ is i.i.d. in $L_1$.
For a random field Curie-Weiss model whoseHamiltonian is  defined by \eqref{rfcw-hamiltonian},
\begin{enumerate}
   \item 
In the high temperature regime with $\theta_1<\frac{1}{\bb E[\sech^2(h)]}$,   for $t\in\bb R$ pointwise,  
\small
\begin{align*}
    \bb E\bigg[\exp\bl n^{-1/2}{t\sum_{i=1}^n\sigma_i}\br\bigg]\to \exp\lef({\ca Vt^2}/2\rig)\text{ and }\bigg\Vert n^{-1/2}{\sum_{i=1}^n\sigma_i}\bigg\Vert_{\psi_2}<\infty.
\end{align*}
\normalsize
with $\ca V:=\frac{1-\theta_1(\bb E[\sech^2(h)])^2}{(1-\theta_1\bb E[\sech^2(h)])^2}$.

 \item 
In the low temperature regime of $\theta_1>\frac{1}{\bb E[\sech^2(h)]}$, $x=\bb E[\tanh(\sqrt{\theta_1} m +h)]$ has two nonzero solutions defined by $x_1<0<x_2$. Define $\ca C_1=(0,\infty)$ and $\ca C_2=\ca C_1^c$. Then   we have for $t\in\bb R$ pointwise,  for $\ell\in\{1,2\}$
\small
\begin{align}\label{gtr0clt}
    &\bb E\lef[\exp\lef(t\frac{\sum^n_{i=1}(\sigma_i -\sqrt{\theta_1}x_\ell)}{\sqrt n}\rig)\bigg |{\frac{\sum_{i=1}^n\sigma_i}{n}\in \ca C_\ell}\rig]{\to} \exp\lef(\frac{\ca V(m_1)t^2}{2}\rig),
\end{align}
\normalsize
and
$
   \lef\Vert n^{-1/2}\sum_{i=1}^n(\sigma_i-\sqrt{\theta_1}x_{\ell})\bigg|n^{-1}\sum_{i=1}^n\sigma_i\in\ca C_{\ell}\rig\Vert_{\psi_2}<\infty,
$ with\\ $\ca V(x):=\frac{(1-\theta_1(\bb E[\sech^2(\sqrt{\theta_1}x+h)])^2-\bb E[\tanh(\sqrt{\theta_1}x+h)]^2)}{(1-\theta_1\bb E[\sech^2(\sqrt{\theta_1}x+h)])^2}$.
 \item 
 At the critical temperature $\theta_1=\frac{1}{\bb E[\sech^2(h)]}$, assume that the critical value defined by \ref{flatnessdefine} is $\tau$, then for $t\in\bb R$ pointwise,  
 \small
\begin{align}\label{supergaussian}
        \bb E\lef[\exp\lef(\frac{t\sum_{i=1}^n\sigma_i}{n^{\frac{4\tau-3}{4\tau-2}}}\rig)\rig]\to \int_{\bb R}\frac{(2\tau-1)x^{2\tau-2}}{\sqrt{2\pi v(0)}}\exp\lef(-\frac{x^{4\tau-2}}{2v(0)}+tx\rig)dx,
    \end{align}
    \normalsize
    and
$        \lef\Vert n^{-\frac{4\tau-3}{4\tau-2}}\sum_{i=1}^n\sigma_i\rig\Vert_{\psi_{4\tau-2}}<\infty,
$ with \small\begin{align}\label{stirlingv}
    v(x):=&((2\tau)!)^2\bb V(\tanh(\sqrt{\theta_1}x+h))(\bb E[\sech^2(\sqrt{\theta_1}x+h)])^{4\tau-2}\nnb\\
    &\cdot\bb E\bigg[(1+\tanh(\sqrt{\theta_1}x+h))\sum_{k=0}^{2\tau-1}\frac{k!}{2^k}S(2\tau-1,k)(\tanh(\sqrt{\theta_1}x+h)-1)^{k}\bigg]^{-2}.
\end{align}\normalsize
     and $S$ is the second type of Stirling number. And if we are at the second case of \eqref{taylorcond},   \eqref{gtr0clt} holds.
\end{enumerate}
\end{theorem}
Then we can also get the following corollary, which gives a limit theorem for the mean magnetization when $h$ is asymmetric.
\begin{corollary}\label{corocltfcr}
    Assume that $h_i\sim \mu(h)$ is i.i.d. in $L_1$ and asymmetric. When the function $\frac{x^2}{2}-\bb E[\log\cosh(\sqrt{\theta_1}x+h)]$ has a single optimum $x^*$ of flatness $\tau=1$, then for $t\in\bb R$ ,  
    \begin{align*}
        \bb E\lef[\exp\bl\frac{t\sum_{i=1}^n(\sigma_i-\sqrt{\theta_1}x^*)}{\sqrt n}\br\rig]{\to}\exp\lef(\frac{\ca V(x^*)t^2}{2}\rig),\quad\bigg\Vert n^{-1/2}\sum_{i=1}^n(\sigma_i-\sqrt{\theta_1}x^*) \bigg\Vert_{\psi_2}<\infty.
    \end{align*}
    with $\ca V(x):=\frac{(1-\theta_1(\bb E[\sech^2(\sqrt{\theta_1}x+h)])^2-\bb E[\tanh(\sqrt{\theta_1}x+h)]^2)}{(1-\theta_1\bb E[\sech^2(\sqrt{\theta_1}x+h)])^2}$.
    
    When there exists a single optimum $x^*$ of flatness $\tau\geq 2$, for $t\in\bb R$  ,  
   \begin{align*}
        &\bb E\lef[\exp\lef(\frac{t\sum_{i=1}^n(\sigma_i-\sqrt{\theta_1}x^*)}{n^{\frac{4\tau-3}{4\tau-2}}}\rig)\rig]\to \int_{\bb R}\frac{(2\tau-1)x^{2\tau-2}}{\sqrt{2\pi v(x^*)}}\exp\lef(-\frac{x^{4\tau-2}}{2v(x^*)}+tx\rig)dx,\\
        &\bigg\Vert n^{-\frac{4\tau-3}{4\tau-2}}\sum_{i=1}^n(\sigma_i-\sqrt{\theta_1}x^*) \bigg\Vert_{\psi_{4\tau-2}}<\infty,
    \end{align*}
with $v(x)$ defined in \eqref{stirlingv}.
    \end{corollary}

\section{Discussions}\label{sect6}
In the previous sections we establish the statistical minimax rates for this hidden clique problem. It is noted that the scanning algorithm and the exact recovery algorithms both need exponential time complexity. In the previous literature, the statistical computational gaps are always observed for the Gaussian sparse PCA problem. However, some new insights are given in this work for the pRFCW model. This section discusses the upper bounds given by computationally efficient algorithms. We provide three polynomial algorithms here, along with their sample complexity results, to achieve powerful testing and exact recovery. 
\begin{enumerate}
    \item Algorithm \ref{alg:sdp} follows the ideas in  \citep{berthet2013optimal} in using the objective value of semi-definite programs to perform statistical tests. Its guarantees are given by proposition \ref{testingcomputationalupperbound}. We show that the results in the high-temperature regime correspond to the order in the Gaussian sparse PCA problem. However, for the critical temperature and critical temperature regimes, the results are much different from all of the existing literature. Surprisingly, we show that the low-temperature regime demonstrates no statistical/computational gaps. It is already known in \citep{brennan2018reducibility} that conditional on the planted clique conjecture being true, the computationally efficient rate of testing the Gaussian sparse PCA problem cannot be improved. However, it is also immediately checked that it is nontrivial to derive the same guarantee for the pRFCW model even at the high-temperature regime since the model is a countable mixture of different discrete sparse PCA model and the original reduction algorithms in \citep{brennan2018reducibility,brennan2019optimal} strongly rely on the Gaussianity of the design. The critical temperature regime demonstrates much different rates from the Gaussian sparse PCA problem. It is further unlikely that the identical method in \citep{brennan2019optimal} can provide a valid reduction from the planted clique problem.
    
    \item Algorithm \ref{sdprecovery} is the same semidefinite program as \citep{amini2009high,d2004direct}. The results obtained by this algorithm are presented in proposition \ref{propexactrecov}. Recall that in the sparse PCA literature, the simplest method that reaches the best computational upper bounds of recovery is the covariance diagonal thresholding algorithm proposed by \citep{johnstone2009sparse}. However, due to the fact that the pRFCW model only has $0$ on the diagonal entries, the thresholding algorithm does not work here. Instead, we extend the results in \citep{amini2009high} towards the complete $k=o(\sqrt n)$ regime using a matrix perturbation results given by \citep{yu2015useful}. For the high and critical temperature regimes, we show that the rates match exactly with the ones obtained by the testing algorithm. However, its optimality continues to be open. In the past, \citep{brennan2018reducibility} used a simple reduction from recovery to testing to obtain the computational lower bounds for recovery. Therefore, if we manage to find the testing lower bounds for the testing in the future, it naturally implies the recovery lower bounds.

    \item Algorithm \ref{alg:exactrecov} works for the large clique recovery that wins over the semi-definite program \ref{sdprecovery} when $k=\omega(\sqrt n)$. If compared with the sparse PCA problem, the best algorithm for this region is the spectral method \citep{krauthgamer2013semidefinite,paul2007asymptotics}. However, the spectral properties of Ising models remain a much more challenging task compared with the Gaussian case (or the Wishart matrices), which has long been established in \citep{paul2007asymptotics}. A simple intuition comes from the fact that the Gaussian vectors can always be seen as the permutation of isotropic multivariate Gaussians, whereas the mean-field Ising empirical covariance matrix is much more complicated to study. However, we show that a simple sum-recovery algorithm can obtain the same rate up to a log factor for the high-temperature regime (this regime sees similar statistical rates as the Gaussian sparse PCA). For the critical temperature regime, we show that the same algorithm undergoes a phase transition where it matches the minimax lower bound given $k=\omega( k^{\frac{2\tau-1}{4\tau-3}})$. This is surprising in the sense that the Gaussian sparse PCA \citep{krauthgamer2013semidefinite} never sees a computationally efficient algorithm that matches with minimax lower bounds given $k=\omega(\sqrt n)$. Moreover, it remains a long-standing open problem \citep{brennan2018reducibility} on whether the PC conjecture implies the sharp computational lower bounds of recovery given $k=\omega(\sqrt n)$. This also implies that much more is yet to be studied for this pRFCW model's computational phases, especially at the critical temperature regime.
\end{enumerate}
The following three propositions give formal guarantees of algorithms \ref{alg:sdp}, \ref{sdprecovery}, and \ref{alg:exactrecov} respectively.
\begin{proposition}\label{testingcomputationalupperbound}
    The sample complexity required for powerful testing using algorithm \ref{alg:sdp} at the different temperature regimes are summarized as follows
    \begin{align*}
        m\gtrsim \begin{cases}
        k^2\log n& \theta_1\in\lef(0,\frac{1}{\bb E[\sech^2(h)]}\rig)\\
         k^{\frac{2}{2\tau-1}}\log n& \theta_1=\frac{1}{\bb E[\sech^2(h)]}\text{ and }0\text{ is the global minimum of \eqref{taylorcond} with flatness }\tau\\
        \log n& \theta_1\in\lef(\frac{1}{\bb E[\sech^2(h)]},\infty\rig)
    \end{cases}.
    \end{align*}
\end{proposition}
\begin{algorithm}[htbp]
\caption{SDP Test of $k$-Clique}\label{alg:sdp}
\KwData{$\{\bfa\sigma^{(i)}\}_{i\in[m]}$ with $\bfa\sigma\in\{-1,1\}^n$}
 Let $\wh\Sigma:=\sum_{j=1}^m\frac{1}{m}\bfa\sigma\bfa\sigma^\top-\bfa I$\;
     Construct the following Semi-definite Program 
    \begin{center}
        Maximize T:=$\tr(\wh\Sigma X)$ \;
        subject to $\tr(X)=1$, $\Vert X\Vert_{1,1}\leq k$, $X\succeq 0$\;
    \end{center}
     
     \uIf{$0<\theta_1<\frac{1}{\bb E[\sech^2(h)]}$}{
        Reject $H_0$ if the value given by above program $SDP_k(\wh\Sigma)>\tau_\delta$\;
      }
      \uElseIf{$\theta_1=\frac{1}{\bb E[\sech^2(h)]}$}{
          \uIf{$0$ is the global minimum of \eqref{taylorcond} with flatness $\tau$}{Reject $H_0$ if the value given by above program $ k^{\frac{-2\tau+2}{2\tau-1}}SDP_k(\wh\Sigma)>\tau_\delta$\;}
          \Else{
             Reject $H_0$ if the value given by above program $\frac{1}{k}SDP_k(\wh\Sigma)>\tau_\delta$\;
          }}
      \Else{
        Reject $H_0$ if the value given by above program $\frac{1}{k}SDP_k(\wh\Sigma)>\tau_\delta$ \;}
        
\end{algorithm}

\begin{proposition}\label{propexactrecov}
     The sample complexity required for exact recovery using algorithm \ref{sdprecovery} at the different temperature regimes are summarized as follows
    \begin{align*}
        m\gtrsim \begin{cases}
        k^2\log n& \theta_1\in\lef(0,\frac{1}{\bb E[\sech^2(h)]}\rig)\\
         k^{\frac{2}{2\tau-1}}\log n& \theta_1=\frac{1}{\bb E[\sech^2(h)]}\text{ and }0\text{ is the global minimum of \eqref{taylorcond} with flatness }\tau\\
        \log n& \theta_1\in\lef(\frac{1}{\bb E[\sech^2(h)]},\infty\rig)
    \end{cases}
    \end{align*}
\end{proposition}

\begin{proposition}\label{propexactrecov2}
    The sample complexity required for exact recovery using algorithm \ref{alg:exactrecov} at the different temperature regimes are summarized as follows
    \tny{
\begin{align*}
        m\gtrsim \begin{cases}
        n\log n& \theta_1\in\lef(0,\frac{1}{\bb E[\sech^2(h)]}\rig)\\
         k^{\frac{1}{2\tau-1}}\log k\vee \frac{n}{k^{\frac{4\tau-4}{2\tau-1}}}\log (n-k)& \theta_1=\frac{1}{\bb E[\sech^2(h)]}\text{ and }0\text{ is the global minimum of \eqref{taylorcond} with flatness }\tau
    \end{cases}.
    \end{align*}}
\end{proposition}
An even more interesting question to be asked is to justify the rate obtained by the computational efficient algorithms are indeed not improvable. In the literature, many attempts have been made in the theoretical computer science community to study the proper certificate for the computational lower bounds. These attempts include \citep{barak2019nearly,hopkins2018statistical,feldman2017statistical,razborov1994natural,zdeborova2016statistical}. In the statistics community a widely accepted criteria is to compare the problem with that of the planted clique \citep{berthet2013complexity}. However, it remains a challenging task to make this comparision since we consider the uncountable mixture model whereas all the literature only consider fixed effect models.
\begin{algorithm}[htp!]
\caption{Computational Efficient Exact Recovery by SDP}\label{sdprecovery}
\KwData{$\{\bfa\sigma^{(i)}\}_{i\in[m]}$ with $\bfa\sigma\in\{-1,1\}^n$}
 Let $\wh\Sigma:=\sum_{j=1}^m\frac{1}{m}\bfa\sigma\bfa\sigma^\top$\;
 Choose $\rho$ and compute
\begin{align}
     \wh Z=
        \argmax\tr(\wh\Sigma Z)-\rho\Vert Z\Vert_{1,1},\quad s.t.\quad tr(Z)=1,\quad  Z\succeq 1,
\end{align}  
Factorize $\wh Z=\sum_{i=1}^n\lambda_i z_i z_i^\top$ with $\lambda_1\geq\lambda_2\geq\ldots\geq\lambda_n$\;
Return the index set of the positive/negative entries in $ z_1$ when $\theta_1>0$ and the positive/negative entries in $z_n$ when $\theta_1<0$.
\end{algorithm}
\begin{algorithm}
\caption{Computationally Efficient Exact Recovery by Row Sums}\label{alg:exactrecov}
\KwData{$\{\bfa\sigma^{(i)}\}_{i\in[m]}$ with $\bfa\sigma\in\{-1,1\}^n$}
 Let $\wh\Sigma:=\sum_{j=1}^m\frac{1}{m}\bfa\sigma\bfa\sigma^\top$\;

\uIf{$0<\theta_1<\frac{1}{\bb E[\sech^2(h)]}$}{
        For $i\in[m]$, compute $\delta_i=\sum_{j\in[n]}\wh\Sigma_{ij}$\;
      }
      \uElseIf{$\theta_1=\frac{1}{\bb E[\sech^2(h)]}$}{
          \uIf{$0$ is the global minimum of \eqref{taylorcond} with flatness $\tau$}{For $i\in[m]$, compute $\delta_i=k^{-\frac{2\tau-2}{2\tau-1}}\sum_{j\in[n]}\wh\Sigma_{ij}$\;}
          \Else{
             Use algorithm \ref{sdprecovery} to do the recovery\;
          }}
      \Else{
       Use algorithm \ref{sdprecovery} to do the recovery\;}

  Rank $\delta_i$ according to the order of magnitude and return the set of the top $k$ indices\;
\end{algorithm}
Then we summarize a few open problems and future working directions in this work beyond the computational aspects listed above.

\smallbreak
\emph{Connections to models with larger support.}
Despite that the major discussion of this work is based on the $\{-1,1\}$ spins, the same methods and techniques can also generalize to graphical models defined on $\bfa X\in \bb X^n$ for some $\bb X\subset\bb R$.
\begin{align*}
    \bb P(\bfa X):=\frac{1}{Z_{\theta_1}(\bfa\rho)}\exp\bl \frac{\theta_1}{n}\bl\sum_{i=1}^n X_i\br^2\br\prod_{i=1}^n\rho(X_i,\bfa h).
\end{align*}
for some random measure $\rho(\cdot,\bfa h)$ depending on r.v.s. $\bfa h$. This model is more diverse in mixed effect models than the $\{-1,1\}$ model considered in this work. However, the analysis of $\rho$ might require certain tail assumptions. 
\smallbreak
\emph{The Missing $\log k$ in the lower bounds.}
We note that the lower bounds for the small clique testing at part of the high-temperature regime and the critical temperature regime miss a $\log k$ factor from the upper bound. Although we believe this is improvable by other methods to give sharper control, there is little progress in the effort to drop it.

\begin{supplement}
\stitle{Supplement to “Hidden Clique Inference in Random Ising Model I: the planted random field Curie-Weiss model” }
\sdescription{In this supplementary material we provide the complete proof of the theorems in this work.}
\end{supplement}


\bibliographystyle{imsart-number} 
\bibliography{bib}       





\newpage
\setcounter{page}{1}
\setcounter{section}{0}

\begin{frontmatter}
\title{Supplement to "Hidden Clique Inference in Random Ising Model I:
the planted random field Curie-Weiss model"}






\end{frontmatter}
\startcontents[supplementary]
\printcontents[supplementary]{l}{1}{\setcounter{tocdepth}{2}}

\renewcommand{\thesection}{\Roman{section}}
\renewcommand{\thesubsection}{\thesection.\roman{subsection}}
\renewcommand{\thesubsubsection}{\thesubsection.\roman{subsubsection}}
\section{Additional Proof Sketches}
This section presents a few additional proof sketches for both theorem \ref{minimaxlb} and theorem \ref{cltrfcr}.
\subsection{Proof Sketch of Theorem \ref{minimaxlb}}\label{proofsketch3.1}
    The proof of the theorem can be separated into two different parts. The first part is when $\theta<\frac{1}{2\bb E[\sech^2(h)]}$ and the second part is the rest of the high temperature regime. The underlying reason for giving separate treatment to the two parts is boundedness of a typical statistics. When the temperature gets close to the critical, we need the additional \emph{the fake measure method}. The key idea underlying the lower bound is the Le Cam's lemma, which translates the lower bound of testing to that of upperbounding the TV distance. However, the direct computation of TV distance is not tractable and we resort to alternative divergence metrics given by the chi-square divergence. In particular, for $\bb P_{\bar S}:=\sum_{S\subset[n],|S|=k}\frac{1}{\binom{n}{k}}\bb P_{S}$ as the composed measure with uniform prior on the $k$-subsets of $[n]$, it is simple to check that the following general conditional Le Cam's method holds
     \begin{lemma}[Conditional Le Cam]
Assume that $P$, $Q$ are two probability measure dependent on random variables $\bfa h$ with $P\ll Q$ almost surely w.r.t. $\mu(\bfa h)$, we define
the conditional chi-squared divergence as $D_{\chi^2}(P,Q|\bfa h)=\int \lef((\frac{P(d\bfa\sigma|\bfa h)}{Q(d\bfa\sigma|\bfa h)})^2-1\rig)Q(d\bfa\sigma|\bfa h)$. Denote the $m$-fold product measure of alternative as $\bb P_{S,m}$. Denote $\bar{\bb P}_m=\frac{1}{\binom{n}{k}}\sum_{S:|S|=k}\bb P_{S,m}$ be the $m$-fold product measure under the arbitrary prior over the $k$-cardinal subsets of $[n]$ and $\bfa h^m=(\bfa h_1,\ldots,\bfa h_m)$ contains $m$ independent copies of $\bfa h$. Then, we have for all $\{0,1\}$-valued test statistics $\psi$ measurable with sigma-algebra generated with $\{\bfa\sigma^{(i)}\}_{i\in[m]}$, the following holds
\begin{align*}
    \inf_{\psi}\lef[\bb P_{0,m}(\psi=1)+\sup_{S:|S|=k}\bb P_{S,m}\lef(\psi=0\rig)\rig]\geq 1-\frac{1}{2}\sqrt{\bb E\lef[D_{\chi^2}(\bar{\bb P}_{m},\bb P_{0,m}|\bfa h^m)\rig]}
.\end{align*}
\end{lemma}

    Define $m_S:=\frac{1}{k}\sum_{i=1}^k\sigma_i$ and $m_{S^\prime}=\frac{1}{k}\sum_{i=r+1}^{r+k}$.
    By the decomposibility of $D_{\chi^2}$, our strategy of proof is to bound each term in $D_{\chi^2}$ separately, in the form of 
    \small
    \begin{align*}
        \bb E\lef[\frac{\bb P_{S}\bb P_{S^\prime}}{\bb P_{0}}\rig]=\bb E\lef[\frac{\sum_{\bfa\sigma}\exp\lef(\frac{\theta_1 k}{2}(m^2_{S}+m^2_{S^\prime})+\sum_{i\in[k+r]}h_i\sigma_i\rig)\sum_{\bfa\sigma}\exp\lef(\sum_{i\in[k+r]}h_i\sigma_i\rig)}{\sum_{\bfa\sigma}\exp\lef(\frac{\theta_1k}{2}m_S^2+\sum_{i\in[k+r]}h_i\sigma_i\rig)\lef(\sum_{\bfa\sigma}\exp\lef(\frac{\theta_1k}{2}m_{S^\prime}^2+\sum_{i\in[k+r]}h_i\sigma_i\rig)\rig)}\rig].
    \end{align*}
    \normalsize
    The above quantity can be studied using the Hubbard–Stratonovich transform, which turns the r.h.s. in the above equation into the following form
    \begin{align*}
        \bb E\lef[\frac{\int_{\bb R\times\bb R}\exp(-k\mca H_{0,k}(x,y,\bfa h))dxdy}{\int_{\bb R}\exp(-k\mca H_{1,k}(x,\bfa h))dx\int_{\bb R}\exp(-k\mca H_{2,k}(y,\bfa h))dy}\rig].
    \end{align*}
    Then we are ready to utilize the multivariate Laplace approximation method to transform the quantity within the expectation to a fixed point analysis problem, where the integral becomes in the form of
    \begin{align}\label{fixedpoint version}
        \bb E[\exp(-k\mca H_{0,k}(x^*,y^*,\bfa h)+k\mca H_{1,k}(x_1^*,\bfa h)+k\mca H_{2,k}(x_2^*,\bfa h))(1+o(1))],
    \end{align}
    where $x^*,y^*,x_1^*,x_2^*$ are the global maximum of their respective function. Notice that they are random given the filtration generated by $\bfa h$.
    Unsurprisingly, these fixed point is governed by the `mean field equation' under the physics terminology, in the form of
    \begin{align}\label{mfdequa}
        x_1^*= \bb E[\tanh(\sqrt{\theta_1}x_1^*+h)].
    \end{align}
  Then we analyze \eqref{fixedpoint version} using the canonical analysis of $Z$-estimators. Another technical question is how do we obtain the error terms in \eqref{fixedpoint version} since we observe that the majority of the terms indexed by the set $\ca S:=\{S:|S|=k,S\subset [n]\}$ is very close to $1$. This problem can be explicitly computed using a complicated iteration procedure given by lemma \ref{diverge}. Luckily, the error term in the Laplace method is negligible compared with the fluctuation in the exponential. Finally, we achieve the following upper bound after a few algebraic analyses.
  \begin{align*}
      \bb E\lef[\frac{\bb P_{S}\bb P_{S^\prime}}{\bb P_0}\rig]\leq\exp\lef(C\frac{(k-r)^2}{k^2}\rig)\quad\text{ when }r<k;\quad \bb E\lef[\frac{\bb P_{S}\bb P_{S^\prime}}{\bb P_0}\rig]=1\quad\text{ when }S\cap S^\prime=\emptyset.
  \end{align*}
    Then it comes to the composition procedure. This procedure is alternatively seen as an expectation w.r.t. the uniform measure on the $k$-subset of $[n]$:
    \begin{align}\label{combsum}
        \sum_{S:|S|=k,S\subset[n]}\frac{1}{\binom{n}{k}}\bb E\lef[\frac{\bb P_{S}\bb P_{[k]}}{\bb P_0}\rig]&=\sum_{r:r\leq k}\frac{\binom{k}{r}\binom{n-k}{k-r}}{\binom{n}{k}}\bb E\lef[\frac{\bb P_{[r:k+r]}\bb P_{[k]}}{\bb P_0}\rig]\nnb\\
        &\leq\sum_{r<k}f(n,r,k)\exp\lef(-C\frac{(k-r)^2}{k^2}\rig)+\frac{\binom{n-k}{k}}{\binom{n}{k}}.
    \end{align}
    And we can approximate the combinatorial quantity $f(n,r,k):=\frac{\binom{k}{r}\binom{n-k}{k-r}}{\binom{n}{k}}$ on the R.H.S. through Stirling's approximation. The final step is to first approximate the sum \eqref{combsum} through the Riemannian integral and use Laplace approximation again to get back to the estimate of chi-square.

    However, the reader might notice that this procedure requires a crucial regularity condition in \eqref{fixedpoint version}. It turns out that this question is closely related to the fixed point equation \ref{mfdequa}. Consider the simplest case when $r=0$, we see that $x^*=y^*$ both satisfy
    \begin{align*}
        x^*=\bb E[\tanh(2\sqrt{\theta_1}x^*+h)].
    \end{align*}
    Hence when $\theta_1>\frac{1}{2\bb E[\sech^2(h)]}$, $x_1^*$ does not converge to the same limit as $x^*$ and the chi-square divergence is of the order of $\exp(k)$ ( In particular, consider the case of $h$ being symmetric, the numerator converges to two points of order $O(1)$ symmetric with $0$ whereas the denominator converges to $0$. ). This problem implies that this method does not work for temperature that is close to critical. Instead, an alternative idea is to construct an intermediate `fake measure' interpolating into the TV distance. The name `fake measure' is due to the fact that the interpolating measure is not a probability measure but we still study the information divergence of it. However, to compute the chi-square divergence of fake measure is still non-trivial, which involves a more delicate transfer principle to complete the analysis. This method can overcome the problem of diverging maximum but necessarily induce a loss of logarithmic factor in $k$. These technicalities can be traced in the formal proof of this theorem and we deferred it to the proof sketch of lemma \ref{minimaxlbct}.

\subsection{Proof Sketch of Theorem \ref{cltrfcr}}
        Due to space limitations, we give only a heuristic proof of the above two results. Our results are proved using the idea of asymptotic integral expansion and transfer principle.
        After the H-S transformation. We directly work on the mgf, for $\beta\in(0,1)$, the mgf is given by 
        \begin{align}\label{mgfintegralform}
            \bb E\bigg[\exp\bl\frac{t}{n^{\beta}}\sum_{i=1}^n\sigma_i\br\bigg]=\bb E\bigg[\frac{\int_{\bb R}\exp(-n\mca H_{0,n}(x))dx}{\int_{\bb R}\exp(-n\mca H_{1,n}(x))dx}\bigg].
        \end{align}
        with 
        \begin{align*}
    \mca H_{0,n}(x,\bfa h):&=\frac{1}{2}x^2-\frac{1}{n}\sum_{i=1}^n\log\cosh\lef(\sqrt{\theta_1} x+h_i+\frac{t}{ n^\beta}\rig),\\ \mca H_{1,n}(x,\bfa h):&=\frac{1}{2}x^2-\frac{1}{n}\sum_{i=1}^n\log\cosh\lef(\sqrt{\theta_1} x+h_i\rig).
    \end{align*}
    Therefore, the convergence rate of average magnetization depends  on the local landscape w.r.t. $t$ around the global minimum of $\mca H_{0,n}(x,\bfa h)$ and $\mca H_{1,n}(x,\bfa h)$, given by
    \begin{align*}
        x_0=\frac{1}{n}\sum_{i=1}^n\tanh\lef(\sqrt{\theta_1}x+h_i+\frac{t}{n^{\beta}}\rig),\qquad x_1=\frac{1}{n}\sum_{i=1}^n\tanh\lef(\sqrt{\theta_1}x+h_i\rig).
    \end{align*}
    At the high and critical temperature regimes, we need to study the fluctuations given by the principle term in \eqref{mgfintegralform}, which are given by
$        \bb E[\exp( -n\mca H_{0,n}(x_0)+n\mca H_{1,n}(x_1))].
$ And the main technical tool is a strong $Z$-estimator linearization lemma appears in the supplementary material and a fixed point analysis for the limiting value of $x_0$ and $x_1$, governed by the uniform convergence of $\mca H_{0,n}$ and $\mca H_{1,n}$ to their respective limits.  For the low temperature regime ($\theta_1>\frac{1}{\bb E[\sech^2(h)]}$) this problem becomes more complicated since we have multiple global minima that take almost the same function value of $\mca H_{0,n}$ and $\mca H_{1,n}$. To overcome this barrier, we prove a new transfer principle for the random integral, which implies the following
    \begin{align*}
         \bb E\bigg[\exp\bl\frac{t}{n^{\beta}}\sum_{i=1}^n\sigma_i\br\bigg]=\bb E\bigg[\frac{\int_{\bb R^+}\exp(-n\mca H_{0,n}(x))dx}{\int_{\bb R^+}\exp(-n\mca H_{1,n}(x))dx}\bigg]+O(\exp(-n\delta)),\quad\delta>0.
    \end{align*}
    And the rest of the proof follows through.

\section{Proof of Results in Section \ref{sect3}}
\subsection{Proof of Theorem \ref{minimaxlb}}\label{proofminimaxlb}
    We note that the $\chi$-square divergence admits an alternative form. We denote $\ca S:=\{S:|S|=k\}$. Let $\bfa h=\{ h_{ij}\}_{i,j\in[n]}$ denote the random field, then
    \begin{align}\label{chisqlb}
        \bb E\lef[D_{\chi^2}\lef(\bar{\bb P}_m,\bb P_{0,m}|\bfa h^m\rig)\rig]:=\frac{1}{\binom{n}{k}^2}\sum_{S,S^\prime\in\ca S}\bb E_{0,m}\lef[\frac{\bb P_{S,m}(\bfa\sigma|\bfa h^m)\bb P_{S^\prime,m}(\bfa\sigma|\bfa h^m)}{\bb P_{0,m}^2(\bfa\sigma|\bfa h^m)}\rig]-1
    .\end{align}

In the following proof we discuss over the two possible regimes: (1) When there is no overlap between $S$ and $S^\prime$ or $r=1$. (3) When the overlap is $\frac{r}{k}:=c$

\begin{center}
    \textbf{1. When $r=k$}
\end{center}
We start with analyzing non-overlapped $S$ and $S^\prime$ since according to lemma \ref{combprob} when $k=o\lef(\sqrt n\rig)$  we have almost surely  $S\cap S^\prime=\emptyset$. It is then noticed that 
\ttny{
\begin{align}\label{firstcase}
   \bb E\lef[\frac{\sum_{\bfa\sigma}\exp\bigg(\frac{\theta_1}{2k}\lef(\lef(\sum_{i\in [k]}\sigma_i\rig)^2+\lef(\sum_{k\in[k+1:2k]}\sigma_i\rig)^2\rig)+\sum_{i\in[2k]}h_i\sigma_i\bigg)\sum_{\bfa\sigma}\exp\bigg(\sum_{i\in[2k]}\sigma_ih_i\bigg)}{\lef(\sum_{\bfa\sigma}\exp\lef(\frac{\theta_1}{2k}\lef(\sum_{i\in[k]}\sigma_i\rig)^2+\sum_{i\in[2k]}h_i\sigma_i\rig)\rig)\lef(\sum_{\bfa\sigma}\exp\lef(\frac{\theta_1}{2k}\lef(\sum_{i\in[k+1:2k]}\sigma_i\rig)^2+\sum_{i\in[2k]}h_i\sigma_i\rig)\rig)}\rig]=1
.\end{align}}
\begin{center}
    \textbf{2. When $r<k$ }
\end{center}
Using $\bfa h:=(h_1^\prime,\ldots,h_{k+r}^\prime)^\top$ to denote a random vector. We show that there exists $G_{0,k}(x,y,\bfa h):\bb R^2\times\bb R^{k+r}\to\bb R$, and $G_{1,k}(x,\bfa h),G_{2,k}(x,\bfa h):\bb R\times\bb R^{k+r}\to\bb R$ such that
\begin{align}\label{vphi0g}
\bb E\lef[\frac{\bb P_{S}(\bfa\sigma)\bb P_{S^\prime}(\bfa\sigma)}{\bb P_0(\bfa\sigma)}\rig]=\bb E\lef[\frac{\prod_{i=r+1}^k\cosh(h_i)\int\exp\lef(-kG_{0,k}(x,y,\bfa h)\rig)dxdy }{\int\exp\lef(-kG_{1,k}(x,\bfa h)\rig)dx\int\exp\lef(-kG_{2,k}(y,\bfa h)\rig)dy}\rig]
.\end{align}
Using Gaussian integration identity, we observe that the numerator can be written as
\begin{align*}
    &\sum_{\bfa\sigma}\exp\bl\frac{\theta_1k}{2}(m_S^2+m_{S^\prime}^2)+\sum_{i\in[k+r]}\sigma_ih_i\br\\
    &=\frac{1}{2\pi}\sum_{\bfa\sigma}\int\int\exp\bl-\frac{x^2+y^2}{2}+\sqrt{\theta_1k}(m_Sx+m_{S^\prime}y)+\sum_{i\in[k+r]}\sigma_ih_i\br dxdy\\
   &=\frac{1}{2\pi}\int\int\sum_{\bfa\sigma}\exp\bl-\frac{x^2+y^2}{2}+\sqrt{\frac{\theta_1}{k}}\bl\sum_{i\in[r+1:k]}\sigma_i(x+y)+\sum_{i\in[r]}\sigma_ix+\sum_{i\in[k+1:k+r]}\sigma_iy\br\\&+\sum_{i\in[k+r]}\sigma_ih_i\br dxdy\\
    &=\frac{k2^{k+r}}{2\pi}\int\int\exp\lef(-k G_{0,k}(x,y,\bfa h)\rig)dxdy
.\end{align*}
with  
\begin{align*}
G_{0,k}:&=\frac{x^2+y^2}{2}-\frac{1}{k}\bl\sum_{i=1}^{r}\log\cosh\lef(\sqrt{\theta_1}x+h_i\rig)\\
&+\sum_{i=r+1}^k\log\cosh(\sqrt{\theta_1}(x+y)+h_i)+\sum_{i=k+1}^{k+r}\log\cosh(\sqrt{\theta_1}y+h_i)\br
.\end{align*}
And analogously we can check that
\begin{align*}
    \sum_{\bfa\sigma}\exp\bl\frac{\theta_1k}{2}m_S^2+\sum_{i=1}^{k+r}\sigma_ih_i\br=\sqrt{\frac{k}{2\pi}}2^{k+r}\prod_{i=k+1}^{k+r}\cosh(h_i)\int\exp(-kG_{1,k}(x,\bfa h))dx
.\end{align*}
with
\begin{align*}
    G_{1,k}(x,\bfa h):=\frac{x^2}{2}-\frac{1}{k}\sum_{i=1}^k\log\cosh\lef(\sqrt{\theta_1}x+h_i\rig)
.\end{align*}
And using the symmetry between $S$ and $S^\prime$  we have
\begin{align*}
    G_{2,k}(x,\bfa h):=\frac{x^2}{2}-\frac{1}{k}\sum_{i=r+1}^{k+r}\log\cosh\lef(\sqrt{\theta_1}x+h_i\rig)
.\end{align*}

Then we define the following population varieties
\begin{align*}
    G_0(x,y)&=\frac{x^2+y^2}{2}-c\bb E[\log\cosh(\sqrt{\theta_1}x+h)\cosh(\sqrt{\theta_1}y+h)]\\&-(1-c)\bb E[\log\cosh(\sqrt{\theta_1}(x+y)+h)],\\
    G_1(x)&=\frac{x^2}{2}-\bb E[\log\cosh(\sqrt{\theta_1}x+h)]
.\end{align*}
The following regularity conditions is important in the uniform convergence criteria.
\begin{lemma}[Regularity Conditions]\label{convergeasG} When $h$ is in $L_1$.
    Almost surely in $\mu(\bfa h)$ and uniformly on $(x,y)$  we have
    \begin{align*}
        G_{0,k}^{(j_1,j_2)}(x,y,\bfa h):=\frac{\pta^{j_1+j_2}G_{0,k}(x,y,\bfa h)}{\pta x^{j_1}\pta y^{j_2}}\to G_0^{(j_1,j_2)}(x,y)
    .\end{align*}
    with $G_{0,k}^{(0,0)}:=G_{0,k}$. Similar argument holds for $G_{1,k},G_{2,k}\to G_1$. And condition \eqref{laplace2}, \eqref{laplace3} in lemma \ref{laplace} holds for $G_{0,k}, G_{1,k}, G_{2,k}$.
\end{lemma}
\begin{proof}
For the first condition, we define (note that $h_i$ depends on $k$ and we add it as subscript )
\begin{align*}
    \varphi_k(x,y,\bfa h)&:=-\frac{1}{k}\bl\sum_{i\in[r]}\log \cosh(\sqrt{\theta_1} x+h_{i}) +\sum_{i\in[r+1:k]}\log \cosh(\sqrt{\theta_1} x+\sqrt{\theta_1}y+h_{i}) \\
    &+\sum_{i\in[k+1:k+r]}\log\cosh(\sqrt{\theta_1}y+h_i)\br
\end{align*}
and
\begin{align*}
    \varphi(x,y)&:=-c\bb E[(\log \cosh(\sqrt{\theta_1} x+h)]-(1-c)\bb E[\log \cosh(\sqrt{\theta_1} x+\sqrt{\theta_1}y+h)]\\&-c\bb E[\log \cosh(\sqrt{\theta_1}y+h)]
.\end{align*}
It is not hard to see that by SLLN almost surely we have $\varphi_k(x,y,\bfa h)\to\varphi(x,y)$ point-wise. 
Then we can check that for $(x_1,y_1)$ and $(x_2,y_2)\in\bb R^2$  we have
\begin{align*}
    \varphi_k(x_1,y_1,\bfa h)-\varphi_k(x_2,y_2,\bfa h)\leq 2\lef(\sqrt{\theta_1}|x_1-x_2|+\sqrt{\theta_1}|y_1-y_2|\rig),\forall k
.\end{align*}
implies that $\varphi_k$ form an uniformly equicontinuous sequence. Since countable intersection of sets with measure $1$ has also measure $1$ we conclude that it is possible to choose $A\subset\Omega$ such that $\mu(A)=1$ such that $\forall \bfa h\in A$, $\varphi_k(x,y,\bfa h)\to\varphi(x,y)$. This implies that $G_{0,k}\to G_0$ uniformly almost surely (A simple exercise using Arzelà–Ascoli theorem). Similar argument can be  verified to hold for $G_{1,k}, G_{2,k}$ and we omit it here.

Then we move toward the discussion over the derivatives. Since we verified that the derivatives of $\varphi^{(i,j)}(x,y,\bfa h):=\frac{\pta^{i+j}\varphi(x,y,\bfa h)}{\pta x^i\pta y^j}$ is bounded. Therefore, we conclude that $\varphi^{(i,j)}(x,y,\bfa h):=\frac{\pta^{i+j}\varphi(x,y,\bfa h)}{\pta x^i\pta y^j}$ is equicontinuous and hence uniformly almost surely converging to $\varphi^{(i,j)}(x,y):=\frac{\pta^{i+j}\varphi(x,y)}{\pta x^i\pta y^j}$. This implies that the derivatives also converges unformly almost surely. Similar arguments can be analogously applied to $\mca H_1$.

For the second condition, noticing that $\log \cosh(x+y)\leq 2\log2 +|x|+|y|$ we see that:
\begin{align*}
-\varphi_k(x,y,\bfa h)&\leq \frac{1}{k}\sum_{i=1}^{k+r}2|h_i|+\sqrt{\theta_1}|x|+\sqrt{\theta_1}|y|+4\log 2\\
&\leq 2|x|+2|y|+4\log 2+\frac{2}{k+r}\sum_{i\in[k+r]}|h_i|
.\end{align*}
which consequently shows that
\begin{align*}
    G_{0,k}(x,y,\bfa h)\geq\frac{x^2+y^2}{2}-2|x|-2|y|-4\log 2-\frac{2}{k+r}\sum_{i\in[k+r]}|h_i|
.\end{align*}
We denote  $C(\bfa h)=16\exp\lef(\frac{2}{k+r}\sum_{i\in[k+r]}|h_i|\rig)$ it is checked that by dominated convergence theorem and $h$ is in $L_1$  we have
\begin{align*}
    \int\exp(-G_0(x,y))dxdy&=\lim_{k\to\infty}\int\exp\lef(-G_{0,k}(x,y,\bfa h)\rig)dxdy\\
    &\leq \exp\lef(\int_{\bb R}2|h|d\mu(h)\rig)\int_C\exp\lef(-\frac{x^2+y^2}{2}+2|x|+2|y|\rig)dxdy\\
    &\leq A\exp(2|u|)<\infty
.\end{align*}
for some constant $A$ not dependent on $x,y,u$
Similar argument also holds for $\mca H_{1}$ and we complete the proof.
\end{proof}

Then we define $(x_k,y_k)=\argmin_{(x,y)\in\bb R^2}G_{0,k}(x,y,\bfa h)$ that is almost surely unique by the fact that $G_0(x,y)$ has unique minimum  at $\theta_1<1$ regime denoted by $(x^*,y^*)$ with $(x_k,y_k)\to (x^*,y^*)$ by uniform convergence of all order of derivatives. Analogously we can define $x_{1,k}=\argmin_{x\in\bb R}G_{1,k}(x,\bfa h)$ and $x_{2,k}=\argmin_{x\in\bb R}G_{2,k}(x,\bfa h)$ that both converge to $x_1^*=0=\argmin_{x\in\bb R}G_1(x,\bfa h)$. Note that here we do not take derivative w.r.t. $\bfa h$ as in lemma \ref{laplace} and the following holds
\begin{align*}
    \int\exp\lef(-k G_{0,k}(x,y,\bfa h)\rig)dxdy &= \exp\lef(-k G_{0,k}(x_k,y_k,\bfa h)\rig)\det\bl\frac{k\nabla^2G_{0,k}(x_k,y_k,\bfa h)}{2\pi}\br^{1/2}\\
    &\cdot\lef(1+{a_0(\bfa h)}{k^{-1}}+O\lef({k^{-2}}\rig)\rig)
.\end{align*}
And for $i\in\{1,2\}$  we have
\begin{align*}
   \int\exp (-kG_{i,k}(x,\bfa h))dx &=\exp\lef(-kG_{i,k}(x_{i,k},\bfa h)\rig)\det\bl\frac{k G_{i,k}^{(2)}(x_{i,k},\bfa h)}{2\pi}\br^{1/2}\\
   &\cdot\lef(1+{a_i(\bfa h)}{k^{-1}} + O\lef(k^{-2}\rig )\rig)
.\end{align*}
Thus we can rewrite \eqref{vphi0g} as follows:
\begin{align}\label{varphi0rep}
    \bb E\lef[\frac{\bb P_{S}(\bfa\sigma)\bb P_{S^\prime}(\bfa\sigma)}{\bb P_0(\bfa\sigma)}\rig]&=\frac{\lef(G_{1,k}^{(2)}(x_{1,k},\bfa h)G_{2,k}^{(2)}(x_{2,k},\bfa h)\rig)^{1/2}}{\det(\nabla^2 G_{0,k}(x_k,y_k,\bfa h))^{1/2}}\nnb\\&\cdot\exp\lef(-kG_{0,k}(x_k,y_k,\bfa h)+kG_{1,k}(x_{1,k},\bfa h)+kG_{2,k}(x_{2,k},\bfa h)\rig)\nnb\\ &\cdot\prod_{i=r+1}^k\cosh(h_i)\cdot\lef(1+\frac{1}{2}\frac{a_1(\bfa h)+a_2(\bfa h)-a_0(\bfa h)}{k}+O\lef(\frac{1}{k^2}\rig)\rig)
.\end{align}

For the first term we introduce $\sum_{1}:=\sum_{i=1}^r$, $\sum_{2}:=\sum_{i=r+1}^k$, and $\sum_{3}:=\sum_{i=k+1}^{k+r}$ to simplify notations and get the following:
\tny{
\begin{align}\label{hessian}
    &G^{(2)}_{1,k}(x,\bfa h)=1-\theta_1+\frac{\theta_1}{k}\sum_{i=1}^k\tanh^2(\sqrt{\theta_1}x+h_i),\thinspace G^{(2)}_{2,k}(x,\bfa h)=1-\theta_1+\frac{\theta_1}{k}\sum_{i=r+1}^{k+r}\tanh^2(\sqrt{\theta_1}x+h_i),\nnb\\
    & \nabla^2G_{0,k}(x,y,\bfa h)=\\
    &\tiny\begin{bmatrix} 
    1-\frac{\theta_1}{k}\lef(\sum_1\sech^2(\sqrt{\theta_1}x+h_i)+\sum_2\sech^2(\sqrt{\theta_1}(x+y)+h_i)\rig)&-\frac{\theta_1}{k}\sum_{2}\sech^2(\sqrt{\theta_1}(x+y)+h_i))\nnb\\-\frac{\theta_1}{k}\sum_{2}\sech^2(\sqrt{\theta_1}(x+y)+h_i))
    &1-\frac{\theta_1}{k}\lef(\sum_3\sech^2(\sqrt{\theta_1}y+h_i)+\sum_2\sech^2(\sqrt{\theta_1}(x+y)+h_i)\rig)
    \end{bmatrix}
.\end{align}}
Here we introduce a few quantities as $\det\lef(\nabla^2 G_{0,k}(0,0,\bfa h)\rig)$ and the derivatives of \\$\det\lef(\nabla^2 G_{0,k}(x_k,y_k)\rig)$ at $(x,y)=(0,0)$, given by
\sm{\begin{align*}
A_0(\bfa h)&:=\bl 1-\frac{\theta_1}{k}\sum_{i=1}^k\sech^2(h_i)\br\bl 1-\frac{\theta_1}{k}\sum_{i=r+1}^{k+r}\sech^2(h_i)\br-\bl\frac{\theta_1}{k}\sum_{i=r+1}^k\sech^2(h_i)\br^2,\\
    A^k_{1}(\bfa h)&:=\frac{2\sqrt{\theta_1}}{k}\sum_{i=1}^k\tanh(h_i)\sech^2(h_i),\thinspace A^k_{2}(\bfa h) :=\frac{\theta_1}{k}\sum_{i=1}^k\lef(\sech^4(h_i) - 2 \sech^2(h_i)\tanh^2(h_i)\rig),\\
    A^k_3(\bfa h)&:=\frac{2\sqrt{\theta_1}}{k}\sum_{i=1}^k\tanh(h_i)\sech^2(h_i),\thinspace A^k_4(\bfa h):=\frac{\theta_1}{k}\sum_{i=r+1}^{k+r}\lef(\sech^4(h_i) - 2 \sech^2(h_i)\tanh^2(h_i)\rig).
\end{align*}}
Therefore we note that $A_0,A_1,A_2,A_3,A_4$ are bounded from below as well as above. Further we note that $\sqrt kA_1$ and $\sqrt k A_3$ converges to Gaussian with constant variance.
\begin{gather*}
    A_5^k(\bfa h) :=-\frac{\pta G_{0,k}(0,0,\bfa h))}{\pta x}=\frac{\sqrt{\theta_1}}{k}\sum_{i=1}^k\tanh(h_i),\\ A_6^k(\bfa h):=-\frac{\pta^2G_{0,k}(0,0,\bfa h)}{\pta x^2}=\frac{\theta_1}{k}\sum_{i=1}^k\sech^2(h_i)-1,\\
     A_7^k(\bfa h) :=-\frac{\pta G_{0,k}(0,0,\bfa h)}{\pta y}=\frac{\sqrt{\theta_1}}{k}\sum_{i=r+1}^{r+k}\tanh(h_i),\\
    A_8^k(\bfa h):=-\frac{\pta^2 G_{0,k}(0,0,\bfa h)}{\pta y^2}=\frac{\theta_1}{k}\sum_{i=r+1}^{r+k}\sech^2(h_i)-1,\\
    A^k_9(\bfa h):=-\frac{\pta^2 G_{0,k}(0,0,\bfa h)}{\pta y\pta x}=\frac{\theta_1}{k}\sum_{i=r+1}^k\sech^2(h_i)
.\end{gather*} 
where it is   checked that $\sqrt kA_5,\sqrt k A_7$ are Gaussian with constant variance. It is also checked that $A_6^k$ and $A_8^k$ are even function.

Then doing Taylor expansion of $\nabla G_{0,k}$ around $(0,0)$ we see that
\begin{align}\label{cross}
\bfa 0=-\begin{bmatrix}
A_5^k(\bfa h)\\
A_7^k(\bfa h)
\end{bmatrix}-\begin{bmatrix}
    A_6^k(\bfa h) & A_9^k(\bfa h)\\
    A_9^k(\bfa h) & A_8^k(\bfa h)
\end{bmatrix}
\begin{bmatrix}
    x_k\\
    y_k
\end{bmatrix} + O(x_k^2+y_k^2)
.\end{align}
where we note that asking for $\begin{bmatrix}
    A^k_6(\bfa h)&A^k_9(\bfa h)\\
    A^k_9(\bfa h)&A^k_8(\bfa h)
\end{bmatrix}$ to be positive definite in we need to have
\begin{align}\label{hessianhe}
   \det(&\nabla^2G_{0,k}(0,0,\bfa h))= A_6^k(\bfa h)A_8^k(\bfa h)-A_9^{k,2}(\bfa h)\nnb\\
   &=\bl 1-\frac{\theta_1}{k}\sum_{i=1}^k\sech^2(h_i)\br\bl1-\frac{\theta_1}{k}\sum_{i=r+1}^{k+r}\sech^2(h_i)\br-\bl\frac{\theta_1}{k}\sum_{i=r+1}^k\sech^2(h_i)\br^2>0.
\end{align}
Consider its asymptotics  we have
\begin{align*}
    \det\lef(\nabla^2G_{0}(0,0)\rig)&=c\theta_1\lef(\theta_1\bb E[\sech^2(h)]^2-\bb E[\sech^2(h)]+\theta_1\bb V(\sech^2(h))\rig)\\
    &+(1-\theta_1^2\bb E[\sech^2(h)]^2-\theta_1^2\bb V(\sech^2(h))).
\end{align*}
Therefore, when the following holds,  we have the equation \eqref{hessianhe} is positive.
\begin{align}\label{caseI}
     c> \frac{-(1-\theta_1^2\bb E[\sech^2(h)]^2-\theta_1^2\bb V(\sech^2(h)))}{\theta_1(\theta_1\bb E[\sech^2(h)]^2-\bb E[\sech^2(h)]+\theta_1\bb V(\sech^2(h)}\quad\text{ or }\quad\theta_1<\frac{1}{2\bb E[\sech^2(h)]}
\end{align}
\begin{center}
    \textbf{Case I: $c> \frac{-(1-\theta_1^2\bb E[\sech^2(h)]^2-\theta_1^2\bb V(\sech^2(h)))}{\theta_1(\theta_1\bb E[\sech^2(h)]^2-\bb E[\sech^2(h)]+\theta_1\bb V(\sech^2(h)}$ or $\theta_1<\frac{1}{2\bb E[\sech^2(h)]}$}\label{convergencond}
\end{center}

It is   checked that in case I the denominator and the numerator converges together to $(0,0)$, we then have
\begin{align*}
    \begin{bmatrix}
        x_k\\
        y_k
    \end{bmatrix}=\frac{1}{A_9^2-A_6A_8}\begin{bmatrix}
        A_8A_5-A_9A_7\\
        A_6A_7-A_9A_5
    \end{bmatrix}+o\lef(\frac{1}{\sqrt k}\rig).
\end{align*}
And for the other two functions $G_{1,k}^\prime$ and $G_{2,k}^\prime$  we can write their expansion together as:
\begin{align}\label{nocross}
\bfa 0=-\begin{bmatrix}
A_5^k(\bfa h)\\
A_7^k(\bfa h)
\end{bmatrix}-\begin{bmatrix}
    A_6^k(\bfa h) & 0\\
    0 & A_8^k(\bfa h)
\end{bmatrix}
\begin{bmatrix}
    x_{1,k}\\
    x_{2,k}
\end{bmatrix} + O(x_{1,k}^2)
.\end{align}
which also implies that $\sqrt k x_k,\sqrt k y_k,\sqrt kx_{1,k},\sqrt kx_{2,k}$ converging in distribution to Gaussian.
Therefore \eqref{cross} and \eqref{nocross} yield that
\begin{align*}
    \begin{bmatrix}
        x_k-x_{1,k}\\
        y_k-x_{2,k}
    \end{bmatrix}&=\begin{bmatrix}
    \frac{A_9^2}{(A_6A_8-A_9^2)A_6} & \frac{A_9}{A_6A_8-A_9^2}\\
    \frac{A_9}{A_6A_8-A_9^2} & \frac{A_9^2}{(A_6A_8-A_9^2)A_8}
\end{bmatrix}
\begin{bmatrix}
        A_5^k(\bfa h)\\
        A_7^k(\bfa h)
    \end{bmatrix}+o(x_{1,k}-x_k)\\
    &=\begin{bmatrix}
        \frac{A_9(A_9A_5+A_7A_6)}{(A_6A_8-A_9^2)A_6} \\
    \frac{A_9(A_9A_7+A_5A_8)}{(A_6A_8-A_9^2)A_8}
    \end{bmatrix}+o(x_{1,k}-x_k)
.\end{align*}
which implies that $\frac{k\sqrt k}{(k-r)}(x_{1,k}-x_k)$ converges to Gaussian with constant variance.
By definition of $G_{0,k},G_{1,k},$ and $G_{2,k}$ it is   checked that
\begin{align*}
    \frac{\pta^\tau G_{0,k}(0,0,\bfa h)}{\pta x^\tau}=G^{(\tau)}_{1,k}(0,\bfa h),\qquad\qquad\frac{\pta^\tau G_{0,k}(0,0,\bfa h)}{\pta y^\tau}=G_{2,k}^{\tau}(0,\bfa h).
\end{align*}
And,  we have
\begin{align*}
   G_{1,k}^{(2)}(x_{1,k},\bfa h)&=1-\theta_1\frac{1}{k}\sum_{i=1}^k\sech^2(h_i)+A_1^k(\bfa h)x_{1,k}+ A_2^k(\bfa h)x_{1,k}^{2}+O(x_{1,k}^{3}),\\
   G_{2,k}^{(2)}(x_{2,k},\bfa h)&=1-\theta_1\frac{1}{k}\sum_{i=r+1}^{k+r}\sech^2(h_i)+A_3^k(\bfa h)x_{2,k}+ A_4^k(\bfa h)x_{2,k}^{2}+O(x_{2,k}^{3})
.\end{align*}
And analogously, using \eqref{hessian} and expand at $(0,0)$   give
\begin{align*}
    \det\lef(\nabla^2 G_{0,k}(x_k,y_k)\rig)&=A_0(\bfa h)+A_1^k(\bfa h) x_k +A_3^k(\bfa h)y_k+A_1^k(\bfa h)A_3^k(\bfa h)x_ky_k\nnb\\
    &+A_2^k(\bfa h)x_k^2+A_4^k(\bfa h)y_k^2 +O\lef(x_k^{2\tau-1}+y_k^{2\tau-1}\rig)
.\end{align*}
Consequently for the first term in \eqref{varphi0rep}  we have for all $\delta>1$
\begin{align}\label{varphi01}
    \bb E& \lef[\frac{\lef(G_{1,k}^{(2)}(x_{1,k},\bfa h)G_{2,k}^{(2)}(x_{2,k},\bfa h)\rig)^{\delta/2}}{\det(\nabla^2 G_{0,k}(x_k,y_k,\bfa h))^{\delta/2}}\rig]=1+\bb E\bigg[\frac{\delta}{2A_0(\bfa h)}\big(A_1^k(\bfa h)(x_{1,k}-x_k)+A_3^k(\bfa h)(x_{2,k}-y_k)\nnb\\
    &+A_1^k(\bfa h)A_3^k(\bfa h)(x_{1,k}x_{2,k}-x_ky_k)+A_2^k(\bfa h)(x_{1,k}^2-x_k^2)+A_4^k(\bfa h)(x_{2,k}^2-y_k^2)\big) +O\lef(x_k^{3}-x_{1,k}^{3}\rig)\bigg]\nnb\\
    &=1+O\lef(({k-r}){k^{-2}}\rig)
.\end{align}

Then  we have by symmetry of $x_k$ and $y_k$ and expand $G_{0,k}, G_{1,k},G_{2,k}$ at $(0,0)$, and $0$, $0$ respectively to see that there exists $C_1,C_2>0$ such that
\begin{align}\label{varphi02}
\bb E&\exp\bl -k\bl  G_{0,k}(x_k,y_k,\bfa h)-\frac{1}{k}\sum_{i=r+1}^k\cosh(h_i)-G_{1,k}(x_{1,k},\bfa h)-G_{2,k}(x_{2,k},\bfa h)\br \br\nnb\\
&=\bb E\exp\bl -k(A_5^k(\bfa h)(x_{1,k}-x_k)\nnb+A_7^k(\bfa h)(x_{2,k}-y_k))-\frac{k}{2}A^k_6(\bfa h)(x_{1,k}^2-x_k^2) \nnb\\
&-\frac{k}{2}A^k_8(\bfa h)(x_{2,k}^2-y_k^2)+kA_9^k(\bfa h)x_ky_k +O(k-r)\frac{A_5^k(\bfa h)^2A_7^k(\bfa h)}{A_6^k(\bfa h)^2A_8^k(\bfa h)}\br\nnb\\
&=\bb E\lef[\exp \lef(\frac{k}{2}A^k_5(\bfa h)(x_k-x_{1,k})+\frac{k}{2}A^k_7(\bfa h)(x_k-x_{2,k})+O((k-r)A_5^2A_7/(A_6^2A_8))\rig)\rig]\nnb\\
&=\bb E\lef[\exp\lef(\lef(\frac{k}{2}\frac{A_9^2(A^2_5/A_6+A_7^2/A_8)+2A_5A_7A_9}{(A_6A_8-A_9)^2}\rig)+O((k-r)A_5^2A_7/(A_6^2A_8))\rig)\rig]\nnb\\
&\leq\lef(\bb E\exp\lef(\tau_1kC_1(A_5^2+A_7^2)A_9^2\rig)\rig)^{1/\tau_1}\lef(\bb E\exp\lef(\tau_2kC_2A_5A_7A_9\rig)\rig)^{1/\tau_1}\nnb\\
&\cdot\bb E\lef(\exp\lef(O(k-r)\frac{A_5^k(\bfa h)^2A_7^k(\bfa h)}{A_6^k(\bfa h)^2A_8^k(\bfa h)}\rig)\rig)^{1/\tau_3}
.\end{align}
for all $\tau_1,\tau_2,\tau_3>1$ and $1/\tau_1+1/\tau_2=1$ that comes from the condition of H\"older's inequality.
We note that it is checked that $\sqrt kA_5^k(\bfa h)=\frac{1}{\sqrt k}\sum_{i=1}^k\tanh(h_i)$ is Sub-gaussian with constant sub-Gaussian norm and similar arguments hold for $\sqrt k A_7^k(\bfa h)$. Hence we   see that the above term is Sub-exponential according to lemma 2.7.7. in \citep{vershynin2018high}. First notice that
\tny{\begin{align*}
    \bb E&\lef[kC_2(A_5^2+A_7^2)A_9\rig]\leq C\frac{(k-r)^2}{k^2}\bb E\lef[\frac{1}{k}\sum_{i=1}^k\sum_{j=1}^{k}\tanh(h_i)\tanh(h_j^\prime)\rig]=O\lef(\frac{(k-r)^2}{k^2}\rig).
\end{align*}}
And analogously for some $C_3>0$:
\tny{
\begin{align*}
\bb E&\lef[kC_2A_5A_7A_9\rig]\leq C(k-r)\bb E\lef[\frac{1}{k^2}\sum_{i=1}^k\sum_{j=r+1}^{k+r}\tanh(h_i)\tanh(h_j^\prime)\rig]\\
&=\frac{C_3(k-r)}{k^2}\bl\sum_{i=1}^k\sum_{j=r,j\neq i}^{k+r}\bb E\lef[\tanh(h_i)\rig]\bb E[\tanh(h_j^\prime)]+\sum_{i=r+1}^k\bb E[\tanh^2(h_i)]\br
=O\lef(\frac{(k-r)^2}{k^2}\rig)
.\end{align*}}
The sub-exponential norm can then be estimated using the sub-Gaussian properties in \citep{vershynin2018high} as follows
\begin{align*}
    \lef\Vert\frac{1}{k}\tanh(h_i)\rig\Vert_{\psi_2}\leq\frac{1}{k\sqrt{\log 2}}
   \quad \Rightarrow \quad\Vert A_5\Vert_{\psi_2},\Vert A_7\Vert_{\psi_2}=O\lef(\frac{1}{\sqrt{k}}\rig),
\end{align*}
which further implies that $\Vert A_i^2-\bb E[A_i^2]\Vert_{\psi_1}=O\lef(\frac{1}{k}\rig)$ for $i\in\{5,7\}$ and $ \Vert A_5A_7-\bb E[A_5A_7]\Vert_{\psi_1}=O\lef(\frac{1}{k}\rig)$.
Then, by the sub-exponential property,  we have when $\frac{k-r}{k}\leq C_0 $ for some $C_0>0$  we have $C_4,C_5,C_6,C_7>0$ such that
\begin{align}\label{subquad1}
&\bb E\lef[\exp\lef(kC_1(A_5^2+A_7^2)A_9^2)\rig)\rig]\leq\bb E\lef[\exp\lef(C_4\frac{(k-r)^2}{k^2}k(A_5^2+A_7^2)\rig)\rig]= \exp\lef(\frac{(k-r)^2}{k^2}C_5\rig), \nnb\\
 &\bb E\lef[\exp\lef(kC_2A_5A_7A_9\rig)\rig]\leq\bb E[\exp\lef(C_6(k-r)(A_5^kA_7^k-\bb E[A_5^kA_7^k])\rig)]=\exp\lef(\frac{(k-r)^2}{k^2}C_7\rig)
.\end{align}
For some constant $C_8,C_9>0$, the correction term can be expanded as (here we omit $k$ and $\bfa h$ for simplicity)
\begin{align}\label{varphi04}
\bb E\lef[\exp\lef(C_8(k-r)A_5^2A_7\rig)\rig]&=\sum_{j=0}^\infty \frac{C_8^j(k-r)^j}{j!}\bb E[A_5^{2j}A_7^j]\nnb\\
&=\sum_{j\text{ is even}} \frac{C_8^j(k-r)^j}{j!}\bb E[A_5^{2j}A_7^j]+\sum_{j\text{ is odd}} \frac{C_8^j(k-r)^j}{j!}\bb E[A_5^{2j}A_7^j]\nnb\\
&=\sum_{j\text{ is even}} \frac{C_8^j(k-r)^j}{k^{3j}j!}\bb E\bigg[\bl\sum_{i=1}^k\tanh(h_i)\br^{2j}\bl\sum_{i=r+1}^{k+r}\tanh(h_i)\br^j\bigg]\nnb\\
&\leq\sum_{j\text{ is even}} \frac{(8C_8)^j(k-r)^j}{k^{3j}j!}\leq\exp\lef(\frac{C_9(k-j)}{k^3}\rig)
.\end{align}

Then we need to analyze the first order error term in the Laplace approximation. This is a complicated task where we used the results derived in \citep{barndorff1989asymptotic} Section 6. To be concise in the presentation, we need to introduce the following new sets of notations.
$B^0:=\lef(\nabla^2G_{0,k}(x_k,y_k,\bfa h)\rig)^{-1}, B^1:=\begin{bmatrix}
G^{(2)}_{1,k}(x_{1,k},\bfa h)^{-1}&0\\0&G^{(2)}_{2,k}(x_{2,k},\bfa h)^{-1}\end{bmatrix}$ and $G_{3,k}(x,y,\bfa h):=G_{1,k}(x,\bfa h)+G_{2,k}(y,\bfa h)$. Using the new notation of $U_{\mca{pqrs}}:=\frac{\pta^4U}{\pta x_{\mca p}\pta x_{\mca q}\pta x_{\mca r}\pta x_{\mca s}}$ and similar definition for other number of subscripts. We also omit $x_k,y_k,x_{1,k},y_{1,k}$ and the subscript $k$ of $G_{0,k} $ and $G_{3,k}$ here. Then the first order term in \eqref{varphi0rep} can be calculated as

\begin{align}\label{varphi03}
a_0(\bfa h)-a_1(\bfa h)-a_2(\bfa h)&=-\frac{1}{8}\sum_{\mca{p,q,r,s}\in\{1,2\}}\lef(G_{0,\mca{pqrs}}B^0_{\mca{pr}}B^0_{\mca{qs}}-G_{3,\mca{pqrs}}B^1_{\mca{pr}}B^1_{\mca{qs}}\rig)\nnb\\
&+\sum_{\mca{p,q,r,s,t,u}\in\{1,2\}}\bl G_{0,\mca{pqr}}G_{0,\mca{stu}}\lef(\frac{1}{8}B^0_{\mca{ps}}B^0_{\mca{qr}}B^0_{\mca{tu}}+\frac{1}{12}B^0_{\mca{ps}}B^0_{\mca{qt}}B^0_{\mca{ru}}\rig)\nnb\\
&-G_{3,\mca{pqr}}G_{3,\mca{stu}}\lef(\frac{1}{8}B^1_{\mca{ps}}B^1_{\mca{qr}}B^1_{\mca{tu}}+\frac{1}{12}B^1_{\mca{ps}}B^1_{\mca{qt}}B^1_{\mca{ru}}\rig)\br
.\end{align}
By the fact that any term in the above expression is upper bounded by $O\lef(\frac{k-r}{k}\rig)$ (Using the fact that $x_{1,k}-x_k=O(\frac{k-r}{k\sqrt k})\sqrt k x_k$ and the cross term has only $\frac{k-r}{r}$ bounded components) and the total number does not go with $k\to\infty$, we complete that the correction term is $O\lef(\frac{k-r}{k^2}\rig)$.

Therefore, collecting pieces in \eqref{varphi01}, \eqref{varphi02}, \eqref{subquad1}, \eqref{varphi04}, and \eqref{varphi03} and we can conclude by H\"older's inequality over the three terms in \eqref{varphi0rep} that there exists constant $C_0>0$ such that  $\forall \frac{k-r}{k} <C_0$ the following holds for some constant $C>0$
\begin{align}
\bb E\lef[\frac{\bb P_{S}(\bfa\sigma)\bb P_{S^\prime}(\bfa\sigma)}{\bb P_0(\bfa\sigma)}\rig]\leq \exp\lef(C\lef(\frac{k-r}{k}\rig)^2\rig).\label{seconcase}
\end{align}

However, noticing that by \eqref{varphi02}  we have $\bb E\lef[\frac{\bb P_{S}(\bfa\sigma)\bb P_{S^\prime}(\bfa\sigma)}{\bb P_0(\bfa\sigma)}\rig]$ to be a monotonic function of $r$, since $A_9^2$ and $A_9A_5A_7$ are both monotonic decreasing function of $r$. Then it suffices to show that for the extreme case of $r=0$  we have boundedness of $\bb E\lef[\frac{\bb P_{S}(\bfa\sigma)\bb P_{S^\prime}(\bfa\sigma)}{\bb P_0(\bfa\sigma)}\rig]$ to cover the rest of the cases not given by \eqref{seconcase}. Note that in this case  we have $A_7=A_5$ and $A_6=A_8=A_9-1$ and $G_{0,k}$  degenerates to 
\begin{align*}
    G_{3,k}(x,\bfa h):=\frac{x^2}{2}-\sum_{i=1}^k\log\cosh(\sqrt{2\theta_1}x+h_i),\thinspace -G_{3,k}^{\prime}(0,\bfa h)=\sqrt{2}A_5,\thinspace -G^{(2)}_{3,k}(0,\bfa h)=2A_6+1,
\end{align*}
and $G_{1,k}=G_{2,k}$. Denote $x_3^k=\argmin_{x}G_{3,k}$, here we notice that still by Taylor expansion  we have
\begin{align*}
    G_{3,k}^\prime(x_3^k,\bfa h)=0 =G_{3,k}^\prime(0,\bfa h)+ G_{3,k}^{(2)}(0,\bfa h)x_3^k+o_{\psi_2}\lef(\frac{1}{\sqrt k}\rig),
\end{align*}
which implies that
\begin{align*}
    x_3^k=\frac{-G_{3,k}^\prime(0,\bfa h)}{G_{3,k}^{(2)}(0,\bfa h)}+o\lef(1\rig)=-\frac{\sqrt 2A_5}{2A_6+1}+o(1).
\end{align*}
And recalling that $x_{1,k}=x_{2,k}=-\frac{A_5}{A_6}$  we arrive at
\begin{align}\label{maxterm}
    &\exp\bl-k\bl G_{3,k}(x_3^k,\bfa h)-\frac{1}{k}\sum_{i=1}^k\cosh(h_i)-G_{1,k}(x_{1,k},\bfa h)-G_{2,k}(x_{2,k},\bfa h)\br\br\nnb\\
    &=\exp\bl k\bl\frac{2A_5^2}{2A_6+1}-\frac{A_5^2}{A_6} \br+o(1)\br=\exp\lef(\frac{-kA_5^2}{(2A_6+1)A_6}+o(1)\rig).
\end{align}
Note that if $2A_6+1\leq 0$ the above term is always less than $1$ since $A_6$ is negative. This is equivalent to having  $\theta_1<\frac{1}{2\bb E[\sech^2(h)]}$.

Collecting the above pieces, we conclude that there exists constant $C>0$ such that for all $\frac{k-r}{k}\leq 1$:
\begin{align}\label{finalphi1}
\bb E\lef[\frac{\bb P_{S}(\bfa\sigma)\bb P_{S^\prime}(\bfa\sigma)}{\bb P_0(\bfa\sigma)}\rig]\leq\exp\bl C\lef(\frac{k-r}{k}\rig)^2\br.
\end{align}

Since the average can be seen as taking the expectation over overlap between two randomly picked $k$- cardinal subset of $[n]$. We introduce $v:=k-r$ to be the overlap between $S$ and $S^\prime$ and $E_k(v):=\bb E_{0}\lef[\frac{\bb P_{S}(\bfa\sigma|\bfa h)\bb P_{S^\prime}(\bfa\sigma|\bfa h)}{\bb P_{0}^2(\bfa\sigma|\bfa h)}\rig]$.  Introducing a random variable $V:=|S\cap S^\prime|$ when $S$ and $S^\prime$ are uniformly randomly picked $k$-sets in $n$ elements. In what follows we dissect the term \eqref{chisqlb} according to $v$ as:
\begin{align}\label{hkk}
    \mca H_k:&=\frac{1}{\binom{n}{k}}\sum_{v=0}^k\sum_{S^\prime:|S^\prime\cap S|=v}E_k^m\lef(v\rig)=\bb P(V=0)E_k^m(0)+\sum_{v=1}^k\bb P(V=v)E_k^m\lef(v\rig)\\&=\frac{\binom{n-k}{k}}{\binom{n}{k}}+\sum_{v=1}^{k}\frac{\binom{n-k}{k-v}\binom{k}{v}}{\binom{n}{k}}E_k^m\lef(v\rig)
.\end{align}
where in the last equality \eqref{firstcase} is used.
For the first term we note that by lemma \ref{combprob},
\begin{align}\label{limratio}
    \lim_{k\to\infty}\frac{\binom{n-k}{k}}{\binom{n}{k}}=\begin{cases}
        1&\text{ if }k=o(\sqrt n)\\
        \exp(-\lambda)&\text{ if }\lim_{k\to\infty}\frac{k^2}{n}=\lambda\\
        0&\text{ if }k=\omega(\sqrt n)
    \end{cases}
.\end{align}
Then we consider the first and the second/third case separately.
\begin{center}
\textbf{1. When $k=o(\sqrt n)$}
\end{center}
For the second term on the R.H.S. of \eqref{hkk}, given some $p=\epsilon k$ for some very small $\epsilon>0$  we have
\begin{align}\label{divdep}
    \sum_{v=1}^k\bb P(V=v)E_k^m(v)=\sum_{v=1}^p\bb P\lef(V=v\rig)E_k^m(v)+\sum_{v=p+1}^{k}\bb P(V=v) E_k^m(v).
\end{align}
We then note that by lemma \ref{combprob} there exists $C_1,C_2>1$ such that:
\tny{\begin{align}\label{therest}
\sum_{v=p}^k\bb P(V=v) E_k^m(v)\leq\sum_{v=p+1}^k\frac{1}{v!}\lef(\frac{k^2}{n}\rig)^vC_1^m\leq\sum_{v=p}^kv\lef(\frac{ek^2}{nv}\rig)^vC_1^m\leq \lef(C_2\frac{k\log k}{n}\rig)^{\epsilon k}C_1^m=o(1)
\end{align}}
when $m=o\lef(k\log n\rig)$.
Based on different limiting settings in \eqref{limratio} we analyze the value of the first term in \eqref{divdep} separately. 

For the first term in \eqref{divdep} we note that by lemma \ref{combprob} and \eqref{finalphi1}, for some constant $C_1,C_2,C_3>0$  we have
\begin{align*}
\sum_{v=1}^{p}\bb P(V=v)E_k^m\lef(v\rig)&\leq\sum_{v=1}^p\frac{1}{v!}\lef(\frac{k^2}{n}\rig)^v\lef(C_1\exp\lef(\frac{v^2}{k^2}\rig)\rig)^m \leq\sum_{v=1}^p\frac{1}{v!}\lef(\frac{k^2}{n}\rig)^v\exp\lef(C_3\frac{mv^2}{k^2}\rig)\\
&\leq\sum_{v=1}^p\frac{1}{v!}\lef(\frac{k^2}{n}\exp\lef(C_3\frac{mv}{k^2}\rig)\rig)^v\leq\exp\lef(\frac{k^2}{n}\exp\lef(C_1\frac{m}{k}\rig)\rig)-1
.\end{align*}
Therefore, when $k=o(\sqrt n)$ picking $m=o\lef( k\log n\rig)$  we have the above term is $o(1)$. Combining with \eqref{limratio} and \eqref{therest}, we show that for all $m=o\lef(k\log n\rig)$ all tests are powerless asymptotically.

\begin{center}
    \textbf{2. When $k=\Omega(\sqrt n)$}
\end{center}
Note that by \eqref{finalphi1}  we have for sufficiently large $k$ there exists $C_1>0$ such that using the standard Laplace method in lemma \ref{laplaceu1} and \eqref{laplaceu2}  we have
\begin{align}\label{riemann}
    \sum_{v=1}^{p}\bb P(V=p) E_{k}^m(v)
    &\leq\sum_{v=1}^{\epsilon k}\frac{1}{(1-\frac{p}{k})\sqrt{2\pi p}}\exp\bigg(\lef(\frac{4k}{n}-\frac{p}{n}-\log\frac{pn}{k^2}-1\rig)p\nnb\\&-\frac{2k^2}{n}-2(k-p)\log\bigg(1-\frac{p}{k}\bigg)-\frac{1}{12p+1}+o(1)\bigg)\nnb\\
    &=\int_{(\frac{1}{k},\epsilon)}\frac{\sqrt k}{(1-x)\sqrt{2\pi x}}\exp\lef(k \mca f(x)\rig)dx(1+o(1)).
\end{align}
We define $\gamma:=\frac{k}{n}$ and
\begin{align*}
    \mca f(x):&=\lef((4-x)\gamma-\log\frac{x}{\gamma}-1\rig)x-2\gamma-2(1-x)\log\lef(1-x\rig) +C_1\frac{mx^2}{k}.
\end{align*}
Applying Laplace method in lemma \ref{laplaceu1}, \ref{laplaceu2} again, we note that the derivatives can be written as 
\sm{\begin{align*}
    &\mca f^\prime(x) =  (4-2x)\gamma-\log\frac{x}{\gamma}+2\log(1-x)+\frac{2C_1mx}{k},\qquad \mca f^{(2)}(x)=-2\gamma-\frac{1}{x}-\frac{2}{1-x}+\frac{2C_1m}{k}.
\end{align*}}
Then we study the maximum of $\mca f(x)$, using Fermat's condition  we have
\sm{\begin{align*}
    \mca f^\prime(x^*)=(4-2x^*)\gamma -\log\frac{x^*}{\gamma}+2\log(1-x^*)+\frac{mx^*}{k}=0\Rightarrow x^*=\frac{2\log(1-x^*)+4\gamma-\log\frac{x^*}{\gamma}}{-\frac{m}{k}+2\gamma}
\end{align*}}
which admits the maximum point $x^*=\gamma$ given $m=O(k)$. Therefore we use the Laplace method to conclude that when $m=o\lef(\frac{n^2}{k^2}\rig)$ for $k\gtrsim n^{\frac{2}{3}}$  we have
\begin{align*}
    \sum_{v=1}^{\epsilon k}\bb P(V=v)E_k^m(v)=\frac{1}{(1-\gamma)}\exp(\mca f(\gamma))\to 1.
\end{align*}
And for $k=o\lef(n^{2/3}\rig)$ we check that the maximum point $x^*$ is obtained at $c_0\in(0,1)$. This implies that we can have $m=o\lef(k\log n\rig)$ to make the sum given by \eqref{riemann} converge to $1$. Collecting pieces the sample complexity lower bound is $o(k\log n)$ for $k=o(n^{2/3})$ and $m=o(\frac{n^2}{k^2})$ for $k=\Omega(n^{2/3})$.

\begin{center}
    \textbf{ Case II: Divergence}
\end{center}
We note that when $\theta_1\in\lef[\frac{1}{2\bb E[\sech^2(h)]},\frac{1}{\bb E[\sech^2(h)]}\rig)$ we rely on the first condition in \eqref{caseI} to guarantee the positivity of \eqref{hessianhe}. This makes the chi-square method overly optimistic. The reason comes from the fact that the small probability event contributes too much to the chi-square. Here we present a new method to sharpen this result. Going back to the TV distance, we define the following event for some $c_0>0$, pick an $\epsilon$ such that $\epsilon=\omega(1)$ and $\epsilon=o(\log k)$,  we have
\begin{align}\label{esgoodset}
    E_S:=\lef\{|m_S|\leq c_0\sqrt{\frac{1}{k}\log( (m\vee k)\epsilon)}\rig\},\qquad \bb P^\prime_S(\bfa\sigma)=\begin{cases}
        \bb P_S(\bfa\sigma)&\text{ if }\bfa\sigma\in E_S\\
        0&\text{ otherwise }
    \end{cases}.
\end{align}
And analogously we define $\bb P^\prime_{\bar S}(\bfa\sigma)=\frac{1}{\binom{n}{k}}\sum_{S\subset [n]:|S|=k}\bb P_{S}^\prime(\bfa\sigma)$.
And it is simply checked that for the mixture measure $\bb P_{\bar S}$  we have by $\Vert\sqrt km_S\Vert_{\psi_2}<\infty$,
\begin{align*}
    \Vert\bb P_{\bar S}-\bb P_{\bar S}^\prime\Vert_{TV}&=\int|\bb P_{\bar S}(\bfa\sigma)-\bb P_{\bar S}^\prime(\bfa\sigma)|d\mu(\bfa\sigma)\leq \frac{1}{\binom{n}{k}}\sum_{S:|S|=k}\Vert \bb P_{S}-\bb P_{S}^\prime\Vert_{TV}\\
    &=\frac{1}{\binom{n}{k}}\sum_{S:|S|=k}\bb P_S(E_S^c)=O\lef(\frac{1}{(m\vee k)\epsilon}\rig).
\end{align*}
Therefore, by the tensorization of TV distance in \citep{donoho1991geometrizing},  we have for all $\delta\in(0,1)$
\begin{align*}
    \Vert\bb P_{\bar S}-\bb P^\prime_{\bar S}\Vert_{TV}\leq 1-\lef(\frac{1-\delta^2}{2}\rig)^{1/m}\qquad\Rightarrow\qquad \Vert \bar{\bb P}_m-\bar{\bb P}_m^\prime\Vert_{TV}\leq\delta.
\end{align*}
Hence  we have
\begin{align*}
    \Vert\bar{\bb P}_{m}-\bar{\bb P}^\prime_{m}\Vert_{TV}=O\lef(\frac{1}{\sqrt{\epsilon}}\rig)=o(1).
\end{align*}

Therefore, introducing $\bar{\bb P}_{m}^\prime=\bb P^{\prime\otimes k}_{\bar S}$ if we manage to upper bound
$    \bb E[D_{\chi^2}(\bar{\bb P}^\prime_m,\bb P_{0,m}|\bfa h^m)]
$, by triangle inequality of TV distance we then have
\begin{align*}
    \Vert\bar{\bb P}_{m} -\bar{\bb P}_{0,m}\Vert_{TV}&\leq \Vert\bar{\bb P}_{m}-\bar{\bb P}^\prime_{m}\Vert_{TV}+\Vert\bar{\bb P}_{m}^\prime-\bb P_{0,m}\Vert_{TV}\\
    &\leq\Vert\bar{\bb P}_{m}-\bar{\bb P}^\prime_{m}\Vert_{TV}+\sqrt{\bb E[D_{\chi^2}(\bb P^\prime_m,\bb P_{0,m}|\bfa h^m)]},\\
    \bb E[D_{\chi^2}(\bb P^\prime_m,\bb P_{0,m}|\bfa h^m)]&=\frac{1}{\binom{n}{k}^2}\sum_{S,S^\prime, S:|S|=|S^\prime|=k}\bb E\lef[\int\frac{\bb P^\prime_{S}\bb P^\prime_{S^\prime}}{\bb P_0}d\mu(\bfa\sigma)\rig]^m-2\int\bar{\bb P}_{m}d\mu(\bfa\sigma)+1\\
    &=\frac{1}{\binom{n}{k}^2}\sum_{S,S^\prime, S:|S|=|S^\prime|=k}\bb E\lef[\int\frac{\bb P^\prime_{S}\bb P^\prime_{S^\prime}}{\bb P_0}d\mu(\bfa\sigma)\rig]^m-1+o(1).
\end{align*}
Then we study the following decoupled quantity:
\tny{\begin{align}\label{psi00}
   \bb E&\lef[\frac{\bb P^\prime_{S}(\bfa\sigma)\bb P^\prime_{S^\prime}(\bfa\sigma)}{\bb P_0(\bfa\sigma)}\rig]:\nnb\\
   &=\bb E\lef[\frac{\sum_{\bfa\sigma:E_S\cap E_{S^\prime}}\exp\lef(\frac{\theta_1k}{2}(m_S^2+m_{S^\prime}^2)+\sum_{i\in[k+r]}\sigma_ih_i\rig)\sum_{\bfa\sigma}\exp\lef(\sum_{i\in[k+r]}\sigma_ih_i\rig)}{\lef(\sum_{\bfa\sigma}\exp\lef(\frac{\theta_1k}{2}m_S^2+\sum_{i\in[k+r]}h_i\sigma_i\rig)\rig)\lef(\sum_{\bfa\sigma}\exp\lef(\frac{\theta_1k}{2}m_{S^\prime}^2+\sum_{i\in[k+r]}h_i\sigma_i\rig)\rig)}\rig]
\end{align}}

Then we analyze the denominator and numerator of $\psi(0)$ separately. Here we introduce a new measure
\begin{align*}
\rho(m_S=a,m_{S^\prime}=b|\bfa h)=\frac{\sum_{\bfa\sigma:m_S=a,m_{S^\prime}=b}\exp(\sum_{i\in[k+r]}h_i\sigma_i)}{\sum_{\bfa\sigma}\exp(\sum_{i\in[k+r]}h_i\sigma_i)}. 
\end{align*}
Then there exists $C>0$ such that for all $t\geq C_1\sqrt{k\log(k\vee m\log k)}$  we have
\begin{align*}
    \rho(|m_S-\bb E[m_S|\bfa h]|&\geq t)\leq\exp(-Ckt^2), \\
    \rho(|m_S-\bb E[m_S|\bfa h]|&\geq t_2,|m_S^\prime-\bb E[m_S^\prime|\bfa h]|\geq t_1)\leq\exp(-Ckt_1^2-Ckt_2^2).
\end{align*}
Moreover we recall that under $\rho$, $\bb E[m_S|\bfa h]=\frac{1}{k}\sum_{i=1}^k\tanh(h_i)$ and $\bb E[m_S^\prime|\bfa h]=\frac{1}{k}\sum_{i=r+1}^{r+k}\tanh(h_i)$. Then  we have by their sub-Gaussian property there exists $C_2>0$ such that
\begin{align*}
    \rho(\bb E[m_S|\bfa h]\vee \bb E[m_{S^\prime}|\bfa h]\geq t)\leq\exp(-C_2kt^2).
\end{align*}
After introducing this measure $\psi(0)$ can be rewritten as
\begin{align*}
   \bb E\lef[\frac{\bb P^\prime_{S}(\bfa\sigma)\bb P^\prime_{S^\prime}(\bfa\sigma)}{\bb P_0(\bfa\sigma)}\rig]:&=\bb E\lef[\frac{\sum_{m_S,m_S^\prime:E_S\cap E_S^\prime}\exp\lef(\frac{\theta_1k}{2}(m_S^2+m_{S^\prime}^2)\rig)\rho(m_S,m_S^\prime|\bfa h)}{\sum_{m_S}\exp\lef(\frac{\theta_1k}{2}m_S^2\rig)\rho(m_S|\bfa h)\sum_{m_S^\prime}\exp\lef(\frac{\theta_1k}{2}m_{S^\prime}^2\rig)\rho(m_S^\prime|\bfa h)}\rig].
\end{align*}
This essentially implies that the regularity conditions holds as in \citep{ellis1980limit}.
Recalling the definition of $m_S$ and $m_S^\prime$,  we have:
\tny{\begin{align*}
    &\sum_{m_S,m_S^\prime:E_S\cap E_S^\prime}\exp\lef(\frac{\theta_1k}{2}(m_S^2+m_{S^\prime}^2)\rig)\rho(m_S,m_S^\prime|\bfa h)\\
   &=\frac{k}{2\pi}\sum_{E_S,E_{S^\prime}}\int_{\bb R}\int_{\bb R}\exp\bl-\frac{k(x^2+y^2)}{2}+k\sqrt{\theta_1}(m_Sx+m_{S^\prime}y)\br \rho(dm_S,dm_{S^\prime}|\bfa h)dxdy\\
   &=\frac{k}{2\pi}\sum_{m_S,m_{S^\prime}\in[-1,1]}\int_{|x|\leq c_1}\int_{|y|\leq c_2}\exp\bl-\frac{k(x^2+y^2)}{2}+k\sqrt{\theta_1}(m_Sx+m_{S^\prime}y)\br \rho(dm_S,dm_{S^\prime}|\bfa h)dxdy\\
   &+\ub{\frac{k}{2\pi}\sum_{E_S,E_{S^\prime}}\int_{|x|\leq c_1}\int_{|y|> c_2}\exp\bl-\frac{k(x^2+y^2)}{2}+k\sqrt{\theta_1}(m_Sx+m_{S^\prime}y)\br \rho(dm_S,dm_{S^\prime}|\bfa h)dxdy}_{T_1}\\
   &+\ub{\frac{k}{2\pi}\sum_{E_S,E_{S^\prime}}\int_{|x|> c_1}\int_{|y|> c_2}\exp\bl-\frac{k(x^2+y^2)}{2}+k\sqrt{\theta_1}(m_Sx+m_{S^\prime}y)\br \rho(dm_S,dm_{S^\prime}|\bfa h)dxdy}_{T_2}\\
   &-\frac{k}{2\pi}\sum_{E_S, E^c_{S^\prime}}\int_{|x|\leq c_1}\int_{|y|\leq c_2}\exp\bl-\frac{k(x^2+y^2)}{2}+k\sqrt{\theta_1}(m_Sx+m_{S^\prime}y)\br \rho(dm_S,dm_{S^\prime}|\bfa h)dxdy \\
   &-\frac{k}{2\pi}\sum_{ E^c_S, E^c_{S^\prime}}\int_{|x|\leq c_1}\int_{|y|\leq c_2}\exp\bl-\frac{k(x^2+y^2)}{2}+k\sqrt{\theta_1}(m_Sx+m_{S^\prime}y)\br \rho(dm_S,dm_{S^\prime}|\bfa h)dxdy.
\end{align*}}
Therefore, when picking $c_1=c_2=C\sqrt{\frac{1}{k}\log( (m\vee k)\epsilon)}$ for some proper $C>0$,  we have
\begin{align*}
   T_1,T_2= O\lef(\frac{1}{(m\vee k)\epsilon}\rig)=o\lef(\frac{1}{m}\rig).
\end{align*}
Then  we have a corresponding form of \eqref{vphi0g}:
\tny{\begin{align}\label{boundingprocedure}
   \bb E\lef[\frac{\bb P^\prime_{S}(\bfa\sigma)\bb P^\prime_{S^\prime}(\bfa\sigma)}{\bb P_0(\bfa\sigma)}\rig]\leq\bb E\lef[\frac{\prod_{i=r+1}^k\cosh(h_i)\int_{|x|\vee|y|\leq c_1}\exp\lef(-kG_{0,k}(x,y,\bfa h)\rig)dxdy }{\int_{\bb R}\exp\lef(-kG_{1,k}(x,\bfa h)\rig)dx\int_{\bb R}\exp\lef(-kG_{2,k}(y,\bfa h)\rig)dy}\rig]+ o\lef(\frac{1}{m}\rig)
\end{align}}
It is   checked that for $c>\frac{-(1-\theta_1^2\bb E[\sech^2(h)]^2-\theta_1^2\bb V(\sech^2(h)))}{\theta_1(\theta_1\bb E[\sech^2(h)]^2-\bb E[\sech^2(h)]+\theta_1\bb V(\sech^2(h)}$ $\bb E\lef[\frac{\bb P_{S}(\bfa\sigma)\bb P_{S^\prime}(\bfa\sigma)}{\bb P_0(\bfa\sigma)}\rig]$ with the new measure is the same as Case I. The only difference is that
when considering $c< \frac{-(1-\theta_1^2\bb E[\sech^2(h)]^2-\theta_1^2\bb V(\sech^2(h)))}{\theta_1(\theta_1\bb E[\sech^2(h)]^2-\bb E[\sech^2(h)]+\theta_1\bb V(\sech^2(h)}$, we notice that here the maximum of $G_{0,k}$ is taken at the boundary points. Then by Laplace method \ref{laplaceu2}, there exists $C>0$ such that
\begin{align*}
   \bb E\lef[\frac{\bb P^\prime_{S}(\bfa\sigma)\bb P^\prime_{S^\prime}(\bfa\sigma)}{\bb P_0(\bfa\sigma)}\rig]\leq \exp\lef(C(1-c)^2\log\lef((m\vee k)\log k\rig)\rig).
\end{align*}
Going back to the chi-square divergence, we notice that the total contribution of this divergence region can be bounded as
\begin{align}\label{reduceepsk}
    \sum_{v=0}^{(1-c)k}\sum_{S^\prime:|S^\prime\cap S|=v}&\frac{1}{\binom{n}{k}}\exp\lef(Cm\lef(\frac{v}{k}\rig)^2\log\lef((m\vee k)\log k\rig)\rig)\nnb\\
    &\leq k\frac{\binom{k}{ck}\binom{n-k}{(1-c)k}}{\binom{n}{k}}\exp\lef(Cmc^2\log\lef((m\vee k)\epsilon \rig)\rig)\nnb\\
    &\leq k\lef(\frac{cn}{k}\rig)^{-ck}\exp\lef(Cmc^2\log\lef((m\vee k)\epsilon\rig)\rig)=o(1).
\end{align}
Therefore if we ask $m=o\lef( \frac{k}{\log k}\log n\rig)$, the above inequality holds. Notice that the loss of $\log k$ term only occured at the previous $k=o(n^{2/3})$ region since the maximum is taken at the $c\asymp 1$ region. The previous region of $m=\Omega\lef(n^{2/3}\rig)$ (where the maximum is taken with $c=\frac{k}{n}$ is not affected by this loss of $\log k$ factor. 

\subsection{Proof of Theorem \ref{thm1}}
Our proof goes by analyzing the two parts of the algorithms separately. The local part proves the guarantee of the $k=o(n^{2/3})$ region of the test. The global part proves the guarantee of the $k\gtrsim n^{2/3}$ region of the test.
\begin{center}
    \textbf{1. Local Part}
\end{center}
The proof goes along by considering the two parts in the inequality separately and bound them individually. Before we start the proof, the following properties are needed, which is a direct result of the central limit theorem \ref{cltrfcr}
\begin{lemma}\label{tailh0}
 Define $W:=\frac{2}{k}\sum_{1\leq i <j\leq k}\sigma_i\sigma_j,$ then we have $\Vert W\Vert_{\psi_1}\asymp 1$ under the null hypothesis.
\end{lemma}
\begin{lemma}\label{h1tail}
    We use the same notations as in lemma \ref{tailh0} and assume that $\theta_1< \theta_c$. Then,  under the alternative hypothesis with the index of clique defined by $[k]$  we have $\Vert W\Vert_{\psi_1}\asymp 1$.
\end{lemma}
With the above preparation, we can obtain the upper bound of the local test. Here we define $W_S:=\frac{1}{k}\sum_{i,j\in S}\sigma_i\sigma_j$ and $\{W_S^{(i)}\}_{i\in[n]}$ to be $n$ independent copies of $W_{S}$. We denote $\bb P_0$ to be the measure under the null hypothesis, given lemma \ref{tailh0}, lemma \ref{regularbernstein} and union bound we check that the type I error can be upper bounded as
\begin{align}\label{h0uni}
    \bb P_0\lef(\sup_{S\subset [n]:|S|=k}\phi_S\geq\tau_\delta\rig) &\leq \sum_{S\subset[n]:|S|=k}\bb P_0\lef(\frac{1}{m}\sum_{i=1}^m W_S^{(i)}\geq\frac{\tau_\delta}{2}\rig)\nnb\\
    &\leq \lef(\frac{en}{k}\rig)^k\exp\lef(-c\lef(m\tau_\delta^2\wedge \tau_\delta m\rig)\rig)
.\end{align}
Asking the R.H.S. to be less than $\delta/2$ we have
\begin{align*}
    \tau_\delta \geq C_0\sqrt{\frac{k\log\frac{en}{k}-\log(\delta/2)}{m}}\vee \frac{k\log\lef(\frac{en}{k}\rig)-\log(\delta/2)}{m}
.\end{align*}

Then we prove upper bound for the Type II error. To simply notation we introduce $\bb P_{S}$ as the probability measure under alternative hypothesis.
Then, using lemma \ref{h1tail} we conclude that there exists $C>0$ such that
\begin{align}\label{alternative}
    \bb P_{S_0}\lef(\sup_{S:|S|=k}\phi_S\leq\tau_\delta\rig)
    &\leq \bb P_\hks\lef(\frac{1}{m}\sum_{i=1}^m (W_{S_0}^{(i)}-\bb E[W_{S_0}])\leq\frac{1}{2}\tau_\delta - \bb E[W_{S_0}]\rig)\nnb\\
    &\leq \exp\lef(-Cm\lef(\lef(\tau_\delta-\bb E[W_{S_0}]\rig)^2\vee (\tau_\delta-\bb E[W_{S_0}])\rig)\rig)
.\end{align}
Let the R.H.S. be less than $\delta/2$ we conclude that
\begin{align*}
    \bb E[W_{S_0}]-\tau_\delta\geq C_1\sqrt{\frac{-\log\delta/2}{m}}\vee \frac{-\log\delta/2}{m}
.\end{align*}
Hence, for some constant $C_0,C_1$ there exists an interval for $\tau_\delta$ controlling the sum of Type-I and Type-II error to be less than $\delta$ defined by
\tny{\begin{align}\label{val_inter}
    \tau_\delta\in\lef(C_0\sqrt{\frac{k\log\frac{en}{k}-\log(\delta/2)}{m}}\vee \frac{k\log\lef(\frac{en}{k}\rig)-\log(\delta/2)}{m} ,\bb E[W_{S_0}]-C_1\sqrt{\frac{-
    \log\delta/2}{m}}\vee \frac{-\log\delta/2}{m}\rig)
.\end{align}}
A crucial quantity that we need to find out is the order of $\bb E[W_{S_0}]$. This can be directly calculated using theorem \ref{cltrfcr} as
\begin{align}\label{expectationofclique}
   \bb E[W_{S_0}] 
   &=\frac{1-\theta_1\bb E[\sech^2(h)]^2}{(1-\theta_1\bb E[\sech^2(h)])^2}+o(1)
.\end{align}
Therefore we ask that the interval \eqref{val_inter} exists as long as we have
$
 m\gtrsim k\log\lef(\frac{n}{k}\rig).
$
\begin{center}
\textbf{2. Global Part}
\end{center}
First we notice that at $\mca H_0$  we have $\sigma_i$ are i.i.d. Rademacher random variables, then by the composition of i.i.d. sub-Gaussian random variables,  we have
$\lef\Vert\frac{1}{\sqrt n}\sum_{i=1}^n\sigma_i\rig\Vert_{\psi_2}<\infty$. This further implies that $\lef\Vert\frac{1}{n}(\sum_{i=1}^n\sigma_i)^2\rig\Vert_{\psi_1}<\infty$ by the lemma \ref{squaresubg}. Then we use Bernstein's inequality ( lemma \ref{regularbernstein} ) to note that there exists $C>0$ such that for all $t>0$:
\begin{align*}
    \bb P_0\lef(\phi_2\geq t\rig )\leq 2\exp\lef(-Ct^2m\frac{k^2}{n^2}\wedge tm\frac{k}{n}\rig).
\end{align*}
Then we consider the alternative. Assume that $S_0$ is the set of indices in the clique, then  we have by \eqref{expectationofclique}:
\sm{\begin{align*}
    \bb E[\phi_2]=\frac{1-\theta_1(\bb E[\sech^2(h)])^2}{(1-\theta_1\bb E[\sech^2(h)])^2}-1+o(1)=\frac{2\theta_1\bb E[\sech^2(h)]-\theta_1(1+\theta_1)(\bb E[\sech^2(h)])^2}{(1-\theta_1\bb E[\sech^2(h)])^2}+o(1).
\end{align*}}
And by lemma \ref{tailh0}, and $\frac{1}{\sqrt n}\sum_{i\in S^c}\sigma_i\perp \frac{1}{\sqrt k}\sum_{i\in S}\sigma_i$,  we have
$    \bigg\Vert\bigg(\frac{1}{\sqrt n}\sum_{i=1}^n\sigma_i\bigg)^2\bigg\Vert_{\psi_1}<\infty.
$
And  we have by Bernstein's inequality, there exists $C>0$ such that:
\begin{align*}
    \bb P_{S_0}\lef(|\phi_2-\bb E[\phi_2]|\geq t\rig)\leq 2\exp\lef(-Ct^2\frac{mk^2}{n^2}\wedge t\frac{mk}{n}\rig).
\end{align*}
Therefore, we can pick $\tau_\delta\in\lef(0,\frac{2\theta_1\bb E[\sech^2(h)]-\theta_1(1+\theta_1)(\bb E[\sech^2(h)])^2}{(1-\theta_1\bb E[\sech^2(h)])^2}\rig)$ and ask $m\gtrsim\frac{n^2}{k^2}$ to complete the proof.

\subsection{Proof of Corollary \ref{hightemperaturerec}}
Using the result in \eqref{mrksmalzero} by setting $a=0$,  we have
\begin{align*}
    \bb E\lef[\exp\lef(t\sqrt{k}(1-c)m_{rk}\rig)\rig]=\exp\lef(\frac{V(c)}{2}t^2\rig)(1+o(1)),
\end{align*}
with $V(c):=(1-c)\frac{1-\theta_1(\bb E[\sech^2(h)])^2}{(1-\theta_1\bb E[\sech^2(h)])^2}$.

And we can  see that $V(c)+c$ is a monotonic decreasing function of $c$. Then we use the fact that $m_r^\prime$ is average of $k-r$ i.i.d. Rademacher r.v.s. to get
\begin{align*}
    \bb E[\exp(t\sqrt km_{S^\prime})]=\exp\lef(\frac{t^2}{2}(V(c)+c)\rig)(1+o(1)).
\end{align*}
    Recall that our test statistics $\phi_S = \frac{1}{k}\lef(\mbbm 1_{S}^\top\wh {\bb E}[\bfa\sigma\bfa\sigma^\top]\mbbm 1_{S} - k\rig)= (k-1)m_{S}^2-1$.
    Then we can use Chi-square Tail bound in \eqref{chisqtail} to conclude that for $x>0$, there exists $C_1,C_2,C_3,C_4>0$ such that for all $t>0$:
    \tny{\begin{align*}
        \bb P\lef(\phi_{S^\prime}-\bb E[\phi_{S^\prime}]\geq t\rig)\leq\exp\lef(-C_1 mt\wedge C_2mt^2\rig),\quad\bb P\lef(\phi_{S^\prime}-\bb E[\phi_{S^\prime}]\leq -t\rig)\leq\exp\lef(-C_3mt\wedge C_4mt^2\rig).
    \end{align*}}
    And moreover we notice that moment convergence implies that
    \begin{align*}
        \bb E\lef[\phi_{S^\prime}\rig]= V(c)+c-1+o(1).
    \end{align*}
    Here we introduce the notation $S_c:=[ck+1:ck+k]$, it is   checked that $S=S_0$.
    Therefore, by union bound and the above discussion over the tail of $\phi_{S_c}$ that the following holds with some constant $C_1, C_2,C_3, C_4>0$:
    \tny{\begin{align}\label{overlapprob}
        \bb P&\lef(|S_0\Delta S_{\max}|\leq k\epsilon\rig)\geq\bb P\lef(\lef|\phi_{S_0}-\bb E[\phi_{S_0}]\rig|\leq\delta_1,\forall S^\prime\Delta S\geq k\epsilon:\lef|\phi_{S^\prime}-\bb E[\phi_{S^\prime}]\rig|\leq \bb E[\phi_{S_0}]-\bb E[\phi_{S_\epsilon}]-\delta_1\rig)\nnb\\
        &=1-\bb P\lef(\lef|\phi_{S_0}-\bb E[\phi_{S_0}]\rig|>\delta_1\text{ or }\exists S^\prime\Delta S\geq k\epsilon:\lef|\phi_{S^\prime}-\bb E[\phi_{S^\prime}]\rig|>\bb E[\phi_{S_0}]-\bb E[\phi_{S_\epsilon}]-\delta_1\rig)\nnb\\
        &\geq 1-\bb P\lef(-C_1m\rig)-\sum_{p=\epsilon k/2}^{k}\binom{n-k}{k-p}\binom{k}{p}\bb P\lef(-C_2m\rig)\nnb\\
        &=1-\frac{1+\sum_{p=\epsilon k/2}^{k}\binom{n-k}{k-p}\binom{k}{p}}{\binom{n}{k}}\exp\lef(-Ck\log n\rig)\nnb\\
        &=1-o(1),
    \end{align}}
    where we used the upper tail bound in lemma \ref{combprob} and the fact that $m\gtrsim k\log n$. 
\subsection{Proof of Theorem  \ref{lbltsmallclique}}
    The proof goes similarly as that of theorem \ref{minimaxlb}.  
    Since the analysis of $r=k$ is identical, we only discuss over the $r<k$ case, where similar integration identity is applied. Recall that
    \begin{align}\label{vphi0glt}
    \bb E\lef[\frac{\bb P_{S}(\bfa\sigma)\bb P_{S^\prime}(\bfa\sigma)}{\bb P_0(\bfa\sigma)}\rig]=\bb E\lef[\frac{\prod_{i=r+1}^k\cosh(h_i)\int\exp\lef(-kG_{0,k}(x,y,\bfa h)\rig)dxdy }{\int\exp\lef(-kG_{1,k}(x,\bfa h)\rig)dx\int\exp\lef(-kG_{2,k}(y,\bfa h)\rig)dy}\rig]
    .\end{align}
    And we recall that uniformly almost surely:
    \begin{align*}
G_{0,k}(x,y,\bfa h):&=-\frac{1}{k}\bl\sum_{i=1}^{r}\log\cosh\lef(\sqrt{\theta_1}x+h_i\rig)+\sum_{i=r+1}^k\log\cosh(\sqrt{\theta_1}(x+y)+h_i)\\
&+\sum_{i=k+1}^{k+r}\log\cosh(\sqrt{\theta_1}y+h_i)\br+\frac{x^2+y^2}{2}\to G_0(x,y),
\end{align*}
and uniformly,
\begin{align*}
    G_{1,k}(x,\bfa h):&=\frac{x^2}{2}-\frac{1}{k}\sum_{i=1}^k\log\cosh\lef(\sqrt{\theta_1}x+h_i\rig)\to G_1(x),\\
G_{2,k}(x,\bfa h):&=\frac{x^2}{2}-\frac{1}{k}\sum_{i=r+1}^{k+r}\log\cosh\lef(\sqrt{\theta_1}x+h_i\rig)\to G_2(x)
.\end{align*}
Note that in the low temperature regime, the function $G_{0}$ has four nonzero stationary point defined by $(x_i^*,y_i^*)$ with $i\in [4]$ and function $G_{1}$, $G_{2}$ both have two nonzero stationary  points $x^*_{1,i}, y^*_{1,i}$ for $i\in[2]$ respectively. Then by uniform convergence we know that there exists a sequence of stationary points of $G_{0,k}$, $G_{1,k}$, $G_{2,k}$ converging towards their respective population variety. These converging stationary points are denoted by $(x^{(k)}_i,y^{(k)}_i)$, $x^{(k)}_{1,i}$ and $x^{(k)}_{2,i}$ respectively with $k\in\bb N$.  Introduce $c:=\lim_{k\to\infty}\frac{r}{k}$ as we did other places. By Fermat's condition  we have
\tny{\begin{align}\label{g0derivalt}
    \nabla G_{0}(x_i^*,y_i^*,\bfa h) = \begin{bmatrix}
        x_i^*-c\sqrt{\theta_1}\bb E\tanh(\sqrt{\theta_1}x_i^*+h)-(1-c)\sqrt{\theta_1}\bb E\tanh(\sqrt{\theta_1}(x_i^*+y_i^*)+h)
        \\ y_i^*-(1-c)\sqrt{\theta_1}\bb E\tanh(\sqrt{\theta_1}(x_i^*+y_i^*)+h)-c\sqrt{\theta_1}\tanh(\sqrt{\theta_1}y_i^*+h)
    \end{bmatrix}=\bfa 0.
\end{align}}
By the derivative being an even function, we see that the four nonzero solutions forms two pair which we denoted by $(x_1^*,y_1^*)>0$, $(x_2^*,y_2^*)=-(x_1^*,y_1^*)$, and $(x_3^*>0,y_3^*<0)$, $(x_4^*,y_4^*)=-(x_3^*,y_3^*)$. Using the fact that $G_{0,k}$ is monotonic in $x+y$ we know that $(x_1,y_1)$ and $(x_2^*,y_2^*)$ are the global minimum. Since the analysis of $(x^*_2,y^*_2)$ (and the converging sequence of it) are identical with the analyze of $(x_1^*,y_1^*)$, we only analyze the sequence converging to $(x^*_1,y^*_1)$. For $G_{1}$ and $G_{2}$, it is analogously seen that  we have two they also have this symmetry and we assume that $x_{1,1}^*,y_{1,1}^*$ are the positive ones. We can also checked that when $c=1$  we have $(x^*_i,y^*_i)=(x_{1,i}^*,y_{1,i}^*)$ for $i\in[2]$.

In what follows, for notation simplicity we omit the $k$ subscript in all intermediate stationary points and for example write $(x_1,y_1):=(x_i^{(k)},y_i^{(k)})$ to simplify notations.  For the second order derivatives, we introduce $\nabla^2(G_{0,k}(x_1^*,y_1^*,\bfa h))=\bfa B:=\begin{bmatrix}
    B_{xx}&B_{xy}\\
    B_{xy}&B_{yy}
\end{bmatrix}$ and:
\ttny{\begin{align*}
    B_{xx}:=&\nabla^2_{xx}(G_{0,k}(x_1^*,y_1^*,\bfa h)) = 1-\frac{\theta_1}{k}\sum_{i=1}^r\sech^2(\sqrt{\theta_1}x_1^*+h_i)-\frac{\theta_1}{k}\sum_{i=r+1}^k\sech^2(\sqrt{\theta_1}(x_1^*+y_1^*)+h_i), \\
    B_{xy}:=&\nabla^2_{xy}(G_{0,k}(x_1^*,y_1^*,\bfa h)) =  -\frac{\theta_1}{k}\sum_{i=r+1}^k\sech^2(\sqrt{\theta_1}(x_1^*+y_1^*)+h_i), \\
    B_{yy}:=&\nabla^2_{yy}(G_{0,k}(x_1^*,y_1^*,\bfa h)) =1-\frac{\theta_1}{k}\sum_{i=r+1}^{k}\sech^2(\sqrt{\theta_1}(x_1^*+y_1^*)+h_i)-\frac{\theta_1}{k}\sum_{i=k+1}^{k+r}\sech^2(\sqrt{\theta_1}y_1^*+h_i),\\
    B_1:=&G_{1,k}^{(2)}(x_{1,1}^*,\bfa h)=1-\frac{\theta_1}{k}\sum_{i=1}^k\sech^2(\sqrt{\theta_1}x_{1,1}^*+h_i),\\ B_2:=&G_{2,k}^{(2)}(x_{2,1}^*,\bfa h)=1-\frac{\theta_1}{k}\sum_{i=r+1}^{k+r}\sech^2(\sqrt{\theta_1}y_{1,1}^*+h_i).
\end{align*}}
And analogously define 
\begin{align*}
   \begin{bmatrix}
       A_x\\
       A_y
   \end{bmatrix} :&=\begin{bmatrix}
       \frac{\pta G_{0,k}(x_1^*,y_1^*,\bfa h)}{\pta x}\\
       \frac{\pta G_{0,k}(x_1^*,y_1^*,\bfa h)}{\pta y}
   \end{bmatrix} \\
   &=\begin{bmatrix}
        x_1^*-\frac{\sqrt{\theta_1}}{k}\sum_{i=1}^r\tanh(\sqrt{\theta_1}x_1^*+h_i)-\frac{\sqrt{\theta_1}}{k}\sum_{i=r+1}^k\tanh(\sqrt{\theta_1}(x_1^*+y_1^*)+h_i)
        \\ y_1^*-\frac{\sqrt{\theta_1}}{k}\sum_{i=r+1}^{k}\tanh(\sqrt{\theta_1}(x_1^*+y_1^*)+h_1^\prime)-\frac{\sqrt{\theta_1}}{k}\sum_{i=k+1}^{k+r}\tanh(\sqrt{\theta_1}y_1^*+h_i)
    \end{bmatrix},\\
    A_1:&= G^\prime_{1,k}(x_{1,1}^*) = x_{1,1}^*-\frac{\sqrt{\theta_1}}{k}\sum_{i=1}^k\tanh(\sqrt{\theta_1}x_{1,1}^*+h_i), \\A_2:&= G^\prime_{2,k}(y_{1,1}^*)=y_{1,1}^*-\frac{\sqrt{\theta_1}}{k}\sum_{i=1}^k\tanh(\sqrt{\theta_1}y_{1,1}^*+h_i).
\end{align*}
Therefore together with \eqref{g0derivalt} we see that $\sqrt k A_x,\sqrt k A_y, \sqrt k A_1,\sqrt{k} A_2$ are converging to Gaussian.
Before we start analyzing the desired quantity $\bb E\lef[\frac{\bb P_{S}(\bfa\sigma)\bb P_{S^\prime}(\bfa\sigma)}{\bb P_0(\bfa\sigma)}\rig]$, we first analyze the asymptotic distribution of intermediate stationary points. By Taylor expansion and the Fermat's condition  we have:
\begin{align*}
    &\begin{bmatrix}
        \frac{\pta G_{0,k}(x_1,y_1,\bfa h)}{\pta x}\\
        \frac{\pta G_{0,k}(x_1,y_1,\bfa h)}{\pta y}
    \end{bmatrix} =\bfa 0 =\begin{bmatrix}
        A_x\\
        A_y
    \end{bmatrix}+\begin{bmatrix}
        B_{xx} & B_{xy}\\
        B_{xy} & B_{yy}
    \end{bmatrix}\cdot\begin{bmatrix}
        x_1-x_1^*\\
        y_1-y_1^*
    \end{bmatrix}+o_{\psi_2}\lef(x_1-x_1^*\rig),\\
    &\begin{bmatrix}
        \frac{\pta G_{1,k}(x_{1,1},\bfa h)}{\pta x}\\
        \frac{\pta G_{2,k}(y_{1,1},\bfa h)}{\pta y}
    \end{bmatrix}=\bfa 0=\begin{bmatrix}
        A_1\\
        A_2
    \end{bmatrix} +\begin{bmatrix}
        B_1&0\\
        0&B_2
    \end{bmatrix}\cdot\begin{bmatrix}
        x_{1,1}-x_{1,1}^*\\
        y_{1,1}-y_{1,1}^*
    \end{bmatrix}+o_{\psi_2}(x_{1,1}-x_{1,1}^*).
\end{align*}
And therefore we note that
\begin{align*}
    \sqrt k\begin{bmatrix}
        x_1-x_1^*\\
        y_1-y_1^*
    \end{bmatrix} &=-\sqrt k\begin{bmatrix}
        B_{xx} & B_{xy}\\
        B_{xy} & B_{yy}
    \end{bmatrix}^{-1}\begin{bmatrix}
        A_x\\
        A_y
    \end{bmatrix}+o_{\psi_2}(1),\\
    \sqrt k\begin{bmatrix}
        x_{1,1}-x_{1,1}^*\\
        y_{1,1}-y_{1,1}^*
    \end{bmatrix}&=-\sqrt k\begin{bmatrix}
        B_1^{-1}A_1\\
        B_2^{-1}A_2
    \end{bmatrix}+ o_{\psi_2}(1).
\end{align*}
And we also noticed that by mean value theorem and the fact that $G_{0,k}$ is every infinitely differentiable in $\bb R^2$, there exists $x_0\in[x_1\wedge x_{1,1}, x_1\vee x_{1,1}]$ such that
\begin{align*}
    \bfa 0=\frac{\pta G_{0,k}(x_{1,1},y_{1})-G_1(x_{1,1})}{\pta x}+\frac{\pta^2 G_{0,k}(x_0,y_1)}{\pta x^2}\lef(x_1-x_{1,1}\rig),
\end{align*}
which implies that $x_1-x_{1,1}=O\lef(\frac{k-r}{k}\rig)$. Analogously we also have $y_1-y_{1,1}=O\lef(\frac{k-r}{k}\rig)$.
Given above preparation, we apply the Laplace method in lemma \ref{laplace} (whose regularity condition is already checked in lemma \ref{convergeasG}), and noticing that by symmetry the minimum values along with the second order derivatives are identical for the two global minimum.
\begin{align*}
    \bb E&\lef[\frac{\bb P_{S}(\bfa\sigma)\bb P_{S^\prime}(\bfa\sigma)}{\bb P_0(\bfa\sigma)}\rig]=\ub{\frac{\lef(B_1B_2\rig)^{1/2}}{\det(\bfa B)^{1/2}}}_{T_1}\\&
    \cdot\ub{\exp\bl-kG_{0,k}(x_1,y_1,\bfa h)+kG_{1,k}(x_{1,1},\bfa h)+kG_{2,k}(y_{1,1},\bfa h)+\sum_{i=r+1}^k\log\cosh(h_i)\br}_{T_2}\nnb\\ &\cdot\lef(1+\frac{a(\bfa h)}{k}+O\lef(\frac{1}{k^2}\rig)\rig).
\end{align*}
where $a_1(\bfa h)$ is the correction term dependent on $\bfa h$. Then we analyze the different term separately. For the first term, we notice that by $x_1-x_{1,1}=O\lef(\frac{k-r}{k}\rig)$ 
 we have:
\begin{align*}
    B_{xx}=B_{1}+O\lef(\frac{k-r}{k}\rig),\qquad B_{yy}=B_2+O\lef(\frac{k-r}{k}\rig),\qquad B_{xy}=O\lef(\frac{k-r}{k}\rig).
\end{align*}
and we finally get:
\begin{align*}
    T_1= \frac{1}{\sqrt{\frac{B_{xx}B_{yy}}{B_1B_2}-\frac{B_{xy}^2}{B_1B_2}}}= 1+O\lef(\frac{k-r}{k}\rig).
\end{align*}
We then study the population version of quantity in exponential, by $1$-Lipschitzness of $\log\cosh$  we have
\begin{align*}
    -kG_{0,k}(x^*,y^*,\bfa h)+kG_{1,k}(x_{11}^*,\bfa h)+k G_{2,k}(y_{11}^*,\bfa h)+\sum_{i=r+1}^k\log\cosh(h_i)= O\lef(k-r\rig).
\end{align*}
And similarly  we have
\begin{align*}
    A_x-A_1=O\lef(\frac{k-r}{k}\rig),\qquad A_y-A_2=O\lef(\frac{k-r}{k}\rig).
\end{align*}
And for the second term  we have by Taylor expansion, for some $C>0$:
\tny{\begin{align*}
    T_2&=\exp\bl -k\bl\begin{bmatrix}
        A_x\\
        A_y
    \end{bmatrix}\begin{bmatrix}
        x_1-x_1^*\\
        y_1-y_1^*
    \end{bmatrix}+\frac{1}{2}\begin{bmatrix}
        x_1-x_1^*\\
        y_1-y_1^*
    \end{bmatrix}^\top\begin{bmatrix}
        B_{xx}&B_{xy}\\
        B_{xy}&B_{yy}
    \end{bmatrix}\begin{bmatrix}
        x_1-x_1^*\\
        y_1-y_1^*
    \end{bmatrix}-\begin{bmatrix}
        A_1\\
        A_2
    \end{bmatrix}\begin{bmatrix}
        x_{1,1}-x_{1,1}^*\\
        y_{1,1}-y_{1,1}^*
    \end{bmatrix}\\
    &+\frac{1}{2}\begin{bmatrix}
        x_{1,1}-x_{1,1}^*\\
        y_{1,1}-y_{1,1}^*
    \end{bmatrix}^\top\begin{bmatrix}
        B_{1}&0\\
        0&B_{2}
    \end{bmatrix}\begin{bmatrix}
        x_{1,1}-x_{1,1}^*\\
        y_{1,1}-y_{1,1}^*
    \end{bmatrix}\br +o(k-r)\br\\
    &=\exp\lef(-\frac{k}{2}\begin{bmatrix}
        A_x\\
        A_y
    \end{bmatrix}^\top\begin{bmatrix}
        B_{xx}&B_{xy}\\
        B_{xy}&B_{yy}
    \end{bmatrix}^{-1}\begin{bmatrix}
        A_x\\
        A_y
    \end{bmatrix} + \frac{k}{2}\begin{bmatrix}
        A_1\\
        A_2
    \end{bmatrix}^\top\begin{bmatrix}
        B_{1}&0\\
        0&B_{2}
    \end{bmatrix}^{-1}\begin{bmatrix}
        A_1\\
        A_2
    \end{bmatrix}+o(k-r)\rig).
\end{align*}}
Combining pieces and noticing that $a(\bfa h)=O(1)$, we conclude that there exists $C>0$ such that
\begin{align*}
    \bb E\lef[\frac{\bb P_{S}(\bfa\sigma)\bb P_{S^\prime}(\bfa\sigma)}{\bb P_0(\bfa\sigma)}\rig]\leq \exp\lef(C(k-r)\rig).
\end{align*}

And we can go back to the decomposition in \eqref{hkk}. Recall that we denote $V$ to be the random overlap to a specified $k$-element subset $S$ $[n]$ when sample uniformly at random from $[n]$ another $k$-element subset $S^\prime$. This corresponds to the quantity of $k-r$. Therefore, in the region of $k=o(\sqrt n)$, the only thing needs to be proved is under some $m$ $\sum_{v=1}^k\bb P(V=v)E_k^m(v)\to 0$ to finish the proof.

In this proof we still utilize the second approximation result given in lemma \ref{combprob} to get:
\begin{align*}
    \sum_{v=1}^{k-1}&\bb P(V=v)E_k^m(v)\leq\sum_{v=1}^{k-1}\frac{1}{\lef(1-\frac{v}{k}\rig)\sqrt{2\pi v}}\exp\lef(-v\log\frac{vn}{k^2}-\frac{2v^2}{k}-\frac{1}{12v+1}+o(1)+mv\rig)\\
    &=\int_{1}^{k-1}\frac{dv}{\lef(1-\frac{v}{k}\rig)\sqrt{2\pi v}}\exp\lef(-v\log\frac{vn}{k^2}-\frac{2v^2}{k}-\frac{1}{12v+1}+o(1)+mv\rig)+o(1)\\
    &\leq\int_{[\frac{1}{k},1-\frac{1}{k}]}\frac{\sqrt kdx}{(1-x)\sqrt{2\pi x}}\exp\lef(g(x)+o(1)\rig)+o(1),
\end{align*}
where we denote $g(x):=-xk\log\frac{nx}{k}-2x^2k+mkx$.
To find the maximum of $g(x)$ in the interval of integral, we found that there exists $c>0$ such that for all $m\leq c\log n$ its first derivative satisfies:
\begin{align*}
    g^\prime(x)=-k\log\frac{nx}{k}-k+4xk +mk<0 \qquad\text{ for sufficiently large }k\text{ and }\forall x\in\lef[\frac{1}{k},1-\frac{1}{k}\rig].
\end{align*}
Therefore we utilize lemma \ref{laplaceu2} to get 
\begin{align*}
    \sum_{v=1}^{k-1}\bb P(V=v) E_k^m(v)\leq \exp\lef(-\log\frac{n}{k^2}+o(1)+m\rig)\frac{k}{k}+o(1)=o(1).
\end{align*}
Therefore we finish the proof that there exists $c>0$ such that for all $m< c\log n$ all tests are powerless.

\subsection{Proof of Theorem \ref{thmlowtempfindclique}}
This proof is also divided by the local and the global parts, separated by $k\asymp n^{\frac{1}{2}}$.
\begin{center}
    \textbf{Local Part}
\end{center}
    We first analyze the null. Here we make use of the following facts and lemmas.
\begin{fact}[\citep{vershynin2018high}]\label{squaresubg}
    A random variable is sub-Gaussian if and only if $X^2$ is sub-exponential. Moreover  we have $\Vert X^2\Vert_{\psi_1}=\Vert X\Vert_{\psi_2}^2$.
\end{fact}
\begin{fact}[\citep{vershynin2018high}]\label{abssubg}
        When $X$ is sub-Gaussian r.v. $X$  we have $\Vert X\Vert_1\leq C_1$ for some $C_1>0$. And $Z=|X|$ satisfy $\Vert Z\Vert_{\psi_2}\leq C_2$ for some $C_2>0$.
\end{fact}
\begin{fact}[\citep{vershynin2018high}]\label{centeredsubgaussian}
    For sub-Gaussian r.v. $X$  we have $\Vert X-\bb E[X]\Vert_{\psi_2}\leq C\Vert X\Vert_{\psi_1}$ for some constant $C>0$.
\end{fact}
\begin{lemma}\label{lowtempaltertail}
        Under the high temperature regime with $\theta_1>\frac{1}{\bb E[\sech^2(h^\prime)]}$, with $S$ being the clique,  we have 
        \begin{align*}
             \bigg\Vert\bigg|\frac{1}{\sqrt k}\sum_{i\in S}\sigma_i\bigg|-\bb E\bigg[\bigg|\frac{1}{\sqrt k}\sum_{i\in S}\sigma_i\bigg|\bigg]\bigg\Vert_{\psi_2}\asymp 1,
        \end{align*}
        and $\bb E\lef[\lef|\frac{1}{k}\sum_{i\in S}\sigma_i\rig|\rig]\asymp 1$.
    \end{lemma}
    Therefore under the $\mca H_0$  we have by lemma \ref{tailh0} and lemma \ref{squaresubg}, for all $S$ such that $|S|=k$:
    \begin{align*}
       \bigg\Vert\frac{1}{\sqrt k}\sum_{i\in S}\sigma_i^{(j)}\bigg\Vert_{\psi_2}\asymp 1\quad\Rightarrow \bigg\Vert\frac{1}{m\sqrt k}\sum_{j=1}^m\sum_{i\in S}\sigma_i^{(j)}\bigg\Vert_{\psi_2}\asymp \frac{1}{\sqrt m}.
    \end{align*}
    And by lemma \ref{abssubg}  we have for some constant $C>0$, for sufficently large $k$  we have for all $t>0$
    \begin{align*}
        \bb P_0\bl\frac{1}{m}\sum_{j=1}^m\bigg|\frac{1}{k}\sum_{i\in S}\sigma_i^{(j)}\bigg| -\bb E\bigg[\bigg|\frac{1}{k}\sum_{i\in S}\sigma_i\bigg|\bigg]\geq t\br\leq \exp\lef(-Cmkt^2\rig).
    \end{align*}
    And by union bound  we have for some constant $C>0$ the Type I error can be upper bounded by
    \begin{align*}
        \bb P_0(\phi^\prime\geq\tau_\delta)&\leq\bb P_0\bigg(\phi^\prime-\bb E[\phi^\prime]\geq \tau_\delta-\bb E[\phi^\prime]\bigg)\leq \binom{n}{k}\exp\lef(-Cmk(\tau_\delta-\bb E[\phi^\prime])^2\rig)\\
        &\leq\lef(\frac{en}{k}\rig)^k\exp\lef(-Cmk(\tau_\delta-\bb E[\phi^\prime])^2\rig).
    \end{align*}
     Here we make use of lemma \ref{lowtempaltertail} to get that for some constant $C>0$ the Type II error can be upper bounded by 
    \begin{align*}
         \bb P_{S_0}\lef(\phi^\prime\leq\tau_\delta\rig)&\leq\bb P_{S_0}\lef(\phi_{S_0}\leq\tau_\delta\rig)\leq \bb P_{S_0}\bl \bb E\bigg[\bigg|\frac{1}{k}\sum_{i\in S_0}\sigma_i\bigg|\bigg]-\phi_{S_0} \geq \bb E\bigg[\bigg|\frac{1}{k}\sum_{i\in S_0}\sigma_i\bigg|\bigg] - \tau_\delta \br\\
         &\leq \exp\bl-Cmk\bl\bb E\bigg[\bigg|\frac{1}{k}\sum_{i\in S_0}\sigma_i\bigg|\bigg] - \tau_\delta\br^2\br.
    \end{align*}
    Therefore there exists $m\asymp \log n$ such that both Type I and II error are upper bounded by $\frac{\delta}{2}$.
    \begin{center}
        \textbf{Global Part}
    \end{center}
    First we consider controling the Type II error:
    We define $m_S:=\frac{1}{k}\sum_{i\in S}\sigma_i$ and $m_{S^c}:=\frac{1}{n-k}\sum_{i\in S^c}\sigma_i$. By lemma \ref{cltrfcr} we check that the asymptotic value of $\bb E[|m_S|]$ is the positive solution to the following equation:
    \begin{align*}
        x=\bb E[\tanh(\theta_1x+h)].
    \end{align*}
     We easily checked that $m_S\perp m_{S^c}$. It is not hard to see that by central limit theorem of i.i.d. Rademacher r.v.s. for all $t\in\bb R$,
     \begin{align*}
         \sqrt{n-k}m_{S^c}\overset{d}{\to} N(0,1).
     \end{align*}
     And analogously by lemma \ref{cltrfcr}  we have
     \begin{align*}
         \lef(\sqrt k(m_S-\bb E[|m_S|])|m_S>0\rig)\overset{d}{\to}N\lef(0,\frac{1-\theta_1(\bb E[\sech^2(\theta_1x+h)])^2}{\lef(1-\theta_1(\bb E[\sech^2(\theta_1x+h)])\rig)^2}\rig).
     \end{align*}
     Therefore,  we have
     \begin{align*}
         \lef(m_S+\frac{n-k}{k}m_{S^c}\bigg| m_{S}>0\rig)\overset{d}{\to} N\lef(\bb E[|m_S|],\frac{n}{k^2}\rig)\overset{d}{=} N\lef(x,\frac{n}{k^2}\rig).
     \end{align*}
     And analogously  we have $\lef(m_S+\frac{n-k}{k}m_{S^c}\bigg| m_{S}<0\rig)\overset{d}{\to} N\lef(-\bb E[|m_S|],\frac{n}{k^2}\rig)$. Then  we have
     \begin{align*}
         \lef(\lef|m_S+\frac{n-k}{k}m_{S^c}\rig|\bigg|m_S>0\rig)&\overset{d}{=}\lef(\lef|m_S+\frac{n-k}{k}m_{S^c}\rig|\bigg|m_S<0\rig)\overset{d}{=}\lef|m_S+\frac{n-k}{k}m_{S^c}\rig|\\
         &\overset{d}{\to}\ca {FN}\lef(x,\frac{n}{k^2}\rig),
     \end{align*}
     where $\ca{FN}$ is the short hand of folded normal distribution.
     Hence, by symmetry  we have $\bb P(m_S>0)=\bb P(m_S<0)$. By the property of folded Gaussian and the fact that $m_S+\frac{n-k}{k}m_{S^c}$ is uniformly integrable random variable  we have \citep{billingsley2013convergence}:
     \begin{align*}
         \bb E[\phi_4]&=\frac{1}{2}\bb E\lef[\lef|m_S+\frac{n-k}{k}m_{S^c}\rig|\bigg|m_{S}>0\rig]+\frac{1}{2}\bb E\lef[\lef|m_S+\frac{n-k}{k}m_{S^c}\rig|\bigg|m_{S}<0\rig]+o(1)\\
         &= \frac{\sqrt n}{k}\sqrt{\frac{2}{\pi}}\exp\lef(-\frac{x^2k^2}{2n}\rig)+x\lef[1-2\Phi\lef(-\frac{xk}{\sqrt n}\rig)\rig]+o(1).
     \end{align*}
     Then, we observe that $\bb E\lef[\exp\lef(\lef|m_S+\frac{n-k}{k}m_{S^c}\rig|\rig)\rig]<\infty$ implies 
     $\Vert\lef|m_S+\frac{n-k}{k}m_{S^c}\rig|\Vert_{\psi_1}<\infty$. Hence, by Bernstein inequality, there exists $C>0$ such that
     \begin{align*}
         \bb P(\phi_4\leq\bb E[\phi_4]-t)\leq\exp\lef(-Cmt^2\wedge mt\rig).
     \end{align*}
     Then we move to the null. For the Type I error, we note that $\Vert \frac{1}{\sqrt n}\sum_{i=1}^k\sigma_i\Vert_{\psi_2}\asymp\Vert\frac{1}{k}\sum_{i=1}^k\Vert_{\psi_2}<\infty$. Similar to the Type II error  we have $|\frac{1}{k}\sum_{i=1}^n\sigma_i|\overset{L_1}{\to}\ca {FN}(0,\frac{n}{k^2})$. Therefore  we have
     \begin{align*}
         \bb E_0[\phi_4]=\sqrt{\frac{2k^2}{n\pi}}+o(1).
     \end{align*}
     By Hoeffding's inequality, there exists $C>0$ such that
     \begin{align*}
         \bb P_0(\phi_4\geq\bb E_0[\phi_4]+t)\leq\exp(-Cmt^2).
     \end{align*}
     Therefore we can control the Type I+ Type II error by $\delta$ with $m\asymp 1$ and choosing $$\tau_\delta\in\lef(\sqrt{\frac{2n}{\pi k^2}},\frac{\sqrt n}{k}\sqrt{\frac{2}{\pi}}\exp\lef(-\frac{x^2k^2}{2n}\rig)+x\lef[1-2\Phi\lef(-\frac{xk}{\sqrt n}\rig)\rig]\rig).$$
     Then we move to the case of $k=\omega(\sqrt n)$. It is checked that when $m=1$  we have
     \begin{align*}
         \text{Under the Null},\quad\phi_4\overset{p}{\to}x,\qquad \text{ Under the Alternative },\phi_4\overset{p}{\to}0.
     \end{align*}
     And one sample is already be enough.

\subsection{Proof of Corollary \ref{recoveryguaranteelowtemp}}
    The underlying idea of this proof is to extend the method we used to derive the central limit theorem to a more general quantity of the form $m_{S^\prime}=\frac{1}{k}\lef(\sum_{i\in S\cap S^\prime}\sigma_i+\sum_{u\in S^\prime\setminus(S^\prime\cap S)}\sigma_i\rig)$ for some set $S^\prime$ with $|S^\prime\cap S|=(1-c)k$ with $S$ being the index set of the hidden clique. Without loss of generality we assume $S=[k]$ and $S^\prime=[r+1:r+k]$ with $r=ck$ for some $c\in[0,1]$. 

    We  define $m_{rk}:=\frac{1}{k-r}\sum_{i=r+1}^k\sigma_i$, $m_{r}:=\frac{1}{r}\sum_{i=1}^r\sigma_i$, and $m_r^\prime:=\frac{1}{r}\sum_{i=k+1}^{k+r}\sigma_i$. Then  we have $\frac{1}{k}\sum_{i\in S^\prime}\sigma_i=(1-c)m_{rk}+cm_{r}^\prime$ and we analyze each term separately. Note that the second part is i.i.d. and  we have
    \begin{align}\label{mrclt}
        \sqrt{r} m_{r}^\prime\overset{d}{\to} N(0, 1)\quad\text{ and }\quad\bb E\lef[\exp\lef(t\sqrt rm_{r^\prime}\rig)\rig]\to\exp\lef(-\frac{t^2}{2}\rig).
    \end{align}
    Then we move toward the study of $m_{rk}$, here we continue make use of the Laplace method and the transfer principle in \citep{ellis1980limit} as in the proof of theorem \ref{cltrfcr}. For some $a\in\bb R$, there exists some $C,\delta>0$ such that
    \tny{\begin{align*}
        \bb E\lef[\exp\lef(t\sqrt k((1-c)m_{rk}-a)\rig)\bigg| m_{rk}>0\rig]&=\bb E\lef[\frac{\int_{\Vert x-x_{0,k}\Vert\leq C}\exp\lef(-k\mca H_{0,k}(x,\bfa h)-ta\sqrt k\rig)dx}{\int_{\Vert x-x_{1,k}\Vert\leq C}\exp\lef(-k\mca H_{1,k}(x,\bfa h)\rig)dx}\rig]\\
        &\cdot\lef(1+O\lef(\exp\lef(-k\delta\rig)\rig)\rig),
    \end{align*}}
    with \begin{align*}
        \mca H_{0,k}(x,\bfa h):&=\frac{x^2}{2}-\frac{1}{k}\lef(\sum_{i=1}^{r}\log\cosh(h_i)+\sum_{i=r+1}^k\log\cosh\lef(h_i+\frac{t}{\sqrt k}\rig)\rig),\\
        \mca H_{1,k}(x,\bfa h):&=\frac{x^2}{2}-\frac{1}{k}\sum_{i=1}^k\log\cosh\lef(h_i\rig).
    \end{align*} for some $C>0$ and $x_{0,k}$, $x_{1,k}$ to be the positive local minimum of $\mca H_{0,k}, \mca H_{1,k}$ respectively. Further assume that $x_0^*:\in\argmin_{x}\mca H_{0}$ and $x_1^*:\in\argmin_{x}\mca H_1$ to be the positive root respectively. By uniform convergence we can get $x_{0,k}\to x_0^*$ and $x_{1,k}\to x_1^*$ with $x_1^*=x_0^*$. By Fermat's condition  we have
    \begin{align*}
       \mca H_{1.k}^\prime(x_{1,k},\bfa h)=0,
    \end{align*}
    which further implies that
    \begin{align*}
        0= \mca H_{1,k}^\prime(x_{1,k},\bfa h)=\mca H_{1,k}^\prime(x_1^*,\bfa h)+\mca H_1^{(2)}(x_1^*,\bfa h)(x_{1,k}-x_1^*) +O\lef((x_{1,k}-x_1^*)^2\rig).
    \end{align*}
    Noticing that $x_1^*=\sqrt{\theta_1}\bb E[\tanh(\sqrt{\theta_1}x_1^*+h)]$  we use the linearization of $Z$ estimator in lemma \ref{linearZ} to get
    \begin{align}
       \sqrt{k}(x_{1,k}-x_1^*)&=\frac{\sqrt{\theta_1}\sum_{i=1}^k(\tanh(\sqrt{\theta_1}x_1^*+h_i)-\bb E[\tanh(\sqrt{\theta_1}x_1^*+h_i)])}{\sqrt {k}(1-\theta_1\bb E[\sech^2(\sqrt{\theta_1}x_1^*+h)])}+o_{\psi_2}(1).\label{linearasymp}
    \end{align}
    We define \\$F(y,x,\bfa h):=\frac{x^2}{2}-\frac{1}{k}\sum_{i=1}^{r}\log\cosh\lef(\sqrt{\theta_1}x+h_i\rig)-\frac{1}{k}\sum_{i=r+1}^{k}\log\cosh(\sqrt{\theta_1}x+h_i+y)$. Then, $F(t/\sqrt k,x,\bfa h)= \mca H_{0,k}(x,\bfa h)$ and $F(0,x,\bfa h)=\mca H_{1,k}(x,\bfa h)$. Further using the fact that
    \begin{align*}
        \frac{\pta \mca H_{0,k}(x_{0,k},\bfa h)}{\pta x}= \frac{\pta \mca H_{1,k}(x_{1,k},\bfa h)}{\pta x}=0,
    \end{align*}
    we subsequently get
    \begin{align*}
        k(\mca H_{0,k}(x_{0,k},\bfa h)&-\mca H_{1,k}(x_{1,k},\bfa h))=\frac{\pta F(0,x_{1,k},\bfa h)}{\pta y}\sqrt k t +\frac{1}{2}\frac{\pta^2 F(0,x_{1,k},\bfa h)}{\pta y^2}t^2 +o_{\psi_2}(1)\\
        &= -\frac{t}{k^{1/2}}\sum_{i=r+1}^k\tanh(\sqrt{\theta_1}x_{1,k}+h_i)-\frac{t^2}{2k}\sum_{i=r+1}^k\sech^2(\sqrt{\theta_1}x_{1,k}+h_i)+o_{\psi_2}(1)\\
        &=-\frac{t}{k^{1/2}}\sum_{i=r+1}^k\tanh(\sqrt{\theta_1}x_{1}^*+h_i)-\frac{t^2}{2k}\sum_{i={r+1}}^k\sech^2(\sqrt{\theta_1}x_{1}^*+h_i)\\
        &-\frac{\sqrt{\theta_1}t}{k}\sum_{i=r+1}^k\sech^2(\sqrt{\theta_1}x_1^*+h_i)\sqrt k(x_{1,k}-x_1^*) +o_{\psi_2}(1).
    \end{align*}
    And we also notice that $\frac{\mca H_{0,k}^{(2)}(x_{0,k},\bfa h)}{\mca H_{1,k}^{(2)}(x_{1,k},\bfa h)}=1+o(1)$. Picking $a=(1-c)\bb E[\tanh(\sqrt{\theta_1}x_1^*+h)]$ and applying lemma \ref{laplace} we get 
    \begin{align*}
       \bb E\lef[\frac{\int_{\Vert x-x_{0,k}\Vert\leq C}\exp\lef(-k\mca H_{0,k}(x,\bfa h)-ta\sqrt k\rig)dx}{\int_{\Vert x-x_{1,k}\Vert\leq C}\exp\lef(-k\mca H_{1,k}(x,\bfa h)\rig)dx}\rig]
        &=\exp\lef(\frac{V(c)}{2}t^2\rig)\lef(1+o(1)\rig),
    \end{align*}
    with $V(c):=(1-c)\frac{(1-\theta_1(\bb E[\sech^2(\sqrt{\theta_1}x_1^*+h)])^2-\bb E[\tanh(\sqrt{\theta_1}x_1^*+h)]^2)}{(1-\theta_1\bb E[\sech^2(\sqrt{\theta_1}x_1^*+h)])^2}$.

    Combining pieces, for the current $a$, we get
    \begin{align*}
        \bb E\lef[\exp\lef(t\sqrt k((1-c)m_{rk}- a)\rig)\bigg |m_{rk}>0\rig]&=\exp\lef(\frac{V(c)}{2}t^2\rig)(1+o(1)).
    \end{align*}
   Combining the above result with \eqref{mrclt}, we can see that
   \tny{\begin{align}
       \bb E\lef[\exp\lef(t\sqrt k((1-c)m_{rk}+cm_r^\prime-a)\rig)\bigg | m_{rk}>0\rig]=\exp\lef(\frac{V(c)+c}{2}t^2\rig)(1+o(1)).\label{mrksmalzero}
   \end{align}}
   Noticing the fact that for $Z\sim N(\mu,\sigma^2)$  we have $\bb E[|Z|]= \sqrt{\frac{2}{\pi}}\sigma\exp\lef(-\frac{\mu^2}{2\sigma^2}\rig)+\mu\lef(1-2\Phi(-\frac{\mu}{\sigma})\rig)$, and
    \begin{align*}
        \bb E[|mk||m_{rk}>0]=\bb E[|(1-c)m_{rk}+cm^\prime_r||m_{rk}>0]= (1-c)\bb E[\tanh(\sqrt{\theta_1}x_1^*+h)] +o(1).
    \end{align*}
    Similar result holds for $\bb E[|mk||m_{rk}<0]$ and we finally conclude that
    \begin{align*}
        \bb E[\phi_{S_c}]=\bb E[|m_S|]=(1-c)\bb E[\tanh(\sqrt{\theta_1}x_1^*+h)].
    \end{align*}
    And it is also seen that by \eqref{mrksmalzero} and the same quantity given $m_{rk}<0$,  we have
    \sm{\begin{align*}
        \bb E\lef[\exp\lef(t\sqrt k(|(1-c)m_{rk}+cm_r^\prime|-\bb E[\tanh(\sqrt{\theta_1}x_1^*+h)])\rig)\rig]\leq 2\exp\lef(\frac{V(c)+c}{2}t^2\rig)(1+o(1)).
    \end{align*}}
    Then we conclude that for $S_c$ such that $|S_c\cap S|=(1-c)k$ there exists $C>0$ such that $\Vert\phi_{S_c}\Vert_{\psi_2}\leq C$ and $\phi_{S_c}=(1-c)\bb E[\tanh(\sqrt{\theta_1}x_1^*+h)]+o(1)$. Therefore one can follow a similar path as in the proof of \eqref{hightemperaturerec} upon noticing that there exists $C_1>0$ and
    \begin{align*}
        \bb P\lef(\phi_{S_c}-\bb E[\phi_{S_c}]\geq x\rig)&\leq \exp\lef(-C_1mx^2\rig),\quad \bb P\lef(\phi_{S_c}-\bb E[\phi_{S_c}]\leq -x \rig)\leq\exp\lef(-C_1mx^2\rig).
    \end{align*}
    And we therefore notice that by the same method as \eqref{overlapprob}, when $m\gtrsim \log n$ we have
    \begin{align*}
        \bb P\lef(|S\Delta S_\max|\leq k\epsilon\rig)\geq 1-\exp\lef(-C_4\log\lef(\frac{\epsilon n}{k}\rig)\epsilon k\rig),
    \end{align*}
    and therefore we complete the proof.

\subsection{Proof of Theorem \ref{minimaxlbct}}
Before we start the proof, we note that in this proof we only need to consider the $k=o\lef(n^{\frac{2\tau-1}{4\tau-3}}\rig)$ case. The rest of the region can be accomplished by a $O(1)$ samples and has no need to elaborate on lower bounds.
Here we reuse all the notations in the proof of theorem \ref{minimaxlb} at appendix \ref{proofminimaxlb}. Here, instead of \eqref{esgoodset} we use the following good sets for some $C>0$:
\begin{align*}
    E_S^*:=\lef\{|k^{\frac{1}{4\tau-2}}m_S|\leq C\lef(\log (m\vee k)\log k\rig)^{\frac{1}{4\tau-2}}\rig\},\qquad \bb P_{S}^*(\bfa\sigma)=\begin{cases}
        \bb P_S(\bfa\sigma) &\text{ if }\bfa\sigma\in E_S^*\\
        0&\text{otherwise}
    \end{cases}.
\end{align*}
Then we analogously have
\begin{align*}
    \Vert\bb P_{\bar S}-\bb P^*_{\bar S}\Vert_{TV}=\int|d\bb P_{\bar S}(\bfa\sigma)-d\bb P^*_{\bar S}(\bfa\sigma)|\leq\frac{1}{\binom{n}{k}}\sum_{S:|S|\leq k}\bb P_{S}(E_S^c)=O\lef(\frac{1}{(m\vee k)\log k}\rig).
\end{align*}
Here we analogously have
\tiny{\begin{align}\label{psi00}
&\bb E\lef[\frac{\bb P^*_{S}\bb P_{S^\prime}^*}{\bb P_0}\rig]=\bb E\lef[\frac{\sum_{\bfa\sigma:E^*_S\cap E^*_{S^\prime}}\exp\bigg(\frac{\theta_1}{2k}\lef(\lef(\sum_{i\in [k]}\sigma_i\rig)^2+\lef(\sum_{i\in[k+r]}\sigma_i\rig)^2\rig)+\sum_{i\in[k+r]}h_i\sigma_i\bigg)\sum_{\bfa\sigma}\exp\bigg(\sum_{i\in[k+r]}\sigma_ih_i\bigg)}{\lef(\sum_{\bfa\sigma}\exp\lef(\frac{\theta_1}{2k}\lef(\sum_{i\in[k]}\sigma_i\rig)^2+\sum_{i\in[k+r]}h_i\sigma_i\rig)\rig)\lef(\sum_{\bfa\sigma}\exp\lef(\frac{\theta_1}{2k}\lef(\sum_{i\in[r+1:k+r]}\sigma_i\rig)^2+\sum_{i\in[k+r]}h_i\sigma_i\rig)\rig)}\rig].
\end{align}}\normalsize
Analogous to \eqref{boundingprocedure},  we have
\tny{\begin{align}\label{vphi0gc}
\bb E\lef[\frac{\bb P^*_{S}(\bfa\sigma)\bb P^*_{S^\prime}(\bfa\sigma)}{\bb P_0(\bfa\sigma)}\rig]\leq\bb E\lef[\frac{\prod_{i=r+1}^k\cosh(h_i)\int_{|x|\vee|y|\leq c_1}\exp\lef(-kG_{0,k}(x,y,\bfa h)\rig)dxdy }{\int_{\bb R}\exp\lef(-kG_{1,k}(x,\bfa h)\rig)dx\int_{\bb R}\exp\lef(-kG_{2,k}(y,\bfa h)\rig)dy}\rig]+o\lef(\frac{1}{m}\rig),
\end{align}}
where we define $c_1= C\bl\log (m\vee k)\log k\br^{\frac{1}{4\tau-2}}$ for some $C>0$. In \eqref{vphi0gc}, we recall the definition
\begin{align*}
G_{0,k}(x,y):&=\frac{x^2+y^2}{2}-\frac{1}{k}\bl\sum_{i=1}^{r}\log\cosh\lef(\sqrt{\theta_1}x+h_i\rig)\\
&+\sum_{i=r+1}^k\log\cosh(\sqrt{\theta_1}(x+y)+h_i)+\sum_{i=k+1}^{k+r}\log\cosh(\sqrt{\theta_1}y+h_i)\br
.\end{align*}
and
\begin{align*}
    G_{1,k}(x,\bfa h):&=\frac{x^2}{2}-\frac{1}{k}\sum_{i=1}^k\log\cosh\lef(\sqrt{\theta_1}x+h_i\rig),\\ G_{2,k}(x,\bfa h):&=\frac{x^2}{2}-\frac{1}{k}\sum_{i=r+1}^{k+r}\log\cosh\lef(\sqrt{\theta_1}x+h_i\rig).
\end{align*}
And  we have their population version
\begin{align*}
    G_0(x,y)&=\frac{x^2+y^2}{2}-c\bb E[\log\cosh(\sqrt{\theta_1}x+h)\cosh(\sqrt{\theta_1}y+h)]\\
    &-(1-c)\bb E[\log\cosh(\sqrt{\theta_1}(x+y)+h)],\\
    G_1(x)&=\frac{x^2}{2}-\bb E[\log\cosh(\sqrt{\theta_1}x+h)]
.\end{align*}
Then we recall that $(x_{k},y_k)$ is the minimum of $G_{0,k}$. To find this maximum we notice that for all $c\in[0,1)$ we consider the Fermat's condition
\begin{align*}
    \begin{bmatrix}
        x_k\\
        y_k
    \end{bmatrix}=\frac{1}{k}\begin{bmatrix}
\sum_{i=1}^r\tanh(\sqrt{\theta_1}x_k+h_i)+\sum_{i=r+1}^k\tanh(\sqrt{\theta_1}(x_k+y_k)+h_i)\\
        \sum_{i=r+1}^k\tanh(\sqrt{\theta_1}(x_k+y_k)+h_i)+\sum_{i=k+1}^{k+r}\tanh(\sqrt{\theta_1}y_k+h_i)
    \end{bmatrix}.
\end{align*}
And we also consider the Hessian, recall that $\sum_{1}:=\sum_{i=1}^r$, $\sum_{2}:=\sum_{i=r+1}^k$, and $\sum_{3}:=\sum_{i=k+1}^{k+r}$ we write the Hessian as
\begin{align*}
    &\nabla^2G_{0,k}(x,y,\bfa h)=\\
    &\tiny\begin{bmatrix} 
    1-\frac{\theta_1}{k}\lef(\sum_1\sech^2(\sqrt{\theta_1}x+h_i)+\sum_2\sech^2(\sqrt{\theta_1}(x+y)+h_i)\rig)&-\frac{\theta_1}{k}\sum_{2}\sech^2(\sqrt{\theta_1}(x+y)+h_i))\nnb\\-\frac{\theta_1}{k}\sum_{2}\sech^2(\sqrt{\theta_1}(x+y)+h_i))
    &1-\frac{\theta_1}{k}\lef(\sum_3\sech^2(\sqrt{\theta_1}y+h_i)+\sum_2\sech^2(\sqrt{\theta_1}(x+y)+h_i)\rig)
    \end{bmatrix}\normalsize.
\end{align*}
From here we omit $\bfa h$ in the notations for the purpose of clarity.

Consider the population version  we have for all $c\in[0,1)$:
\begin{align*}
    \det(\nabla^2G_0(0,0))&=(1-c\theta_1\bb E[\sech^2(h+\sqrt{\theta_1}x)])(1-c\theta_1\bb E[\sech^2(h+\sqrt{\theta_1}y)])\\
    &- (1-c)\theta_1\bb V[\sech^2(h+\sqrt{\theta_1}(x+y))]<0.
\end{align*}
which implies that $(0,0)$ is a local maximum. This further implies that the sequence of global minimum $x_k,y_k$ do not converge to $(0,0)$ when $(1-c)$ is non-vanishing. Then we aim to decide the scale at which it converges to $(0,0)$.

Then we consider when $c=\frac{r}{k}\to 1$. Introducing the notation $\mca H_1(x)=\frac{1}{r}\sum_{1}\sech^2(\sqrt{\theta_1}x+h_i)$, $\mca H_2(x,y)=\frac{1}{k-r}\sum_{2}\sech^2(\sqrt{\theta_1}(x+y)+h_i)$, and $\mca H_3(y):=\frac{1}{r}\sum_3\sech^2(\sqrt{\theta_1}y+h_i)$.
And using the result of lemma \eqref{cltrfcr} we note that $(0,0)$ is a global minimum for $c=1$. This implies that when $1-c$ vanishing, the global minimum sequence of $(x_k,y_k)$ converge to $(0,0)$ by boundedness on closed interval and uniform convergence. We then notice that $\nabla^iG_{0,k}(0,0)\overset{a.s.}{\to}\nabla^i G_{0}(0,0)$ using regularity condition given by lemma \ref{convergeasG}. Introduce $\times$ as the notation for $k-$mode tensor product. Introducing the notation $\bfa x=(x_k,y_k)^\top$ and using the Fermat's condition  we have
\ttny{\begin{align*}
    \nabla G_{0,k}(x_k,y_k)=0=\nabla G_{0,k}(0,0)+\nabla^2 G_{0,k}(0,0)\times\bfa x+\ldots+\frac{1}{(2\tau-1) !}\nabla^{2\tau}G_{0,k}(0,0)\times\bfa x^{2\tau-1}+ O(\Vert\bfa x\Vert_2^{2\tau}).
\end{align*}}
Notice that $\nabla^{(2\tau)}G_{0,k}(0,0) $ is diagonally dominated tensor with positive diagonal values exactly equivalent to $G_{1,k}^{(2\tau)}(0)$ and $G_{2,k}^{(2\tau)}(0)$. Hence when $1-c=o(k^{-\frac{\tau-1}{2\tau-1}})$, the term $\nabla^2 G_{0,k}(0,0)\times\bfa x$ is dominated by the term $\frac{1}{(2\tau-1)!}\nabla^{2\tau} G_{0,k}(0,0)\times\bfa x^{2\tau-1}$. We then have
\tny{\begin{align}\label{crossterms}
   \sqrt{k}\nabla^2G_{0,k}(0,0)\times\bfa x+ \sqrt k\nabla^{2\tau} G_{0,k}(0,0)\times\bfa x^{(2\tau-1)}&=\sqrt k\begin{bmatrix}
        G_{1,k}^{(2\tau)}(0)x_{k}^{2\tau-1}+O(1-c)y_k\\
        G_{2,k}^{(2\tau)}(0)y_{k}^{2\tau-1}+O(1-c)x_k
    \end{bmatrix}\nnb\\
    &=-\begin{bmatrix}
        \frac{\sqrt{\theta_1}}{\sqrt k}\sum_{i=1}^k\tanh(h_i)\\
         \frac{\sqrt{\theta_1}}{\sqrt k}\sum_{i=r+1}^{k+r}\tanh(h_i)
    \end{bmatrix}+o_{\psi_2}(1).
\end{align}}


Similarly for the rest of two functions  $G_{1,k}$ and $G_{2,k}$  we have
\begin{align*}
    G^\prime_{1,k}(x_{1,k})&=0=G_{1,k}^\prime(0)+G_{1,k}^{(2)}(0)x_{1,k}+\ldots+ \frac{1}{(2\tau-1)!}G_{1,k}^{(2\tau)}(0)x_{1,k}^{2\tau-1}+O(x_{1,k}^{2\tau}),\\
    G^\prime_{2,k}(x_{2,k})&=0=G_{2,k}^\prime(0)+G^{(2)}_{2,k}(0)x_{2,k}+\ldots+\frac{1}{(2\tau-1)!}G^{(2\tau)}_{2,k}(0)x_{2,k}^{2\tau-1}+O(x_{2,k}^{2\tau}).
\end{align*}
Hence, using the linearization results given by lemma \ref{linearZ}  we have
\begin{align*}
    \sqrt kx_{1,k}^{2\tau-1}&=\frac{-(2\tau-1)!\sum_{i=1}^k\tanh(h_i)}{\sqrt kG_{1,k}^{(2\tau)}(0)} + o_{\psi_2}(1), \\
    \sqrt kx_{2,k}^{2\tau-1}&=\frac{-(2\tau-1)!\sum_{i=r+1}^{k+r}\tanh(h_i)}{\sqrt kG_{2,k}^{(2\tau)}(0)} + o_{\psi_2}(1).
\end{align*}
Therefore,  we have
\tny{\begin{align}\label{boundedorlicz}
    \Vert k^{1/2}x_{1,k}^{2\tau-1}\Vert_{\psi_2}<\infty,\quad \Vert k^{1/2}x_{2,k}^{2\tau-1}\Vert_{\psi_2}<\infty.\quad\Rightarrow\quad\Vert k^{\frac{1}{4\tau-2}}x_{1,k}\Vert_{\psi_{2\tau-1}},\Vert k^{\frac{1}{4\tau-2}}x_{2,k}\Vert_{\psi_{2\tau-1}}<\infty.
\end{align} }
And otherwise when $1-c=\omega( k^{-\frac{2\tau-2}{4\tau-2}})$,  we have
\begin{align}\label{normalxy}
    \begin{bmatrix}
        \sum_{i=r+1}^k\frac{\theta_1}{\sqrt{k}}\sech^2(h_i)y_{k}\\
        \sum_{i=r+1}^k\frac{\theta_1}{\sqrt{k}}\sech^2(h_i)x_{k} 
    \end{bmatrix} =\begin{bmatrix}
        \frac{\sqrt{\theta_1}}{\sqrt k}\sum_{i=1}^k\tanh(h_i)\\
         \frac{\sqrt{\theta_1}}{\sqrt k}\sum_{i=r+1}^{k+r}\tanh(h_i)
    \end{bmatrix}+o_{\psi_2}(1).
\end{align}
And similarly it is not hard to verify that
\begin{align*}
    \lef\Vert(1-c)\frac{1}{\sqrt k} x_{k}\rig\Vert_{\psi_2}<\infty,\qquad\lef\Vert(1-c)\frac{1}{\sqrt k} y_k\rig\Vert_{\psi_2}<\infty.
\end{align*}
Then we consider the region of $1-c=o(k^{-\frac{2\tau-1}{4\tau-2}})$, using \ref{crossterms}  we have
\begin{align}\label{centeredrelation}
    \sqrt k\begin{bmatrix}
        G_{1,k}^{(2\tau)}(0)x_{k}^{2\tau-1}+O(1-c)y_k\\
        G_{2,k}^{(2\tau)}(0)y_{k}^{2\tau-1}+O(1-c)x_k
    \end{bmatrix}=\begin{bmatrix}
        x_{1,k}^{2\tau-1}\\
        x_{2,k}^{2\tau-1}
    \end{bmatrix}(1+o_{\psi_2}(1)).
\end{align}
And  we have when $1-c=o\lef( k^{-(2\tau-2)/(4\tau-2)}\rig)$, let $A_k:=\frac{\sqrt{\theta_1}}{\sqrt k}\sum_{i=1}^k\tanh(h_i)$, $B_k:=\frac{\sqrt{\theta_1}}{\sqrt k}\sum_{i=r+1}^{r+k}\tanh(h_i)$, 
\tny{\begin{align*}
    k^{\frac{1}{4\tau-2}}x_k=\sign(A_k)(| A_k|+O(1-c)y_k\sqrt k)^{\frac{1}{2\tau-1}}=\sign(A_k)|A_k|^{\frac{1}{2\tau-1}}+ {|A_k|^{-\frac{2\tau}{2\tau-1}}}\sqrt k O(1-c)y_k(1+o(1)),\\
   k^{\frac{1}{4\tau-2}}y_k=\sign(B_k)(| B_k|+O(1-c)x_k\sqrt k)^{\frac{1}{2\tau-1}}=\sign(B_k)|B_k|^{\frac{1}{2\tau-1}}+ {|B_k|^{-\frac{2\tau}{2\tau-1}}}\sqrt k O(1-c)x_k(1+o(1)).
\end{align*}}
Solving the principle terms in the above equation it is not hard to conclude that  we have $x_{1,k}-x_{k}= O(1-c)k^{\frac{2\tau-2}{4\tau-2}}x_{2,k}$. $x_{2,k}-y_k= O(1-c)k^{\frac{2\tau-2}{4\tau-2}}x_{1,k}$. 
And using \eqref{varphi01}  we have for $\delta >1$, $\bb E[|x_k^{2\tau-1}|]$, $\bb E[|y_k^{2\tau-1}|]$,$\bb E[|x_{1,k}^{2\tau-1}|]$, $\bb E[|x_{2,k}^{2\tau-1}|]=O\lef(\frac{1}{k^{1/2}}\rig)$.
Going back to \eqref{vphi0gc} we notice that
\tny{\begin{align*}
    \bb E\bigg[\frac{\bb P_{S}(\bfa\sigma)\bb P_{S^\prime}(\bfa\sigma)}{\bb P_0(\bfa\sigma)}\bigg]=\bb E\bigg[\exp\bigg(-k(G_{0,k}(x_k,y_k)+G_{1,k}(x_{1,k})+G_{2,k}(x_{2,k}))+\sum_{i=r+1}^k\log\cosh(h_i)\bigg)\ca A\bigg],
\end{align*}}
with
\sm{\begin{align*}
    \ca A:&=\frac{\int_{|x|\vee|y|\leq c_1}\exp\lef(-k(G_{0,k}(x,y)-G_{0,k}(x_k,y_k)) \rig)dxdy}{\int_{\bb R}\exp(-k(G_{1,k}(x,\bfa h)-G_{1,k}(x_{1,k},\bfa h)))dx\int_{\bb R}\exp(-k(G_{2,k}(y,\bfa h)-G_{2,k}(x_{2,k},\bfa h)))dy}.
\end{align*}}
First we consider the numerator of $\ca A$, define $\delta\bfa x:=(x-x_k,y-y_k)$ and
\begin{align*}
    G_{0,k}(x,y,\bfa h)&-G_{0,k}(x_k,y_k,\bfa h)=\sum_{i=2}^{2\tau}\frac{1}{i!}\nabla^{i}G_{0,k}(x_k,y_k)\times\delta\bfa x^i+O(\Vert\delta\bfa x\Vert_2^{2\tau+1})\\
    &=\frac{1}{(2\tau)!}\nabla_x^{2\tau}G_{0,k}(x_k,y_k)(x-x_k)^{2\tau}+\frac{1}{(2\tau)!}\nabla_y^{2\tau}G_{0,k}(x_k,y_k)(x-x_k)^{2\tau}\\
    &+\nabla^2_{xy}G_{0,k}(x_k,y_k)(x-x_k)(y-y_k)+O(\Vert\delta\bfa x\Vert_2^{2\tau+1}).
\end{align*}
The underlying idea of the above proof is to decouple the cross terms and analyze them separately,
\tny{\begin{align}\label{separationofcrossterms}
    &\int_{|x|\vee|y|\leq c_1}\exp\lef(-k(G_{0,k}(x,y)-G_{0,k}(x_k,y_k)) \rig)dxdy\leq \exp\lef(C(1-c)^2k^{\frac{4\tau-4}{2\tau-1}}\rig)\nnb\\
    &\cdot\int_{|x|\leq c_1}\exp\lef(\sum_{i=2}^\infty\frac{1}{i!}\nabla_x^{i}G_{0,k}(x_k,y_k)(x-x_k)^{i}\rig)dx\int_{|y|\leq c_1}\exp\lef(\sum_{i=2}^\infty\frac{1}{i!}\nabla_x^{i}G_{0,k}(x_k,y_k)(y-y_k)^{i}\rig)dy.
\end{align}}
Then the rest of the integral is a product and we can apply the higher order $1$-dimensional Laplace Approximation (Also see the \citep{bolthausen1986laplace}) to analyze them. Then we can use the Laplace approximation of integral of the denominator and H\"older's inequality to obtain that for all $\delta>0$ there exists $\tau_1>1$ such that,
\begin{align}\label{controlofsmallterms}
    \bb E[\ca A^{1+\delta}]\leq\exp(C(1-c)^2k^{\frac{4\tau-4}{2\tau-1}})\bb E\lef[\lef(\frac{G^{(2\tau)}(x_{1,k})G^{(2\tau)}(x_{2,k})}{\nabla_x^{2\tau}G(x_k,y_k)\nabla_y^{2\tau}G(x_k,y_k)}\rig)^{\frac{\tau_1}{2\tau}}\rig]^{\frac{1}{\tau_1}}.
\end{align}
Then we notice that
\tny{\begin{align*}
    \nabla_x^{2\tau}G(x_k,y_k)-G^{(2\tau)}(x_{1,k})=G^{(2\tau+1)}(x_{1,k})(x_k-x_{1,k})+O((x_k-x_{1,k})^2)=O((1-c)^2k^{\frac{4\tau-4}{2\tau-1}}),\\
    \nabla_y^{2\tau}G(x_k,y_k)-G^{(2\tau)}(x_{2,k})=G^{(2\tau+1)}(x_{1,k})(y_k-x_{2,k})+O((y_k-x_{2,k})^2)=O((1-c)^2k^{\frac{4\tau-4}{2\tau-1}}).
\end{align*}}
Therefore  we have
\begin{align*}
    \bb E[\ca A^{1+\delta}]
    &\leq\exp(C_1(1-c)^2k^{\frac{4\tau-4}{2\tau-1}}).
\end{align*}
Then, we analyze the first term in $\bb E\lef[\frac{\bb P_{S}(\bfa\sigma)\bb P_{S^\prime}(\bfa\sigma)}{\bb P_0(\bfa\sigma)}\rig]$, there exists $C>0$ such that
\tny{\begin{align*}
    &\ca B:=-kG_{0,k}(x_k,y_k)+\sum_{i=r+1}^k\log\cosh(h_i)+kG_{1,k}(x_{1,k})+kG_{2,k}(x_{2,k})\\
    &=-k\bl\nabla G_{0,k}(0,0)\times\bfa x+\frac{1}{2}\nabla^2 G_{0,k}(0,0)\times\bfa x^2+\cdots+\frac{1}{(2\tau)!}\nabla^{2\tau}G_{0,k}(0,0)\times\bfa x^{2\tau}-G^\prime_{1,k}(0)x_{1,k}\\
    &-\frac{1}{2}G^{(2)}_{1,k}(0)x_{1,k}^2-\cdots-\frac{1}{(2\tau)!}G_{1,k}^{(2\tau)}(0)x_{1,k}^{2\tau}-G^\prime_{2,k}(0)x_{2,k}-\frac{1}{2}G^{(2)}_{2,k}(0)x_{2,k}^2-\cdots-\frac{1}{(2\tau)!}G_{2,k}^{(2\tau)}(0)x_{2,k}^{2\tau}
    \\
    &+(1-c)O(x_{1,k}^{2\tau+1}+x_{2,k}^{2\tau+1})\br\\
    &=-\bl 1-\frac{1}{(2\tau)!}\br kG_{1,k}^\prime(0)(x_k-x_{1,k})-\bl 1-\frac{1}{(2\tau)!}\br kG^\prime_{2,k}(0)(y_k-x_{2,k})+k\nabla^2_{xy}G_{0,k}(0,0)x_{k}y_{k}\\
    &+(1-c)O(x_{1,k}^{2\tau+1}+x_{2,k}^{2\tau+1}).
\end{align*}}
Therefore, when $1-c=o\lef( k^{-\frac{2\tau-2}{4\tau-2}}\rig)$ the cross term in $\ca B$ dominates the higher order terms, and  we have
\tny{\begin{align*}
    \bb E[\ca B]&\asymp\bb E[kG_{1,k}^\prime(0)(x_k-x_{1,k})+kG_{2,k}^\prime(0)(y_k-x_{2,k})+k\nabla^2_{1,2}G_{0,k}(0,0)x_ky_k]=O\lef(k^{\frac{4\tau-4}{2\tau-1}}(1-c)^2\rig).
\end{align*}}
Then we use the fact that $\Vert k^{\frac{1}{2}}G_{1,k}^\prime(0)\Vert_{\psi_{2}}<\infty$ and $\Vert k^{\frac{1}{2}}G_{2,k}^\prime(0)\Vert_{\psi_{2}}<\infty$ to get that by H\"older's inequality there exists $\tau_1,\tau_2,\tau_3>1$ and $\frac{1}{\tau_1}+\frac{1}{\tau_2}+\frac{1}{\tau_3}=1$ such that for some $\lambda>1,C>0$ such that
\begin{align*}
     \bb E[\exp\lef(\lambda(\ca B-\bb E[\ca B])\rig)]&\leq  \bb E[\exp\lef(\tau_1 Ck(1-c)x_{1,k}x_{2,k}\rig)]^{\frac{1}{\tau_1}}\bb E[\exp(CkG^\prime_{1,k}(0)(x_k-x_{1,k}))]^{\frac{1}{\tau_2}}\\
     &\cdot\bb E[\exp(CkG^\prime_{2,k}(0)(y_k-x_{2,k}))]^{\frac{1}{\tau_2}}\\
     &\leq \exp(C(1-c)^2k^{\frac{4\tau-4}{2\tau-1}}\wedge (1-c)k^{\frac{2\tau-2}{2\tau-1}}\log^{\frac{1}{2\tau-1}}(m\vee k)).
\end{align*}
Then we consider $1-c=\omega(k^{-\frac{2\tau-2}{4\tau-2}})$ and $1-c=o(1)$. In this case  we have by \eqref{normalxy}:
\tny{\begin{align*}
    \ca B&=-k\frac{1}{2}\nabla G_{0,k}(0,0)\times \bfa x+k\bl 1-\frac{1}{(2\tau)!}\br G_{1,k}^\prime(0,0)x_{1,k}(1+o(1))+k\bl 1-\frac{1}{(2\tau)!}\br G_{2,k}^\prime(0,0)x_{2,k}\\
    &=O\bl\sum_{i=r+1}^k\theta_1\sech^2(h_i)x_ky_k\br.
\end{align*}}
Noticing that $\bb E[\ca B]=O\lef((1-c)^2k^{\frac{2\tau-2}{2\tau-1}}\log^{\frac{1}{2\tau-1}}(m\vee k)\rig)$ and the fact that $\ca B$ is also bounded and hence Sub-Gaussian, we conclude that for some $C>0,\lambda >1$  we have
\begin{align*}
    \bb E[\exp(\lambda(\ca B-\bb E[\ca B])])]\leq\exp((1-c)^2k^{\frac{4\tau-4}{2\tau-1}}\wedge (1-c)k^{\frac{2\tau-2}{2\tau-1}}\log^{\frac{1}{2\tau-1}}(m\vee k)).
\end{align*}
Then finally we consider $(1-c)=\Theta(1)$, it is   checked that for some constant $C,C_2>0$:
\begin{align*}
    \bb E\lef[\exp(\lambda\ca B)\rig]\leq\exp\lef(Ck^{\frac{2\tau-2}{2\tau-1}}\log^{\frac{1}{2\tau-1}}(m\vee k)\rig)\leq \exp\lef(C_2(1-c)k^{\frac{2\tau-2}{2\tau-1}}\log^{\frac{1}{2\tau-1}}(m\vee k)\rig).
\end{align*}
Similar to \eqref{boundingprocedure}, collecting the above pieces, there exists $C>0$ such that
\begin{align*}
    \bb E\lef[\frac{\bb P^*_{S}(\bfa\sigma)\bb P^*_{S^\prime}(\bfa\sigma)}{\bb P^*_0(\bfa\sigma)}\rig]=\exp(C(1-c)^2k^{\frac{4\tau-4}{2\tau-1}}\wedge (1-c)k^{\frac{2\tau-2}{2\tau-1}}\log^{\frac{1}{2\tau-1}}(m\vee k)).
\end{align*}

Then we consider the case with $k=\Omega(\sqrt n)$. It suffices to check $m=o\lef(\frac{n}{k^{5/3}}\rig)$ for the two competing terms $\frac{mk^{5/3}}{n}$ and $m(1-c)^2k^{\frac{2\tau-2}{2\tau-1}}\log^{\frac{1}{2\tau-1}}(k)$ separately. Essentially for the first one  we have $\ca G=1+o(1)$.  We notice that for some small $\epsilon>0$, there exists $C,C_1$ such that for some $\epsilon\in(0,1)$ the following holds
\tny{\begin{align*}
    \ca G&=\sum_{v=1}^{ k}\bb P(V=v)E_k^m(v)=\ub{\sum_{v=1}^{\epsilon k}\bb P(V=v)E_k^m(v)}_{T_1}+\ub{\sum_{v=\epsilon k+1}^k\bb P(V=v)E_k^m(v)}_{T_2}\\
    &\leq\sum_{v=\epsilon k+1}^k\bb P(V=v)\exp\lef(\frac{mv}{k^{\frac{1}{2\tau-1}}}\log^{\frac{1}{2\tau-1}}k\rig)
    +\sum_{v=1}^{\epsilon k}\frac{1}{(1-\frac{p}{k})\sqrt{2\pi p}}\exp\bigg(\lef(\frac{4k}{n}-\frac{p}{n}-\log\frac{pn}{k^2}-1\rig)p\\
    &-\frac{2k^2}{n}-2(k-p)\log\bigg(1-\frac{p}{k}\bigg)-\frac{1}{12p+1}+C(1-c)^2k^{\frac{4\tau-4}{2\tau-1}}\wedge(1-c)k^{\frac{2\tau-2}{2\tau-1}}\log^{\frac{1}{2\tau-1}}(m\vee k)\bigg)\\
    &+o(1).
\end{align*}}
 
And we treat the term $T_1$ using similar argument as \eqref{reduceepsk}.
 The strategy is to approximate address $T_2$ is to approximate it with Riemannian integral. Introducing $\gamma=\frac{k}{n}$ and similar to \eqref{riemann},  we have
\begin{align*}
    T_2 =\int_{(\frac{1}{k},\epsilon)}\frac{\sqrt k}{(1-x)\sqrt{2\pi x}}\exp\lef(k \mca f(x)\rig)dx(1+o(1)).
\end{align*}
where $\mca f$ is defined by
\begin{align*}
    \mca f(x):&=\lef((4-x)\gamma-\log\frac{x}{\gamma}-1\rig)x-2\gamma-2(1-x)\log\lef(1-x\rig) \\
    &+C_1x{m}{k^{-\frac{1}{2\tau-1}}}\log^{\frac{1}{2\tau-1}}(m\vee k)\wedge x^2mk^{\frac{2\tau-3}{2\tau-1}}.
\end{align*}
Applying Laplace method in lemma \ref{laplaceu1}, \ref{laplaceu2} again, we note that the derivative can be written as 
\tny{\begin{align*}
&\mca f^\prime(x) =  (4-2x)\gamma-\log\frac{x}{\gamma}+2\log(1-x)+2xC_1{m}{k^{\frac{2\tau-3}{2\tau-1}}},&&\text{ when }1-c=o(k^{-\frac{2\tau-2}{2\tau-1}}),\\
    &\mca f^{(2)}(x)=-2\gamma-\frac{1}{x}-\frac{2}{1-x}+C_1mk^{\frac{2\tau-3}{2\tau-1}}.\\
    &\mca f^\prime(x) =  (4-2x)\gamma-\log\frac{x}{\gamma}+2\log(1-x)+C_1{m}{k^{-\frac{1}{2\tau-1}}}\log^{\frac{1}{2\tau-1}}k,&&\text{ when }1-c=\Omega(k^{-\frac{2\tau-2}{2\tau-1}}),\\
    &\mca f^{(2)}(x)=-2\gamma-\frac{1}{x}-\frac{2}{1-x}<0.
\end{align*}}
And the stationary point $x^*$ (which is also maximum) satisfies
$f^\prime(x^*)=0$,
which admits the only solution $x^*=\gamma(1+o(1))=\omega\lef(\frac{1}{k}\rig)$ when $n^{\frac{4\tau-2}{8\tau-5}}\lesssim k\lesssim n^{\frac{2\tau-1}{4\tau-3}}$ and $m=o\lef(n^2k^{-\frac{2(4\tau-3)}{2\tau-1}}\rig)$. Then we use Laplace method for interior point in lemma \ref{laplaceu1} to get that
\begin{align*}
    T_2=\frac{1}{(1-\gamma)}\exp(\mca f(x^*))\to 1.
\end{align*}
Similarly  we have when $k=o(n^{\frac{4\tau-2}{8\tau-5}})$, $m=o\lef(\lef(\frac{k}{\log k}\rig)^{\frac{1}{2\tau-1}} \log n\rig)$ and complete the proof.

\subsection{Proof of Theorem \ref{guaranteecriticaltest}}
The proof is, as usual, be divided by the local part and the global part
\begin{center}
    \textbf{1. The Local Part}
\end{center}
    We first recall from theorem \ref{cltrfcr}, for sufficently large $k$, there exists $C>0$ such that for all $t>0$ not dependent on $k$:
\begin{align}\label{tailboundsupergaussiancon}
    \bb P_{S_0}\lef(\lef|k^{\frac{1}{4\tau-2}}m_{S_0}\rig|\geq t\rig)\leq 2\exp\lef(-Ct^{4\tau-2}\rig).
\end{align}
And alternatively,  we have for all $t>0$ and large $k$, for some $C>0$ we have
\begin{align*}
    \bb P_{S_0}\lef( \lef|k^{\frac{1}{2\tau-1}}m_{S_0}^2-\bb E\lef[ k^{\frac{1}{2\tau-1}}m_{S_0}^2\rig]\rig|\geq t\rig)\leq \exp(-Ct^{2\tau-1}),
\end{align*}
which implies that $\Vert k^{\frac{1}{2\tau-1}}m_{S_0}^2\Vert_{\psi_{2\tau-1}}<\infty$.
And using the uniform integrability of measure, under the alternative hypothesis,  we have
\begin{align}\label{complicate1}
    \bb E_{S_0}\lef[k^{\frac{1}{2\tau-1}}m_{S_0}^2\rig]\to \frac{\int_{\bb R}x^{2\tau}\exp\lef(-\frac{x^{4\tau-2}}{2\ca V}\rig)dx}{\int_{\bb R}x^{2\tau-2}\exp\lef(-\frac{x^{4\tau-2}}{2\ca V}\rig)dx}=\frac{(2\ca V(\tau))^{\frac{1}{2\tau-1}}\Gamma(\frac{2\tau+1}{4\tau-2})}{\sqrt{\pi}}(1+o(1)), 
\end{align}
where $\ca V(\tau)$ is defined by \eqref{stirlingv}.
And we also notice that under the null hypothesis  we have $$\bb E_0[k^{\frac{1}{2\tau-1}}m_{S_0}^2]=k^{-\frac{2\tau-2}{2\tau-1}}\to 0.$$
The following lemma gives concentration of i.i.d. sub-Weibull r.v.s.
\begin{lemma}[\citep{zhang2022sharper}]\label{weibullnormsubadditivity}
    For i.i.d. centered random variables $X_1,\ldots X_n$ such that $\Vert X_i\Vert_{\psi_{\theta}}<\infty$ for some $\theta>2$, there exists $C_1,C_2>0$ such that
    \begin{align*}
        \bb P\bl\bigg|\frac{1}{n}\sum_{i=1}^nX_i\bigg|\geq t\br\leq 2\exp\bl-\frac{C_1nt^\theta}{\Vert X\Vert^\theta_{\psi_\theta}}\wedge\frac{C_2nt^2}{\Vert X\Vert_{\psi_\theta}^2}\br.
    \end{align*}
\end{lemma}

Recall from \eqref{h0uni}  we have under the null hypothesis there exists $C>0$ such that
\begin{align}\label{phi411}
    \bb P_0\lef(\phi_5-\bb E_{0}[\phi_5]\geq \tau_\delta\rig)\leq\bl\frac{en}{k}\br^k \exp\lef(-Cmk^{\frac{4\tau-4}{2\tau-1}}\tau_\delta^2\wedge k^{\frac{2\tau-2}{2\tau-1}}\tau_\delta m\rig).
\end{align}
And under the alternative hypothesis there exists $C>0$ such that for $t>0$
\begin{align}\label{phi421}
    \bb P_{S_0}(\phi_5\leq \bb E[\phi_5]-t)\leq \exp(-Cmt^{2\tau-1}\wedge mt^2 ).
\end{align}
Therefore, let both \eqref{phi411} and \eqref{phi421} less than $\delta/2$, we let $\tau_\delta\in\lef(0,\frac{(2\ca V(\tau))^{\frac{1}{2\tau-1}}\Gamma(\frac{2\tau+1}{4\tau-2})}{\sqrt{\pi}}\rig)$ and let $m\gtrsim k^{\frac{1}{2\tau-1}}\log n$.
\begin{center}
    \textbf{2. The Global Part}
\end{center}
The proof utilizes the fact that under the null hypothesis the spins are i.i.d. Rademacher random variables. Then  we have $\lef\Vert\frac{1}{n}\lef(\sum_{i=1}^n\sigma_i\rig)^2\rig\Vert_{\psi_1}<\infty$, $\bb E_0[\phi_6]=0$. Then by Bernstein inequality, there exists $C_1,C_2>0$ such that
  \tny{  \begin{align*}
        \bb P_0(\phi_6\geq\tau_\delta)&=\bb P_0\bl\frac{1}{mn}\sum_{j=1}^m\bl\sum_{i=1}^n\sigma_i^{(j)}\br^2-1\geq\frac{k^{\frac{4\tau-3}{2\tau-1}}}{n}\tau_\delta\br\\
        &\leq C_1\exp\bl-C_2\tau_\delta^2\frac{mk^{\frac{2(4\tau-3)}{2\tau-1}}}{n^2}\wedge \tau_\delta\frac{mk^{\frac{4\tau-3}{2\tau-1}}}{n}\br.
    \end{align*}}
    Then, we study the Type II error. Recall that \eqref{complicate1},  we have
    \begin{align*}
        \bb E_{S}[\phi_6]=\pi^{-\frac{1}{2}}(2\ca V(\tau))^{\frac{1}{2\tau-1}}\Gamma\lef(\frac{2\tau+1}{4\tau-2}\rig)+o(1).
    \end{align*}
    And we also use the sub-additivity of Orlicz norm
    \begin{align*}
        \bigg\Vert k^{-\frac{4\tau-3}{4\tau-2}}\wedge n^{-1/2}\sum_{i=1}^n\sigma_i\bigg\Vert_{\psi_2}\leq \bigg\Vert k^{-\frac{4\tau-3}{4\tau-2}}\sum_{i\in S}\sigma_i\bigg\Vert_{\psi_2}+\bigg\Vert n^{-1/2}\sum_{i\in S^c}\sigma_i\bigg\Vert_{\psi_2}<\infty.
    \end{align*}
    And we then have
    \begin{align*}
        \bigg\Vert k^{-\frac{4\tau-3}{2\tau-1}}\wedge n^{-1}\bl\bl\sum_{i=1}^n\sigma_i\br^2-\bb E\bigg[\bl\sum_{i=1}^n\sigma_i\br^2\bigg]\br\bigg\Vert_{\psi_1}<\infty.
    \end{align*}
    Therefore, when $k=O\lef(n^{\frac{2\tau-1}{4\tau-3}}\rig)$,  we have
    \begin{align*}
        \bb P_{S}(\phi_6-\bb E_S[\phi_6]\leq -\tau_\delta)&=\bb P_S\bl m^{-1}k^{-\frac{4\tau-3}{2\tau-1}}\sum_{j=1}^m\bl\bl\sum_{i=1}^n\sigma_i\br^2-\bb E\bigg[\bl\sum_{i=1}^n\sigma_i\br^2\bigg]\br\leq -\tau_\delta\br\\
        &\leq \exp\lef(-C_2\tau_\delta^2n^{-2}mk^{\frac{2(4\tau-3)}{2\tau-1}}\wedge \tau_\delta n^{-1}mk^{\frac{4\tau-3}{2\tau-1}}\rig).
    \end{align*}
    And when $k=\Omega\lef(n^{\frac{2\tau-1}{4\tau-3}}\rig)$,  we have
    \begin{align*}
         \bb P_{S}(\phi_6-\bb E_S[\phi_6]\geq\tau_\delta)&\leq \exp\lef(-C_2m\tau_\delta\wedge m\tau_\delta^2\rig).
    \end{align*}
    Therefore, for $k=\omega(n^{\frac{2\tau-1}{4\tau-3}})$, $m\asymp 1$ is enough. For $k\lesssim n^{\frac{2\tau-1}{4\tau-3}}$ we need $m\gtrsim n^2k^{-\frac{2(4\tau-3)}{2\tau-1}}$.

\subsection{Proof of Corollary \ref{recoverycrtical}}
    The proof goes by first analyzing the general $m_{S^\prime}:=\frac{1}{ k}\sum_{i\in S^\prime}\sigma_i$ with $S^\prime:=[r+1:r+k]$ with $r=ck$. 
    We can decompose $m_{S^\prime}=m_{rk}+m_r$ as $m_{rk}:=\frac{1}{k}\sum_{i=r+1}^k\sigma_i$ and $ m_{r}^\prime:=\frac{1}{k}\sum_{i=k+1}^{k+r}\sigma_i$. Notice that $m_r^\prime\perp m_{rk}$. We first analyze the mgf of $k^{\frac{1}{4\tau-2}}m_{rk}$ under the RFCW model using similar method as \eqref{integraltransform},
    \tny{\begin{align*}
        \bb E\bigg[\exp\bl t\frac{\sum_{i=r+1}^k\sigma_i}{k^{\frac{4\tau-3}{4\tau-2}}}\br\bigg]=\bb E\bigg[\frac{\int_{\bb R}\exp(-n\mca H_{0,k}(x))dx}{\int_{\bb R}\exp(-n \mca H_{1,k}(x))dx}\bigg]=\bb E[\exp(-n\mca H_{0,k}(x_0)+n\mca H_{1,k}(x_1))]\lef(1+o(1)\rig),
    \end{align*}}
     where we define the following
    \begin{align*}
        \mca H_{0,k}(x):&=\frac{1}{2}x^2-\frac{1}{k}\sum_{i=1}^{r}\log\cosh\lef(\sqrt{\theta_1}x+h_i\rig)-\frac{1}{k}\sum_{r+1}^k\log\cosh\lef(\sqrt{\theta_1}x+h_i+\frac{t}{k^{\frac{4\tau-3}{4\tau-2}}}\rig),\\
        \mca H_{1,k}(x):&=\frac{1}{2}x^2-\frac{1}{k}\sum_{i=1}^{k}\log\cosh\lef(\sqrt{\theta_1}x+h_i\rig),
    \end{align*} and $x_0,x_1$ to be the global minimum of $\mca H_{0,k}$ and $\mca H_{1,k}$ respectively. We notice that by the equicontinuity of $\mca H_0,\mca H_1$, uniformly  we have
    \tny{\begin{align*}
        &\mca H_{0,k}(x)\to \mca H_0(x):=\frac{1}{2}x^2-c\bb E[\log\cosh(\sqrt{\theta_1}x+h)]-(1-c)\bb E\lef[\log\cosh\lef (\sqrt{\theta_1}x+h+\frac{t}{k^{\frac{4\tau-3}{4\tau-2}}}\rig)\rig],\\
        &\mca H_{1,k}(x)\to \mca H_1(x):=\frac{1}{2}x^2-\bb E[\log\cosh(\sqrt{\theta_1}x+h)].
    \end{align*}}
    Denote $x_0^*, x_{1}^*$ to be the global minimum of $\mca H_{0,k}$ and $\mca H_{1,k}$ respectively, notice that similarly to \eqref{weakconvergencex1}, by the inearization lemma \ref{linearZ} we have
   \tny{\begin{align*}
        \sqrt k(x_1-x_1^*)^{2\tau-1}
        =\frac{-(2\tau)!\sqrt{\theta_1}}{\sqrt nH^{(2\tau)}_{1,n}(x_1^*)}\sum_{i=1}^n\lef(\tanh(\sqrt{\theta_1}x_1^*+h_i)-\bb E[\tanh(\sqrt{\theta_1}x_1^*+h_i)]\rig)+o_{\psi_2}(1).
    \end{align*}}
    Therefore, one will get the following 
    \tny{\begin{align*}
        n\mca H_{1,k}(x_1)-n\mca H_{0,k}(x_0)=\frac{(1-c)^{\frac{4\tau-3}{4\tau-2}}t\sqrt{\theta_1}}{k-r}\bl\sum_{i=r+1}^k\sech^2(h_i)\br(k-r)^{\frac{1}{4\tau-2}}(x_1^*-x_1)+o_{\psi_2}(1).
    \end{align*}}
    Therefore, following similar path as the derivation of \eqref{fromweaktomgf} we get
    \tny{\begin{align*}
        \bb E\bigg[\exp\bl t\frac{\sum_{i=r+1}^k\sigma_i}{k^{\frac{4\tau-3}{4\tau-2}}}\br\bigg]\to\int_{\bb R}\frac{(2\tau-1)x^{2\tau-2}}{\sqrt{2\pi v}}\exp\bl-\frac{x^{4\tau-2}}{2v}+t(1-c)^{\frac{4\tau-3}{4\tau-2}}\sqrt{\bb E[\sech^2(h)]}x\br dx.
    \end{align*}}
    with $v:=\frac{((2\tau)!)^2\theta_1^{2\tau}\bb V(\tanh(\sqrt{\theta_1}x_1^*+h))(\bb E[\sech^2(\sqrt{\theta_1}x_1^*+h)])^{4\tau-2}}{(\mca H_{1}^{(2\tau)}(x_1^*))^2}$.
    Therefore, we also have
\begin{align*}
    \bb E[k^{\frac{1}{2\tau-1}}m_{rk}^2]\to \frac{\int_{\bb R}x^{2\tau}\exp\lef(-\frac{x^{4\tau-2}}{2\ca V}\rig)dx}{\int_{\bb R}x^{2\tau-2}\exp\lef(-\frac{x^{4\tau-2}}{2\ca V}\rig)dx}=\frac{(2V(c))^{\frac{1}{2\tau-1}}\Gamma(\frac{2\tau+1}{4\tau-2})}{\sqrt{\pi}}(1+o(1)),
\end{align*}
where we define $V(c):=(1-c)\ca V(\tau)$ for $\ca V(\tau)$ defined by \eqref{stirlingv}.
Moreover, it is also checked that 
$    \lef\Vert k^{\frac{1}{4\tau-2}}m_{rk}\rig\Vert_{\psi_{4\tau-2}}<\infty,
$ which implies that $\lef\Vert k^{\frac{1}{2\tau-1}}m^2_{rk}\rig\Vert_{\psi_{2}}<\infty$.
Then for the $m_r^\prime$,  we have by \eqref{mrksmalzero}:
\begin{align*}
    \frac{1}{c}\sqrt km_r^\prime\overset{d}{\to} N(0,1)\quad\text{ and }\quad \Vert\sqrt km_r^\prime\Vert_{\psi_2}<\infty, \quad \bb E[k^{\frac{1}{2\tau-1}}m_r^{\prime 2}]\to 0.
\end{align*}
Therefore  we have by the independence between $m_{rk}$ and $m_r^\prime$:
\begin{align*}
    \bb E[k^{\frac{1}{2\tau-1}}m_{S^\prime}^2]=\bb E[k^{\frac{1}{2\tau-1}}m_{rk}^2]+\bb E[k^{\frac{1}{2\tau-1}}m_r^{\prime 2}]<\bb E[k^{\frac{1}{2\tau-1}}m_{S}^2].
\end{align*}
And using the Sub-additivity of Orlicz norm,  we have \\$\Vert k^{\frac{1}{4\tau-2}}m_{S^\prime}\Vert_{\psi_2}\leq\Vert k^{\frac{1}{4\tau-2}}m_{rk}\Vert_{\psi_2}+\Vert k^{\frac{1}{4\tau-2}}m_{r}^\prime\Vert_{\psi_2}<\infty$. Then, $\Vert k^{\frac{1}{2\tau-1}}m^2_{S^\prime}\Vert_{\psi_1}<\infty$ we also have the following for some constant $C>0$,
\begin{align*}
    \bb P(|\phi_{5,S^\prime}-\bb E[\phi_{5,S^\prime}]|\geq t)\leq\exp(-Cmt^2\wedge mt).
\end{align*}
Finally following similar procedure of \eqref{overlapprob} we obtain the final result.

\subsection{Proof of Theorem \ref{lowerboundforrecovery}}
    We notice that the tensorization property of the mutual information  we have
    \begin{align*}
        I(S;\{\bfa\sigma^{(i)}\}_{i\in[m]})&= mI(\bfa\sigma;S)= m D_{kl}(\bb P(\bfa\sigma|S)\Vert\bb P(\bfa\sigma)|S).
    \end{align*}
    And then we use the fact that $\log x\leq x-1$ and by convexity,

    \tny{\begin{align*}
        D_{kl}(\bb P^*(\bfa\sigma)\Vert\bb P(\bfa\sigma|S)|S)&=\bb E\bigg[\sum_{S\in\ca S}\frac{1}{|\ca S|}\sum_{\bfa\sigma}\bb P^*(\bfa\sigma|\bfa h)\log\frac{\bb P^*(\bfa\sigma|\bfa h)}{\bb P(\bfa\sigma |S,\bfa h)}\bigg]\\
        &\leq \bb E\bigg[\sum_{S\in\ca S}\frac{1}{|\ca S|}\sum_{\bfa\sigma}\bb P^*(\bfa\sigma|\bfa h)\bl\frac{\bb P^*(\bfa\sigma|\bfa h)}{\bb P(\bfa\sigma|S,\bfa h)}-1\br\bigg]\\
        &=\frac{|\ca S|-1}{|\ca S|}\ub{\bb E\bigg[\sum_{\bfa\sigma}\frac{\bb P(\bfa\sigma|S_1,\bfa h)\bb P(\bfa\sigma|S_2,\bfa h)}{\bb P(\bfa\sigma|S_3,\bfa h)}\bigg]}_{T_1}+\frac{1}{|\ca S|}\ub{\bb E\bigg[\sum_{\bfa\sigma}\frac{\bb P^2(\bfa\sigma|S_1,\bfa h)}{\bb P(\bfa\sigma|S_2,\bfa h)}\bigg]}_{T_2}-1.
    \end{align*}}
    From here on we analyze the above quantities ($T_1$ and $T_2$) according to their temperature regimes. Note that without loss of generality we assume that $S_1=[k]$, $S_2=[k-1]\cup\{k+1\}$, and $S_3=[k-1]\cup\{k+2\}$.
    
    We use the H-S quantity $T_1$ can be analyzed as
    \tny{\begin{align}\label{targetkl}
        &\bb E\bigg[\frac{\sum_{\bfa\sigma}\exp\lef(\frac{\theta_1k}{2}(m_{S_1}^2+m_{S_2}^2-m_{S_3}^2)+\sum_{i=1}^{k+1}h_i\sigma_i-h_{k+2}\sigma_{k+2}\rig)\sum_{\bfa\sigma}\exp\lef(\frac{\theta_1k}{2}m_{S_3}^2+\sum_{i\in S_3}\sigma_ih_i\rig)}{\sum_{\bfa\sigma}\exp\lef(\frac{\theta_1k}{2}m_{S_1}^2+\sum_{i\in S_1}\sigma_ih_i\rig)\sum_{\bfa\sigma}\exp\lef(\frac{\theta_1 k}{2}m_{S_2}^2+\sum_{i\in S_2}\sigma_ih_i\rig)}\bigg]\nnb\\
        &=\frac{k}{2\pi}\bb E\bigg[\frac{\int_{\bb R^3}\exp(-kG_{0,k}(x,y,z,\bfa h))dxdydz\int_{\bb R}\exp(-kG_{3,k}(x,\bfa h))dx}{\int_{\bb R}\exp(-kG_{1,k}(x,\bfa h))dx\int_{\bb R}\exp(-kG_{2,k}(x,\bfa h))dx}\bigg],
    \end{align}}
    with uniformly in $x,y,z$,
    \tny{\begin{align*}
        G_{0,k}(x,y,z,\bfa h):&=\frac{x^2+y^2+z^2}{2}-\frac{1}{k}\sum_{i=1}^{k-1}\log\cosh(\sqrt{\theta_1}(x+y+\mca iz)+h_i)-\frac{1}{k}\log\cosh(\sqrt{\theta_1}x+h_{k})\\
        &-\frac{1}{k}\log\cosh(\sqrt{\theta_1}y+h_{k+1})-\frac{1}{k}\log\cosh(\mca i\sqrt{\theta_1}z-h_{k+2}),\\
        G_{0,k}(x,y,z,\bfa h)&\overset{a.s.}{\to}\frac{x^2+y^2+z^2}{2}-\bb E[\log\cosh(\sqrt{\theta_1}(x+y+\mca iz)+h)]=:G_0(x,y,z).\\
        G_{j,k}(x,\bfa h):&=\frac{x^2}{2}-\frac{1}{k}\sum_{i\in S_j}\log\cosh(\sqrt{\theta_1}x+h_i)\\
        &\overset{a.s.}{\to} G_2(x):=\frac{x^2}{2}-\bb E[\log\cosh(\sqrt{\theta_1}x+h)]\quad\text{ for }j\in\{1,2,3\}.
    \end{align*}}
    Since here we are involved in the discussion of complex integral in $G_{(0,k)}$, a more natural method to use is the method of the steepest descent (see, for example in \citep{bleistein1975asymptotic}).  Introducing
$ 
         f_i=\frac{\sqrt{\theta_1}}{k}\tanh(\sqrt{\theta_1}(x_k^*+y_k^*+\mca iz_k^*)+h_i)
$, $ 
         f^\prime_i=\frac{\sqrt{\theta_1}}{k}\tanh(\sqrt{\theta_1}(x_k^*+y_k^*+\mca iz_k^*)-h_i)
$  and the complex stationary points $(x_k^*,y_k^*,z_k^*)$ of $G_{0,k}$ is given by
    \begin{align*}
    \begin{bmatrix}
        \frac{\pta G_{0,k}}{\pta x}\\
        \frac{\pta G_{0,k}}{\pta y}\\
        \frac{\pta G_{0,k}}{\pta z}
    \end{bmatrix}
        =\begin{bmatrix}
            x_k^*-\sum_{i\in[k]}f_i\\
            y_k^*-\sum_{i\in[k-1]\cup\{k+1\}}f_i\\
            z_k^*-\mca i\lef(\sum_{i\in[k-1]}f_i+f_{k+2}^\prime\rig)
        \end{bmatrix}=\bfa 0.
    \end{align*}
    And consider the population version of equation, and the limit of $(x_k^*,y_k^*,z_k^*)\to(x^*,y^*,z^*)$,  we have
    \begin{align}\label{populationstationarypoints}
        \begin{bmatrix}
        \frac{\pta G_{0}}{\pta x}\\
        \frac{\pta G_{0}}{\pta y}\\
        \frac{\pta G_{0}}{\pta z}
    \end{bmatrix}
        =\begin{bmatrix}
            x^*-\sqrt{\theta_1}\bb E[\sech^2(\sqrt{\theta_1}(x^*+y^*+\mca i z^*)+h)]\\
            y^*-\sqrt{\theta_1}\bb E[\sech^2(\sqrt{\theta_1}(x^*+y^*+\mca i z^*)+h)]\\
            z^*-\mca i\sqrt{\theta_1}\bb E[\sech^2(\sqrt{\theta_1}(x^*+y^*+\mca i z^*)+h)]
        \end{bmatrix}=\bfa 0.
    \end{align}
    which further implies that
    \begin{align*}
        x^*+y^*+\mca iz^*=\sqrt{\theta_1}\bb E[\tanh(\sqrt{\theta_1}(x^*+y^*+\mca iz^*)+h)].
    \end{align*}
    And at high/critical temperature it is not hard to see that in the real domain, the only solution to the above equation is is $x^*+y^*+\mca iz^*=0$. And at the low temperature  we have $x^*+y^*+\mca iz^*$ can take two values being symmetric w.r.t. $0$. The reason for taking only real valued ones is due to the admissibility of stationary points, which stays on the (distorted) integral path using the method of the steepest descent.
    \begin{center}
        \textbf{1. High Temperature}
    \end{center}
    And we check that the admissible stationary point is $\bfa 0$ here, which gives the global minimum.
   Then we introduce $ 
         \tde g_i=\frac{\sqrt{\theta_1}}{k}\tanh(h_i)
$ and $ 
          g_i=\frac{\theta_1}{k}\sech^2(h_i)
$. Then the Hessian can be computed as
    \begin{align*}
        \nabla G_{0,k}(0,0,0)&=-\begin{bmatrix}
            \sum_{i\in S_1}\tde g_i\\
            \sum_{i\in S_2}\tde g_i\\
            \mca i\sum_{i\in [k-1]}\tde g_i-\mca i\tde g_{k+2}
        \end{bmatrix},\\
        \nabla^2G_{0,k}(0, 0, 0)&=\begin{bmatrix}
            1-\sum_{i=1}^k g_i&-\sum_{i=1}^{k-1} g_i&\mca -i\sum_{i=1}^{k-1} g_i\\
            -\sum_{i=1}^{k-1} g_i& 1-\sum_{i\in[k-1]\cup\{k+1\}} g_i&-\mca i\sum_{i=1}^{k-1} g_i\\
            -\mca i\sum_{i=1}^{k-1} g_i& -\mca i\sum_{i=1}^{k-1} g_i&1+\sum_{i\in[k-1]\cup\{k+2\}} g_i
        \end{bmatrix}.
    \end{align*}
    Then we check that,
    \tny{\begin{align*}
       \nabla^2G_{0,k}(0,0,0)\overset{a.s.}{\to}\nabla^2G_0(0,0,0)=\begin{bmatrix}
            1-\theta_1\bb E[\sech^2(h)]&-\theta_1\bb E[\sech^2(h)] &-\mca i\theta_1\bb E[\sech^2(h)]\\
            -\theta_1\bb E[\sech^2(h)]&1-\theta_1\bb E[\sech^2(h)]&-\mca i\theta_1\bb E[\sech^2(h)]\\
            -\mca i\theta_1\bb E[\sech^2(h)]&-\mca i\theta_1\bb E[\sech^2(h)]&1+\theta_1\bb E[\sech^2(h)]
        \end{bmatrix}.
    \end{align*}}
    And we also have
    \tny{\begin{align*}
       \Vert\nabla^2 G_{0,k}(\bfa 0)-\nabla^2G_0(\bfa 0) \Vert_2\leq\Vert\nabla^2 G_{0,k}(\bfa 0)-\nabla^2G_0(\bfa 0) \Vert_F=O_{\psi_2}(k^{-\frac{1}{2}}).
    \end{align*}}
    Taking its inverse  we have
    \tny{\begin{align}\label{inverse}
        \lef(\nabla^2 G_0(0,0,0)\rig)^{-1}=\frac{1}{1-\theta_1\bb E[\sech^2(h)]}\begin{bmatrix}
            -1&-\theta_1\bb E[\sech^2(h)] &-\mca i\theta_1\bb E[\sech^2(h)]\\
            -\theta_1\bb E[\sech^2(h)]&-1&-\mca i\theta_1\bb E[\sech^2(h)]\\
            -\mca i\theta_1\bb E[\sech^2(h)]&-\mca i\theta_1\bb E[\sech^2(h)]&2\theta_1\bb E[\sech^2(h)]-1
        \end{bmatrix}.
    \end{align}}
    And then by the a.s. boundedness of $\Vert(\nabla^2 G_{0}(\bfa 0))^{-1}\Vert_2$, $\Vert(\nabla^2 G_{0,k}(\bfa 0))^{-1}\Vert_2$  we have
    \tny{\begin{align}\label{inverseapprox}
       \Vert(\nabla^2 G_{0,k}(\bfa 0))^{-1}- (\nabla^2G_0(\bfa 0))^{-1}\Vert_2&\leq\Vert \nabla^2G_{0,k}(\bfa 0)-\nabla^2G_{0}(\bfa 0)\Vert_2\Vert (\nabla^2 G_{0,k}(\bfa 0))^{-1}\Vert_2\Vert (\nabla^2 G_{0}(\bfa 0))^{-1}\Vert_2\nnb\\
       &=O_{\psi_2}\lef(k^{-\frac{1}{2}}\rig).
    \end{align}}
    Then we Taylor expand at $(x,y,z)=\bfa 0$ and introduce $\bfa\delta:=(x^*,y^*,z^*)^\top$ to get that
    \begin{align*}
      \bfa 0=\nabla G_{0,k}(x^*,y^*,z^*)=\nabla G_{0,k}(0,0,0)+ \nabla^2G_{0,k}(0,0,0)\bfa\delta+o(\Vert\bfa\delta\Vert_2).
    \end{align*}
    And  we have $\det(\nabla^2G_{0,k})=(1-\sum_{i=1}^kg_i)(1+O(k^{-1}))$, this implies the invertibility a.s. at high temperature.
    And we further have 
    \begin{align*}
        \sqrt k\bfa\delta=-\lef(\nabla^2G_{0,k}(0,0,0)\rig)^{-1}\sqrt k\nabla G_{0,k}(0,0,0)(1+o(1)).
    \end{align*}
    And analogously we, define the global minimum point of $G_{j,k}$ to be $x_j^*$ ( It is also easily checked that there exists only $1$ global minimum a.a.s. And we can pick any sequence of stationary points converging to it. )
    \begin{align*}
        G_{j,k}^{(2)}(0)=1-\sum_{i\in S_j}g_j,\qquad \sqrt kx_j^*=-\frac{\sqrt kG_{j,k}^\prime(0)}{G_{1}^{(2)}(0)}(1+o(1)).
    \end{align*}
    Similar to \ref{inverseapprox} we also have for all $j\in\{1,2,3\}$,
    \begin{align}\label{replacehess}
       |G_{1}^{(2)}(0))^{-1}-(G_{j,k}^{(2)}(0))^{-1}|\leq|G_{1}^{(2)}(0)-G_{j,k}^{(2)}(0)|= O_{\psi_2}\lef(k^{-\frac{1}{2}}\rig).
    \end{align}
    Then by the method of the steepest descent, H\"older's inequality, \eqref{inverse}, and \eqref{inverseapprox}, there exists $\delta_1,\delta_2>0$ such that \eqref{targetkl} can be written as 
    \tny{\begin{align}\label{klsecond}
        \bb E\bigg[&\frac{\int_{\bb R^3}\exp(-kG_{0,k}(x,y,z,\bfa h))dxdydz\int_{\bb R}\exp(-kG_{3,k}(x,\bfa h))dx}{\int_{\bb R}\exp(-kG_{1,k}(x,\bfa h))dx\int_{\bb R}\exp(-kG_{2,k}(x,\bfa h))dx}\bigg]\nnb\\
        &=\bb E\bigg[\bl{G_{1,k}^{(2)}(0)G_{2,k}^{(2)}(0)}{\lef(G_{3,k}^{(2)}(0)\rig)^{-1}\lef(\det(\nabla^2G_{0,k}(\bfa 0))\rig)^{-1}}\br^{\frac{1}{2}}\nnb\\
        &\cdot\exp\bl -kG_{0,k}(x^*,y^*,z^*,\bfa h)-k G_{3,k}(x_3^*,\bfa h)+ kG_{1,k}(x_1^*,\bfa h)+kG_{2,k}(x_2^*,\bfa h)\br\bigg]( 1+O(k^{-1}))\nnb\\
        &\leq \bb E\bigg[\bl{G_{1,k}^{(2)}(0)G_{2,k}^{(2)}(0)}{\lef(G_{3,k}^{(2)}(0)\rig)^{-1}\lef(\det(\nabla^2G_{0,k}(\bfa 0))\rig)^{-1}}\br^{\frac{1}{2}(1+\delta_1)}\bigg]^{\frac{1}{1+\delta_1}}(1+O(k^{-1}))\nnb\\&
        \bb E\bigg[\exp\bl\frac{k}{2}(1+\delta_2)\bl -(G_{1,k}^{(2)}(0))^{-1}{(G^\prime_{1,k}(0))^2}-(G_{2,k}^{(2)}(0))^{-1}{(G^\prime_{2,k}(0))^2}\nnb\\&+\nabla G^\top_{0,k}(\bfa 0)\lef(\nabla^2G_{0,k}(\bfa 0)\rig)^{-1}\nabla G_{0,k}(\bfa 0)+(G_{3,k}^{(2)}(0))^{-1}{(G^\prime_{3,k}(0))^2}\br\bigg]^{\frac{1}{1+\delta_2}}.
    \end{align}}
    The term in the exponential is analyzed by noticing that after some algebraic manipulations,
    \begin{align*}
        \nabla G^\top_{0,k}&(\bfa 0)\lef(\nabla^2G_{0}(\bfa 0)\rig)^{-1}\nabla G_{0,k}(\bfa 0)+(G_{1}^{(2)}(0))^{-1}(G^\prime_{3,k}(0))^2-(G_{1}^{(2)}(0))^{-1}{(G^\prime_{1,k}(0))^2}\\
        &-(G_{1}^{(2)}(0))^{-1}{(G^\prime_{2,k}(0))^2}\\
        &=O\lef(\frac{1}{k}\rig) \lef(G_1^{(2)}(0)\rig)^{-1}G_{1,k}^\prime(0)+O_{\psi_2}\lef(\frac{1}{k^2}\rig).
    \end{align*}
    And by $\sqrt k G^\prime_{1,k}(0)$ being sub-Gaussian and centered,  we have for all $t>0$,
    \begin{align}\label{klpart3}
        \bb E\lef[\exp\lef( tG_1^{(2)}G_{1,k}^\prime(0)\rig)\rig]=\exp\lef(O(k^{-1})\rig).
    \end{align}
    And for the first term it is checked that
    \tny{\begin{align}\label{klpart1}
        {G_{1,k}^{(2)}(0)G_{2,k}^{(2)}(0)}{\lef(G_{3,k}^{(2)}(0)\rig)^{-1}\lef(\det(\nabla^2G_{0,k}(\bfa 0))\rig)^{-1}}=\frac{\lef(1-\sum_{i\in S_1}g_i\rig)\lef(1-\sum_{i\in S_2}g_i\rig)}{\lef( 1-\sum_{i\in S_3}g_i\rig)\lef(1-\sum_{i\in S_1}g_i\rig)}\lef(1+O(k^{-1})\rig).
    \end{align}}
    Therefore, collecting \eqref{klpart1} and \eqref{klpart3} we finally conclude that
    \begin{align*}
        T_1=1+O(k^{-1}).
    \end{align*}
    And we similarly, we derive that $T_0=1+O(k^{-1})$. This conclude that
    \begin{align*}
        D_{kl}(\bb P(\bfa\sigma)\Vert\bb P(\bfa\sigma|S)|S)= O(k^{-1})\quad\Rightarrow\quad \bb P(\wh S\neq S)\geq 1- O\lef(\frac{m}{k\log n}\rig)\vee 1.
    \end{align*}
    \begin{center}
        \textbf{ 2. Low Temperature }
    \end{center}
    We notice that at the low temperature regime, the admissible stationary point is $x^*=y^*=x_0$, $z^*=\mca ix_0$ with $x_0$ be the positive/negative root of  $x=\sqrt{\theta_1}\bb E[\tanh(\sqrt{\theta_1}x+h)]$. Since by symmetry, the two stationary points achieve the same function value, we discuss the $x_0>0$ case without loss of generality. Introducing $\tde q_i:=\frac{\sqrt{\theta_1}}{k}\tanh(\sqrt{\theta_1}x_0+h_i)$ and $q_i:=\frac{{\theta_1}}{k}\sech^2(\sqrt{\theta_1}x_0+h_i)$, we rewrite the gradient at $(x_0,x_0,\mca ix_0)$ as:
    \ttny{\begin{align*}
        \nabla G_{0,k}(x_0,x_0,\mca ix_0)&=\begin{bmatrix}x_0-\sum_{i\in S_1}\tde q_i\\
        x_0-\sum_{i\in S_2}\tde q_i\\
        \mca ix_0-\mca i\lef(\sum_{i\in [k-1]}\tde q_i-\tde q_{k+2}\rig)
        \end{bmatrix},\\
        \nabla^2G_{0,k}(x_0,x_0,\mca ix_0)&=\begin{bmatrix}
            1-\sum_{i=1}^k q_i&-\sum_{i=1}^{k-1}q_i&\mca -i\sum_{i=1}^{k-1} q_i\\
            -\sum_{i=1}^{k-1} q_i& 1-\sum_{i\in[k-1]\cup\{k+1\}} q_i&-\mca i\sum_{i=1}^{k-1} q_i\\
            -\mca i\sum_{i=1}^{k-1} q_i& -\mca i\sum_{i=1}^{k-1} q_i&1+\sum_{i\in[k-1]\cup\{k+2\}} q_i
        \end{bmatrix},\\
        \nabla^2G_0(x_0.x_0,\mca ix_0)&=\begin{bmatrix}
            1-\theta_1\bb E[\sech^2(\sqrt{\theta_1}x_0+h)]&-\theta_1\bb E[\sech^2(\sqrt{\theta_1}x_0+h)] &-\mca i\theta_1\bb E[\sech^2(\sqrt{\theta_1}x_0+h)]\\
            -\theta_1\bb E[\sech^2(\sqrt{\theta_1}x_0+h)]&1-\theta_1\bb E[\sech^2(\sqrt{\theta_1}x_0+h)]&-\mca i\theta_1\bb E[\sech^2(\sqrt{\theta_1}x_0+h)]\\
            -\mca i\theta_1\bb E[\sech^2(\sqrt{\theta_1}x_0+h)]&-\mca i\theta_1\bb E[\sech^2(\sqrt{\theta_1}x_0+h)]&1+\theta_1\bb E[\sech^2(\sqrt{\theta_1}x_0+h)]
        \end{bmatrix}.
    \end{align*}}
    And by Taylor expansion at $(x_0,x_0,\mca ix_0)$, defining $\bfa\delta=(x^*-x_0,y^*-x_0,z^*-\mca ix_0)^\top$,  we have
    \begin{align*}
      \bfa 0=\nabla G_{0,k}(x^*,y^*,z^*)=\nabla G_{0,k}(x_0,x_0,\mca ix_0)+ \nabla^2G_{0,k}(x_0,x_0,\mca ix_0)\bfa\delta+o(\Vert\bfa\delta\Vert_2).
    \end{align*}
    And then  we have
    \begin{align}\label{delta2}
        \sqrt k\bfa\delta=-\lef(\nabla^2G_{0,k}(x_0,x_0,\mca ix_0)\rig)^{-1}\sqrt k\nabla G_{0,k}(x_0,x_0,\mca ix_0)(1+o(1)).
    \end{align}
    Similarly, we check that the global minimum of $G_{j,k}$ for all $j\in\{1,2,3\}$ will all converges to $x_0/-x_0$. By similar symmetry, without loss of generality, we consider the $x_0$ case. We rewrite their derivatives as
    \begin{align*}
        G_{j,k}^\prime(x_0)=x_0-\sum_{i\in S_j}\tde q_i,\quad G_{j,k}^{(2)}(x_0)=1-\sum_{i\in S_j}q_i,\quad G_{1}^{(2)}(x_0)=1-\theta_1\bb E[\sech^2(\sqrt{\theta_1}x_0+h)].
    \end{align*}
    And by a similar way of Taylor expansion,  we have
    \begin{align}\label{xjdiff2}
        \sqrt k(x_j^*-x_0)=-\sqrt k(G_{j,k}^{(2)}(x_0))^{-1} G^\prime_{j,k}(x_0)\lef(1+o(1)\rig).
    \end{align}
    Therefore, using \eqref{delta2}, \eqref{xjdiff2}, \eqref{inverseapprox}, \eqref{replacehess}, we can rewrite \eqref{klsecond} as : (Notice that the only difference from the high temperature case is the non-zero mean of $\tanh(x_0+h)$.)
     \ttny{\begin{align}\label{klsecond2}
        \bb E\bigg[&\frac{\int_{\bb R^3}\exp(-kG_{0,k}(x,y,z,\bfa h))dxdydz\int_{\bb R}\exp(-kG_{3,k}(x,\bfa h))dx}{\int_{\bb R}\exp(-kG_{1,k}(x,\bfa h))dx\int_{\bb R}\exp(-kG_{2,k}(x,\bfa h))dx}\bigg]\nnb\\
        &=\bb E\bigg[\bl{G_{1,k}^{(2)}(x_0)G_{2,k}^{(2)}(x_0)}{\lef(G_{3,k}^{(2)}(x_0)\rig)^{-1}\lef(\det(\nabla^2G_{0,k}(x_0,x_0,\mca ix_0))\rig)^{-1}}\br^{\frac{1}{2}}\nnb\\
        &\cdot\exp\bl -kG_{0,k}(x^*,y^*,z^*,\bfa h)-k G_{3,k}(x_3^*,\bfa h)+ kG_{1,k}(x_1^*,\bfa h)+kG_{2,k}(x_2^*,\bfa h)\br\bigg]( 1+O(k^{-1}))\nnb\\
        &\leq \bb E\bigg[\bl{G_{1,k}^{(2)}(x_0)G_{2,k}^{(2)}(x_0)}{\lef(G_{3,k}^{(2)}(x_0)\rig)^{-1}\lef(\det(\nabla^2G_{0,k}(x_0,x_0,\mca ix_0))\rig)^{-1}}\br^{\frac{1}{2}(1+\delta_1)}\bigg]^{\frac{1}{1+\delta_1}}(1+O(k^{-1}))\nnb\\&
        \bb E\bigg[\exp\bl\frac{k}{2}(1+\delta_2)\bl -(G_{1,k}^{(2)}(x_0))^{-1}{(G^\prime_{1,k}(x_0))^2}-(G_{2,k}^{(2)}(x_0))^{-1}{(G^\prime_{2,k}(x_0))^2}\nnb\\
        &+\nabla G^\top_{0,k}(x_0,x_0,\mca ix_0)\lef(\nabla^2G_{0,k}(x_0,x_0,\mca ix_0)\rig)^{-1}\nabla G_{0,k}(x_0,x_0,\mca ix_0)+(G_{3,k}^{(2)}(x_0))^{-1}{(G^\prime_{3,k}(x_0))^2}\br\bigg]^{\frac{1}{1+\delta_2}}\nnb\\
        &=O(1).
    \end{align}}
    And we similarly derive that $T_2=O(1)$ and conclude that
    \begin{align*}
        D_{kl}(\bb P(\bfa\sigma)\Vert\bb P(\bfa\sigma|S)|S)=O(1)\quad\Rightarrow\quad \bb P(\wh S\neq S)\geq 1-O\lef(\frac{m}{\log n}\rig)\vee 1.
    \end{align*}
    \begin{center}
        \textbf{3. Critical Temperature}
    \end{center}
    At the critical temperature regime, the solution to \ref{populationstationarypoints} remains to be $\bfa 0$, but the Taylor expansion of stationary point will involve higher order derivatives. Introduce $\times$ as the notation of tensor product,  we have by Taylor expansion
    \begin{align*}
        \nabla G_{0,k}(x^*,y^*,z^*)=0=\nabla G_{0,k}(\bfa 0)+\sum_{\ell=2}^{2\tau}\frac{1}{(\ell-1)!}\nabla^{2\tau}G_{0,k}(\bfa 0)\bfa\times\bfa\delta^{\ell-1}+O(\Vert\bfa\delta\Vert_2^{\ell}).
    \end{align*}
    And we notice that the first $2\tau-1$ th order derivatives have 0 on the diagonals. This will cause the related term to be smaller in order than the first derivative having non-zeros on the diagonal.  Hence,  we have
    \begin{align*}
       -\sqrt k\nabla G_{0,k}(\bfa 0)=\frac{1}{(2\tau-1)!}\nabla^{2\tau}G_{0,k}(\bfa 0)\times(k^{\frac{1}{4\tau-2}}\bfa\delta)^{2\tau-1}(1+o(1)).
    \end{align*}
    Therefore  we have, noticing that the cross term will cancel aganist each other and is of order $O(\frac{1}{k})$,
    \ttny{\begin{align}\label{criticalxyzstar}
        k^{\frac{1}{2}}\begin{bmatrix}
            (x^*)^{2\tau-1}\\
            (y^*)^{2\tau-1}\\
            (z^*)^{2\tau-1}
        \end{bmatrix}= -\sqrt k(2\tau-1)!\frac{1}{G_{0,k}^{2\tau}(0)}\begin{bmatrix}
             G^{\prime}_{1,k}(0) \\
             G^{\prime}_{2,k}(0)\\
            \mca i G^{\prime}_{3,k}(0)
        \end{bmatrix}+o_{\psi_2}(1)\quad \Rightarrow\quad \Vert x^*\Vert_{\psi_2}\vee\Vert y^*\Vert_{\psi_2}\vee \Vert z^*\Vert_{\psi_2}\lesssim k^{-\frac{1}{4\tau-2}}.
    \end{align}}
    And analogously,  we have
    \begin{align*}
G_{j,k}^\prime(x_{j}^*)=0=G_{j,k}^\prime(0)+\sum_{\ell=2}^{2\tau}\frac{1}{(\ell-1)!}G_{j,k}^{(\ell)}(0)(x_j^*)^{\ell-1}+O((x_j^*)^{\ell}).
    \end{align*}
    And  we have
$        -\sqrt kG_{j,k}^\prime(0)=\frac{1}{(2\tau-1)!}G_{j,k}^{(2\tau)}(0)(k^{\frac{1}{4\tau-2}}x_j^*)^{2\tau-1}(1+o(1))
$. Then we go back to \eqref{klsecond} and notice that by H\"older's inequality, there exists $\delta_1>1,\delta_2>0$ with $\frac{1}{1+\delta_1}+\frac{1}{1+\delta_2}=1$ such that
\tiny{\begin{align}\label{criticalseparate}
    \bb E\bigg[&\frac{\int_{\bb R^3}\exp(-kG_{0,k}(x,y,z,\bfa h))dxdydz\int_{\bb R}\exp(-kG_{3,k}(x,\bfa h))dx}{\int_{\bb R}\exp(-kG_{1,k}(x,\bfa h))dx\int_{\bb R}\exp(-kG_{2,k}(x,\bfa h))dx}\bigg]\nnb\\
    &\leq\ub{\bb E\bigg[\exp\bl(1+\delta_1)\bl\ub{-kG_{0,k}(x^*,y^*,z^*,\bfa h)-kG_{3,k}(x_3^*,\bfa h)+k G_{1,k}(x_1^*,\bfa h)+kG_{2,k}(x_2^*,\bfa h)}_{\ca T_1}\br\br\bigg ]^{\frac{1}{1+\delta_1}}}_{\ca T_0}\nnb\\
    &\cdot\bb E\bigg[\bl\ub{\frac{\int_{\bb R^3}\exp(-k(G_{0,k}(x,y,z,\bfa h)-G_{0,k}(x^*,y^*,z^*,\bfa h)))dxdydz\int_{\bb R}\exp(-k(G_{3,k}(x,\bfa h)-G_{3,k}(x^*_3,\bfa h)))dx}{\int_{\bb R}\exp(-k(G_{1,k}(x,\bfa h)-G_{1,k}(x_1^*,\bfa h)))dx\int_{\bb R}\exp(-k(G_{2,k}(x,\bfa h)-G_{2,k}(x_2^*,\bfa h)))dx}}_{\ca T_2}\br^{1+\delta_2}\bigg]^{\frac{1}{1+\delta_2}}.
\end{align}}\normalsize
For the first term within the expectation can be analyzed as
\tny{\begin{align*}
    \ca T_1 &= -k\nabla G_{0,k}(\bfa 0)\times \bfa\delta -\frac{1}{(2\tau)!}k \nabla^{2\tau}G_{0,k}(\bfa 0)\times\bfa\delta^{2\tau}-k G^{(1)}_{3,k}(\bfa 0) x_{3}^*-\frac{1}{(2\tau)!}kG_{3,k}^{(2\tau)}(0) (x_{3}^*)^{2\tau}\\
    &+k G^{(1)}_{1,k}(\bfa 0) x_{1}^*+\frac{1}{(2\tau)!}kG_{1,k}^{(2\tau)}(0) (x_{1}^*)^{2\tau}+k G^{(0)}_{2,k}(\bfa 0) x_{2}^*+\frac{1}{(2\tau)!}kG_{2,k}^{(2\tau)}(0) (x_{2}^*)^{2\tau}+O_{\psi_2}\lef(k^{-\frac{1}{2\tau-1}}\rig)\\
    &=-\frac{2\tau-1}{2\tau}k\ub{\lef(\nabla G_{0,k}(\bfa 0)\times\bfa\delta+G_{3,k}^\prime(0)x_{3}^*-G^\prime_{1,k}(0)x_1^*-G^\prime_{2,k}(0)x_2^*\rig)}_{\ca T_3}+O_{\psi_2}\lef(k^{-\frac{1}{2\tau-1}}\rig).
\end{align*}}
It is also checked that by the symmetry of $h$'s measure, $\ca T_3$ also has symmetric measure w.r.t. $0$ and satisfies $\bb E[T^{2r+1}_3]=0$. And we notice that by the factorization of difference of powers,  we have 
\tny{\begin{align*}
    k\lef( G^\prime_{1,k}(0)x_1^*-G^\prime_{3,k}(0)x_3^*\rig)&\asymp k\lef(G_{1}^{(2)}(0)\rig)^{-1}\lef((G_{1,k}^\prime(0))^{1+\frac{1}{2\tau-1}}-(G_{3,k}^\prime(0))^{1+\frac{1}{2\tau-1}}\rig)\\
    &\asymp k\lef(G_{1}^{(2)}(0)\rig)^{-1}\lef(G_{1,k}^\prime(0)-G_{3,k}^\prime(0)\rig)\lef((G_{1,k}^\prime(0))^{\frac{1}{2\tau-1}}+(G_{3,k}^\prime(0))^{\frac{1}{2\tau-1}}\rig)\\
    &\asymp \lef(G_{1}^{(2)}(0)\rig)^{-1}\lef(|G_{1,k}^\prime(0)|^{\frac{1}{2\tau-1}}+|G_{3,k}^\prime(0)|^{\frac{1}{2\tau-1}}\rig).
\end{align*}}
And analogously, using \eqref{criticalxyzstar}  we have,
\tny{\begin{align*}
     k\lef( \nabla G_{0,k}(\bfa 0)\times\bfa\delta-G^\prime_{2,k}(0)x_2^*\rig)&\asymp \lef(G_{1}^{(2)}(0)\rig)^{-1}\lef(|G_{1,k}^\prime(0)|^{\frac{1}{2\tau-1}}+|G_{2,k}^\prime(0)|^{\frac{1}{2\tau-1}}+|G_{3,k}^\prime(0)|^{\frac{1}{2\tau-1}}\rig).
\end{align*}}
Then we use the fact that there exists constant $C>1$ such that for all $r\in\bb N$, by H\"older's on $\ca T_0$ to separate out $\exp(\ca T_1)$,
\begin{align}\label{t0bound}
    \bb E[T_1^{2r}]\lesssim C^{2r}\bb E[(G_{1,k}^\prime(0))^{\frac{2r}{2\tau-1}}]\leq C^{2r}k^{-\frac{r}{2\tau-1}}\quad\Rightarrow \quad \ca T_0=1+O(k^{-\frac{1}{2\tau-1}}).
\end{align}
Then the next step is to analyze each individual terms in $T_2$. Using Taylor expansion and defining $\bfa\delta_0:=(x-x^*,y-y^*,z-z^*)$, the first term is given by
\begin{align*}
    \ca T_{11}:&=\int_{\bb R^3}\exp(-k (G_{0,k}(x,y,z,\bfa h)-G_{0,k}(x^*,y^*,z^*,\bfa h)))dxdydz\\
    &=\int_{\bb R^3}\exp\bl -k\bl\sum_{\ell=2}^\infty \nabla^{\ell}G_{0,k}(x^*,y^*,z^*)\times\bfa\delta_0^{\ell}\br\br dxdydz.
\end{align*}
Notice that when $\ell>2$  we have
\ttny{\begin{align*}
    \nabla^\ell G_{0,k}(x^*,y^*,z^*)\times \bfa\delta^\ell_0=
        -G_{1,k}^{(\ell)}(x^*)(x+y+\mca iz-x^*-y^*-\mca iz^*)^\ell+O(k^{-1})((x-x^*)^{\ell}+(z-z^*)^{\ell}+(y-y^*)^{\ell}.
\end{align*}}
And when $\ell=2$  we have
\ttny{\begin{align*}
    \nabla^2 G_{0,k}(x^*,y^*,z^*)\times \bfa\delta^2_0&=
        G_{1,k}^{(2)}(x^*)(x+y+\mca iz-x^*-y^*-\mca iz^*)^2+O(k^{-1})((x-x^*)^{2}-(z-z^*)^{2}+(y-y^*)^{2})\\
        &+(x+\mca iz-x^*-\mca iz^*)^2+(x-y-x^*+y^*)^2.
\end{align*}}
And if we do the change of variables
\begin{align*}
    (x,y,z)\quad\to\quad(-\mca ix-\mca iy+ z,x+\mca i z,x-y):=(u,v,r).
\end{align*}
And it is easy to verify that the determinant of the Jacobian $\det(\ca J):=\det\lef(\frac{\pta (u,v,r)}{\pta (x,y,z)}\rig)=1$. Then we check that
\tny{\begin{align*}
   \ca T_{11}=\int_{\bb R}\exp\bl-k(G_{1,k}(x)-G_{1,k}(x^*))-k\frac{y^2+z^2}{2}\br\bl 1+O((x-x^*)^2+y^2+z^2) \br dxdydz.
\end{align*}}
And by Laplace method in lemma \ref{laplace} it is not hard to see that
\begin{align*}
   \ca  T_{11}=\frac{2\pi}{k}\int_{\bb R}\exp\bl -k (G_{1,k}(x)-G_{1,k}(x^*))\br dx(1+O(k^{-\frac{1}{2\tau}})).
\end{align*}
Then we combine pieces together to check that
\begin{align*}
    \ca T_2=\frac{2\pi}{k}\bl\frac{G_{2,k}^{(2\tau)}(x_2^*)}{G_{3,k}^{(2\tau)}(x_3^*)} \br^{\frac{1}{2\tau}}(1+O(k^{-\frac{1}{2\tau}}))=\frac{2\pi}{k}(1+O(k^{-\frac{1}{2\tau}})).
\end{align*}
And finally we collect the above result, \eqref{t0bound}, and \eqref{criticalseparate} to conclude that $ T_1=1+O(k^{-\frac{1}{2\tau-1}})$. Similarly we can derive that $ T_2=1+O(k^{-\frac{1}{2\tau-1}})$. Then we can get
\begin{align*}
    D_{kl}(\bb P(\bfa\sigma)\Vert\bb P(\bfa\sigma|S)|S)=O(k^{-\frac{1}{2\tau-1}})\quad\Rightarrow\quad \bb P(\wh S\neq S)\geq 1-O\bl\frac{m}{k^{\frac{1}{2\tau-1}}\log n}\br\vee 1.
\end{align*}
\subsection{Proof of Theorem \ref{asscreening}}
    First we notice that the set $S^\prime= S_1\cup S_2$ with $S_1\subset S$ and $S_2\notin S$ with $|S_2|=o(k)$ with probability $1-o(1)$. Then we  prove  under the condition of $|S_2|=o(k)$ since it is an asymptotically almost sure event. Consider all the spins outside $S$, at the high temperature we  have $\Vert m^{\frac{1}{2}} k^{-\frac{1}{2}}(\phi_i-\bb E[\phi_i])\Vert_{\psi_2}<\infty$ and $\bb E[\phi_i]=0$ by the independence. Then the following tail bound holds
    \begin{align*}
        \bb P(\phi_i-\bb E[\phi_i]\geq t)\leq\exp\bl-C\frac{m}{k}t^2\br.
    \end{align*}
    And by union bound we  have (noticing that here $\bb E[\phi_i]$ is the same for all $i\in S^c$)
    \begin{align*}
        \bb P\bl\sup_{i\in S^c}\phi_i\geq \bb E[\phi_i]+t \br\leq n\exp\bl-C\frac{mt^2}{k}\br=o(1),
    \end{align*}
    when picking $m\geq Ck\log n$ for large $C>0$.
    Then, for the critical temperature, it is analogously shown that 
    \begin{align*}
        \bb P\bl\sup_{i\in S^c}\phi_i\geq\bb E[\phi_i]+t\br\leq n\exp\bl -Ck^{-\frac{1}{2\tau-1}}{mt^2}{}\br=o(1),
    \end{align*}
    when picking $m\geq Ck^{\frac{1}{2\tau-1}}\log n$ for large $C$. And at low temperature, we  have $\bb E[\phi_1]=o(1)$ and
    \begin{align*}
        \bb P\bl\sup_{i\in S^c}\phi_i\geq\bb E[\phi_i]+t\br\leq n\exp( -Cmt^2)=o(1),
    \end{align*}
    when picking $m\geq C\log n$ for large $C$. And finally we notice that by the limit theorem in \ref{cltrfcr}, for all $i\in S$ we  have $\bb E[\phi_i]=C>0$ for all $i\in S$ at the high/low and critical temperatures. Furthermore we define the average magnetization by $M_{-i}=\sum_{j\neq i,j\in S}\sigma_i$, then it is immediate to see that for all $r\in\bb N$, by convexity,
    \begin{align*}
        \bb E[(\phi_i-\bb E[\phi_i])^{2r}]\leq 2^{2r}\bb E[\phi_i^{2r}]=\begin{cases}
            m^{-2r}M_{-i}^{2r}&\text{ at the high temperature}\\
            m^{-2r}k^{\frac{4\tau-4}{2\tau-1}r}M_{-i}^{2r}&\text{ at the critical temperature}\\
            m^{-2r}k^{-2r}M_{-i}^{2r}&\text{ at the low temperature}
        \end{cases}.
    \end{align*}
    Then, by the results of $M_{-i}$ given in the proof of corollary \ref{hightemperaturerec}, \ref{recoveryguaranteelowtemp}, and \ref{recoverycrtical}, at high temperature we  have $\Vert m^{\frac{1}{2}}k^{-\frac{1}{2}}(\phi_i-\bb E[\phi_i])\Vert_{\psi_2}<\infty$, at critical temperature we  have $\Vert m^{\frac{1}{2}}k^{-\frac{1}{4\tau-2}}(\phi_i-\bb E[\phi_i])\Vert_{\psi_{4\tau-2}}<\infty$, and at low temperature we  have $\Vert m^{\frac{1}{2}}k^{\frac{1}{2}}(|\phi_i|-\bb E[|\phi_i|])\Vert_{\psi_2}<\infty$.
    
    Then we  use the union bound to get at high temperature for $t>0$
    \begin{align*}
        \bb P\bl\inf_{i\in S}(\phi_i-\bb E[\phi_i])\leq -t\br\leq k\exp\lef(-C\frac{m}{k}t^2\rig)=o(1),
    \end{align*}
    when $m\geq C k\log k$ for some large $C>0$. And at critical temperature we  have
    \begin{align*}
        \bb P\bl\inf_{i\in S}(\phi_i-\bb E[\phi_i])\leq -t\br\leq k\exp\bl-C\frac{m}{k^{\frac{1}{2\tau-1}}}t^2\br=o(1).
    \end{align*}
    And at critical temperature we  have
    \begin{align*}
        \bb P\bl\inf_{i\in S}(\phi_i-\bb E[\phi_i])\leq-t\br\leq k\exp\bl-Cmkt^2\br=o(1).
    \end{align*}
    Therefore, collecting the above pieces it is not hard to check that picking $t$ within the region of $(0,\bb E[\phi_i])$ we  have at all temperature regimes,
    \begin{align*}
        \bb P(S^{\prime\prime}=S)=1-o(1).
    \end{align*}
\section{Proof of Results in Section \ref{sect4} and \ref{sect5}}
\subsection{Proof of Theorem \ref{lbnoncentered}}
    We follow similar path as section \ref{proofminimaxlb}. Notice that here we not be using the \emph{fake measure trick}. Recall from \eqref{vphi0g}, the decomposed Chi-square is computed as
    \begin{align*}
        \bb E\lef[\frac{\bb P_{S}(\bfa\sigma)\bb P_{S^\prime}(\bfa\sigma)}{\bb P_0(\bfa\sigma)}\rig]=\bb E\lef[\frac{\prod_{i=r+1}^k\cosh(h_i)\int\exp\lef(-kG_{0,k}(x,y,\bfa h)\rig)dxdy }{\int\exp\lef(-kG_{1,k}(x,\bfa h)\rig)dx\int\exp\lef(-kG_{2,k}(y,\bfa h)\rig)dy}\rig].
    \end{align*}
    Consider two sets $S$ $S^\prime$ with $S$ being the index set of clique and $S^\prime$ being the overlapping set such that $|S\cap S^\prime|=k-r$. We also define $c:=\frac{r}{k}$. Then, $x_{k},y_{k}$ satisfies the following optimality conditions.
    \begin{align*}
    (x_{k},y_{k})=\argmin_{x,y}G_{0,k}(x,y).
    \end{align*}
    Using the Fermat's first order condition,  we have
    \begin{align*}
        x_k&=\frac{\sqrt{\theta_1}}{k}\sum_{i=1}^r\tanh(\sqrt{\theta_1}x_k+h_i)+\frac{\sqrt{\theta_1}}{k}\sum_{i=r+1}^k\tanh(\sqrt{\theta_1}x_k+\sqrt{\theta_1}y_k+h_i),\\
        y_k&=\frac{\sqrt{\theta_1}}{k}\sum_{i=r+1}^k\tanh(\sqrt{\theta_1}x_k+\sqrt{\theta_1}y_k+h_i)+\sum_{i=k+1}^{k+r}\tanh(\sqrt{\theta_1}y_k+h_i).
    \end{align*}
    Then we consider the following $x^*,y^*$ that is the stationary point of $G_{0}(x,y)$,  we have
    \begin{align*}
        x^*&= c\sqrt{\theta_1}\bb E[\tanh(\sqrt{\theta_1}x^*+h)]+(1-c)\sqrt{\theta_1}\bb E[\tanh(\sqrt{\theta_1}x^*+\sqrt{\theta_1}y^*+h)],\\
        y^*&=c\sqrt{\theta_1}\bb E[\tanh(\sqrt{\theta_1}y^*+h)]+(1-c)\sqrt{\theta_1}\bb E[\tanh(\sqrt{\theta_1}x^*+\sqrt{\theta_1}y^*+h)].
    \end{align*}
    And analogously  we have $x_{1,k},x_{2,k}\to x_{1}^*$, almost surely with
$        x_1^*=\sqrt{\theta_1}\bb E[\tanh(\sqrt{\theta_1}x_1^*+h)].
$ Considering their difference, there exists $x_1^\prime\in(x^*\wedge x_1^*,x^*\vee x_1^*),y_1^\prime\in(y^*\wedge x_1^*,y^*\vee x_1^*)$ such that
\sm{\begin{align*}
    x^*-x_1^*&=\theta_1\bb E[\sech^2(\sqrt{\theta_1}x_1^*+h)]((x^*-x^*_1)+(1-c)(y^*-x^*_1)) + G_{1,k}^{(2\tau)}(x_1^\prime)(x^*-x_1^*)^{2\tau-1},\\
    y^*-x_1^*&=\theta_1\bb E[\sech^2(\sqrt{\theta_1}x_1^*+h)]((y^*-x_1^*)+(1-c)(x^*-x_1^*))+G^{(2\tau)}_{1,k}(y_1^\prime)(y^*-x_1^*)^{2\tau-1}.
\end{align*}}
Noticing that $\theta_1\bb E[\sech^2(\sqrt{\theta_1}x_1^*+h)]=1$,  we have when $\tau\neq 1$,
\begin{align*}
    |x^*-x_1^*|=O\lef((1-c)^{\frac{1}{2\tau-2}}\rig),\qquad |y^*-x_1^*|=O\lef((1-c)^{\frac{1}{2\tau-2}}\rig).
\end{align*}
And when $\tau=1$, $|x^*-x_1^*|\vee |y^*-x_1^*|=O(1-c)$.
Define $\bfa\delta=(x_k-x^*,y_k-x^*)^\top$ and use $\times$ as the notation for tensor product, then  we have
   \ttny{\begin{align*}
       \nabla G_{0,k}(x_k, y_k)&=0=\nabla G_{0,k}(x^*,y^*)+\nabla^2G_{0,k}(x^*, x^*)\times\bfa\delta +\sum_{i=3}^{2\tau}\nabla^{(i)}G_{0,k}(x^*, x^*)\times\bfa\delta^{i-1}+O(\Vert\bfa\delta\Vert_2^{2\tau}).\\
       G_{1,k}^\prime(x_{1,k})&=0 =G^\prime_{1,k}(x^*)+G^{(2)}_{1,k}(x^*)(x_{1,k}-x^*)+\sum_{i=3}^{2\tau}G^{(i)}_{1,k}(x^*)(x_{1,k}-x^*)^{i-1}+O((x_{1,k}-x^*)^{2\tau}).
   \end{align*}}
   Therefore, using similar arguments with the proof of theorem \ref{minimaxlbct},  we have
   \ttny{\begin{align*}
       -kG_{0,k}&(x_{k},y_k)+kG_{1,k}(x_{1,k})+kG_{2,k}(x_{2,k})+\sum_{i=r+1}^k\log\cosh(h_i)=-kG_{0,k}(x^*,y^*)+kG_{1,k}(x_{1}^*)+kG_{2,k}(x_1^*)\\
       &+\sum_{i=r+1}^k\log\cosh(h_i)-k\bl \nabla_x G_{0,k}(x^*,y^*)(x_k-x^*)+\nabla_yG_{0,k}(x^*,y^*)(y_k-y^*)-G^\prime_{1,k}(x_1^*)(x_k- x^*)\\
       &-G_{2,k}^\prime(x_2^*)(x_k-x^*)+\nabla^2_{xy}G_{0,k}(0,0)(x_k-x^*)(y_k-y^*)+O((x_k-x^*)^{2\tau+1} \vee (y_k-y^*)^{2\tau+1})\\
       &+O(1-c)(x_k-x^*)(y_k-y^*)^2\br.
   \end{align*}}
   Then we notice that by the Lipschitzness of $\log\cosh$ and the fact that the random terms can be bounded by H\"older's inequality to be upperbounded by $\exp(C(1-c)k)$ for some $C>0$. Then  we have
   \ttny{\begin{align*}
       \bb E\bigg[\exp\bl-kG_{0,k}(x_{k},y_k)+kG_{1,k}(x_{1,k})+kG_{2,k}(x_{2,k})+\sum_{i=r+1}^k\log\cosh(h_i)\br\bigg]\leq\exp\lef(Ck(1-c)^{\frac{1}{2\tau-2}}\rig).
   \end{align*}}
   Similar to the analysis for the rest of the terms, similar to the derivation of \ref{separationofcrossterms} and \ref{controlofsmallterms}, we finally arrive at for some constant $C>0$,
   \begin{align*}
       \bb E\lef[\frac{\bb P_{S}\bb P_{S^\prime}}{\bb P_0}\rig]\leq\begin{cases}
           \exp(C(1-c)^{\frac{1}{2\tau-2}}k) &\text{ when }\tau\neq 1,\\
           \exp(C(1-c)k) &\text{ when }\tau=1.
       \end{cases}
   \end{align*}
   First we consider when $\tau\neq 1$. We go back to  \eqref{divdep} and \eqref{therest} to get for $p=\epsilon k$, $C_1>1$.
   \begin{align*}
    \sum_{v=p}^k\bb P(V=v) E_k^m(v)&\leq\sum_{v=p+1}^k\frac{1}{v!}\lef(\frac{k^2}{n}\rig)^vC_1^{km}\leq\sum_{v=p}^kv\lef(\frac{ek^2}{nv}\rig)^vC_1^{km}\\
    &\leq \lef(C_2\frac{k\log k}{n}\rig)^{\epsilon k}C_1^{km}=o(1).
\end{align*}
Then we analyze the rest through integral approximation for some $C>1$ and $\mca f(x)$ defined in \eqref{riemann}:
\tny{\begin{align*}
   \sum_{v=0}^{p}\bb P(V=v)E_k^m(v)\leq 1+\sum_{v=1}^k\bb P(V=v)E_k^m(v)\leq 1+\int_{(\frac{1}{k},\epsilon]}\frac{\sqrt k}{(1-x)\sqrt{2\pi x}}\exp\lef(Ck \mca f(x)\rig)dx.
\end{align*}}
   Therefore we use Laplace approximation again, recalling $\gamma:=\frac{k}{n}$  we have
   \begin{align*}
       \mca f^\prime(x)=(4-2x)\gamma-\log\frac{x}{\gamma}+2\log(1-x)+C_1mx^{\frac{1-2\tau}{2\tau-2}}.
   \end{align*}
   Therefore it is not hard to see that the maximum is taking at $x=\epsilon$, which implies that when $k=o(\sqrt n)$ and $m=o(\log n)$  we have $\sum_{v=1}^k\bb P(V=v) E_k^m(v)=o(1)$. 

   Then we consider when $\tau=1$, this directly corresponds to the low temperature case of the symmetric random field case. And we finish the proof.

\subsection{Proof of Theorem \ref{thm4.2}}
The proof follows by proving local part and the global part separately.
\begin{center}
    \textbf{1. Local Part}
\end{center}
   To upper bound the Type I error, we first notice that under the null, $\bb E[\frac{1}{n}\sum_{i=1}^n\sigma_i]=\bb E[\tanh(h)]$ and by Hoeffding's inequality, there exists $C>0$ such that
    \begin{align*}
        \bb P_0\lef(|\xi-\bb E[\tanh(h)]|\geq t\rig)\leq\exp\lef(-Cnmt^2\rig).
    \end{align*}
    And going back to the analysis of $\phi_S$,  we have by union bound, using also the fact that $\lef\Vert\frac{1}{\sqrt{m k}}\sum_{i=1}^m\mbbm 1_{S}\bfa\sigma^{(i)}\rig\Vert_{\psi_2}<\infty$ there exists $C>0$ such that
    \begin{align*}
        \bb P_{0}(\phi_7^{\max}-\bb E[\tanh(h)]\geq t)&\leq\binom{n}{k}\bb P_{0}(\phi_S-\bb E[\tanh(h)]\geq t)\leq\binom{n}{k}\exp\lef(-Cmkt^2\rig),\\
        \bb P_{0}(\bb E[\tanh(h)]-\phi_7^{\min}\geq t)&\leq\binom{n}{k}\bb P_{0}(\bb E[\tanh(h)]-\phi_S\geq t)\leq\binom{n}{k}\exp\lef(-Cmkt^2\rig).
    \end{align*}
    Therefore combing the above two inequalities, there exists $C_1,C_2,C_3,C_4>0$ such that
    \begin{align*}
        \bb P_0(\phi_7^\max-\xi>t)&\leq \bb P_0(\phi_7^{\max}-\bb E[\tanh(h)]+\bb E[\tanh(h)]-\xi>t)\\
        &\leq\bb P_0\lef(\phi_7^{\max}-\bb E[\tanh(h)]>\frac{t}{2}\rig)+\bb P_0\lef(\bb E[\tanh(h)]-\xi>\frac{t}{2}\rig)\\
        &\leq\binom{n}{k}\exp(-C_1mkt^2)+\exp(-C_2nmt^2)\leq \binom{n}{k}\exp\lef(-C_3mkt^2\rig),\\
        \bb P_0(\xi-\phi_7^{\min}>t)&\leq\bb P_0(\xi-\bb E[\tanh(h)]+\bb E[\tanh(h)]-\phi_7^{\min}>t)\leq\binom{n}{k}\exp(-C_4mkt^2).
    \end{align*}
    And we   conclude by union bound there exists $C_5>0$ such that
    \begin{align}\label{type1localh}
        \bb P_0(\xi-\phi_7^{\min}>t\text{ or }\phi_7^\max-\xi>t)\leq\binom{n}{k}\exp(-C_5mkt^2).
    \end{align}
    Then we analyze the Type II error, noticing that when the hidden clique has a index set of $S_0$   by corollary \ref{corocltfcr}, $\Vert k^{-\frac{4\tau-3}{4\tau-2}}\sum_{i\in S_0}(\sigma_i-\bb E[\sigma_i])\Vert_{\psi_{4\tau-2}}<\infty$ and $\bb E[\phi_{S_0}]=\sqrt{\theta_1}x^*+o(1)$. We consider when $\bb E[\phi_{S_0}]>\bb E[\tanh(h)]$ as the other side is   achieveable analogously. We notice that under the alternative hypothesis, by the concentration inequality given by the sum of i.i.d. sub-Weibull r.v.s. in lemma \ref{weibullnormsubadditivity},   for $t>0$ there exists $C>0$ such that
    \begin{align*}
        \bb P(\phi_{S_0}-\bb E[\phi_{S_0}]>t)&=\bb P\bl m^{-1}k^{-\frac{4\tau-3}{4\tau-2}}\sum_{j=1}^m\sum_{i\in S_0}(\sigma_i^{(j)}-\bb E[\sigma_i^{(j)}])>k^{\frac{1}{4\tau-2}}t\br\\
        &\leq\exp\lef(-Cmkt^{4\tau-2}\wedge m k^{\frac{1}{2\tau-1}}t^2\rig),
    \end{align*}
    which also implies that
    \begin{align*}
        \bigg \Vert m^{-\frac{1}{2}}k^{-\frac{4\tau-3}{4\tau-2}}\sum_{j=1}^m\sum_{i\in S}(\sigma^{(j)}_i-\bb E[\sigma_i])\bigg\Vert_{\psi_2}<\infty.
    \end{align*}
    Then, using the sub-additivity of sub-Weibull norm  we have
    \begin{align*}
        \bigg\Vert(m)^{-\frac{1}{2}}(n^{-\frac{1}{2}}\wedge k^{-\frac{4\tau-3}{4\tau-2}})&\sum_{j=1}^m\sum_{i=1}^n(\sigma_i^{(j)}-\bb E[\sigma_i^{(j)}])\bigg\Vert_{\psi_2}\leq\bigg\Vert(nm)^{-\frac{1}{2}}\sum_{j=1}^m\sum_{i\in S^c}(\sigma_i^{(j)}-\bb E[\sigma_i])\bigg\Vert_{\psi_2}\\
        &+\bigg\Vert m^{-\frac{1}{2}}k^{-\frac{4\tau-3}{4\tau-2}}\sum_{j=1}^m\sum_{i\in S}(\sigma_i^{(j)}-\bb E[\sigma_i])\bigg\Vert_{\psi_2}<\infty.
    \end{align*}
    And regarding the expectation, under the alternative, we have
$        \bb E[\xi]=\bb E[\tanh(h)]\frac{n-k}{n}+\frac{k\bb E[\tanh(x^*+h)]}{n}+o(1).
$    Therefore, it is not hard to see that for $t>0$:
    \begin{align*}
        \bb P_{S}(|\xi-\bb E[\xi]|\geq t)\leq\exp\lef(-Cnmt^2\wedge n^2mk^{-\frac{4\tau-3}{2\tau-1}}t^2\rig).
    \end{align*}
    Then we consider $\phi_7^{\max}$ and $\phi_7^{\min}$. We denote $\bb P_{S_0}$ as the measure under the alternative hypothesis with the hidden clique indexed by $S_0$, then   for $t>0$, introducing $\Delta:=\bb E[\phi_{S_0}]-\bb E[\xi]-t$:
    \tny{\begin{align}\label{type2localh}
        \bb P_{S}&(\phi_7^{\max}-\xi< -t)\leq \bb P_{S_0}(\phi_{S_0}-\xi<-t)\leq \bb P_{S_0}(\phi_{S_0}-\bb E[\phi_{S_0}]+\bb E[\xi]-\xi <t-\bb E[\phi_{S_0}]+\bb E[\xi])\nnb\\
        &\leq\bb P_{S_0}\lef(\phi_{S_0}-\bb E[\phi_{S_0}]\leq\frac{1}{2}(-t-\bb E[\phi_{S_0}]+\bb E[\xi])\rig)+\bb P_{S_0}\lef(\bb E[\xi]-\xi\leq\frac{1}{2}(-t-\bb E[\phi_{S_0}]+\bb E[\xi])\rig)\nnb\\
        &\leq\exp\lef(-Cmk\Delta^{4\tau-2}\wedge m k^{\frac{1}{2\tau-1}}\Delta^2\rig) + \exp\lef(-Cnm\Delta^2\wedge n^2mk^{-\frac{4\tau-3}{2\tau-1}}\Delta^2\rig),\nnb\\
        \bb P_{S}&(\xi-\phi_7^{\min}< -t)\leq\exp\lef(-Cmk\Delta^{4\tau-2}\wedge m k^{\frac{1}{2\tau-1}}\Delta^2\rig) + \exp\lef(-Cnm\Delta^2\wedge n^2mk^{-\frac{4\tau-3}{2\tau-1}}\Delta^2\rig).
    \end{align}}
    We notice that $\Delta >0$ for sufficient large $k$ and it is   checked that for $m\gtrsim \log n$ one will have the quantity in \eqref{type1localh} and \eqref{type2localh} to be arbitrarily small asymptotically.
    \begin{center}
        \textbf{2. Global Part}
    \end{center}
    To study the independent copy, we first notice that by independence, the following holds
    \begin{align*}
        \bb E[\sigma_i^{(1)}-\sigma_i^{(2)}]=0,\qquad \bb E\bigg[\bl \sum_{i=1}^n(\sigma_i^{(1)}-\sigma_i^{(2)})\br^2\bigg]=2\bb V\bigg[\sum_{i=1}^n\sigma_i\bigg].
     \end{align*}
     By Jensen's inequality, for all $t\in\bb R$,  we have
     \begin{align*}
     \bb E\bigg[\exp\bl \frac{t}{k}\sum_{i=1}^n(\sigma_i^{(1)}-\sigma_i^{(2)})\br\bigg]\leq\bb E\bigg[\exp\bl\frac{t}{k}\sum_{i=1}^n(\sigma_i-\bb E[\sigma_i])\br\bigg].
    \end{align*}
    Then, we start analyzing the Type I error when $k=O\lef(n^{\frac{2\tau-1}{4\tau-3}}\rig)$, we first notice that
$        \frac{1}{n}\bb V\lef[\sum_{i=1}^n\sigma_i\rig]\asymp 1.
$ and $\lef\Vert\frac{1}{\sqrt n}\sum_{i=1}^n(\sigma_i^{(1)}-\sigma_i^{(2)})\rig\Vert_{\psi_2}<\infty$. Therefore
$
    \lef\Vert\frac{1}{n}\lef(\sum_{i=1}^n(\sigma_i^{(1)}-\sigma_i^{(2)})\rig)^2\rig\Vert_{\psi_1}<\infty
$, using Bernstein's inequality, there exists $C>0$ such that
\tny{\begin{align*}
    \bb P_0\bl& k^{-\frac{4\tau-3}{2\tau-1}}m^{-1}\bigg|\sum_{j=1}^m\frac{1}{m}\bl\bl\sum_{i=1}^n(\sigma_i^{(2j-1)}-\sigma_i^{(2j)})\br^2-\bb E\bigg[\bl\sum_{i=1}^n(\sigma_i^{(2j-1)}-\sigma_i^{(2j)})\br^2\bigg]\br\bigg|\geq t\br\\
    &=\bb P_0\bl \bigg|\sum_{j=1}^m\frac{1}{mn}\bl\bl\sum_{i=1}^n(\sigma_i^{(2j-1)}-\sigma_i^{(2j)})\br^2-\frac{2}{n}\bb V\bigg[\sum_{i=1}^n\sigma_i\bigg]\br\bigg|\geq n^{-1}k^{\frac{4\tau-3}{2\tau-1}} t\br\\
    &\leq 2\exp(-Cmn^{-1}k^{\frac{4\tau-3}{2\tau-1}} t\wedge mn^{-2}k^{\frac{2(4\tau-3)}{2\tau-1}} t^2).
\end{align*}}
And therefore noticing that under the null $\bb E[\phi_8]=o(1)$, there exists $C>0$ such that for $t>0$,
\begin{align}\label{TypeI10}
    \bb P(\phi_8>t)\leq \exp(-Cmn^{-1}k^{\frac{4\tau-3}{2\tau-1}}t\wedge mn^{-2}k^{\frac{2(4\tau-3)}{2\tau-1}}t^2).
\end{align}
    And under the alternative when $k=O\lef(n^{\frac{2\tau-1}{4\tau-3}}\rig)$, noticing that by the exponential inequality in corollary \ref{corocltfcr}, $\lef\Vert k^{-\frac{4\tau-3}{4\tau-2}}\sum_{i\in S}(\sigma^{(1)}_i-\sigma^{(2)}_i)\rig\Vert_{\psi_{4\tau-2}}<\infty$. By the sub-additivity of sub-Weibull norm  we have
\ttny{\begin{align*}
    \bigg\Vert n^{-1/2}\wedge k^{-\frac{4\tau-3}{4\tau-2}}\sum_{i=1}^n(\sigma^{(1)}_i-\sigma^{(2)}_i)\bigg\Vert_{\psi_{2}}\leq \bigg\Vert k^{-\frac{4\tau-3}{4\tau-2}}\sum_{i\in S}(\sigma^{(1)}_i-\sigma^{(2)}_i)\bigg\Vert_{\psi_{2}}+\bigg\Vert n^{-1/2}\sum_{i\in S^c}(\sigma^{(1)}_i-\sigma^{(2)}_i)\bigg\Vert_{\psi_{2}}<\infty.
\end{align*}}
And we also notice that by independence between $\bfa\sigma_S$ and $\bfa\sigma_{S^c}$  we have when $k=\omega\lef(n^{\frac{2\tau-1}{4\tau-3}}\rig)$, using corollary \ref{corocltfcr}, under the alternative  we have
\ttny{\begin{align*}
    \bb E[\phi_8]=2k^{-\frac{4\tau-3}{2\tau-1}}\bb V\bigg[\sum_{i=1}^n\sigma_i\bigg]=2k^{-\frac{4\tau-3}{2\tau-1}}\bb V\bigg[\sum_{i\in S}\sigma_i\bigg]+2k^{-\frac{4\tau-3}{2\tau-1}}\bb V\bigg[\sum_{i\in S^c}\sigma_i\bigg]=2k^{-\frac{4\tau-3}{2\tau-1}}\bb V\bigg[\sum_{i\in S}\sigma_i\bigg]+o(1)\asymp 1.
\end{align*}}
Therefore, we again use the fact that under the alternative $
    \lef\Vert k^{-\frac{4\tau-3}{2\tau-1}}\lef(\sum_{i=1}^n(\sigma_i^{(1)}-\sigma_i^{(2)})\rig)^2\rig\Vert_{\psi_1}<\infty
$, there exists $C_1,C_2>0$ such that for all $t>0$,
\begin{align}\label{typeII10}
    \bb P\lef(\phi_8\leq\bb E[\phi_8] - t\rig)\leq\exp(-C_1mt^2\wedge C_2mt).
\end{align}
Therefore, collecting \ref{TypeI10} and \ref{typeII10} if we pick $m=C_0$ for some $C_0>0$, we can control the Type I + Type II error to be arbitarily small.


\subsection{Proof of Theorem \ref{oracletestguarantee}}
    We first control the Type I error, noticing that $\bb E[\phi_9]=0$ and by the fact that when $\sigma_i$ are i.i.d. r.v.s.  we have $\Vert\frac{1}{\sqrt n}(\sum_{i=1}^n\sigma_i-\bb E[\sigma_i])\Vert_{\psi_2}<\infty$ and $\Vert\frac{1}{n}(\sum_{i=1}^n\sigma_i-\bb E[\sigma_i])^2\Vert_{\psi_1}<\infty$. By Bernstein's inequality, there exists $C>0$ such that under the null  we have
    \begin{align*}
        \bb P(|\phi_9|\geq t)&=\bb P\bl \bigg |m^{-1}k^{-2}\sum_{j=1}^m\bl\bigg (\sum_{i=1}^n\sigma_i^{(j)} -n\bb E[\sigma_i] \bigg )^2-n(1-\bb E[\tanh(h)]^2)\br\bigg|\geq t\br\\
        &\leq \exp(-Cmk^4n^{-2}t^2\wedge mk^{2}n^{-1}t).
    \end{align*}
    And then we study the Type II error. It is first checked that under the alternative  we have 
    \tny{\begin{align*}
        &\bb E[\phi_9]=\frac{1}{k^2}\bl\bb E\bigg[\bl\sum_{i=1}^n\sigma_i-n\bb E[\tanh(h)]\br^2\bigg]-n(1-\bb E[\tanh(h)]^2)\br\\
        &=\frac{1}{k^2}\bb E\bigg[\bl\sum_{i=1}^n\sigma_i-\bb E\bigg[\sum_{i=1}^n\sigma_i\bigg]\br^2\bigg]+\frac{1}{k^2}\bl\bb E\bigg[\sum_{i=1}^n\sigma_i\bigg]-n\bb E[\tanh(h)]\br^2-\frac{n}{k^2}(1-\bb E[\tanh(h)]^2)\\
        &=\frac{1}{k^2}\bb E\bigg[\bl\sum_{i\in S^c}\sigma_i-\bb E\bigg[\sum_{i\in S^c}\sigma_i\bigg]\br^2\bigg]-\frac{1}{k}(1-\bb E[\tanh(h)]^2)+\frac{1}{k^2}\bl\bb E\bigg[\sum_{i=1}^n\sigma_i \bigg] -n\bb E[\tanh(h)]\br^2\\
        &=(\bb E[\tanh(x^*+h)]-\bb E[\tanh(h)])^2+O\lef(\frac{1}{k}\rig).
    \end{align*}}
    Then we notice that by the sub-additivity of Orlicz norm  we have when $k\gtrsim \sqrt n$,
    \begin{align*}
        \bigg\Vert k^{-1}\bl \sum_{i=1}^n\sigma_i-n\bb E[\tanh(h)]\br\bigg\Vert_{\psi_2}&\leq\bigg\Vert k^{-1}\bl\sum_{i\in S}\sigma_i-k\bb E[\tanh(h)]\br\bigg\Vert_{\psi_2}\\
        &+\bigg\Vert k^{-1}\bl\sum_{i\in S^c}\sigma_i-\bb E\bigg[\sum_{i\in S^c}\sigma_i\bigg]\br\bigg\Vert_{\psi_2}<\infty.
    \end{align*}
    Hence $\bigg\Vert k^{-2}\bl \sum_{i=1}^n\sigma_i-n\bb E[\tanh(h)]\br^2\bigg\Vert_{\psi_1}<\infty$ and there exists $C>0$ such that by Bernstein's inequality,
    \begin{align*}
        \bb P\bl \frac{1}{m}\sum_{j=1}^m\bl k^{-2}\bl \sum_{i=1}^n\sigma_i-n\bb E[\tanh(h)]\br^2&-\bb E\bigg[ k^{-2}\bl \sum_{i=1}^n\sigma_i-n\bb E[\tanh(h)]\br^2\bigg]\br
        \leq -t\br\\
        &\leq\exp(-Cmt^2\wedge Cmt).
    \end{align*}
    Therefore, for $m\asymp 1$ and picking some small $\tau_\delta>0$ algorithm \ref{alg:nine} can perform asymptotic powerful test.

\subsection{Proof of Theorem \ref{cltrfcr} and Corollary \ref{corocltfcr}}
Here we present the proof of theorem \ref{cltrfcr} and corollary \ref{corocltfcr}.
Our proof goes by three parts according to the high, low, and critical temperature.
\begin{center}
    \textbf{1. The High Temperature Case}
\end{center}
     we have by the Laplace method for constant $a$ and note that the mgf can be written as:
    \tny{\begin{align}\label{integraltransform}
        \bb E&\lef[\exp\lef(t\frac{\sum_{i=1}^n(\sigma_i-a)}{\sqrt n}\rig)\rig]=\bb E\lef[\frac{\sum_{\bfa\sigma}\exp\lef(\frac{\theta_1}{2}\lef(\sum_{i=1}^{n}\sigma_i\rig)^2+\sum_{i=1}^n\lef(\frac{t}{\sqrt n}+h_i\rig)\sigma_i\rig)\exp\lef(-\sqrt n ta\rig)}{\sum_{\bfa\sigma}\exp\lef(\frac{\theta_1}{2}\lef(\sum_{i=1}^n\sigma_i\rig)^2+\sum_{i=1}^nh_i\sigma_i\rig)}\rig]\nnb\\
        &=\bb E\lef[\frac{\int_{\bb R}\exp\lef(-n\mca H_{0,n}(x)-\sqrt n ta\rig)dx}{\int_{\bb R}\exp\lef(-n\mca H_{1,n}(x)\rig)dx}\rig]\nnb\\
        &=\bb E\bigg[\sqrt{\frac{\mca H_{0,n}^{(2)}(x_0)}{\mca H_{1,n}^{(2)}(x_1)}}\exp\lef(-n\mca H_{0,n}(x_0)+n\mca H_{1,n}(x_1)-\sqrt nta\rig)\bl 1+\sum_{i=1}^\infty\frac{b_i(\bfa h)}{n^i}\br\bigg].
    \end{align}}
    with $b_i$ bounded.
    Then it is easily checked that uniformly the following holds:
    \begin{align*}
    \mca H_{0,n}(x,\bfa h):&=\frac{1}{2}x^2-\frac{1}{n}\sum_{i=1}^n\log\cosh\lef(\sqrt{\theta_1} x+h_i+\frac{t}{\sqrt n}\rig),\\
    \mca H_{1,n}(x,\bfa h):&=\frac{1}{2}x^2-\frac{1}{n}\sum_{i=1}^n\log\cosh\lef(\sqrt{\theta_1} x+h_i\rig).
    \end{align*}
    And we define the population version as
    \begin{align*}
        \mca H_{0}(x,\bfa h):&=\frac{1}{2}x^2-\bb E\lef[\log\cosh\lef(\sqrt{\theta_1} x+h_i+\frac{t}{\sqrt n}\rig)\rig],\\
    \mca H_{1}(x,\bfa h):&=\frac{1}{2}x^2-\bb E\lef[\log\cosh\lef(\sqrt{\theta_1} x+h_i\rig)\rig].
    \end{align*}
    Then we propose the following lemma guarantees the uniform convergence whose proof is delayed to the appendix.
   \begin{lemma}[Regularity Conditions] When $h$ is in $L_1$. Almost surely in $\mu(\bfa h)$ and uniformly on $(x,y)$  we have
    \begin{align*}
        \mca H_{0,n}^{(j_1,j_2)}(x,y,\bfa h):=\frac{\pta^{j_1+j_2}G_{0,k}(x,y,\bfa h)}{\pta x^{j_1}\pta y^{j_2}}\to \mca H_0^{(j_1,j_2)}(x,y)
    .\end{align*}
    with $\mca H_{0,n}^{(0,0)}:=\mca H_{0,n}$. Similar argument holds for $\mca H_{1,n}$ and $\mca H_1$. And condition \eqref{laplace2}, \eqref{laplace3} in lemma \ref{laplace} holds for $\mca H_{0,n}, \mca H_{1,n}$, implies the validity of Laplace integral approximation.
\end{lemma}
    Here we denote $x_1:=\argmin_{x}\mca H_{1,n}(x,\bfa h)$, $x_0:=\argmin_{x}\mca H_{0,n}(x,\bfa h)$. Further denote that
    \begin{align*}
        x_0^*=\sqrt{\theta_1}\bb E\lef[\tanh\lef(\sqrt{\theta_1}x_0^*+h_i+\frac{t}{\sqrt n}\rig)\rig],\quad
        x_1^*=\sqrt{\theta_1}\bb E[\tanh(\sqrt{\theta_1}x_1^*+h_i)].
    \end{align*}
    Therefore  we have by Fermat's condition:
    \begin{align*}
        \mca H_{1,n}^\prime(x_1,\bfa h)= 0 =\mca H_{1,n}^\prime(x_1^*,\bfa h)+\mca H_{1,n}^{(2)}(x_1^*,\bfa h)\lef(x_1-x_1^*\rig)+O\lef((x_1-x_1^*)^2\rig),\\
        \mca H_{0,n}^\prime(x_0,\bfa h)= 0 =\mca H_{0,n}^\prime(x_0^*,\bfa h)+\mca H_{0,n}^{(2)}(x_0^*,\bfa h)\lef(x_0-x_0^*\rig)+O\lef(\frac{1}{\sqrt n}\rig).
    \end{align*}
    And further  we have by the linearization of the $Z$ estimator given by lemma \ref{linearZ},
    \tny{\begin{align*}
        \sqrt{n}\lef(x_0-x_0^*\rig)&=\frac{\sqrt{\theta_1}\sum_{i=1}^n\lef(\tanh\lef(\sqrt{\theta_1}x_0^*+h_i+\frac{t}{\sqrt n}\rig)-\bb E\lef[\tanh\lef(\sqrt{\theta_1} x_0^*+h_i+\frac{t}{\sqrt n}\rig)\rig]\rig)}{\sqrt n\lef(1-\theta_1\bb E[\sech^2(\sqrt{\theta_1}x_0^*+h+\frac{t}{\sqrt n})]\rig)}+o_{\psi_2}(1),\\
         \sqrt n \lef(x_1-x_1^*\rig)&=\frac{\sqrt{\theta_1}\sum_{i=1}^n\tanh\lef(\sqrt{\theta_1}x_1^*+h_i\rig)-\bb E\lef[\tanh\lef(\sqrt{\theta_1} x_1^*+h_i\rig)\rig]}{\sqrt n\lef(1-\theta_1\bb E[\sech^2(h+\sqrt{\theta_1}x_1^*)]\rig)}+o_{\psi_2}(1).
    \end{align*}}
    And we also have
    \begin{align*}
        x_1^*-x_0^*&=\sqrt{\theta_1}\bb E[\sech^2(h+\sqrt{\theta_1}x_1^*)]\lef(\sqrt{\theta_1}(x_1^*-x_0^*)-\frac{t}{\sqrt n}\rig)+o_{\psi_2}\lef(\frac{1}{\sqrt n}\rig).
    \end{align*}
    Hence, by the boundedness of $\sech $ and the law of large numbers  we have 
    \tny{\begin{align*}
        \Delta:&=\sqrt{\theta_1}(x_1^*-x_0^*)-\frac{t}{\sqrt n}=\frac{-t}{\sqrt n(1-\theta_1\bb E[\sech^2(h+\sqrt{\theta_1}x_1^*)])}.\\
        \sqrt{\theta_1}(x_1&-x_0)=\sqrt{\theta_1}(x_1-x_1^*+x_0^*-x_0)+\Delta+\frac{t}{\sqrt n}\\
        &=\frac{\theta_1\sum_{i=1}^n(\sech^2(\sqrt{\theta_1}x^*_1+h_i)-\bb E[\sech^2(\sqrt{\theta_1}x^*_1+h_i)])}{n\lef(1-\theta_1\bb E[\sech^2(\sqrt{\theta_1}x_1^*+h)]\rig)}\Delta+\Delta+\frac{t}{\sqrt n}\\
        &-\frac{\theta_1^{2}\bb E[\sech^2(\sqrt{\theta_1}x^*_1+h_i)\tanh(\sqrt{\theta_1}x^*_1+h_i)]}{n\lef(1-\theta_1\bb E[\sech^2(\sqrt{\theta_1}x_1^*+h)\rig)^2}\sum_{i=1}^n\tanh\lef(\sqrt{\theta_1}x_1^*+h_i\rig)\Delta+o_{\psi_2}\lef(\frac{1}{\sqrt n}\rig)\\
        &=-\alpha\frac{t}{\sqrt n}\lef(1+o_{\psi_2}(1)\rig),
    \end{align*}}
    where we let $\alpha:=\frac{\theta_1\bb E[\sech^2(\sqrt{\theta_1}x_1^*+h)]}{1-\theta_1\bb E[\sech^2(\sqrt{\theta_1}x_1^*+h)]}$. 
    Define $F(x,y):=-\frac{1}{2}x^2+\frac{1}{n}\sum_{i=1}^n\log\cosh\lef(\sqrt{\theta_1}x+h_i+y\rig)$, and  we have 
    \ttny{\begin{align*}
        n\mca H_{1,n}(x_1,\bfa h) - n\mca H_{0,n}(x_0,\bfa h)&=n\bl\frac{\pta F(x_1,0)}{\pta x_1}(x_0-x_1)+\frac{\pta F(x_1,0)}{\pta y}\frac{t}{\sqrt n}+\frac{\pta^2 F(x_1,0)}{\pta x_1^2}\frac{1}{2}(x_0-x_1)^2\\
        &+\frac{\pta^2 F(x_1,0)}{\pta x_1\pta y}\frac{t}{\sqrt n}(x_0-x_1)+\frac{\pta^2 F(x_1,0)}{\pta y^2}\frac{t^2}{2n}\br+o_{\psi_2}(1)\\
        &=\frac{t}{\sqrt n}\sum_{i=1}^n\lef(\tanh\lef(\sqrt{\theta_1} x_1 +h_i\rig)-\bb E\lef[\tanh\lef(\sqrt{\theta_1} x_1 +h_i\rig)\rig]\rig)\\
        &+\frac{t^2}{2n}\sum_{i=1}^n\sech^2(\sqrt{\theta_1} x^*_1+h_i)(1 +2\alpha+\alpha^2)-\frac{\alpha^2t^2}{2\theta_1}+o_{\psi_2}\lef(1\rig)\\
        &=\ub{\frac{t}{\sqrt n}\sum_{i=1}^n\lef(\tanh\lef(\sqrt{\theta_1} x_1 +h_i\rig)-\bb E\lef[\tanh\lef(\sqrt{\theta_1} x_1 +h_i\rig)\rig]\rig)}_{T_1}\\
        &+\ub{\frac{t^2}{2n}\sum_{i=1}^n\sech^2(\sqrt{\theta_1} x_1^*+h_i)(1+\alpha)^2}_{T_2}-\frac{\alpha^2t^2}{2\theta_1}\\
        &+\ub{\frac{t\sqrt{\theta_1}}{n}\lef(\sum_{i=1}^n\sech^2\lef(\sqrt{\theta_1}x_1^*+h_i\rig) \rig)\sqrt{n}(x_1-x_1^*)}_{T_3}+o_{\psi_2}\lef(1\rig).
    \end{align*}}
    Moreover, it is checked that 
$        \frac{\mca H_{0,n}^{(2)}(x_0)}{\mca H_{1,n}^{(2)}(x_1)}=1+o_{\psi_2}\lef(\frac{1}{\sqrt n}\rig).
$    And we finally see that picking $a=\bb E[\tanh\lef(\sqrt{\theta_1}x_1^*+h\rig)]$, for all $t\in\bb R$ not dependent on $n$ we can check the boundedness of mgf. By theorem 2 in \citep{chareka2008converse} we get
    \begin{align*}
       \bb E\lef[\exp\lef(t\frac{\sum_{i=1}^n(\sigma_i-a)}{\sqrt n}\rig)\rig] \to \exp\lef(\frac{\ca Vt^2}{2}\rig),
    \end{align*}
    with $\ca V:=\frac{1-\theta_1(\bb E[\sech^2(\sqrt{\theta_1} x_1^*+h)])^2-\bb E[\tanh(\sqrt{\theta_1}x_1^*+h)]^2}{(1-\theta_1\bb E[\sech^2(\sqrt{\theta_1}x_1^*+h)])^2}$. To prove the sub-Gaussian norm being bounded involve two steps. First we show that by H\"older's inequality, for all $t=o(n^{1/2})$, there exists $\delta_1,\delta_2,\delta_3,\delta_4>0$ such that $\frac{1}{1+\delta_1}+\frac{1}{1+\delta_2}+\frac{1}{1+\delta_3}+\frac{1}{1+\delta_4}<1$ and
    \tny{\begin{align*}
        \bb E[\exp(n\mca H_{1,n}(x_1,\bfa h)-n\mca H_{0,n}(x_0,\bfa h))](1+o(1))&\leq\bb E[\exp((1+\delta_1)A_1)]^{\frac{1}{1+\delta_1}}\bb E[\exp((1+\delta_2)A_2)]^{\frac{1}{1+\delta_2}}\\
        &\cdot\bb E[\exp((1+\delta_3)A_3)]^{\frac{1}{1+\delta_3}}\bb E[(1+\delta_4)o_{\psi_2}(1)]^{\frac{1}{1+\delta_4}}\\
        &=(1+o(1))\exp(C(1+o(1))t^2).
    \end{align*}}
    Let $X:=n^{-1/2}\sum_{i=1}^n\sigma_i$.
    By standard Chernoff bound, we obtain that for all $t=o(\sqrt n)$,
    \begin{align*}
        \bb P(X\geq t)\leq C\exp(-Ct^2).
    \end{align*}
    Then we use moment method to extend the above results to $t\in\bb R$ case, this is done by the proposition 2.5.2 in  \citep{vershynin2018high}, which requires us to obtain uniform control on the moment of order $p\in\bb N$. It is checked that by $|X|\leq\sqrt n$, using the property of Gamma functions, \emph{ for all $p\in\bb N$ }, there exists universal $C_0$ such that
    \begin{align*}
        \bb E|X|^p&=\int_{\bb R^+}\bb P(|X|^p\geq u)du=\int_{\bb R^+}\bb P(|X|\geq t)pt^{p-1}dt\\
        &\leq\int_{0}^{\sqrt n/\log n}C\exp(-Ct^2)pt^{p-1}dt+\bb P(|X|\geq \sqrt n/\log n)(\sqrt n)^p\\
        &\leq C^p\Gamma(p/2)+n^{\frac{p}{2}}\exp(-Cn/\log^2 n)\leq C_0(C_0p)^{p/2}.
    \end{align*}
    And we complete the proof by the equivalent definition of sub-Gaussian random variables given by proposition 2.5.2. in \citep{vershynin2018high}, which implies that $\Vert X\Vert_{\psi_2}<\infty$.
    \begin{center}
        \textbf{2. The Low Temperature Case}
    \end{center}
    And we come to prove when the solution to $\sqrt{\theta_1}\bb E[\tanh(\sqrt{\theta_1}x+h)]=x$ is more than $1$. In particular, when $h$ is symmetric, asymptotically with $n$,  we have two global maximum with similar value. Here we we extend the transfer principle initially proposed by \citep{ellis1980limit} to the random measure case. Without loss of generality we assume that $h$ has a distribution symmetric w.r.t. $0$.
    
    The proof then goes as follows
    \tny{\begin{align}\label{transferprincipa}
        \bb E&\lef[\exp\lef(t\frac{\sum_{i=1}^n(\sigma_i-b)}{\sqrt n}\rig)\bigg|m<0\rig]=\bb E\lef[\frac{\sum_{\bfa\sigma:m<0}\exp\lef(\theta_1m^2n/2+\sum_{i=1}^n(\frac{t}{\sqrt n}+h_i)\sigma_i\rig)\exp\lef(-\sqrt n tb\rig)}{\sum_{\bfa\sigma:m<0}\exp\lef(\theta_1m^2n/2+\sum_{i=1}^nh_i\sigma_i\rig)}\rig]\nnb\\
        &=\bb E\lef[\frac{\int_{\bb R}\sum_{m<0}\exp\lef(-\frac{n}{2}x^2+\sum_{i=1}^n(\sqrt{\theta_1}x+\frac{t}{\sqrt n}+h_i)\sigma_i)\rig)\exp\lef(-\sqrt ntb\rig)dx}{\int_{\bb R}\sum_{m<0}\exp\lef(-\frac{n}{2}x^2+\sum_{i=1}^n(\sqrt{\theta_1}x+h_i)\sigma_i)\rig)dx}\rig]
    \end{align}}
    Then we introduce a probability measure $\rho(m):=\frac{\sum_{\bfa\sigma:m}\exp\lef(\sum_{i=1}^n( h_i+\sqrt{\theta_1}x)\sigma_i\rig)}{2^n\prod_{i\in[n]}\cosh\lef(\sqrt{\theta_1}x+h_i\rig)}$ is the product measure of independent Bernoulli r.v. It is the not hard to see that under $\rho$, there exists $C>0$ such that by classical result of Large Deviation Principle:
    \tny{\begin{align*}
        &\bb E_{\rho}[m]=\frac{1}{n}\sum_{i=1}^n\tanh(\sqrt{\theta_1}x+h_i),\quad\rho\bl \bigg|m-\frac{1}{n}\sum_{i=1}^n\tanh(\sqrt{\theta_1}x+h_i)\bigg|>t\br\leq 2\exp( -Cnt^2).
    \end{align*} }
    Further we notice that for $t>0$ there exists $C>0$ such that
    \begin{align*}
        \bb P\bigg(\bigg|\frac{1}{n}\sum_{i=1}^n\tanh(\sqrt{\theta_1}x+h_i)-\bb E[\tanh(\sqrt{\theta_1}x+h_i)]\bigg|>t\bigg)\leq\exp(-Cnt^2).
    \end{align*}
    Therefore, for $\frac{1}{n}\sum_{i=1}^n\tanh(\sqrt{\theta_1}x+h_i)<0$, there exists $C_1,C_2,C_3>0$ such that for all $x=o(1)$,
   \tny{ \begin{align}\label{tailbd}
        \rho(m>0)&\leq2\exp\bl-\frac{C_1}{n}\bl\sum_{i=1}^n\tanh(\sqrt{\theta_1}x+h_i)\br^2\br\leq 2\exp\lef(-C_2n\bb E[\tanh(\sqrt{\theta_1}x+h)]^2\rig)\nnb\\
        &\leq2\exp\lef(-C_3nx^2\rig).
    \end{align}}
    And furthermore,  we have 
    to get that for some $C>0$  we have
    \begin{align*}
        \int_{\bb R}\sum_{m<0}&\exp\bl-\frac{n}{2}x^2+\sum_{i=1}^n\log\cosh(\sqrt{\theta_1}x+h_i)\br\rho(m) dx\\
        &=\ub{\int_{x< 0}\exp\bl-\frac{n}{2}x^2+\sum_{i=1}^n\log\cosh(\sqrt{\theta_1}x+h_i)\br dx\sum_{m\in[-1,1]}\rho(m)}_{T_0}\\
        &-\ub{\int_{x< 0}\exp\bl-\frac{n}{2}x^2+\sum_{i=1}^n\log\cosh(\sqrt{\theta_1}x+h_i)+\log(\rho(m>0))\br dx}_{T_1}\\
        &+\ub{\int_{x\geq 0}\exp\bl-\frac{n}{2}x^2+\sum_{i=1}^n\log\cosh(\sqrt{\theta_1}x+h_i)+\log(\rho(m<0))\br dx}_{T_2}.
    \end{align*}
    We note the fact that at low temperature $-\frac{n}{2}x^2+\sum_{i=1}^n\log\cosh(\sqrt{\theta_1}x+h_i)$ take its local maximum at $x_1=\frac{1}{n}\sum_{i=1}^n\sqrt{\theta_1}\tanh(\sqrt{\theta_1}x_1+h_i)<0,x_2=\frac{1}{n}\sum_{i=1}^n\sqrt{\theta_1}\tanh(\sqrt{\theta_1}x_2+h_i)>0$.
    Hence  we have by \eqref{tailbd}, there exists $C>0$ such that
    \begin{align*}
        T_1&\leq\int_{x< -\frac{1}{2}x_1}\exp\bl-\frac{n}{2}x^2+\sum_{i=1}^n\log\cosh(\sqrt{\theta_1}x+h_i)-Cnx^2\br dx\\
        &+\int_{x\in(-\frac{1}{2}x_1,0)}\exp\bl-\frac{n}{2}x^2+\sum_{i=1}^n\log\cosh(\sqrt{\theta_1}x+h_i)\br dx\leq T_0\exp(-C_1n).
    \end{align*}
    And analogously  we have $T_2\leq T_0\exp(-C_2n)$. 
    Hence
    \begin{align*}
        \int_{\bb R}\sum_{m<0}\exp\bl-\frac{n}{2}x^2+&\sum_{i=1}^n\log\cosh(\sqrt{\theta_1}x+h_i)\br\rho(m) dx=T_0(1+O\lef(\exp(-C_1n)\rig)).
    \end{align*}
    Noticing that the effect given by $\frac{t}{\sqrt n}$ is at most $\exp(\sqrt nC)$, which implies that there exists $C_2$ such that
    \begin{align*}
        \int_{\bb R}&\sum_{m<0}\exp\bl-\frac{n}{2}x^2+\sum_{i=1}^n\log\cosh\lef(\sqrt{\theta_1}x+h_i+\frac{t}{\sqrt n}\rig)\br\rho(m) dx\\
        &=\int_{x<0}\exp\bl-\frac{n}{2}x^2+\sum_{i=1}^n\log\cosh\lef(\sqrt{\theta_1}x+h_i+\frac{t}{\sqrt n}\rig)\br dx\lef(1+O(\exp(-C_2n))\rig).
    \end{align*}
    And we collect pieces to conclude that \eqref{transferprincipa} become the following for some $\delta>0$
    \tny{\begin{align*}
        \bb E\lef[\exp\lef(t\frac{\sum_{i=1}^n(\sigma_i-b)}{\sqrt n}\rig)\bigg|m<0\rig]&=\bb E\lef[\frac{\int_{x<0}\exp\lef(-n\mca H_{0,n}(x)-\sqrt ntb\rig)dx}{\int_{x<0}\exp\lef(-n\mca H_{1,n}(x)\rig)dx}\rig]\lef(1+O\lef(\exp(-n\delta)\rig)\rig).
    \end{align*}}
    And the sub-Gaussian properties can be analogously derived by similar argument as the high temperature regime.
    And then we can use similar method as the high temperature case to derive that at the low temperature we concentrate
    on $m_1<0,m_2>0$ with
    \begin{align*}
        m_1=\bb E[\tanh(\theta_1m_1+h)]>0,\qquad m_2=\bb E[\tanh(\theta_1m_2+h)]<0.
    \end{align*}
    \begin{center}
        \textbf{3. The Critical Temperature Case}
    \end{center}
    Then we consider the general situation at the critical temperature $\theta_1=\frac{1}{\bb E[\sech^2(h)]}$. 
    Recall that we define the critical value $\tau$ for the critical temperature as
    \begin{align*}
        \mca H_{1}(x)=\mca H_{1}(x_1^*)+\frac{1}{(2\tau)!}\mca H_{1}^{(2\tau)}(x_1^*)(x-x_1^*)^{2\tau}+O((x-x_1^*)^{2\tau+1}),
    \end{align*}
    with $\mca H_{1}^{(2\tau)}(x_1^*)>0$ and $x_1^*$ is the unique minimum of $\mca H_1$.
    Then we study the fluctuation of $x_1$. Using the Fermat's condition  we have
    \begin{align*}
    0=\mca H_{1,n}^\prime(x_{1})&=\mca H_{1,n}^\prime(x_1^*)+\ub{\sum_{i=2}^{2\tau-1}\frac{1}{i!}\mca H_{1,n}^{(i)}(x_1^*)(x_1-x_1^*)^{i-1}}_{T_0}\\
    &+\frac{1}{(2\tau)!}\mca H_{1,n}^{(2\tau)}(x_1^*)(x_1-x_1^*)^{2\tau-1}+O((x_1-x_1^*)^{i+1}).
    \end{align*}
    Hence, noticing that the middle term $T_0=o_{\psi_2}\lef(\frac{1}{\sqrt n}\rig)$  we have
$
        (x_1-x_1^*)^{2\tau-1}=\frac{-\mca H_{1,n}^\prime(x_1^*)}{\mca H_{1,n}^{(2\tau)}(x_1^*)}+o_{\psi_2}\lef(\frac{1}{\sqrt n}\rig),
$
        which implies that
        \tny{\begin{align}\label{weakconvergencex1}
        \sqrt n(x_1-x_1^*)^{2\tau-1}&=\frac{-(2\tau)!\sqrt{\theta_1}}{\sqrt nH^{(2\tau)}_{1,n}(x_1^*)}\sum_{i=1}^n\lef(\tanh(\sqrt{\theta_1}x_1^*+h_i)-\bb E[\tanh(\sqrt{\theta_1}x_1^*+h_i)]\rig)+o_{\psi_2}(1)\nnb\\
        &\overset{d}{\to} N\lef(0,\frac{((2\tau)!)^2\theta_1\bb V(\tanh(\sqrt{\theta_1}x_1^*+h))}{(\mca H_{1}^{(2\tau)}(x_1^*))^2}\rig).
    \end{align}}
    where we already use the fact that by the law of large numbers $\mca H_{1,n}^{(2\tau)}(x_1^*)\overset{a.s.}{\to} \mca H_{1}^{(2\tau)}(x_1^*)$. 
    Similarly, we notice that
    \begin{align*}
    0=\mca H_{0,n}^\prime(x_{0})&=\mca H_{0,n}^\prime(x_0^*)+\ub{\sum_{i=2}^{2\tau-1}\frac{1}{i!}\mca H_{1,n}^{(i)}(x_0^*)(x_1-x_0^*)^{i-1}}_{T_0}\\
    &+\frac{1}{(2\tau)!}\mca H_{1,n}^{(2\tau)}(x_0^*)(x_0-x_0^*)^{2\tau-1}+O((x_0-x_0^*)^{i+1}).
    \end{align*}
    When $\beta>\frac{1}{2}$,  we have $T_0=o_{\psi_2}\lef(\frac{1}{\sqrt n}\rig)$. Then we apply similar arguments as $\mca H_{1,n}$ to get
    \begin{align*}
         \sqrt n(x_0-x_0^*)^{2\tau-1}&=\frac{-\sqrt{\theta_1}}{\sqrt nH^{(2\tau)}_{0,n}(x_0^*)}\sum_{i=1}^n\bl\tanh\lef(\sqrt{\theta_1}x_0^*+h_i+\frac{t}{n^{\beta}}\rig)\\
         &-\bb E\lef[\tanh\lef(\sqrt{\theta_1}x_1^*+h_i+\frac{t}{n^{\beta}}\rig)\rig]\br+o_{\psi_2}(1).
    \end{align*}
    Denote $\Delta:=\sqrt{\theta_1}(x_1^*-x_0^*)-\frac{t}{n^{\beta}}$,  we have by the Fermat's condition
    \sm{\begin{align*}
        \Delta +\frac{t}{n^{\beta}}&=\theta_1\lef(\bb E[\tanh(\sqrt{\theta_1}x_1^*+h)]-\bb E\lef[\tanh\lef(\sqrt{\theta_1}x_0^*+\frac{t}{n^{\beta}}+h\rig)\rig]\rig)\\
        &=\theta_1\lef(\bb E[\sech^2(\sqrt{\theta_1}x_1^*+h)]\Delta +\sum_{i=3}^{2\tau-1}\frac{\mca H_{1}^{(i)}(x_1^*)}{(i-2)!}\Delta^{i-1}+\frac{\mca H_1^{(2\tau)}(x_1^*)}{(2\tau-2)!}\Delta^{2\tau-1}\rig)+O(\Delta^{2\tau}).
    \end{align*}}
    By similar argument the middle term is small in order  we have
$        \frac{t}{n^{\beta}}=\frac{\mca H_1^{(2\tau)}(x_1^*)}{(2\tau-2)!}\Delta^{2\tau-1}+O(\Delta^{2\tau}).
$
    And
    \begin{align*}
        \sqrt{\theta_1}(x_1^*-x_0^*)=-\lef(\frac{(2\tau-2)!t}{\theta_1 \mca H_1^{(2\tau)}(x_1^*)n^{\beta}}\rig)^{\frac{1}{2\tau-1}}(1+o_{\psi_2}(1)).
    \end{align*}
    Notice that by binomial expansion there exists $C(\tau)$ such that $|A-B|\asymp C(\tau)\frac{|A^{2\tau-1}-B^{2\tau-1}|}{(A^{2\tau-2}+B^{2\tau-2})}$,
    \tny{\begin{align*}
        x_1-x_1^*-x_0+x_0^*&=O\lef(\frac{t}{n^{\beta}}\frac{\sum_{i=1}^n(\sech^2(h_i+\sqrt{\theta_1}x_1^*)-\bb E[\sech^2(h_i+\sqrt{\theta_1}x_1^*)])}{n(x_1-x_1^*)^{2\tau-2}}\rig)(1+o_{\psi_2}(1))\\
        &=O\lef(\frac{t}{n^{\beta}}\rig)(1+o_{\psi_2}(1)).
    \end{align*}}
    Therefore, we conclude that
    \tny{\begin{align*}
        x_1-x_0=x_1-x_1^*+x_0^*-x_0+(x_1^*-x_0^*)=-\sign(t)\frac{1}{\sqrt{\theta_1}}\lef(\frac{(2\tau-2)!|t|}{\theta_1 \mca H_1^{(2\tau)}(x_1^*)n^{\beta}}\rig)^{\frac{1}{2\tau-1}}(1+o_{\psi_{4\tau}}(1)).
    \end{align*}}
    Denote $\times$ as the tensor product and $\delta:=(x_0-x_1,t/n^{\beta})$, we pick proper $\beta$ such that when $t=o(n^{\beta})$,
    \begin{align*}
        n\mca H_{1,n}(x_1,\bfa h)&-n\mca H_{0,n}(x_0,\bfa h)=n\lef(\sum_{i=1}^\infty\frac{1}{i!}F^{(i)}(x_1,0)\times\delta^i\rig)\\
        &=\ub{n\frac{1}{(2\tau)!}F_x^{(2\tau)}(x_1,0)(x_0-x_1)^{2\tau}(1+o_{\psi_2}(1))}_{A_1}+\ub{n\frac{\pta F(x_1,0)}{\pta y}\frac{t}{n^{\beta}}(1+o_{\psi_2}(1))}_{A_2}.
    \end{align*}
    Then is to decide the magnitude of the two terms in the bracket, we notice that the first term is in the order of $O_{\psi_2}\lef(n^{-\beta\lef(\frac{2\tau}{2\tau-1}\rig)}\rig)$ and the second order is
    \begin{align*}
        \frac{\pta F(x_1,0)}{\pta y}\frac{t}{n^{\beta}}&=\frac{1}{n}\sum_{i=1}^n\tanh(\sqrt{\theta_1}x_1+h_i)\frac{t}{n^{\beta}}\\
        &=\frac{1}{n}\sum_{i=1}^n\lef(\tanh(\sqrt{\theta_1}x_1^*+h_i)+\sech^2(\sqrt{\theta_1}x_1^*+h_i)\sqrt{\theta_1}x_1\rig)\frac{t}{n^{\beta}}+o_{\psi_2}(1)\\
        &=\frac{t}{n}\sum_{i=1}^n\sech^2(\sqrt{\theta_1}x_1^*+h_i)\sqrt{\theta_1}\frac{1}{n^{\beta}}x_1+o_{\psi_2}\lef(1\rig).
    \end{align*}
    Therefore $A_2$ is more significant in order. Recall the weak convergence result given by \eqref{weakconvergencex1}. Finally we pick $\beta=\frac{4\tau-3}{4\tau-2}$ and check that for all $t\in\bb R$ not dependent on $n$  we have boundedness of mgf. By theorem 2 in \citep{chareka2008converse} we get
    \begin{align}\label{fromweaktomgf}
        \bb E\lef[\exp\lef(\frac{t\sum_{i=1}^n(\sigma_i-\bb E[\sigma_i])}{n^{\frac{4\tau-3}{4\tau-2}}}\rig)\rig]\to \int_{\bb R}\frac{(2\tau-1)x^{2\tau-2}}{\sqrt{2\pi v}}\exp\lef(-\frac{x^{4\tau-2}}{2v}+tx\rig)dx,
    \end{align}
    with $v:=\frac{((2\tau)!)^2\theta_1^{2\tau}\bb V(\tanh(\sqrt{\theta_1}x_1^*+h))(\bb E[\sech^2(\sqrt{\theta_1}x_1^*+h)])^{4\tau-2}}{(\mca H_{1}^{(2\tau)}(x_1^*))^2}$. In particular, using the result in \citep{boyadzhiev2009derivative}  we have
    \begin{align*}
        \mca H_{1}^{(2\tau)}(x_1^*)=\theta_1^{\tau}\bb E\lef[(1+\tanh(\sqrt{\theta_1}x_1^*+h))\sum_{k=0}^{2\tau-1}\frac{k!}{2^k}S(2\tau-1,k)(\tanh(\sqrt{\theta_1}x_1^*+h)-1)^{k}\rig],
    \end{align*}
    where $S(2\tau-1,k)$ is the Stirling number of the second kind. Take $x_1^*=0$ and we complete the proof of theorem \ref{cltrfcr}.
    To derive the sub-Weibull properties, we first use H\"older's inequality to upper bound the m.g.f. for all $t=o(n^{\frac{4\tau-3}{4\tau-2}})$ where there exists $\delta_1,\delta_2>0$ such that $\frac{1}{1+\delta_1}+\frac{1}{1+\delta_2}=1$ and
    \begin{align}\label{mgfub}
        \bb E[\exp(tX)]&=\bb E[\exp(n\mca H_{1,n}(x_1,\bfa h)-n\mca H_{0,n}(x_0,\bfa h))](1+o(1))\nnb\\
        &\leq\bb E[\exp((1+\delta_1)A_1)]^{\frac{1}{1+\delta_1}}\bb E[\exp((1+\delta_2)A_2)]^{\frac{1}{1+\delta_2}}(1+o(1)).
    \end{align}
    We define $Y_i=\tanh(\sqrt{\theta_1}x_1^*+h_i)-\bb E[\tanh(\sqrt{\theta_1}x_1^*+h_i)]$, using the results by \citep{packwood2011moments}  we have for $\eta\in\bb N\cap[1,\infty),\xi\in(1,\infty)$, by $Y_i\in[-1,1]$, introducing $\{\epsilon_i\}_{i\in[n]}$ i.i.d. Rademacher random variables,
    \begin{align*}
        \bb E[Y_i^{2\eta}]\leq\bb E[Y_i^{2}]\quad\Rightarrow\quad\bb E\bigg[\bl\frac{1}{\sqrt n}\sum_{i=1}^n\epsilon_iY_i\br^{\xi}\bigg]\leq C^\xi\bb E[Z^\xi]\text{ with }Z\sim N(0,1).
    \end{align*}
    Then by Taylor expansion and symmetrization arguments in \citep{wainwright2019high} Proposition 4.11, introducing $\{\epsilon_{i}\}_{i\in[n]}$ i.i.d. Rademacher random variables, consider the principle term in $A_2$ we have,
    \begin{align*}
        \bb E\bigg[\exp\bl &Ct\sign\bl\sum_{i=1}^nY_i\br\bl n^{-1/2}|\sum_{i=1}^nY_i|\br^{\frac{1}{2\tau-1}} \br\bigg]=\sum_{i=1}^\infty\frac{C^it^i}{i!n^{\frac{i}{4\tau-2}}}\bb E\bigg[\bl\sum_{i=1}^nY_i\br^{\frac{i}{2\tau-1}}\bigg]\\
        &\leq\sum_{i=1}^{2\tau-1}\frac{C^it^i}{i!n^{\frac{i}{4\tau-2}}}\bb E\bigg[\bl\sum_{i=1}^nY_i\br^{\frac{i}{2\tau-1}}\bigg]+\sum_{i=2\tau-1}^\infty\frac{C^it^i}{i!n^{\frac{i}{4\tau-2}}}\bb E\bigg[\bl\sum_{i=1}^n\epsilon_iY_i\br^{\frac{i}{2\tau-1}}\bigg]\\
        &\leq C+\sum_{i=0}^\infty\frac{C^it^i}{i!}\bb E[Z^{\frac{i}{2\tau-1}}]=C+C\int_{\bb R}x^{2\tau-2}\exp(-Cx^{4\tau-2}+tx)dx\\
        &\leq C+ C\exp(Ct^{\frac{4\tau-1}{4\tau-2}}).
    \end{align*}
    And going back to \eqref{mgfub} we conclude that for all $t=o(n^{\frac{4\tau-3}{2\tau-2}})$,
    \begin{align*}
        \bb E[\exp(tX)]\leq C+C\exp(Ct^{\frac{4\tau-1}{4\tau-2}}).
    \end{align*}
    Therefore by Chernoff bound, one will get for all $t=o(n^{\frac{4\tau-2}{4\tau-3}})$,
    \begin{align*}
        \bb P(|X|\geq t)\leq C\inf_{\lambda}\exp(C\lambda^{\frac{4\tau-1}{4\tau-2}}-t\lambda)+C\exp(-t\lambda)\leq C\exp(-C t^{4\tau-2}).
    \end{align*}
    Notice that $|X|\leq n^{\frac{1}{4\tau-2}}$.
    Then we apply similar procedure to boost from $t=o(n^{\frac{4\tau-2}{4\tau-3}})$ to complete $t\in\bb R^+$ by computing the moments, which are given by, for all $p\in\bb N$, 
    \begin{align*}
       \bb E[|X|^p] &=\int_{\bb R^+}\bb P(|X|^p\geq u)du=\int_{\bb R^+}\bb P(|X|\geq t)pt^{p-1}dt\\
       &\leq \int_{0}^{n^{\frac{4\tau-3}{4\tau-2}}/\log n}C\exp(-Ct^{4\tau-2})pt^{p-1}dt+n^{\frac{p}{4\tau-2}}\bb P(|X|\geq n^{\frac{4\tau-3}{4\tau-2}}/\log n)\\
       &\leq C\Gamma\lef(1+\frac{p}{4\tau-2}\rig)+n^{\frac{p}{4\tau-2}}C\exp(-Cn^{-{4\tau-3}}/\log^{2\tau-2}n)\\
       &\leq (Cp)^{\frac{p}{4\tau-2}}.
    \end{align*}
    And finally  we have for all $p\in\bb N$, $\bb E[|X|^p]^{\frac{1}{p}}\leq Cp^{\frac{1}{4\tau-2}}$.
    Therefore we use theorem 2.1 in \citep{vladimirova2020sub} to complete the proof of sub-Weibull property of $X$ with $\Vert X\Vert_{\psi_{4\tau-2}}<\infty$.

\section{Proof of Results in Section \ref{sect6}}
Here we provide the detailed proof of the propositions in section 6. 
\subsection{Proof of Proposition \ref{testingcomputationalupperbound}}

To prove that the SDP relaxation recovers the hidden clique, we will first utilize the dual program of \ref{alg:sdp}, which is also studied in \citep{bach2010convex} as
\begin{align*}
    SDP_k(\wh\Sigma)\leq\min_{Y\in\bfa S_p}\lambda_{\max}(\wh\Sigma+Y)+k\Vert Y\Vert_{\infty}.
\end{align*}
Essentially, we may write the program of recovering the $k$-set $S$ of indices with the largest $F$-norm of sub-covariance matrix $\Sigma_{SS}$ as the following program:
\begin{center}
Maximize $\tr(\wh\Sigma X)$  \\
subject to $\tr(X)=1$, $\Vert X\Vert_0\leq k^2$, $\rank(X)=1$, $X\succeq 0$.
\end{center}
However, one will check that this program has the rank constraint and will not be computational efficiently solvable. Therefore we just drop the rank constraint to obtain the polynomially solvable program.
We denote the optimal objective value of the above program as $\lambda_{\max}^k(\Sigma)$ since it is also the maximum biggest eigenvalue of the $k\times k$ sub-matrix of $\Sigma$.
And one can easily relax it toward the program given by \ref{alg:sdp}. This implies that for all $Y\in \bfa S_p$ we will have
\begin{align*}
    \lambda^k_{\max}(\wh\Sigma)\leq SDP_k(\wh\Sigma)\leq \lambda_{\max}(\wh\Sigma+Y)+k\Vert Y\Vert_{\infty,\infty}.
\end{align*}
We first analyze the upper bound under the alternative.
We construct a soft thresholding matrix $Y$ as $Y_{ij}=-\wh\Sigma_{ij}\mbbm 1_{|\wh\Sigma_{ij}|\leq z}-\sign(\wh\Sigma_{ij})z\mbbm 1_{|\wh\Sigma_{ij}|> z}$. We will also denote $\bb P_{0}$ be the meaure under the null and $\bb P_{S}$ be the measure under the alternative. We define the event $\ca C:=\lef\{\sup_{i,j:i\neq j}\lef|\frac{1}{m}\sum_{k=1}^m\sigma^{(k)}_i\sigma^{(k)}_j\rig|\leq t\rig\}$ and
\begin{align*}
    \bb P_0(\ca C^c):=\bb P_S\bl\sup_{i,j:i\neq j,}\bigg|\frac{1}{m}\sum_{k=1}^m\sigma^{(k)}_i\sigma^{(k)}_j\bigg|\geq t\br\leq n^2 C\exp(-Cmt^2).
\end{align*}
Therefore, if we pick $t=C_1\sqrt{\frac{\log n}{m}}$, then we will have $\bb P_0(\ca C^c)=o(1)$. We also notice that 
conditional on this event $\ca C$, we will check that $Y+\wh\Sigma $ is a matrix with zero off-diagonal entries and we will have
\begin{align*}
    \lambda_{\max}(\wh\Sigma+Y) = 1.
\end{align*}

Define $f(C_1)=1+C_1k\sqrt{\frac{\log n}{m}}$.
First we control the Type I error and see that when $m\gtrsim k^2\log n$,
\sm{\begin{align*}
   \bb P_0\bl SDP_k(\wh\Sigma)\geq f(C_1)\br\leq\bb P_0(\lambda_{\max}(\wh\Sigma+Y)+k\Vert Y\Vert_{\infty,\infty}\geq f(C_1) |\ca C)\bb P_0(\ca C)+\bb P_0(\ca C^c)=o(1).
\end{align*}}
And then we consider the Type II error, noticing that for $v\in\bb R^n$ such that $v_i=\frac{1}{\sqrt k}, \forall i\in S$ and $v_i=0, \forall i\in S^c$,
\begin{align*}
    SDP_k(\wh\Sigma)&\geq\lambda_{\max}^k(\wh\Sigma)\geq v^\top\wh\Sigma v=\frac{1}{m}\sum_{j=1}^mv^\top\bfa\sigma^{(i)}\bfa\sigma^{(i)\top} v-\Vert v\Vert_2^2\\
    &=\frac{1}{m}\sum_{j=1}^m(\bfa\sigma^{(i)\top}v)^2-1=\frac{1}{mk}\sum_{j=1}^m\bl\sum_{i\in S}\sigma_i^{(j)}\br ^2-1.
\end{align*}
And here we will discuss the different scenarios based on their temperature regimes.
\begin{center}
    \textbf{I. High Temperature}
\end{center}
Using the result given by lemma \ref{cltrfcr}, we will have $\lef\Vert\frac{1}{\sqrt k}\sum_{i\in S}\sigma_i\rig\Vert_{\psi_{2}}<\infty$ and 
\begin{align*}
    \bb E\bigg[\frac{1}{k}\bl\sum_{i\in S}\sigma_i\br^2\bigg]-1&=\frac{1-\theta_1(\bb E[\sech^2(h)])^2}{(1-\theta_1\bb E[\sech^2(h)])^2}-1+o(1)\\
    &=\frac{2\theta_1\bb E[\sech^2(h)]-\theta_1(1+\theta_1)(\bb E[\sech^2(h)])^2}{(1-\theta_1\bb E[\sech^2(h)])^2}+o(1)=:C_0.
\end{align*}
Therefore, noticing that $\Vert\frac{1}{k}(\sum_{i\in S}\sigma_i)^2\Vert_{\psi_1}<\infty$ and use Bernstein inequality we will get
\begin{align*}
    \bb P_S(SDP_k(\wh\Sigma)&<f(C_1))\leq \bb P_S\bl\frac{1}{mk}\sum_{j=1}^m\bl\sum_{i\in S}\sigma_i^{(j)}\br^2-1<f(C_1)\br\\
    &=\bb P_S\bl\frac{1}{mk}\sum_{j=1}^m\bl\sum_{i\in S}\sigma_i^{(j)}\br^2-\frac{1}{k}\bb E\bigg[\bl\sum_{i\in S}\sigma_i^{(j)}\br^2\bigg]< f(C_1)-C_0\br\\
    &\leq C\exp\lef(-Cm(f(C_1)-C_0)^2\wedge m(C_0-f(C_1))\rig),
\end{align*}
Therefore, when $f(C_1)<C_0$ and we pick $m\gtrsim k^2\log n$, the test will be asymptotically powerful.
\begin{center}
    \textbf{II. Low Temperature}
\end{center}
For the low temperature regime, we set $\frac{1}{k}f(C_1)$ as the threshold and the Type I error is controlled by
\begin{align}\label{contTypeIerror}
    \bb P_0\bl \frac{1}{k}SDP_k(\wh\Sigma)>\frac{1}{k}f(C_1)\br=o(1),
\end{align}
when $m\gtrsim \log n$.
And notice that by lemma \ref{cltrfcr}, $
    \bb E\lef[\frac{1}{k^2}(\sum_{i\in S}\sigma_i)^2\rig]=C_2\asymp 1.
$

there exists $\mu>0$ such that $\Vert \sqrt k(\frac{1}{k}|\sum_{i\in S}\sigma_i|-\mu)\Vert_{\psi_2}<\infty$, and one will get by the centering argument of sub-Gaussian r.v.s given by \citep{vershynin2018high}:
\begin{align*}
    \bigg\Vert &\sqrt k\bl\bl\frac{1}{k}\sum_{i\in S}\sigma_i\br^2-\bb E\bigg[\bl\frac{1}{k}\sum_{i\in S}\sigma_i\br^2\bigg]\br\bigg\Vert_{\psi_2}\leq\bigg\Vert \sqrt k\bl\bl\frac{1}{k}\sum_{i\in S}\sigma_i\br^2-\mu^2\br\bigg\Vert_{\psi_2}\\
    &\lesssim\bigg\Vert \sqrt k\bl\frac{1}{k}\bigg|\sum_{i\in S}\sigma_i\bigg|-\mu\br\bigg\Vert_{\psi_2}<\infty.
\end{align*}
And, noticing that $\frac{1}{k}f(C_1)\asymp 1\asymp C_2>0$, there exists $C_3>0$ such that with some small choice of $C_1$,
\sm{\begin{align*}
    \bb P_S\bl\frac{1}{k}SDP_k(\wh\Sigma)&< \frac{1}{k}f(C_1)\br\leq \bb P_S\bl\frac{1}{mk^2}\sum_{j=1}^m\bl\sum_{i\in S}\sigma_i^{(j)}\br^2-\frac{1}{k}\leq \frac{1}{k}f(C_1)\br\\
    &\leq\bb P_S\bl\frac{1}{mk^2}\sum_{j=1}^m\bl\sum_{i\in S}\sigma_i^{(j)}\br^2-\frac{1}{k^2}\bb E\bigg[\bl\sum_{i\in S}\sigma_i^{(j)}\br^2\bigg]<\frac{1}{k}f(C_1)+\frac{1}{k}-C_2\br\\
    &\leq C\exp\lef(-Cm kC_3^2\rig)=o(1),
\end{align*}}
when $m\gtrsim \log n$.
\begin{center}
    \textbf{ III. Critical Temperature}
\end{center}
Before we start analyzing the critical temperature, it is not hard to see that the Type I error can always be upper bounded as for all $\beta>0$,
\begin{align*}
    \bb P(k^{-\beta}SDP_k(\wh\Sigma)\geq k^{-\beta}f(C_1))=o(1).
\end{align*}
However, to understand the concentration of measure for the critical temperature case involve more effort. 
Also notice that using standard results of sub-Weibull norm, (see theorem 2.1 in \citep{vladimirova2020sub}) one will have for all $\theta>1$:
\begin{align*}
    \Vert X\Vert_{\psi_\theta}<\infty\quad\Leftrightarrow\quad \bb P(|X|\geq t)\leq C\exp(-t^{\theta}),\quad\text{ for all }t\in\bb R.
\end{align*}
Then we can conclude by applying the above equation for all the m.g.f. given by lemma \ref{cltrfcr}.
When $0$ is the global minimum of \eqref{taylorcond} with flatness $\tau$ we will have $\Vert k^{-\frac{4\tau-3}{4\tau-2}}\sum_{i\in S}\sigma_i\Vert_{\psi_{4\tau-2}}<\infty$.
The more detailed analysis goes to upper bounding the Type II error. Here we will address the different scenarios separately. Given lemma \ref{cltrfcr}, there exists $\gamma>0$ such that
\begin{align*}
    \bigg\Vert k^{-\gamma}\sum_{i\in S}\sigma_i \bigg\Vert_{\psi_{\frac{1}{1-\gamma}}}<\infty,\quad\Rightarrow\quad \bigg\Vert k^{-2\gamma}\bl\bl\sum_{i\in S}\sigma_i\br^2-\bb E\bigg[\bl\sum_{i\in S}\sigma_i\br^2\bigg]\br\bigg\Vert_{\psi_{\frac{2}{1-\gamma}}}<\infty.
\end{align*}
And therefore we will use the Sub-Weibull concentration inequality \citep{zhang2022sharper} to get
\tny{\begin{align*}
    \bb P_{S}\bl k^{-2\gamma+1}SDP_k(\wh\Sigma)&\leq k^{-2\gamma+1}f(C_1) \br\leq \bb P_S\bl-\frac{1}{mk^{2\gamma}}\sum_{j=1}^m\bl\sum_{i\in S}\sigma_i\br^2-k^{-2\gamma+1}\leq k^{-2\gamma+1}f(C_1)\br\\
    &\leq \bb P_S\bl m^{-1}k^{-2\gamma}\sum_{j=1}^m\bl\sum_{i\in S}\sigma_i^{(j)}\br^2-k^{-2\gamma}\bb E\bigg[\bl\sum_{i\in S}\sigma_i^{(j)}\br^2\bigg]\leq C_3\br\\
    &\leq C\exp(-Cm)=o(1).
\end{align*}}
when $m\gtrsim k^{-4\gamma+4}\log n$. Replacing $\gamma=\frac{4\tau-3}{4\tau-2}$ we complete the proof.

\subsection{Proof of Proposition \ref{propexactrecov}}
Here we use $\gamma=\frac{4\tau-3}{4\tau-2}$ as the shorthand notation.
    The proof idea follows by analyzing the rank 1 programs given by 
    \begin{align}\label{non-convexsdp}
        \wh Z=\argmax_{\rank(Z)=1}\tr(\wh\Sigma Z).
    \end{align}
    And show that the algorithm \ref{sdprecovery} can return the same solution as the program given sufficiently large sample size.
     The proof idea is by first constructing a pair of primal and dual feasible matrices, then prove that they yield the desired selection property. This is also referred to as the \emph{primal dual certificate} in \citep{amini2009high}. The following lemma states the sufficient condition of primal and dual feasibility, which appears in \citep{amini2009high}, and we include here for completeness
    \begin{lemma}[\citep{amini2009high}]\label{primaldualfeasi}
        Suppose that, for each $x\in\bb R^p$ with $\Vert x\Vert_2=1$, there exists a sign matrix $\wh U(x)$ that satisfies the dual feasibility condition given by
        \begin{align*}
            \wh U_{ij}(x)=\begin{cases}
                \sign(\wh z_i)\sign(\wh z_j),&\text{ if }\wh z_i\wh z_j\neq 0\\
                \in[-1,1],&\text{ otherwise}
            \end{cases}.
        \end{align*}
        And the primal variable $\wh z\wh z^\top$ satisfies the primal feasibility condition given by
        \begin{align}\label{ineq}
            \forall x\in\bb R,\qquad & x^\top(\wh \Sigma-\rho\wh U(x))x\leq \wh z^\top(\wh \Sigma-\rho\wh U(x))\wh z. \quad\text{ when }\theta_1>0,\\
            & x^\top(\wh \Sigma-\rho\wh U(x))x\geq \wh z^\top(\wh \Sigma-\rho\wh U(x))\wh z. \quad\text{ when }\theta_1<0.
        \end{align}
        Then $\wh z$ is the solution to \eqref{non-convexsdp}.
    \end{lemma}
    Then, based on the above lemma, we will construct the primal variable $\wh Z$ and dual variable $\wh U$ and parameter $\rho$ as follows:
    \begin{align}\label{designsdp}
        \wh z&=(\wh z_S,\bfa 0),\qquad \wh z_S=\argmax_{z\in\bb R^S}z\lef(\wh\Sigma_{SS}-\rho\wh U_{SS}\rig)z^\top,\nnb\\ \wh U_{ij}(x)&=\begin{cases}
        -\frac{1}{\rho}\wh\Sigma_{ij}\mbbm 1_{i\neq j}   &\text{ when } (i,j)\notin S\times S\\
        1&\text{ otherwise } 
        \end{cases}.
    \end{align}
    Then we will choose $\rho$ based on the variance estimates. This gives us the following (here the constants $C$ might differ across temperature regimes.) 
    \tny{\begin{align*}
        \rho:&=\frac{1}{2}\bb E[\sigma_i\sigma_j]\\
        &=\begin{cases}
            Ck^{-1}&\text{ when }\theta_1\in\lef(0,\frac{1}{\bb E[\sech^2(h)]}\rig)\\
            C k^{-\frac{1}{2\tau-1}}&\text{ when }\theta_1=\frac{1}{\bb E[\sech^2(h)]}\text{ and } 0 \text{ is a global minimum of \eqref{taylorcond} has flatness }\tau \\
            C&\text{ when } \theta_1\in\lef(\frac{1}{\bb E[\sech^2(h)]},\infty\rig)\text{ or }\theta_1=\frac{1}{\bb E[\sech^2(h)]}\text{ and }$0$ \text{ is a global maximum of \eqref{taylorcond}}
        \end{cases}
    \end{align*}}
    Denote $\Sigma$ as the expected covariance matrix, the matrix $\wh\Sigma$ can be extended as \[\wh\Sigma-\rho \wh U(x)=\begin{bmatrix} \wh\Sigma_{SS}-\rho\wh U_{SS}&\wh\Sigma_{SS^c}-\rho\wh U_{SS^c}\\
    \wh\Sigma_{S^cS}-\rho \wh U_{S^cS}&\wh\Sigma_{S^cS^c}-\rho\wh U_{S^cS^c}\end{bmatrix}=\begin{bmatrix}
        \frac{1}{2}\Sigma_{SS}+(\frac{1}{2}-\rho) I_k+\Delta_{SS}& 0_{k\times (n-k)}\\
         0_{(n-k)\times k}& (1-\rho)I_{n-k}
    \end{bmatrix}.\]    
    where $\Delta_{SS}=\wh\Sigma_{SS}-\Sigma_{SS}$.
    The proof will be presented according to the followings four steps under the design of \eqref{designsdp}
    \begin{enumerate}
        \item The feasibility of $\wh U(x)$ being a sign matrix hold a.a.s.
        \item The estimate $\wh z_S$ has the same sign for each entry a.a.s.
        \item The solution returned by the oracle program \eqref{non-convexsdp} matches with the program \ref{sdprecovery} a.a.s.
    \end{enumerate}
    Then the rest of the proof we will establish the different aguments sequentially. 
    \begin{center}
        \textbf{(1) The a.a.s. feasibility of $\wh U(x)$}
    \end{center}
    We first study a single entry in $\wh\Sigma_{SS}-\Sigma_{SS}$. Notice that according to the different temperature regimes (alternatively this implies the different expectation of correlation), we will have different convergence guarantees. Assume that $\bb E[\sigma_i\sigma_j]=k^{-\alpha}$ for some $\alpha=0$ without loss of generality, then it is checked that by the boundedness of $\sigma_i\sigma_j$ (which alternatively implies sub-Gaussianity) we will have
    \begin{align*}
        \bb P\bl\inf_{i,j:(i,j)\in S\times S}\wh\Sigma_{ij}-\rho\wh U_{ij}\leq t\br&=\bb P\bl\inf_{i,j:(i,j)\in S\times S}\frac{1}{m}\sum_{\ell\in[m]}\sigma^{(\ell)}_i\sigma^{(\ell)}_j-\bb E[\sigma_i\sigma_j]\leq t-\frac{1}{2}\bb E[\sigma_i\sigma_j]\br\\
        &\leq k^2\exp\bl-m\bl t-\frac{1}{2}\bb E[\sigma_i\sigma_j]\br^2\br.
    \end{align*}
    And for the entries outside $S\times S$ we will apply the union bound again (notice that here the mean is $0$)
    \begin{align*}
        \bb P\bl\sup_{i,j:(i,j)\notin S\times S}\frac{1}{m}\sum_{\ell\in[m]}\sigma_i^{(\ell)}\sigma_j^{(\ell)}\geq t\br\leq n^2\exp(-mt^2).
    \end{align*}
    Then we let $t=Ck^{-\alpha}$ for some $0<C< \frac{1}{2}k^{\alpha}\bb E[\sigma_i\sigma_j]$ one will have the above two terms to be both $o(1)$ when $m=\omega\lef(k^{2\alpha}\log n\rig)$. Therefore the choice of $\rho$ guarantees that when $m\gtrsim k^{2\alpha}\log n$ all elements in $\wh U_{S,S^c}$, $\wh U_{S^c,S^c}$, and $\wh U_{S,S}$ are all with absolute values smaller than $1$. Therefore, the feasibility of $\wh U$ is satisfied for the $\wh U_{ij}$ with $(i,j)\notin S\times S$. For the feasibility for $\wh U_{ij}$ with $(i,j)\in S$ we will need to verify the next step.
    \begin{center}
        \textbf{(2) The a.a.s. correctness of sign of $\wh z_S$ }
    \end{center}
    The proof of correct sign of $\wh z_S$ relies on the eigenvector perturbation bound. We first notice that the population matrix $\Sigma_{SS}-\rho U_{SS}$ will have the principle eigenvector $z_S=\pm\lef(\frac{1}{\sqrt k},\ldots,\frac{1}{\sqrt k}\rig)^\top$. And by a variant of the Davis-Kahan Theorem in \citep{yu2015useful} theorem 2, we will have (assume that for a matrix $A$ with rank $k$, we have $\lambda_1(A)\geq \lambda_2(A)\ldots\geq\lambda_k(A)$)
    \begin{align}\label{daviskahanbound}
        \Vert\wh z_S-z_S \Vert_{2}\leq\begin{cases}
            \frac{ \Vert\Delta_{SS} \Vert_{2}}{\lambda_1(\Sigma_{SS}-\rho U_{SS})-\lambda_2(\Sigma_{SS}-\rho U_{SS})}&\text{ when }\theta_1>0\\
            \frac{ \Vert\Delta_{SS} \Vert_{2}}{\lambda_{k-1}(\Sigma_{SS}-\rho U_{SS})-\lambda_k(\Sigma_{SS}-\rho U_{SS})}&\text{ when }\theta_1<0
        \end{cases}.
    \end{align}
    To analyze the term on the R.H.S. we will need to derive a concentration bounds for the covariance estimation under sub-Weibull concentration. Use $\Delta$ to denote $\Delta_{SS}$ from now on without loss of generality. We use the variational representation $\Vert \Delta\Vert_{2}=\max_{v\in \bb S^{n-1}}| \bfa v^\top\Delta\bfa v|$. 
        Then we can construct an $\epsilon$-net denoted by $\ca N$ covering the $L_2$-sphere $S_2^{k}\subset\bb R^k$ and notice that by the volume argument, the covering number $N(S_2^k,\epsilon)=|\ca N|$ satisfies
    \begin{align*}
        \log N(S_2^k,\epsilon)\lesssim -k\log\epsilon,
    \end{align*}
    And one one will always be able to find $\wha v=\argmin_{\wha v\in \ca N}\Vert\wha v-\bfa v\Vert_2$. We denote $\bfa\delta:=\bfa v-\wha v$ and it is checked that $\Vert\bfa\delta\Vert_2\leq\epsilon$. Then we will apply the decomposition to get
    \begin{align*}
        |\bfa v^\top\Delta\bfa v|&\leq |\wha v^\top\Delta\wha v|+2\Vert\bfa\delta\Vert_2\Vert\Delta\Vert_2\Vert\wha v\Vert_2+\Vert\bfa\delta\Vert_2^2\Vert\Delta\Vert_2\leq|\wha v^\top\Delta\wha v|+(2\epsilon+\epsilon^2)\Vert\Delta\Vert_2.
    \end{align*}
    And we will be able to discretize the operator norm as
    \begin{align*}
        \Vert\Delta\Vert_2\leq\frac{1}{1-2\epsilon-\epsilon^2}\sup_{\bfa v\in\ca N}|\bfa v^\top\Delta\bfa v|.
    \end{align*}
    Then one can upper bound the mgf by
    \begin{align*}
        \bb E[\exp(\lambda\Vert\Delta\Vert_2)]&\leq\bb E\bigg[\exp\bl\frac{\lambda\sup_{\bfa v\in\ca N}|\bfa v^\top\Delta\bfa v|}{1-2\epsilon-\epsilon^2}\br\bigg]\leq\sum_{\bfa v\in\ca N}\bb E\bigg[\exp\bl\frac{\lambda|\bfa v^\top\Delta\bfa v|}{1-2\epsilon-\epsilon^2}\br\bigg].
    \end{align*}
    Then we will notice that by symmetrization argument, introducing i.i.d. Rademacher random variables $\{\epsilon_{i}\}_{i\in[m]}$ we will have for all $t\in\bb R^+$
    \begin{align*}
        \bb E\lef[\exp\lef( t|\bfa v^\top\Delta\bfa v|\rig)\rig]&=\bb E\bigg[\exp\bl t\bigg|\bfa v^\top\bl\frac{1}{m}\sum_{i=1}^m\bfa\sigma^{(i)}\bfa\sigma^{(i)\top}-\Sigma\br\bfa v\bigg|\br\bigg]\\
        &\leq\bb E\bigg[\exp\bl \frac{1}{m}\sum_{i=1}^mt\bigg|\epsilon_i\bfa v^\top\bl\bfa\sigma^{(i)}\bfa\sigma^{(i)\top}-\Sigma\br\bfa v\bigg|\br\bigg]\\
        &\leq \bb E\bigg[\exp\bl 2t\bfa v^\top\bl\frac{1}{m}\sum_{i=1}^m\epsilon_i\bfa\sigma^{(i)}\bfa\sigma^{(i)\top}\br\bfa v\br\bigg]=\bb E\bigg[\exp\bl \frac{2t\epsilon_i}{m}(\bfa\sigma\bfa v^\top)^2\br\bigg]^m.
    \end{align*}
    We further obtain that (denote $\otimes$ as the tensor power)
    \tny{\begin{align*}
        \bb E\bigg[\exp\bl \frac{2t\epsilon_i}{m}(\bfa\sigma\bfa v^\top)^2\br\bigg]&=\sum_{j=0}^\infty\frac{1}{(2j)!}\bl\frac{2t}{m}\br^{2j}\bb E[(\bfa\sigma\bfa v^\top)^{4j}]=\sum_{j=0}^\infty\frac{1}{(2j)!}\bl\frac{2t}{m}\br^{2j}\la\bb E[\bfa\sigma^{\otimes 4j}], \bfa v^{\otimes 4j}\ra\\
        &\leq\sum_{j=0}^\infty\frac{1}{(2j)!}\bl\frac{2t}{km}\br^{2j}\la\bb E[\bfa\sigma^{\otimes 4j}], \mbbm 1^{\otimes 4j}\ra=\sum_{j=0}^\infty\frac{1}{(2j)!}\bl\frac{2t}{km}\br^{2j}\bb E\bigg[\bl\sum_{i=1}^n\sigma_i\br^{4j}\bigg].
    \end{align*}}
    Then, depending on the different temperature regimes given by lemma \ref{cltrfcr} we might encounter the following different situations according to the moment inequalities given by the sub-Weibull distribution (see theorem 2.1 in \citep{vladimirova2020sub})
    \tny{\begin{align*}
        \bb E\bigg[\bl\sum_{i=1}^n\sigma_i\br^{2j}\bigg]^{\frac{1}{2j}}\leq
        \begin{cases}
            C(2j)^{\frac{1}{2}}k^{\frac{1}{2}}& \theta_1<0\text{ or }\theta_1\in
            \lef(0,\frac{1}{\bb E[\sech^2(h)]}\rig)\\
            C(2j)^{\frac{1}{4\tau-2}}k^{\frac{4\tau-3}{4\tau-2}}& \theta_1=\frac{1}{\bb E[\sech^2(h)]}\\
            C(2j)^{\frac{1}{6}}k^{\frac{5}{6}-\frac{1}{3}\beta}& \theta_1=\frac{1}{\bb E[\sech^2(h)]}\text{ and }\bb E[k^{2\beta}\tan^{2\beta}(h)]\asymp 1\text{ for }\beta\in(0,\frac{1}{4})\\
            C(2j)^{\frac{1}{4}}k^{\frac{3}{4}}& \theta_1=\frac{1}{\bb E[\sech^2(h)]}\text{ and }\bb E[k^{2\beta}\tan^{2\beta}(h)]\asymp 1\text{ for }\beta\geq\frac{1}{4}\\
            Ck&  \theta_1\in\lef(\frac{1}{\bb E[\sech^2(h)]},\infty\rig)
        \end{cases}.
    \end{align*}}
    And the different cases of $\theta_1$ finally leads to, there exists $C>0$ such that when $t<Cm$, we have $\frac{1}{1-a}\leq C\exp(2a)$ and therefore
    \tny{\begin{align*}
        \bb E\bigg[\exp\bl \frac{2t\epsilon_i}{m}(\bfa\sigma\bfa v^\top)^2\br\bigg]^m&\leq\bl \sum_{j=0}^\infty\frac{C^j(4j)^{2j}}{(2j)^{2j}}\bl\frac{2t}{m}\br^{2j}\br^m\leq \bl\sum_{j=0}^\infty C^j\bl\frac{t}{m}\br^{2j}\br^{m}=\bl\frac{1}{1-C\frac{t^2}{m^2}}\br^m\\
        &\leq \exp\bl\frac{Ct^2}{m}\br.
    \end{align*}}
    Then, when $\theta_1=\frac{1}{\bb E[\sech^2(h)]}$ we will have when $t<Cmk^{-\frac{2\tau-2}{2\tau-1}}$,
    \tny{\begin{align*}
        \bb E\bigg[\exp\bl\frac{2t\epsilon_i}{m}(\bfa\sigma\bfa v^\top)^2\br\bigg]^m&\leq \bl\sum_{j=0}^\infty\frac{C^j(4j)^{\frac{2j}{2\tau-1}}}{(2j)^{2j}}\bl\frac{2t}{mk^{1-\frac{4\tau-3}{2\tau-1}}}\br^{2j}\br^{m}\leq\bl\sum_{j=0}^\infty \frac{C^j}{j!}\bl\frac{t}{mk^{-\frac{2\tau-2}{2\tau-1}}}\br^{2j}\br^m\\
        &\leq C\exp\bl\frac{Ct^2k^{\frac{4\tau-4}{2\tau-1}}}{m}\br.
    \end{align*}}
    And analogously, we can compute similarly for the rest of the $\theta_1$ and $h$ choices and conclude by the following
    \ttny{\begin{align*}
        \bb E\bigg[\exp\bl\frac{2t\epsilon_i}{m}(\bfa\sigma\bfa v^\top)^2\br\bigg]^m\leq\begin{cases}
            \exp({Ct^2}{m^{-1}}) & \theta_1<0\text{ or }\theta_1\in\lef(0,\frac{1}{\bb E[\sech^2(h)]}\rig)\text{ and }t\leq C^\prime m\\
            \exp(Ct^{2}k^{\frac{4\tau-4}{2\tau-1}}m^{-1})& \theta_1=\frac{1}{\bb E[\sech^2(h)]}\text{ and }t\leq C^\prime mk^{-\frac{2\tau-2}{2\tau-1}}\\
            \exp(Ct^2k^2m^{-1})& \theta_1\in\lef(\frac{1}{\bb E[\sech^2(h)]},\infty\rig)
        \end{cases}
    \end{align*}}
    Therefore, we will conclude that 
    \begin{align*}
        \bb E[\exp(\lambda\Vert\Delta\Vert_2)]\leq\bb E\bigg[\exp\bl\frac{2t\epsilon_i}{m}(\bfa\sigma\bfa v^\top)^2\br\bigg]^m\exp(Ck).
    \end{align*}
    Hence, we will use the Chernoff bound to get
    \begin{align*}
        \bb P(\Vert\Delta\Vert_2\geq t)\leq\inf_{\lambda}\bb E[\exp(\lambda\Vert\Delta\Vert_2)]\exp(-\lambda t).
    \end{align*}
    Using the above method, for all $\delta>0$, we will uniformly write the different scenarios of $\theta_1$ by two functions $f_1:\bb R\to\bb R$ and $f_2:\bb R^{2}\to\bb R$,
    \begin{align}\label{eigenvalueperturb}
        \bb P\bl\frac{1}{f_1(k)}\Vert\Delta\Vert_2\geq Cf_2(k,m)+\delta\br\leq C\exp(-m\delta\wedge\delta^2),
    \end{align}
    with
    \ttny{\begin{align*}
    \lef(f_1(k),f_2(k,m)\rig)=\begin{cases}
             \lef(1,\sqrt{\frac{k}{m}}\vee\frac{k}{m}\rig)& \theta_1<0\text{ or }\theta_1\in\lef(0,\frac{1}{\bb E[\sech^2(h)]}\rig)\\
             \lef(k^{\frac{2\tau-2}{2\tau-1}},\sqrt{{m^{-1}k^{-\frac{2\tau-3}{2\tau-1}}}}\rig)& \theta_1=\frac{1}{\bb E[\sech^2(h)]}\\
            \lef(k, \sqrt {m^{-1}k^{-1}}\rig)& \theta_1\in\lef(\frac{1}{\bb E[\sech^2(h)]},\infty\rig)
        \end{cases}.
    \end{align*}}
    We then notice that by \eqref{daviskahanbound}, noticing that when we have $\Vert z_S-\wh z_S\Vert_2<\frac{1}{\sqrt k}$. To upper bound the principle and the second principle eigenvalues, we will check that
    \begin{align*}
        \lambda_1(\Sigma_{SS}-\rho U_{SS})-\lambda_2(\Sigma_{SS}-\rho U_{SS})=f_1(k)\text{ when }\theta_1>0.
    \end{align*}
    And analogously we will have
    \begin{align*}
        \lambda_{k-1}(\Sigma_{SS}-\rho U_{SS})-\lambda_k(\Sigma_{SS}-\rho U_{SS})=f_1(k)\text{ when }\theta_1<0.
    \end{align*}
    To justify the correct sign alignment of $\wh z_S$ with $z_S$, we choose the correct constant to have
    \begin{align*}
        \bb P\lef(\Vert z_S-\wh z_S\Vert_2\geq Cf_2(k,m)\rig)\leq \bb P\lef(\frac{1}{f_1(k)}\Vert\Delta\Vert_2\geq C^\prime f_2(k,m)\rig)=o(1).
    \end{align*}
    Therefore we let $m$ such that $Cf_2(k,m)<\frac{1}{\sqrt k}$ and get
    \begin{align*}
        m\gtrsim \begin{cases}
            k^2& \theta_1<0\text{ or }\theta_1\in\lef(0,\frac{1}{\bb E[\sech^2(h)]}\rig)\\
             k^{\frac{2}{2\tau-1}}& \theta_1=\frac{1}{\bb E[\sech^2(h)]}\text{ and }t\leq C^\prime mk^{-\frac{2\tau-2}{2\tau-1}}\\
            1& \theta_1\in\lef(\frac{1}{\bb E[\sech^2(h)]},\infty\rig)
        \end{cases}.
    \end{align*}
    Notice that all the above scaling is smaller than the one given by theorem \ref{testingcomputationalupperbound}. Therefore the second condition is verified.
    \begin{center}
        \textbf{(3) The feasibility of rank $1$ solution}
    \end{center}
    The final goal is to validate that the constructed solution is the true one returned by the SDP in algorithm \ref{sdprecovery}. We notice that $U(x)$ is not depending on $x$ and therefore the maximum eigenvector of $\wh\Sigma-\rho U(x)$ satisfies the feasibility condition given by \eqref{ineq}.

\subsection{Proof of Proposition \ref{propexactrecov2}}
Here we use $\gamma=\frac{4\tau-3}{4\tau-2}$ as the shorthand notation.
    We will consider the general formulation, for all $j\in S$, there exist $\alpha\geq 1,\gamma>0,\eta>0$ such that
    \begin{align*}
        \bigg\Vert k^{-\gamma}\bl\sum_{i\in S}\sigma_i\sigma_j-\bb E\bigg[\sum_{i\in S}\sigma_i\sigma_j\bigg]\br \bigg\Vert_{\psi_{\alpha}}<\infty.
    \end{align*}
    We also notice that for the noise part we will have
    \begin{align*}
        \bigg\Vert (n-k)^{-\frac{1}{2}}\sum_{i\in S^c}\sigma_i\sigma_j\bigg\Vert_{\psi_{\alpha}}<\infty.
    \end{align*}
    And therefore we will use the subaddtivity of the Orlicz norm to get
    \tny{\begin{align*}
        \bigg\Vert (n-k)^{-\frac{1}{2}}\wedge k^{-\gamma}\bl\sum_{i=1}^n\sigma_j\sigma_i-\bb E\bigg[\sum_{i=1}^n\sigma_i\sigma_j\bigg]\br\bigg\Vert_{\psi_{\alpha}}&=\bigg\Vert (n-k)^{-\frac{1}{2}}\wedge k^{-\gamma}\bl\sum_{i=1}^n\sigma_j\sigma_i-k^{\eta}\bb E[\delta_j]\br\bigg\Vert_{\psi_{\alpha}}\\
        &<\infty.
    \end{align*}}
    And moreover, for all $i\in S,j\in S^c$ our choice of $k^{-\eta}$ scaling in the expression of $\delta_i$ guarantees that
    \begin{align*}
        |\bb E[\delta_i]-\bb E[\delta_j]|\asymp 1.
    \end{align*}
    Therefore the underlying idea is that we try prove that all the $\delta_i$ where $i\in S$ are clustered into a small region around $\bb E[\delta_i]$ while all the $\delta_i$ where $i\in S^c$ are clustered into another small region around $\bb E[\delta_j]$ such that the two regions are not overlapping with each other.
    
   Notice that the above formulation covers all the critical and high temperature regimes, according to lemma \ref{cltrfcr}, let the scaling ahead of the sum be $k^{-\eta}$ for some $\eta>0$, we will have $\bb E[\delta_i]=1$ and by $\sigma^2_i=1$ we have $\lef\Vert\frac{1}{\sqrt n}\sum_{j=1}^n\sigma_i\sigma_j\rig\Vert_{\psi_2}\lesssim\lef\Vert\frac{1}{\sqrt n}\sum_{j\in[n]}\sigma_j\rig\Vert_{\psi_2}<\infty$. Therefore, by union bound and Chernoff inequality, for all $t>0$,
    \begin{align*}
        \bb P\lef(\sup_{i\in S^c}\delta_i\geq \bb E[\delta_i]+t\rig)&\leq\sum_{i\in S^c}\bb P\bl\frac{1}{\sqrt n}\sum_{\ell=1}^m\sum_{j=1}^n\sigma^{(\ell)}_i\sigma^{(\ell)}_j\geq\frac{1}{\sqrt n}+\frac{k^{\eta}}{\sqrt n}t \br\\
        &\leq (n-k)\exp\bl-\frac{mk^{2\eta}}{n}t^2\br=o(1),
    \end{align*}
    when we have $m=\omega\lef(\frac{n}{k^{2\eta}}\log (n-k)\rig)$.
    Then we will check that for spin $i\in S$, the following holds by union bound
    \sm{\begin{align*}
        \bb P\lef(\inf_{i\in S}\delta_i\leq\bb E[\delta_i]-t\rig)&\leq \sum_{j\in S}\bb P\bl\frac{1}{\sqrt{n-k}\vee k^{\gamma}}\bl\sum_{\ell=1}^m\sum_{i=1}^n\sigma^{(\ell)}_j\sigma^{(\ell)}_i-k^{\eta}\bb E[\delta_j]\br\leq -\frac{tk^{\eta}}{\sqrt{n-k}\vee k^{\gamma}}\br\\
        &\leq k\exp\lef( -m\frac{t^2k^{2\eta}}{(n-k)\vee k^{2\gamma}}\wedge m\frac{tk^{\eta}}{\sqrt{n-k}\vee k^{\gamma}}\rig)=o(1),
    \end{align*}}
    when $m\gtrsim\frac{(n-k)\vee k^{2\gamma}}{k^{2\eta}}\log k\vee\frac{\sqrt{n-k}\vee k^{\gamma}}{k^{\eta}}\log k$. Therefore, to meet both requirements requires that we have $m\gtrsim\frac{(n-k)\vee k^{2\gamma}}{k^{2\eta}}\log k\vee\frac{\sqrt{n-k}\vee k^{\gamma}}{k^{\eta}}\log k\vee \frac{n}{k^{2\eta}}\log(n-k)$.
    Then the rest of the proof will be devoted to verifying the two criteria under the different temperature regimes. First, to analyze $\gamma$, we notice that by convexity and subadditivity of the Orlicz norm, for the high and critical temperature regimes,
    \begin{align*}
        \bigg\Vert k^{-\gamma}\bl\sum_{i\in S}\sigma_i\sigma_j-\bb E\bigg[\sum_{i\in S}\sigma_i\sigma_j\bigg]\br\bigg\Vert_{\psi_{2}}\lesssim\bigg\Vert k^{-\gamma}\sum_{i\in S}\sigma_i\sigma_j\bigg\Vert_{\psi_{2}}\lesssim \bigg\Vert k^{-\gamma}\sum_{i\in S}\sigma_i\bigg\Vert_{\psi_{2}}.
    \end{align*}
    It is then checked that 
    \sm{\begin{align*}
        (\gamma,\eta)=\begin{cases}
        \lef(\frac{1}{2},0\rig)& \theta_1<0\text{ or }\theta_1\in\lef(0,\frac{1}{\bb E[\sech^2(h)]}\rig)\\
         (\frac{4\tau-3}{4\tau-2},\frac{2\tau-2}{2\tau-1})& \theta_1=\frac{1}{\bb E[\sech^2(h)]}\text{ and }0\text{ is the global minimum of \eqref{taylorcond} with flatness }\tau
    \end{cases}
    \end{align*}}
    which finally implies that we will need to have (here we consider $k\gtrsim n^{1/2}$,
    \tny{\begin{align*}
        m\gtrsim \begin{cases}
        n\log n& \theta_1<0\text{ or }\theta_1\in\lef(0,\frac{1}{\bb E[\sech^2(h)]}\rig)\\
         k^{\frac{1}{2\tau-1}}\log k\vee \frac{n}{k^{\frac{4\tau-4}{2\tau-1}}}\log (n-k)& \theta_1=\frac{1}{\bb E[\sech^2(h)]}\text{ and }0\text{ is the global minimum of \eqref{taylorcond} with flatness }\tau
    \end{cases}.
    \end{align*}}

\section{Proof of Minor Lemmas}

 \begin{lemma}\label{combprob}
     For a set containing $n$ elements, picking two $k$ sized subset uniformly at random, denote the random variable representing the cardinality of overlap by $V$. Then we have for $v\leq k$:
     \begin{align*}
        \bb P\lef(V=v\rig)\leq\frac{1}{v!}\lef(\frac{k^2}{n}\rig)^v.
     \end{align*}
     Then for $p<k$ and $k\leq n$, introducing $x=\frac{p}{k}$ and $\gamma=\frac{k}{n}$ we have the following:
      \begin{align*}
        \bb P(V=p)&\leq\frac{\sqrt k}{(1-x)\sqrt{2\pi x}}\exp\bl k\lef(\lef((4-x)\gamma-\log\frac{x}{\gamma}-1\rig)x-2\gamma-2(1-x)\log\lef(1-x\rig)\rig)\\
        &-\frac{1}{12xk+1}+o(1)\br.
     \end{align*}
     And for $\frac{ek^2}{n}< p<k$ and $p\in\bb N$ we have the tail bound:
     \begin{align*}
        \bb P(V\geq p)\leq\frac{1}{(1-\frac{1}{e})\sqrt{2\pi p}}\exp\lef(\lef(1-\log\frac{pn}{k^2}\rig)p-\frac{2k^2}{n}+o\lef(\frac{k^2}{n}\rig)\rig)
     .\end{align*}
 \end{lemma}
 Introducing the handy notation of $a^{(b)} = a\cdot (a-1)\cdots(a-b+1)$. First consider the case when $\lim_{k\to\infty}\frac{k^2}{n}=\lambda<\infty$, we have
      \begin{align}\label{overlapvprob}
         \bb P(V =v )&=\frac{\binom{n-k}{k-v}\binom{k}{v}}{\binom{n}{k}}=\frac{1}{v!}\cdot\frac{(k^{(v)})^2}{n^{(v)}}\cdot\frac{(n-k)^{(k-v)}}{(n-v)^{(k-v)}}\leq\frac{1}{v!}\frac{(k^{(v)})^2}{n^{(v)}}
     .\end{align}
     For the middle term it is checked that
     \begin{align*}
         \frac{(k^{(v)})^2}{n^{(v)}}&=\prod_{i=0}^{v-1}\frac{(k-i)^2}{(n-i)}=\lef(\frac{k^2}{n}\rig)^v\prod_{i=0}^{v-1}\frac{(1-\frac{i}{k})^2}{(1-\frac{i}{n})}\leq\lef(\frac{k^2}{n}\rig)^v
     ,\end{align*} and we complete the proof of the first inequality.
     
For the second inequality, we use Stirling's approximation:
   \ttny{ \begin{align*}
    &\bb P\lef(V=p\rig)=\frac{\binom{n-k}{k-p}\binom{k}{p}}{\binom{n}{k}}\\
    &=\frac{((n-k)!k!)^2}{((k-p)!)^2p!(n-2k+p)!n!}\\
    &\leq \frac{(n-k)k}{(k-p)\sqrt {2\pi p(n-2k+p)n}}\frac{(n-k)^{2(n-k)}k^{2k}}{(k-p)^{2(k-p)}p^pn^n(n-2k+p)^{n-2k+p}}\exp\lef(-\frac{1}{12 p+1}+O\lef(\frac{1}{k}\rig)\rig)\\
    &=\frac{(1-\frac{k}{n})}{(1-\frac{p}{k})\sqrt{2\pi p(1-\frac{2k-p}{n})}}\frac{(1-\frac{k}{n})^{2(n-k)}}{(1-\frac{p}{k})^{2(k-p)}\lef(\frac{pn}{k^2}-\frac{2p}{k}+\frac{p^2}{k^2}\rig)^p(1-\frac{2k}{n}+\frac{p}{n})^{n-2k}}\exp\lef(-\frac{1}{12p+1}+O\lef(\frac{1}{k}\rig)\rig)\\
    &=\frac{1}{(1-\frac{p}{k})\sqrt{2\pi p}}\exp\lef(\lef(\frac{4k}{n}-\frac{p}{n}-\log\frac{pn}{k^2}-1\rig)p-\frac{2k^2}{n}-2(k-p)\log\lef(1-\frac{p}{k}\rig)-\frac{1}{12p+1}+o(1)\rig)
    \end{align*}}
    Introducing $x=\frac{p}{k}$, $\gamma=\frac{k}{n}$ we have
    \ttny{\begin{align*}
        &\bb P(V=p)\leq\frac{\sqrt k}{(1-x)\sqrt{2\pi x}}\exp\lef(k\lef(\lef((4-x)\gamma-\log\frac{x}{\gamma}-1\rig)x-2\gamma-2(1-x)\log\lef(1-x\rig)\rig)-\frac{1}{12xk+1}+o(1)\rig).
    \end{align*}}
For the last inequality, we use the binomial estimation $\lef(\frac{n}{k}\rig)^k\leq\binom{n}{k}\leq\lef(\frac{en}{k}\rig)^k$. Assume that $p>\frac{ek^2}{n}$ and $p<k$ we have
    \tny{\begin{align*}
    &\sum_{V=p}^{k}\bb P\lef(V=v\rig)=\sum_{v=p}^k\frac{\binom{n-k}{k-v}\binom{k}{v}}{\binom{n}{k}}\\
    &= \frac{\binom{n-k}{k-p}\binom{k}{p}}{\binom{n}{k}}\bl 1+\sum_{i=1}^{k-p}\frac{1}{(p+i)^{(i)}}\frac{\lef((k-p)^{(i)}\rig)^2}{(n-p)^{(i)}}\frac{(n-p)^{(i)}}{(n-2k+p+i)^{(i)}}\br\\
    &\leq\frac{\binom{n-k}{k-p}\binom{k}{p}}{\binom{n}{k}}\bl 1+\sum_{i=1}^{k-p}\frac{(k-p)^{2i}}{p^i(n-2k+p)^i}\br\\
    &\leq\frac{\binom{n-k}{k-p}\binom{k}{p}}{\binom{n}{k}}\frac{1}{1-\frac{(k-p)^2}{p(n-2k+p)}}\\
    &=\frac{((n-k)!k!)^2}{((k-p)!)^2p!(n-2k+p)!n!}\frac{1}{1-\frac{(k-p)^2}{p(n-2k+p)}}\\
    &= \frac{(n-k)k}{(k-p)\sqrt {2\pi p(n-2k+p)n}}\frac{(n-k)^{2(n-k)}k^{2k}}{(k-p)^{2(k-p)}p^pn^n(n-2k+p)^{n-2k+p}}\frac{\exp\lef(O(\frac{1}{p}+\frac{1}{(k-p)})\rig)}{1-\frac{(k-p)^2}{p(n-2k+p)}}\\
    &=\frac{(1-\frac{k}{n})}{(1-\frac{p}{k})\sqrt{2\pi p(1-\frac{2k-p}{n})}}\frac{(1-\frac{k}{n})^{2(n-k)}}{(1-\frac{p}{k})^{2(k-p)}\lef(\frac{pn}{k^2}-\frac{p}{k}+\frac{p}{k^2}\rig)^p(1-\frac{2k}{n}+\frac{p}{n})^{n-2k}}\frac{\exp\lef(O(\frac{1}{p}+\frac{1}{(k-p)})\rig)}{1-\frac{(k-p)^2}{p(n-2k+p)}}\\
    &=\frac{1}{(1-\frac{1}{e})\sqrt{2\pi p}}\exp\bl\bl 1-\log\frac{pn}{k^2}\br p-\frac{2k^2}{n}+o\bl\frac{k^2}{n}\br\br.
    \end{align*}}
    This completes the proof of the tail bound.

\begin{lemma}[Laplace Method (Multivariate with Randomness)]\label{laplace}
    Suppose we are given r.v.s. $\bfa h\in\bb R^d$, parameters $\bfa s\in\bb S\subset\bb R^d$ and $\{\Gamma_n( \bfa s,\bfa h)\}$ is a family of random variables in $\Omega$ with $\Gamma_n$ infinitely differentiable w.r.t. $\bfa s$. Furthermore, let us assume that $\Gamma_n$ has unique global minimum almost surely for all $n\in\bb N$ within $\bb S$, and the following are satisfied:
    \begin{enumerate}
        \item  There exists $C(\bfa h)>0$, independent of $n$ and real $\tau$ such that almost surely
        \begin{align}\label{laplace2}
            \exp\lef(-\Gamma(\bfa s,\bfa h)\rig)\leq C(\bfa h)\exp\lef(-\Vert \bfa s\Vert_2^2/2+\tau\Vert\bfa s\Vert_1\rig)
        \end{align} uniformly on compact sets in $\bb R$.
        \item We have  almost surely:
        \begin{align}\label{laplace3}
            \int_{\bb S}\exp(-\Gamma(\bfa s,\bfa h))d\bfa s:=\int_{S_1}\cdots\int_{S_n} \exp(-\Gamma(\bfa s,\bfa h))\prod_{i\in[d]}ds_i<\infty
        .\end{align}
    \end{enumerate}
    Then, we have almost surely there exists random variables $a_1(\bfa h),\ldots a_M(\bfa h)$ for all  $M\in\bb N$ such that
    \tny{\begin{align*}
    \int_{\bb S}\exp(-n\Gamma_n(\bfa s,\bfa h))d\bfa s\sim\exp(-n\Gamma_n(\bfa s_n^*,\bfa h))\det\lef(\frac{n\nabla^2\Gamma_n(\bfa s^*_n,\bfa h)}{2\pi}\rig)^{-1/2}\lef(1+\frac{a_1(\bfa h)}{n}+\ldots+\frac{a_M(\bfa h)}{n^M}\rig) 
    ,\end{align*}}
    where $\nabla$ only take derivative w.r.t. $\bfa s$.
\end{lemma}
    The proof goes by first slicing the integral into two parts denoted by $\bfa V_n(\delta):=\{\bfa s:\Vert \bfa s-\bfa s_n^*\Vert_2\leq\delta\}$ that contains $\bfa s_n^*:=\argmin_{\bfa s\in\bb R^d}\Gamma_n(\bfa s,\bfa h)$ and let $\bfa V^c(\delta)$ be its complement. Note that there exists $\epsilon>0$ such that
    \begin{align*}
        \inf_{\bfa s\in\bfa V^c(\delta)}\Gamma_n(\bfa s,\bfa h)-\inf_{\bfa s\in\bb R^d}\Gamma_n(\bfa s,\bfa h)\geq\epsilon
    .\end{align*}
Hence, using \ref{laplace2} and \ref{laplace3} we note that
\sm{\begin{align*}
        \exp(n\Gamma_n(\bfa s_n^*,\bfa h))&\int_{\bfa V^c(\delta)}\exp(-n\Gamma_n(\bfa s,\bfa h))d\bfa s=\exp(n\Gamma_n(\bfa s_n^*,\bfa h))\\
        &\cdot\int_{\bfa V^c(\delta)}\exp\lef(-(n-1)\Gamma_n(\bfa s,\bfa h)\rig)\exp\lef(-\Gamma_n(\bfa s,\bfa h)\rig)d\bfa s\\
        &\leq\exp\lef(n\Gamma_n(\bfa s_n^*,\bfa h)-(n-1)\inf_{\bfa s\in\bfa V^c(\delta)}\Gamma_n(\bfa s^*,\bfa h)\rig)\int_{\bfa V^c(\delta)}\exp\lef(-\Gamma_n(\bfa s,\bfa h) \rig)d\bfa s\\
        &\leq O(\exp(-n\epsilon))
    .\end{align*}}

    Then we review in the following an important fact and the divergence theorem in vector calculus.
    \begin{fact}\label{fact2}
        Let $\bfa 0$ lie in the interior of $D\subset \bb R^d$. Then as $\lambda\to\infty$ we have
        \begin{align*}
            \int_{D}\exp\lef(-\frac{\lambda}{2}\bfa\xi^\top\bfa\xi\rig)d\bfa\xi =\lef(\frac{2\pi}{\lambda}\rig)^{d/2} +o(\lambda^{-m})
        \end{align*}
        for all  $m\in\bb N$.
    \end{fact}
    \begin{theorem}[Divergence Theorem]\label{diverge} Suppose $\bfa D$ is a subset of  $\bb R^d$ with  $\bfa D$ a compact space with piecewise smooth boundary  $\bfa S=\partial \bfa D$. If $\bfa F$ is a continuously differentiable vector field defined on a neighborhood of $\bfa D$ then
    \begin{align*}
        \int_{\bfa D}(\nabla\cdot \bfa F)dV=\oint_{\bfa S}(\bfa F\cdot\bfa n)dS
    .\end{align*}
    where $\bfa n$ is the unit outward normal vector to $\bfa S$ and $d S$ is the differential element on the hypersurface $\bfa S$.
    
    Changing $\bfa F$ to $\bfa Fg$ for some smooth scalar function  $g$ we have
    \begin{align*}
        \int_{\bfa D}\lef(\bfa F\cdot\nabla g+g\nabla\cdot\bfa F\rig)dV=\oint_{\bfa S}g\bfa F\cdot \bfa ndS
    .\end{align*}
    \end{theorem}
    
    The next step is to consider what will lie in $\bfa V_n(\delta_1)$.  The proof strategy follows from \citep{barndorff1989asymptotic} and \citep{bleistein1975asymptotic}. By Taylor expansion there exists $\delta_2>0$ sufficiently small such that for all  $\bfa s\in\bfa V_n(\delta_2)$ we have
    \begin{align*}
        \Gamma_{n}(\bfa s,\bfa h) - \Gamma_n(\bfa s_n^*,\bfa h)&=\frac{1}{2}(\bfa s-\bfa s_n^*)^\top\nabla^2 \Gamma_n(\bfa s_n^*,\bfa h)(\bfa s-\bfa s_n^*)+o\lef(\Vert \bfa s-\bfa s^*_n\Vert_2^2\rig)\\
        &=\frac{1}{2}\bfa z^\top\bfa z +o(\Vert\bfa z\Vert_2^2)
    .\end{align*}
    where  $\bfa z:=\lef(\nabla^2\Gamma_n(\bfa s_n^*,\bfa h)\rig)^{1/2}(\bfa s-\bfa s_n^*)$ .
    Then we can introduce $\bfa m:\bb R^d\to\bb R^d$ such that $m_i(\bfa z)= z_i+o(z_i)$ as $z_i\to 0$ and satisfying
    \begin{align*}
        \Gamma_n(\bfa s,\bfa h)-\Gamma_n(\bfa s^*_n,\bfa h)&=\frac{1}{2}\bfa m^\top(\bfa s)\bfa m(\bfa s)
    .\end{align*}
    Defining the function $ G_0(\bfa m):=\ca J(\bfa m)=\frac{\pta (s_1,\ldots,s_d)}{\pta (m_1,\ldots,m_d)}$ to be the Jacobian at $\bfa s$ and we note that
$  \ca J(\bfa 0)=\lef|\det(\nabla^2\Gamma_n(\bfa s_n^*,\bfa h))\rig|^{-1/2}
$   Introducing  $\bfa D$ to be the image of $\bfa V_n(\delta)$ under the two round of change of variables and $\bfa S=\pta\bfa D$. Therefore with the above preparation we can write the integral as:
\begin{align*}
        &\int_{\bfa V_n(\delta_1)}\exp\lef(-n\Gamma_n(\bfa s,\bfa h)\rig)d\bfa s = \exp(-n\Gamma_n(\bfa s_n^*,\bfa h))\int_{\bfa D} G_0(\bfa m)\exp\lef(-\frac{n}{2}\bfa m^\top\bfa m\rig)d\bfa m
.\end{align*}
Note that there exists a function $\bfa H_0:\bb R^d\to\bb R^d$ such that $ G_0(\bfa m)= G_0(\bfa 0)+\bfa m^\top\bfa H_0(\bfa m)$
We then use theorem \ref{diverge} to get 
\begin{align*}
    &\int_{\bfa V_n(\delta_1)}\exp\lef(-n\Gamma_n(\bfa s,\bfa h)\rig)d\bfa s=\exp(-n\Gamma_n(\bfa s_n^*,\bfa h))\bigg[\int_{\bfa D} G_0(\bfa 0)\exp\lef(-\frac{n}{2}\bfa m^\top\bfa m\rig)d\bfa m\\
    &-\frac{1}{n}\int_{\bfa S}(\bfa H_0(\bfa m)\cdot\bfa n)\exp\lef(-\frac{n}{2}\bfa m^\top\bfa m\rig)dS+\frac{1}{n}\int_{\bfa D} G_1(\bfa m)\exp\lef(-\frac{n}{2}\bfa m^\top\bfa m\rig)d\bfa m\bigg]
.\end{align*}
And we can do the above process recursively and get
\begin{align*}
    &\int_{\bfa V_n(\delta_1)}\exp\lef(-n\Gamma_n(\bfa s,\bfa h)\rig)d\bfa s=\exp(-n\Gamma_n(\bfa s_n^*,\bfa h))\bigg[\sum_{j=0}^{M} G_j(\bfa m)\int_{\bfa D}\exp\lef(-\frac{n}{2}\bfa m^\top\bfa m\rig)d\bfa m\\
    &-\frac{1}{n^M}\int_{\bfa D}G_M(\bfa s)\exp\lef(-\frac{n}{2}\bfa m^\top\bfa m\rig)d\bfa s\bigg]
.\end{align*}
since we note that the boundary integral is exponentially small almost surely according to \ref{laplace2} as $n\to\infty$ and could be ignored here. Note that $ G_j$ is defined recursively as
\begin{align*}
    G_j(\bfa m):=G_j(\bfa 0)+\bfa m^\top\cdot \bfa H_j(\bfa m),\\
    G_{j+1}(\bfa m):=\nabla\cdot \bfa H_j(\bfa m)
.\end{align*}
Further notice that by \ref{laplace2} and \ref{laplace3} we can check that almost surely:
\begin{align*}
    \lef|\frac{1}{n^M}\int_{D}\exp\lef(-\frac{n}{2}\bfa m^\top\bfa m\rig)G_M(\bfa m)d\bfa m\rig|=O\lef(\frac{1}{n^M}\rig)
.\end{align*}
Together with the fact \ref{fact2} we can see that almost surely:
\begin{align*}
    \int_{\bfa V_n(\delta_1)}\exp(-n\Gamma_n(\bfa s,\bfa h))d\bfa s=\exp\lef(-n\Gamma_n(\bfa s_n^*,\bfa h)\rig)\lef(\frac{2\pi}{n}\rig)^{d/2}\bl\sum_{j\in[M-1]}\frac{G_j(\bfa 0)}{n^j}+O\lef(\frac{1}{n^M}\rig)\br
.\end{align*}
Note that $G_j(\bfa 0)$ are functions of $\bfa h$ we complete the proof by defining $a_k(\bfa h)=\frac{G_k(\bfa 0)}{G_0(\bfa 0)}$.
\begin{align*}
.\end{align*}

\section{Additional Standard Arguments}

\subsection{Concentration of Measure}
\begin{lemma}[Chi-square Tail \citep{laurent2000adaptive}]\label{chisqtail}
    Let $(Y_1,\ldots, Y_D)$ be i.i.d. Gaussian variables, with mean $0$ and variance $1$. Let $\alpha_1,\ldots,\alpha_D$ be non-negative. Set 
    \begin{align*}
        |\alpha|_\infty=\sup_{i=1,\ldots, D}|\alpha_i|,\qquad |\alpha|_2^2=\sum_{i=1}^D\alpha_i^2.
    \end{align*}
    Let $Z=\sum_{i=1}^D\alpha_i(Y_i^2-1)$.
    Then the following inequality holds for all $x>0$:
    \begin{align*}
        \bb P\lef(Z\geq 2|\alpha|_2\sqrt x+2|\alpha|_\infty x\rig)&\leq\exp\lef(-x\rig),\\
        \bb P\lef(Z\leq-2|\alpha|_2\sqrt x\rig)&\leq\exp\lef(-x\rig).
    \end{align*}
\end{lemma}
    \begin{lemma}[Bernstein Inequality \citep{ledoux1991probability}]\label{regularbernstein}
Let $X_1,\ldots, X_n$ be independent centered sub-exponential random variables, and $K=\max_{i}\Vert X_i\Vert_{\psi_1}$. Then for every $a=(a_1,\ldots, a_n)\in\bb R^n$ and every $t\geq 0$ we have
\begin{align*}
    \bb P\lef(\lef|\sum_{i=1}^N a_iX_i\rig|\geq t\rig)\leq 2\exp\lef(-c\min\lef(\frac{t^2}{K^2\Vert a\Vert_2^2},\frac{t}{K\Vert a\Vert_\infty}\rig)\rig)
.\end{align*}
where $c>0$ is an absolute constant.
\end{lemma}
\subsection{Laplace Method}
\begin{lemma}[Laplace Method (Univariate with interior maximum)\citep{small2010expansions}] \label{laplaceu1}
    Let $-\infty\leq a<b\leq\infty$. Let $h(x)$ be defined on the open interval $(a,b)$. Suppose the following are satisfied:
    \begin{enumerate}
        \item The function $h(x)$ is differentiable throughout $(a,b)$, is uniquely maximised at some point $x_0\in(a,b)$, and that $h^{(2)}(x_0)$ exists and is strictly negative. $g(x)$ is a continuous function defined on the open interval $(a,b)$ s.t. $g(x_0)\neq 0$.
        \item There exists constant $\eta>0$ and $\delta>0$ such that $h(x)<h(x_0)-\eta$ for all $x\in(a,b)$ such that $|x-x_0|\geq\delta$.
        \item The integral below exists for $n=1$.
        Then we have
        \begin{align*}
            \int_a^bg(x)\exp(nh(x))dx\sim g(x_0)\exp\lef(nh(x_0)\rig)\sqrt{\frac{2\pi}{-nh^{(2)}(x_0)}}
        \end{align*}
        as $n\to\infty$.
    \end{enumerate}
\end{lemma}

\begin{lemma}[Laplace Method (Univeriate with boundary maximum)\citep{small2010expansions}]\label{laplaceu2}
    Suppose that the following condition holds:
    \begin{enumerate}
        \item We have $h(x)<h(a)$ for all $a<x<b$, and for all $\delta>0$ we have
        \begin{align*}
            \inf\{h(a)-h(x):x\in[a+\delta,b)\}>0.
        \end{align*}
        \item The functions $h^\prime(x)$ and $g(x)$ are continuous in a neighborhood of $x=a$.
        \item The following integral converges absolutely for sufficiently large $n$.
    \end{enumerate}
    Then we have
        \begin{align*}
            \int_a^bg(x)\exp\lef(nh(x)\rig)\sim \exp(nh(a))\frac{g(a)}{-nh^\prime(a)}.
        \end{align*}
\end{lemma}
    \begin{lemma}[Linearization of Z estimator, Higher Order \citep{van2000asymptotic}]\label{linearZ} Let $\theta\geq 1$.
      Let $\Psi_n$ and $\Psi$ be a sequence of random maps and a fixed map, respectively, from $\Theta$ to a Banach space such that
      \begin{align}\label{equation34}
          \sqrt n(\Psi_n-\Psi)(\wh\theta_n)-\sqrt n(\Psi_n-\Psi)(\theta_0)=o_{\psi_{\theta}}(1+\sqrt n\Vert \wh\theta_n-\theta_0\Vert^{\gamma}),\qquad \gamma\in\bb N^*
      \end{align}
      and such that the sequence $\sqrt n(\Psi_n-\Psi)(\theta_0)$ converges weakly to a tight random element $Z$. Let $\theta\to\Psi(\theta)$ be $\gamma$-order Frechet-differentiable at $\theta_0$ or
      \begin{align*}
          \bigg\Vert\Psi(\theta)-\Psi(\theta_0)-\nabla^{\gamma}\Psi(\theta_0)\cdot(\theta-\theta_0)^{\otimes\gamma}\bigg\Vert=o(\Vert\theta-\theta_0\Vert^{\gamma}),
      \end{align*}
      with a continuous, invertible, and diagonal derivative $\nabla^{\gamma}\Psi_{\theta_0}$. If $\Psi(\theta_0)=0$ and $\wh\theta_n$ satisfies $\Psi(\wh\theta_n)=o_{\psi_{\theta}}(n^{-1/2})$ and converges in $\psi_{\theta}$ norm to $\theta_0$, then
      \begin{align*}
          \sqrt n\nabla^{\gamma}\Psi(\theta_0)\cdot(\wh\theta_n-\theta_0)^{\otimes\gamma}=-\sqrt n(\Psi_n-\Psi)(\theta_0)+o_{\psi_{\theta}}(1).
      \end{align*}
    \end{lemma}
    \begin{proof}
        First, by the definition of $\wh\theta_n$ and equation \eqref{equation34} we have
        \tny{\begin{align}\label{eqqexpand}
            \sqrt n(\Psi(\wh\theta_n)-\Psi(\theta_0))=\sqrt n(\Psi(\wh\theta_n)-\Psi_n(\wh\theta_n))+o_{\psi_{\theta}}(1)=-\sqrt n(\Psi_n-\Psi)(\theta_0)+o_{\psi_{\theta}}(1+\sqrt n\Vert\wh\theta_n-\theta_0\Vert^{\gamma}).
        \end{align}}
        Using the continuously invertibility and diagonality of $\nabla^{\gamma}\Psi_{\theta_0}$ one obtain that there exists a positive constant $c$ such that
        \begin{align*}
            \Vert\nabla^{\gamma}\Psi(\theta_0)(\theta-\theta_0)^{\otimes\gamma}\Vert\geq c\Vert\theta-\theta_0\Vert^{\gamma},
        \end{align*}
        for every $\theta$ and $\theta_0$. Frechet differentiability condition then yields 
        \begin{align*}
            \Vert\Psi(\theta)-\Psi(\theta_0)\Vert\geq c\Vert\theta-\theta_0\Vert^{\gamma}+o(\Vert\theta-\theta_0\Vert^{\gamma}).
        \end{align*}
        Then we immediately have
        \begin{align*}
            \sqrt n\Vert\wh\theta_n-\theta_0\Vert^{\gamma}(c+o_{\psi_{\theta}}(1))\leq O_{\psi_{\theta}}(1)+o_{\psi_{\theta}}(1+\sqrt n\Vert \wh\theta_n-\theta_0\Vert^{\gamma}).
        \end{align*}
        And we have
        \begin{align*}
            \sqrt n\Vert\wh\theta_n-\theta_0\Vert^{\gamma}=O_{\psi_{\theta}}(1).
        \end{align*}
        Then we apply the Frechet differentiability condition to the l.h.s. of \eqref{eqqexpand} to obtain that
        \begin{align*}
            \sqrt n(\Psi(\wh\theta_n)-\Psi(\theta_0))=\sqrt n\nabla^{\gamma} \Psi({\theta_0})(\wh\theta_n-\theta_0)^{\otimes\gamma}+o_{\psi_{\theta}}(1)=-\sqrt n(\Psi_n-\Psi)(\theta_0)+o_{\psi_{\theta}}(1).
        \end{align*}
    \end{proof}

\end{document}